\documentclass[oneside,11pt]{article}
\usepackage[a4paper, total={6in, 8in}]{geometry}
\usepackage[numbers]{natbib}
\usepackage{multirow,makecell}

\usepackage{bm,amsmath,amsfonts,amsthm,booktabs}
\usepackage{mathtools}
\usepackage[justification=justified]{caption}
\usepackage{url}
\usepackage{enumitem}
\usepackage{xcolor}
\usepackage[ruled,vlined,linesnumbered]{algorithm2e}
\SetArgSty{textnormal}
\usepackage{subcaption}
\usepackage{rotating}
\usepackage{array}
\usepackage{verbatim}
\usepackage{tikz-cd}
\usepackage{graphbox}
\usepackage{float}
\usepackage[export]{adjustbox}

\newcolumntype{P}[1]{>{\centering\arraybackslash}p{#1}}
\newcolumntype{M}[1]{>{\centering\arraybackslash}m{#1}}

\newcommand{\edit}[1]{\textcolor{black}{#1}}

\newcommand{\editt}[1]{\textcolor{black}{#1}}

\newcommand{\edittt}[1]{\textcolor{black}{#1}}

\DeclareSymbolFont{ebgletters}{OML}{EBGaramond-Maths}{m}{it}
\DeclareMathSymbol{w}{\mathalpha}{ebgletters}{`k}

\newtheorem{thm}{Theorem}

\newtheorem{rmk}{Remark}

\newcommand{\argmin}{\mathop{\mathrm{argmin}}\limits}
\newcommand{\argmax}{\mathop{\mathrm{argmax}}\limits}

\begin{document}

\title{LDLE: Low Distortion Local Eigenmaps}


\author{Dhruv Kohli \\  dhkohli@ucsd.edu\\
        Department of Mathematics\\
        University of California San Diego\\
        CA 92093, USA \and 
        Alexander Cloninger\\acloninger@ucsd.edu\\
        Department of Mathematics\\
        University of California San Diego\\
        CA 92093, USA \and
        Gal Mishne\\ gmishne@ucsd.edu\\
        Halicio\u{g}lu Data Science Institute\\
        University of California San Diego\\
        CA 92093, USA}


\maketitle

\begin{abstract}
We present Low Distortion Local Eigenmaps (LDLE), a manifold learning technique which constructs a set of low distortion local views of a dataset in lower dimension and registers them to obtain a global embedding. The local views are constructed using the global eigenvectors of the graph Laplacian and are registered using Procrustes analysis. The choice of these eigenvectors may vary across the regions. In contrast to existing techniques, \edit{LDLE can embed closed and non-orientable manifolds into their intrinsic dimension by tearing them apart. It also provides gluing instruction on the boundary of the torn embedding to help identify the topology of the original manifold. Our experimental results will show that LDLE largely preserved distances up to a constant scale while other techniques produced higher distortion. We also demonstrate that LDLE produces high quality embeddings even when the data is noisy or sparse.}
\end{abstract}


\section{Introduction}
Manifold learning techniques such as Local Linear Embedding \citep{roweis2000nonlinear}, Diffusion maps \citep{coifman2006diffusion}, Laplacian eigenmaps \citep{belkin2003laplacian}, t-SNE \citep{maaten2008visualizing} and UMAP \citep{mcinnes2018umap}, aim at preserving local \edit{information} as they map a manifold embedded in higher dimension into lower (possibly intrinsic) dimension. In particular, UMAP and t-SNE follow a top-down approach as they start with an initial low-dimensional global embedding and then refine it by minimizing a local distortion measure on it. In contrast, similar to \editt{LTSA \citep{zhang2003nonlinear} and} \citep{singer2011orientability}, a bottom-up approach for manifold learning can be imagined to consist of two steps, first obtaining low distortion local views of the manifold in lower dimension and then registering them to obtain a global embedding of the manifold. \editt{In this paper, we take this bottom-up perspective to embed a manifold in low dimension, where the local views are obtained by constructing coordinate charts for the manifold which incur low distortion.}

\subsection{Local Distortion}
\label{subsec:local_dist}

Let $(\mathcal{M},g)$ be a $d$-dimensional Riemannian manifold with finite volume. By definition, for every $x_k$ in $\mathcal{M}$, there exists a coordinate chart $(\mathcal{U}_k,\Phi_k)$ such that $x_k \in \mathcal{U}_k$, $\mathcal{U}_k \subset M$ and $\Phi_k$ maps $\mathcal{U}_k$ into $\mathbb{R}^d$. One can imagine $\mathcal{U}_k$ to be a local view of $\mathcal{M}$ in the ambient space. Using rigid transformations, these local views can be registered to recover $\mathcal{M}$. Similarly, $\Phi_k(\mathcal{U}_k)$ can be imagined to be a local view of $\mathcal{M}$ in the $d$-dimensional embedding space $\mathbb{R}^d$. Again, using rigid transformations, these local views can be registered to obtain the $d$-dimensional embedding of $\mathcal{M}$.

As there may exist multiple mappings which map $\mathcal{U}_k$ into $\mathbb{R}^d$, a natural strategy would be to choose a mapping with low distortion.
\edit{Multiple measures of distortion exist in literature \citep{NEURIPS2018_4c5bcfec}. The measure of distortion used in this work is as follows.} Let $d_g(x,y)$ denote the shortest geodesic distance between $x,y \in \mathcal{M}$. The distortion of $\Phi_k$ on $\mathcal{U}_k$ as defined in \citep{jones2007universal} is given by
\begin{align}
    \text{Distortion}(\Phi_k,\mathcal{U}_k) = \left\|\Phi_k\right\|_\text{Lip}\left\|\Phi_k^{-1}\right\|_\text{Lip} \label{DistortionPhikUk}
\end{align}
where $\left\|\Phi_k\right\|_\text{Lip}$ is the Lipschitz norm of $\Phi_k$ given by
\begin{align}
    \left\|\Phi_k\right\|_\text{Lip} &= \sup_{\substack{x,y \in \mathcal{U}_k\\x\neq y}}\frac{\left\|\Phi_k(x)-\Phi_k(y)\right\|_2}{d_g(x,y)},
\end{align}
and similarly,
\begin{align}
    \left\|\Phi^{-1}_k\right\|_\text{Lip} &= \sup_{\substack{x,y \in \mathcal{U}_k\\x\neq y}}\frac{d_g(x,y)}{\left\|\Phi_k(x)-\Phi_k(y)\right\|_2}.
\end{align}
\edit{Note that $\text{Distortion}(\Phi_k,\mathcal{U}_k)$ is always greater than or equal to $1$.} If $\text{Distortion}(\Phi_k,\mathcal{U}_k) = 1$, then $\Phi_k$ is said to have no distortion on $\mathcal{U}_k$. \edit{This is achieved when the mapping $\Phi_k$ preserves distances between points in $\mathcal{U}_k$ up to a constant scale, that is, when $\Phi_k$ is a similarity on $\mathcal{U}_k$}.  
It is not always possible to obtain a mapping with no distortion. For example, \edit{there does not exist a similarity which maps a locally curved region on a surface into a Euclidean plane. This follows from the fact that the sign of the Gaussian curvature is preserved under similarity transformation which in turn follows from the Gauss’s Theorema Egregium.}

\subsection{Our Contributions}
\label{subsec:contrib}
This paper takes motivation from the work in \citep{jones2007universal} where the authors provide guarantees on the distortion of the coordinate charts of the manifold constructed using \editt{\textit{carefully}} chosen eigenfunctions of the Laplacian. However, this only applies to the charts for small neighborhoods on the manifold and does not provide a global embedding. In this paper, we present an approach to realize their work in the discrete setting and obtain low-dimensional low distortion local views of the given dataset \editt{using the eigenvectors of the graph Laplacian}. Moreover, we piece together these local views to obtain a global embedding of the manifold. The main contributions of our work are as follows:
\begin{enumerate}
    \item We present an algorithmic realization of the construction procedure in \citep{jones2007universal} that applies to the discrete setting and yields low-dimensional low distortion views of small metric balls on the given discretized manifold (See Section~\ref{sec:motivation} for a summary of their procedure).
    \item We present an algorithm to obtain a global embedding of the manifold by registering its local views. The algorithm is designed so as to embed closed as well as non-orientable manifolds into their intrinsic dimension by tearing them apart. It also provides gluing instructions for the boundary of the embedding by coloring it such that the points on the boundary which are adjacent on the manifold have the same color (see Figure~\ref{fig:fig14}).
\end{enumerate}

LDLE consists of three main steps. In the first step, we estimate the inner product of the Laplacian eigenfunctions' gradients using the local correlation between them. These estimates are used to choose eigenfunctions which are in turn used to construct low-dimensional low distortion parameterizations $\Phi_k$ of the small balls $U_k$ on the manifold. The choice of the eigenfunctions depend on the underlying ball. 
A natural next step is to align these local views $\Phi_k(U_k)$ in the embedding space, to obtain a global embedding. One way to align them is to use Generalized Procrustes Analysis (GPA) \citep{fabioprocrustes,gower1975generalized,ten1977orthogonal}. However, we empirically observed that GPA is less efficient and prone to errors due to large number of local views with small overlaps between them.
Therefore, motivated from our experimental observations and computational necessity, in the second step, we develop a clustering algorithm to obtain a small number of intermediate views $\widetilde{\Phi}_m(\widetilde{U}_m)$ with low distortion, from the large number of smaller local views $\Phi_k(U_k)$. This makes the subsequent GPA based registration procedure faster and less prone to errors. 

Finally, in the third step, we register intermediate views $\widetilde{\Phi}_m(\widetilde{U}_m)$ using an adaptation of GPA which enables tearing of closed and non-orientable manifolds so as to embed them into their intrinsic dimension. The results on a 2D rectangular strip and a 3D sphere are presented in Figures~\ref{fig:fig13} and~\ref{fig:fig14}, to motivate our approach.



The paper organization is as follows. Section~\ref{sec:motivation} provides relevant background and motivation. In Section~\ref{sec:biliplocparam} we present the construction of low-dimensional low distortion local parameterizations. Section~\ref{sec:clustering} presents our clustering algorithm to obtain intermediate views. Section~\ref{sec:ge} registers the intermediate views to a global embedding. In Section~\ref{sec:compare} we compare the embeddings produced by our algorithm with existing techniques on multiple datasets. Section~\ref{sec:conclusion} concludes our work and discusses future directions.

\subsection{Related Work}
\label{subsec:related}

\editt{Laplacian eigenfunctions are ubiquitous in manifold learning. A large proportion of the existing manifold learning techniques rely on a fixed set of Laplacian eigenfunctions, specifically, on the first few non-trivial low frequency eigenfunctions, to construct a low-dimensional embedding of a manifold in high dimensional ambient space. These low frequency eigenfunctions not only carry information about the global structure of the manifold but they also exhibit robustness to the noise in the data \citep{coifman2006diffusion}. Laplacian eigenmaps \citep{belkin2003laplacian}, Diffusion maps \citep{coifman2006diffusion} and UMAP \citep{mcinnes2018umap} are examples of such top-down manifold learning techniques. While there are limited bottom-up manifold learning techniques in the literature, to the best of our knowledge, none of them makes use of Laplacian eigenfunctions to construct local views of the manifold in lower dimension}.


\editt{
\paragraph{LTSA} is an example of a bottom-up approach for manifold learning whose local mappings project local neighborhoods onto the respective tangential spaces. A local mapping in LTSA is a linear transformation whose columns are the principal directions obtained by applying PCA on the underlying neighborhood. \textit{These directions form an estimate of the basis for the tangential space}. 
Having constructed low-dimensional local views for each neighborhood, LTSA then aligns all the local views to obtain a global embedding. As discussed in their work and as we will show in our experimental results, LTSA lacks robustness to the noise in the data. This further motivates our approach of using robust low-frequency Laplacian eigenfunctions for the construction of local views. Moreover, due to the specific constraints used in their alignment, LTSA embeddings fail to capture the aspect ratio of the underlying manifold (see Appendix~\ref{subsec:ltsa_stitch} for details).}



\begin{figure}[!h]
    \centering
    \begin{tabular}{|M{0.1\textwidth}|M{0.1\textwidth}|M{0.1\textwidth}|M{0.1\textwidth}|M{0.1\textwidth}|M{0.1\textwidth}|M{0.1\textwidth}|}
        \tiny{Input} & \tiny{LDLE} & \tiny{LDLE with $\partial\mathcal{M}$ known apriori} & \tiny{LTSA} & \tiny{UMAP} & \tiny{t-SNE} & \tiny{Laplacian Eigenmaps}\\
        \hline
        \rotatebox{90}{\includegraphics[width=0.15\textwidth,keepaspectratio]{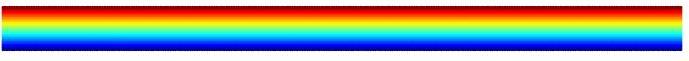}} & \rotatebox{90}{\includegraphics[width=0.15\textwidth,keepaspectratio]{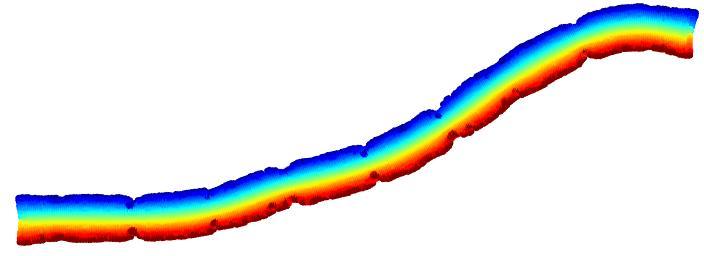}} & \includegraphics[height=0.15\textwidth,keepaspectratio]{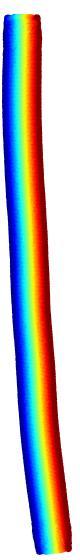} &         \rotatebox{90}{\includegraphics[height=0.1\textwidth,keepaspectratio]{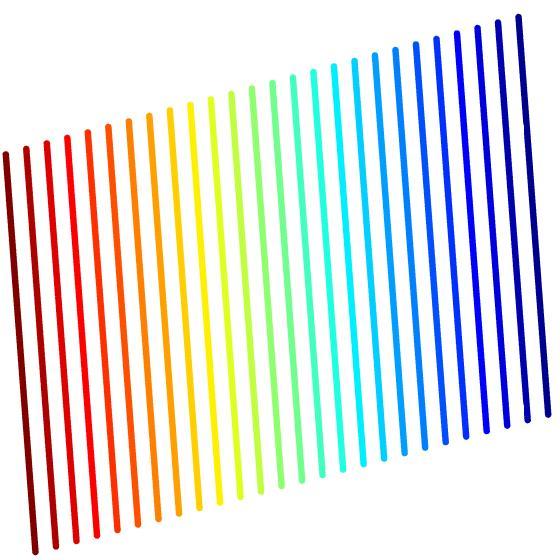}} & \rotatebox{90}{\includegraphics[width=0.135\textwidth,keepaspectratio]{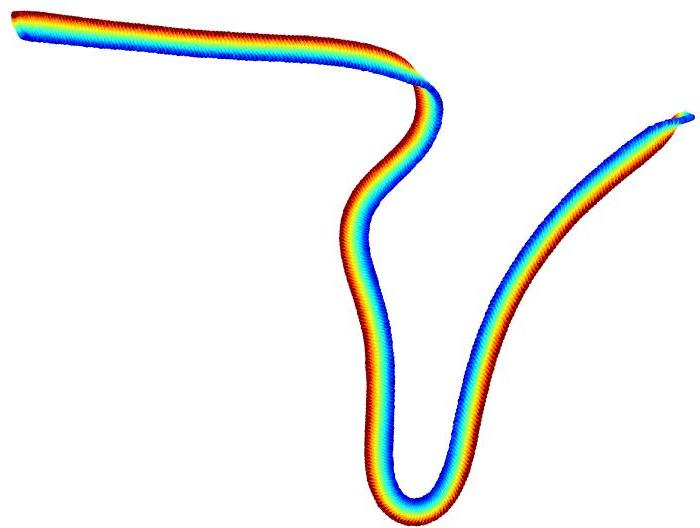}} & \includegraphics[width=0.10\textwidth,height=0.10\textwidth,keepaspectratio]{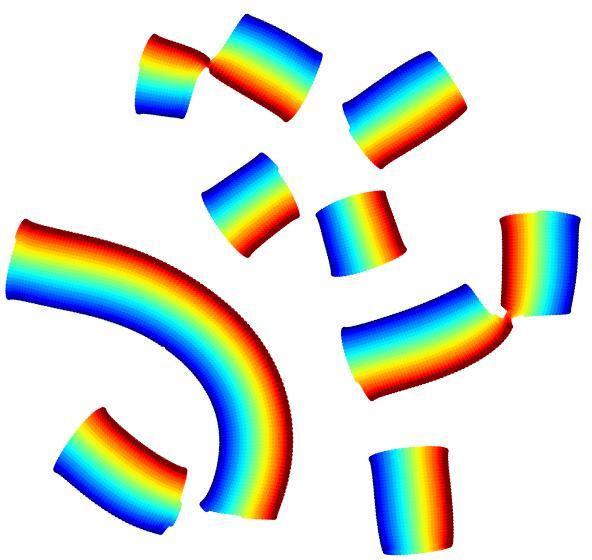} & \includegraphics[width=0.10\textwidth,height=0.10\textwidth,keepaspectratio]{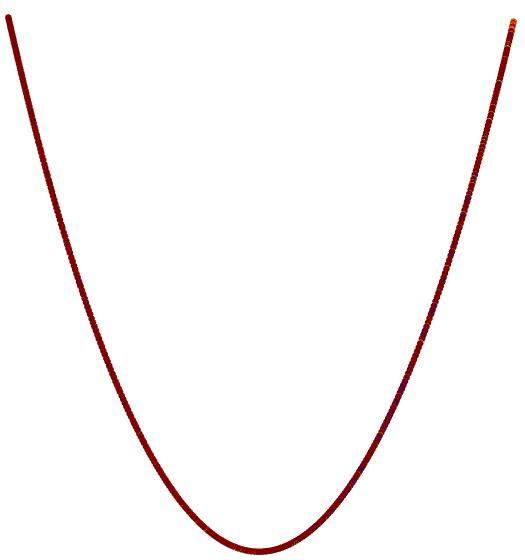} \\
        \hline
    \end{tabular}
    \caption{Embeddings of a rectangle ($4 \times 0.25$) with high aspect ratio in $\mathbb{R}^2$ into $\mathbb{R}^2$.}
    \label{fig:fig13}
\end{figure}

\paragraph{Laplacian eigenmaps} uses the eigenvectors corresponding to the $d$ smallest eigenvalues (excluding zero) of the normalized graph Laplacian to embed the manifold in $\mathbb{R}^d$. \editt{It can also be perceived as a top-down approach which directly obtains a global embedding that minimizes Dirichlet energy under some constraints. For manifolds with high aspect ratio, in the context of Section~\ref{subsec:local_dist}, the distortion of the local parameterizations based on the restriction of these eigenvectors on local neighborhoods, could become extremely high.} For example, as shown in Figure~\ref{fig:fig13}, the Laplacian eigenmaps embedding of a rectangle with an aspect ratio of $16$ looks like a parabola. This issue is explained in detail in \citep{8450808,NEURIPS2019_6a10bbd4,DSILVA2018759,DBLP:journals/corr/BlauM16}.


\paragraph{UMAP,} to a large extent, resolves this issue by first computing an embedding based on the $d$ non-trivial low-frequency eigenvectors of a symmetric normalized Laplacian and then ``sprinkling" white noise in it. It then refines the noisy embedding by minimizing a local distortion measure based on fuzzy set cross entropy. Although UMAP embeddings seem to be topologically correct, they occasionally tend to have twists and sharp turns which may be unwanted (see Figure~\ref{fig:fig13}).


\paragraph{t-SNE} takes a different approach of randomly initializing the global embedding, defining a local t-distribution in the embedding space and local Gaussian distribution in the high dimensional \editt{ambient space}, and finally refining the embedding by minimizing the 
Kullback–Leibler divergence between the two sets of distributions. As shown in Figure~\ref{fig:fig13}, t-SNE tends to output a dissected embedding even when the manifold is connected. \edit{Note that the recent work by \citep{kobak2021initialization} showed that t-SNE with spectral initialization results in a similar embedding as that of UMAP. Therefore, in this work, we display the output of the classic t-SNE construction, with random initialization only.}

\newblock

\begin{figure}[!h]
    \centering
    \begin{tabular}{|M{0.125\textwidth}|M{0.125\textwidth}|M{0.125\textwidth}|M{0.125\textwidth}|M{0.125\textwidth}|M{0.125\textwidth}|}
        \tiny{Input} & \tiny{LDLE} & \tiny{LTSA} & \tiny{UMAP} & \tiny{t-SNE} & \tiny{Laplacian Eigenmaps}\\
        \hline
        \includegraphics[width=0.125\textwidth,height=0.25\textwidth,keepaspectratio]{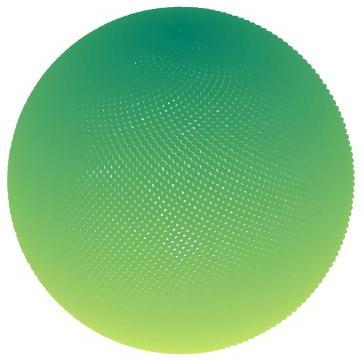} & \includegraphics[width=0.125\textwidth,height=0.25\textwidth,keepaspectratio]{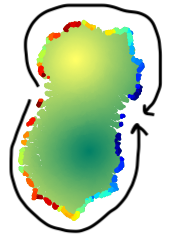} & \includegraphics[width=0.125\textwidth,height=0.25\textwidth,keepaspectratio]{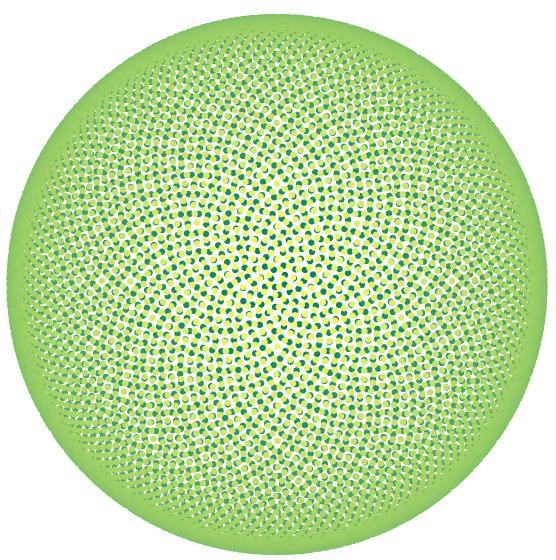} & \includegraphics[width=0.125\textwidth,height=0.25\textwidth,keepaspectratio]{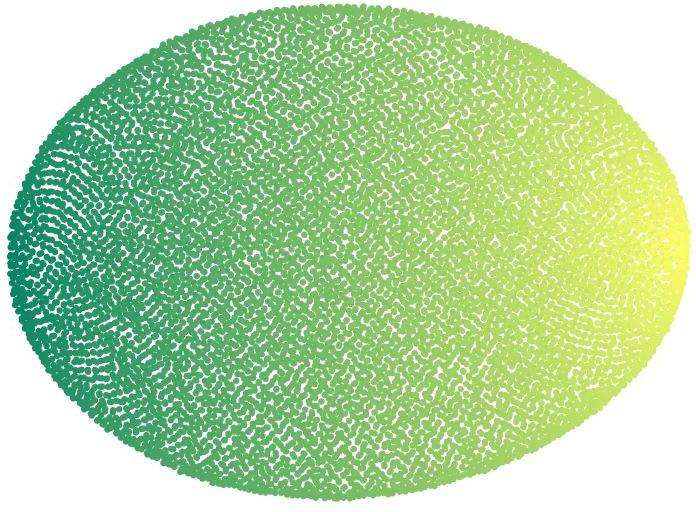} & \includegraphics[width=0.125\textwidth,height=0.25\textwidth,keepaspectratio]{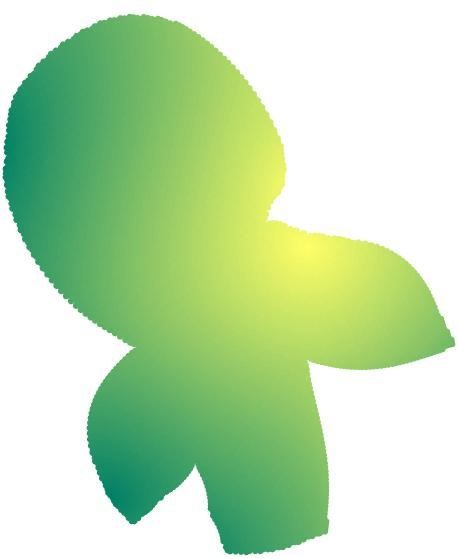} & \includegraphics[width=0.125\textwidth,height=0.25\textwidth,keepaspectratio]{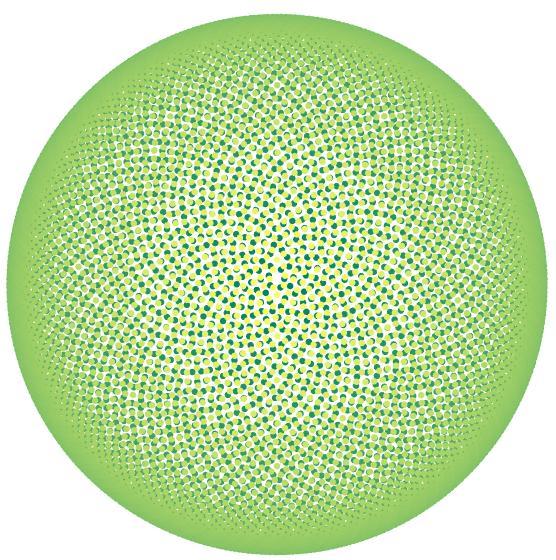} \\
        \hline
        \includegraphics[width=0.125\textwidth,height=0.25\textwidth,keepaspectratio]{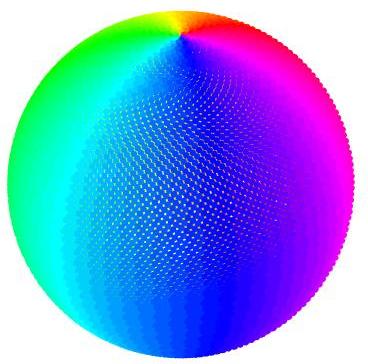} & \includegraphics[width=0.125\textwidth,height=0.25\textwidth,keepaspectratio]{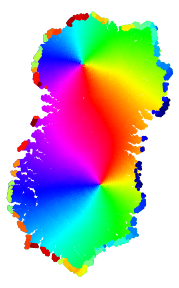} & \includegraphics[width=0.125\textwidth,height=0.25\textwidth,keepaspectratio]{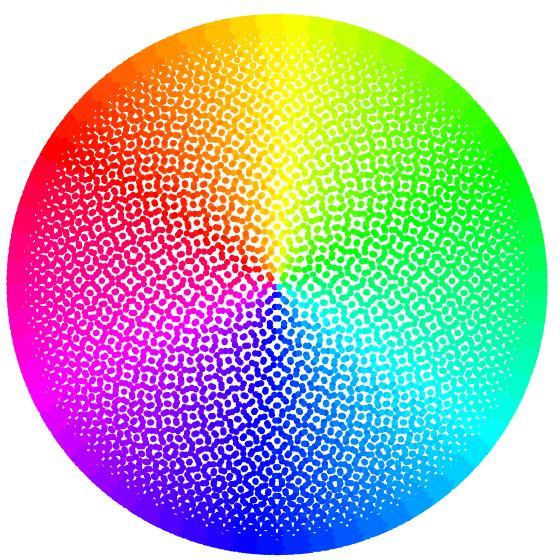} & \includegraphics[width=0.125\textwidth,height=0.25\textwidth,keepaspectratio]{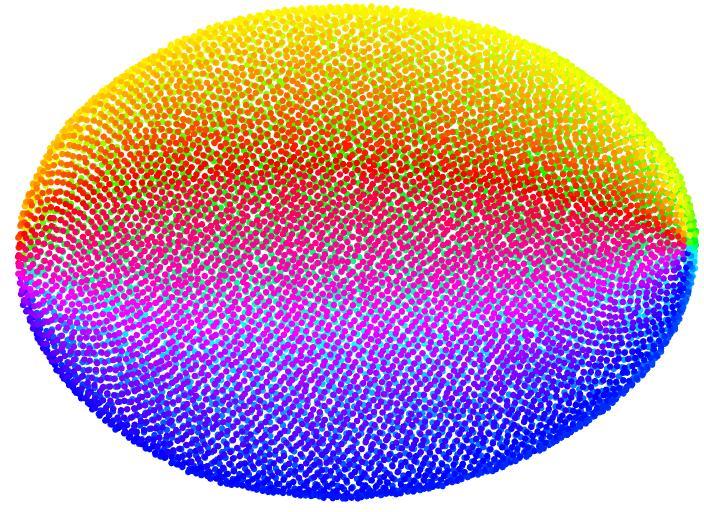} & \includegraphics[width=0.125\textwidth,height=0.25\textwidth,keepaspectratio]{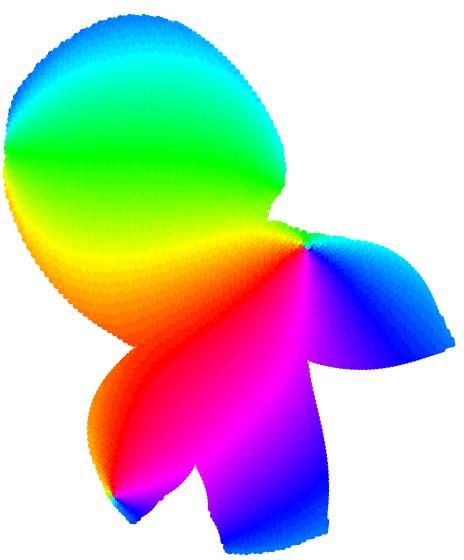} & \includegraphics[width=0.125\textwidth,height=0.25\textwidth,keepaspectratio]{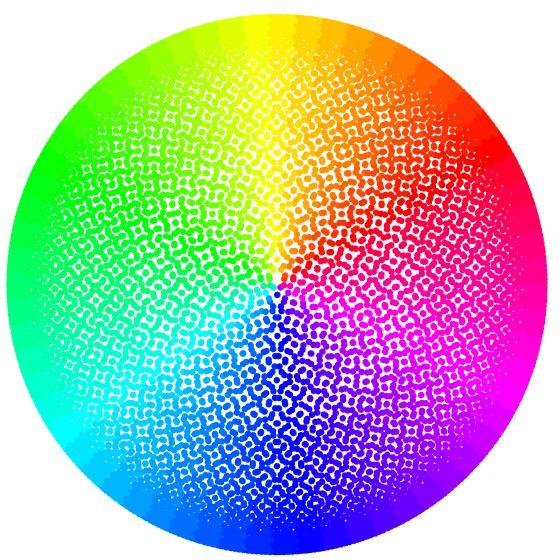} \\
        \hline
    \end{tabular}
    \caption{Embeddings of a sphere in $\mathbb{R}^3$ into $\mathbb{R}^2$. The top and bottom row contain the same plots colored by the height and the azimuthal angle of the sphere ($0-2\pi$), respectively. LDLE automatically colors the boundary so that the points on the boundary which are adjacent on the sphere have the same color. The arrows are manually drawn to help the reader identify the two pieces of the boundary which are to be stitched together to recover the original sphere. LTSA, UMAP and Laplacian eigenmaps squeezed the sphere into different viewpoints of $\mathbb{R}^2$ (side or top view of the sphere). t-SNE also tore apart the sphere but the embedding lacks interpretability as it is ``unaware" of the boundary.}
    \label{fig:fig14}
\end{figure}

\edit{A missing feature in existing manifold learning techniques} is their ability to embed closed manifolds into their intrinsic dimensions. For example, a sphere in $\mathbb{R}^3$ is a $2$-dimensional manifold which can be represented by a connected domain in $\mathbb{R}^2$ with boundary gluing instructions provided in the form of colors. We solve this issue in this paper (see Figure~\ref{fig:fig14}).

\section{Background and Motivation}
\label{sec:motivation}

\editt{Due to their global nature and robustness to noise, in our bottom-up approach for manifold learning,  we propose to construct low distortion (see Eq.~(\ref{DistortionPhikUk})) local mappings using low frequency Laplacian eigenfunctions. A natural way to achieve this is to restrict the eigenfunctions on local neighborhoods. Unfortunately, the common trend of using first $d$ non-trivial low frequency eigenfunctions to construct these local mappings fails to produce low distortion on all neighborhoods. This directly follows from the Laplacian Eigenmaps embedding of a high aspect-ratio rectangle shown in Figure~\ref{fig:fig13}. The following example explains that even in case of unit aspect-ratio, a local mapping based on the same set of eigenfunctions would not incur low distortion on each neighborhood, while mappings based on different sets of eigenfunctions may achieve that.}

\begin{figure}[h]
    \begin{tabular}{cc}
        \includegraphics[width=0.4\textwidth,height=0.4\textwidth,keepaspectratio]{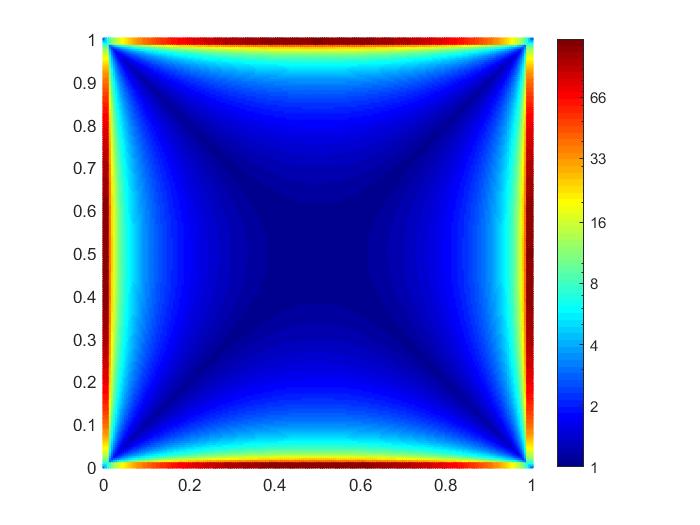}
     & \includegraphics[width=0.4\textwidth,height=0.4\textwidth,keepaspectratio]{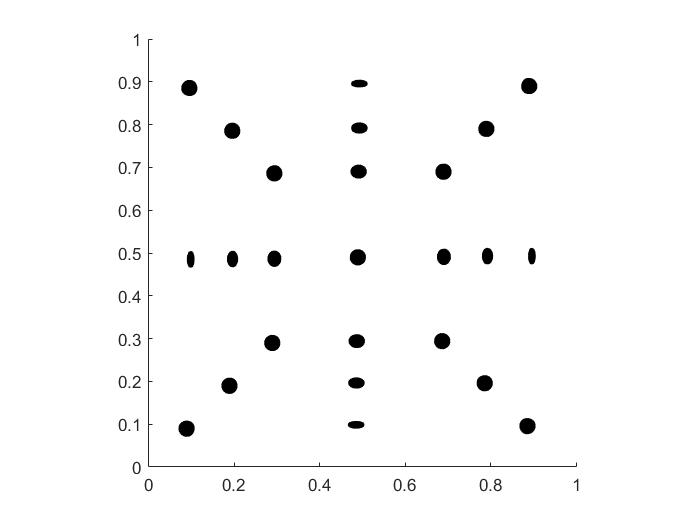} \\
         \includegraphics[width=0.4\textwidth,height=0.4\textwidth,keepaspectratio]{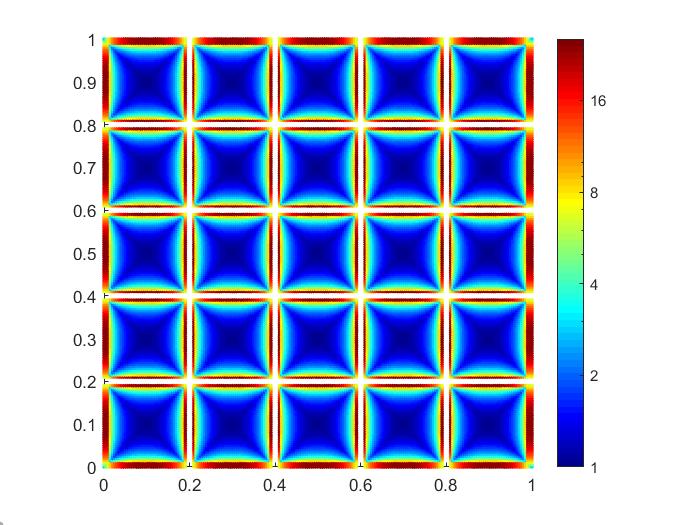}
     & \includegraphics[width=0.4\textwidth,height=0.4\textwidth,keepaspectratio]{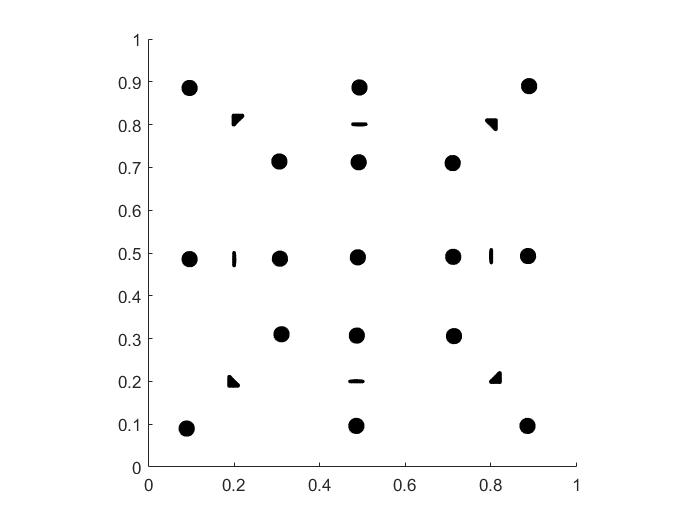}
    \end{tabular}
    \caption{(Left) Distortion of $\Phi_1^*$ (top) and $\Phi_2^*$ (bottom) on discs of radius $0.01$ centered at $(x,y)$ for all $x,y \in [0,1]\times[0,1]$. $\Phi_2^*$ produces close to infinite distortion on the discs located in the white region. (Right) Mapping of the discs at various locations in the square using $\Phi_1^*$ (top) and $\Phi_2^*$ (bottom).}
    \label{fig:fig11}
\end{figure}

Consider a unit square $[0,1]\times [0,1]$ such that for every point $x_k$ in the square, $\mathcal{U}_k$ is the disc of radius $0.01$ centered at $x_k$. Consider a mapping $\Phi_1^*$ based on the first two non-trivial eigenfunctions $\cos(\pi x)$ and  $\cos(\pi y)$ of the Laplace-Beltrami operator on the square with Neumann boundary conditions, that is,
\begin{align}
    \Phi_1^*(x,y) = (\cos(\pi x), \cos(\pi y)).
\end{align}

As shown in Figure~\ref{fig:fig11}, $\Phi_1^*$ maps the discs along the diagonals to other discs. The discs along the horizontal and vertical lines through the center are mapped to ellipses. The skewness of these ellipses increases as we move closer to the middle of the edges of the unit square. Thus, the distortion of $\Phi_1^*$ is low on the discs along the diagonals and high on the discs close to the middle of the edges of the square.



Now, consider a different mapping based on another set of eigenfunctions,
\begin{align}
    \Phi_2^*(x,y) &= (\cos(5\pi x), \cos(5\pi y)).
\end{align}
Compared to $\Phi_1^*$, $\Phi_2^*$ produces almost no distortion on the discs of radius $0.01$ centered at $(0.1,0.5)$ and $(0.9,0.5)$ (see Figure~\ref{fig:fig11}). Therefore, in order to achieve low distortion, it \editt{seem to} make sense to construct local mappings for different regions based on different sets of eigenfunctions. 


\editt{The following result from \citep{jones2007universal} manifests the above claim as it shows that, for a given small neighborhood on a Riemannian manifold, there always exist a subset of Laplacian eigenfunctions such that a local parameterization based on this subset is bilipschitz and has bounded distortion. A more precise statement follows.}

\begin{thm}[\citep{jones2007universal}, Theorem 2.2.1]
\label{thm:thm1}
Let $(\mathcal{M},g)$ be a $d$-dimensional Riemannian manifold. Let $\Delta_g$ be the Laplace-Beltrami operator on it with Dirichlet or Neumann boundary conditions \edit{and let $\phi_i$ be an eigenfunction of $\Delta_g$ with eigenvalue $\lambda_i$}. Assume that $|\mathcal{M}|=1$ where $|\mathcal{M}|$ is the volume of $\mathcal{M}$ and the uniform ellipticity conditions for $\Delta_g$ are satisfied. Let $x_k \in \mathcal{M}$ and $r_k$ be less than the injectivity radius at $x_k$ (the maximum radius where the the exponential map is a diffeomorphism). Then, there exists a constant $\kappa > 1$ which depends on $d$ and the metric tensor $g$ such that the following hold. Let $\rho \leq r_k$ and $B_k \equiv B_{\kappa^{-1}\rho}(x_k)$ where
\begin{align}
    B_{\epsilon}(x) &= \{y \in \mathcal{M} \ | \ d_g(x,y) < \epsilon\}.\label{Bepsx}
\end{align}
Then there exist $i_1,i_2,\ldots,i_d$ such that, if we let
\begin{align}
    \gamma_{ki} = \left(\frac{\int_{B_k} \phi_{i}^2(y)dy}{|B_k|}\right)^{-1/2} \label{gammai}
\end{align}
then the map
\begin{align}
    \Phi_k: B_k &\rightarrow \mathbb{R}^d\nonumber\\
    x &\rightarrow (\gamma_{ki_1}\phi_{i_1}(x),\ldots,\gamma_{ki_d}\phi_{i_d}(x)) \label{Psi2}
\end{align}
is bilipschitz such that for any $y_1,y_2 \in B_k$ it satisfies
\begin{align}
    \frac{\kappa^{-1}}{\rho}d_g(y_1,y_2) \leq \left\|\Phi_k(y_1)-\Phi_k(y_2)\right\| \leq \frac{\kappa}{\rho}d_g(y_1,y_2),
\end{align}
where the associated eigenvalues satisfy
\begin{align}
    \kappa^{-1}\rho^{-2} \leq \lambda_{i_1},\ldots,\lambda_{i_d} \leq \kappa \rho^{-2}, \label{lambdai}
\end{align}
and the distortion is bounded from above by $\kappa^2$ i.e.
\begin{align}
    \sup_{\substack{y_1,y_2 \in B_k\\y_1\neq y_2}}\frac{\left\|\Phi_k(y_1)-\Phi_k(y_2)\right\|}{d_g(y_1,y_2)}\sup_{\substack{y_1,y_2 \in B_k\\y_1\neq y_2}}\frac{d_g(y_1,y_2)}{\left\|\Phi_k(y_1)-\Phi_k(y_2)\right\|} \leq \frac{\kappa}{\rho}\frac{\rho}{\kappa^{-1}} = \kappa^2. \label{distortioneq}
\end{align}
\end{thm}


\editt{Motivated by the above result, we adopt the form of local paramterizations $\Phi_k$ in Eq.~(\ref{Psi2}) as local mappings in our work. The main challenge then is to identify the set of eigenfunctions for a given neighborhood such that the resulting parameterization produces low distortion on it. The existence proof of the above theorem by the authors of \citep{jones2007universal} suggests a procedure to identify this set in the continuous setting. Below, we provide a sketch of their procedure and in Section~\ref{sec:biliplocparam} we describe our discrete realization of it.}

\subsection{Eigenfunction Selection in the Continuous Setting}
\label{subsubsec:jones_method}

\edit{Before describing the procedure used in \citep{jones2007universal} to choose the eigenfunctions, we first provide some intuition about the desired properties for the chosen eigenfunctions $\phi_{i_1}, \ldots, \phi_{i_d}$ so that the resulting parameterization $\Phi_k$ has low distortion on $B_k$.}


\edit{Consider the simple case of $B_k$ representing a small open ball of radius $\kappa^{-1}\rho$ around $x_k$ in $\mathbb{R}^d$ equipped with the standard Euclidean metric. Then the first-order Taylor approximation of $\Phi_k(x)$, $x \in B_k$, about $x_k$ is given by}
\begin{align}
    \Phi_{k}(x) &\approx \Phi_k(x_k) + J(x-x_k) \text{ where } J = [\gamma_{ki_1}\nabla\phi_{i_1}(x_k) \ldots \gamma_{ki_d}\nabla\phi_{i_d}(x_k)]^T.
\end{align}

\edit{Note that $\gamma_{ki_s}$ are positive scalars constant with respect to $x$. Now, $\text{Distortion}(\Phi_k, B_k) = 1$ if and only if $\Phi_k$ preserves distances between points in $B_k$ up to a constant scale (see Eq.~(\ref{DistortionPhikUk})). That is,}
\begin{align}
    \left\|\Phi_k(x) - \Phi_k(y)\right\|_2 = c  \left\|x-y\right\|_2\ \forall x,y \in B_k \text{ and for some constant }c > 0.
\end{align}

\edit{Using the first-order approximation of $\Phi_k$ we get,}
\begin{align}
    \left\|J(x-y)\right\|_2 \approx c \left\|x-y\right\|_2\ \forall x,y \in B_k \text{ and for some constant }c > 0.
\end{align}
\edit{Therefore, for low distortion $\Phi_k$, $J$ must approximately behave like a similarity transformation and therefore, $J$ needs to be approximately orthogonal up to a constant scale. In other words, the chosen eigenfunctions should be such that $\gamma_{ki_1}\nabla\phi_{i_1}(x_k),$ $\ldots,$ $\gamma_{ki_d}\nabla\phi_{i_d}(x_k)$ are close to being orthogonal and have similar lengths.}


\edit{The same intuition holds in the manifold setting too. The construction procedure described in~\citep{jones2007universal} aims to choose eigenfunctions such that}
\begin{enumerate}[label=(\alph*)]
    \item they are close to being locally orthogonal, that is, $\nabla\phi_{i_1}(x_k), \ldots, \nabla\phi_{i_d}(x_k)$ are approximately orthogonal, and
    \item that their local scaling factors $\gamma_{ki_s}\left\|\nabla\phi_{i_s}(x_k)\right\|_2$ are close to each other.
\end{enumerate}


\edit{\textit{Note}. Throughout this paper, we use the convention $\nabla\phi_{i}(x_k) = \nabla(\phi_{i}\ \circ\ \exp_{x_k})(0)$ where $\exp_{x_k}$ is the exponential map at $x_k$. Therefore, $\nabla\phi_{i}(x_k)$ can be represented by a $d$-dimensional vector in a given $d$-dimensional orthonormal basis of $T_{x_k}\mathcal{M}$. Even though the representation of these vectors depend on the choice of the orthonormal basis, the value of the canonical inner product between these vectors, and therefore the $2$-norm of the vectors, are the same across different basis. This follows from the fact that an orthogonal transformation preserves the inner product.}

\begin{rmk}
\label{rmk:ltsa}
\editt{
Based on the above first order approximation, one may take our local mappings $\Phi_k$ to also be projections onto the tangential spaces. However, unlike LTSA \citep{zhang2003nonlinear} where the basis of the tangential space is estimated by the local principal directions, in our case it is estimated by the locally orthogonal gradients of the global eigenfunctions of the Laplacian. Therefore, LTSA relies only on the local structure to estimate the tangential space while, in a sense, our method makes use of both local and global structure of the manifold. 
}
\end{rmk}


\edit{A high level overview of the procedure presented in ~\citep{jones2007universal} to choose eigenfunctions which satisfy the properties in (a) and (b) follows.
\begin{enumerate}
    \item A set $S_k$ of the indices of candidate eigenfunctions is chosen such that $i \in S_k$ if the length of $\gamma_{ki}\nabla\phi_i(x_k)$ is bounded from above by a constant, say $C$.
    \item A direction $p_1 \in T_{x_k}\mathcal{M}$ is selected at random.
    \item Subsequently $i_1 \in S_k$ is selected so that $\gamma_{ki_1}|\nabla\phi_{i_1}(x_k)^Tp_1|$ is sufficiently large. This motivates $\gamma_{ki_1}\nabla\phi_{i_1}(x_k)$ to be approximately in the same direction as $p_1$ and the length of it to be close to the upper bound $C$.
    \item Then, a recursive strategy follows. To find the $s$-th eigenfunction for $s \in \{2,\ldots,d\}$, a direction $p_s \in T_{x_k}\mathcal{M}$ is chosen such that it is orthogonal to $\nabla \phi_{i_1}(x_k), \ldots, \nabla\phi_{i_{s-1}}(x_k)$.
    \item Subsequently, $i_s \in S_k$ is chosen so that $\gamma_{ki_s}|\nabla\phi_{i_s}(x_k)^Tp_s|$ is sufficiently large. Again, this motivates $\gamma_{ki_s}\nabla\phi_{i_s}(x_k)$ to be approximately in the same direction as $p_s$ and the length of it to be close to the upper bound $C$.
\end{enumerate}
Since $p_s$ is orthogonal to $\nabla \phi_{i_1}(x_k), \ldots, \nabla\phi_{i_{s-1}}(x_k)$ and the direction of $\gamma_{ki_s}\nabla\phi_{i_s}$ is approximately the same as $p_s$, therefore $(a)$ is satisfied. Since for all $s \in \{1,\ldots,d\}$, $\gamma_{ki_s}\nabla\phi_{i_s}(x_k)$ has a length close to the upper bound $C$, therefore $(b)$ is also satisfied. The core of their work lies in proving that these $\phi_{i_1},\ldots,\phi_{i_d}$ always exist under the assumptions of the theorem such that the resulting parameterization $\Phi_k$ has bounded distortion (see Eq. (\ref{distortioneq})). This bound depends on the intrinsic dimension $d$ and the natural geometric properties of the manifold. The main challenge in practically realizing the above procedure lies in the estimation of $\nabla\phi_{i_s}(x_k)^Tp_s$. In Section~\ref{sec:biliplocparam}, we overcome this challenge.}

\section{Low-dimensional Low Distortion Local Parameterization}
\label{sec:biliplocparam}

\edit{In the procedure to choose $\phi_{i_1},\ldots,\phi_{i_d}$ to construct $\Phi_k$} as described above, the selection of the first eigenfunction $\phi_{i_1}$ relies on the derivative of the eigenfunctions at $x_k$ along an arbitrary direction \edit{$p_1 \in T_{x_k}\mathcal{M}$}, that is, on $\nabla\phi_i(x_k)^Tp_1$. In our algorithmic realization of the construction procedure, we take $p_1$ to be the gradient of an eigenfunction at $x_k$ itself (say $\nabla\phi_j(x_k)$). 
We relax the unit norm constraint on $p_1$; note that this will neither affect the math nor the output of our algorithm. Then the selection of $\phi_{i_1}$ would depend on the inner products $\nabla\phi_i(x_k)^T\nabla\phi_j(x_k)$. \edit{The value of this inner product does not depend on the choice of the orthonormal basis for $T_{x_k}\mathcal{M}$}. \edit{We discuss several ways to obtain a numerical estimate} of this inner product by making use of the local correlation between the eigenfunctions \citep{steinerberger2017spectral, doi:10.1080/10586458.2018.1538911}. These estimates are used to select the subsequent eigenfunctions too.


\edit{In Section~\ref{subsec:scaledloccorr}, we first review the local correlation between the eigenfunctions of the Laplacian. In Theorem~\ref{thm:thm2} we show that the limiting value of the scaled local correlation between two eigenfunctions equals the inner product of their gradients. We provide two proofs of the theorem where each proof leads to a numerical procedure described in Section~\ref{subsec:Atilde_esimate}, followed by examples to empirically compare the estimates. Finally, in Section~\ref{subsec:biliplocparam}, we use these estimates to obtain low distortion local parameterizations of the underlying manifold.}

\subsection{Inner Product of Eigenfunction Gradients using Local Correlation}
\label{subsec:scaledloccorr}

Let $(\mathcal{M},g)$ be a $d$-dimensional Riemannian manifold with or without boundary, rescaled so that $|\mathcal{M}|\leq 1$. Denote the volume element at $y$ by $\omega_g(y)$. Let $\phi_i$ and $\phi_j$ be the eigenfunctions of the Laplacian operator $\Delta_g$ (see statement of Theorem~\ref{thm:thm1}) with eigenvalues $\lambda_i$ and $\lambda_j$. Let $x_k \in \mathcal{M}$ and define
\begin{align}
    \Psi_{kij}(y)=(\phi_i(y)-\phi_i(x_k))(\phi_j(y)-\phi_j(x_k)).\label{Phikij}
\end{align}
Then the local correlation between the two eigenfunctions $\phi_i$ and $\phi_j$ at the point $x_k$ at scale $t_k^{-1/2}$ as defined in \citep{steinerberger2017spectral, doi:10.1080/10586458.2018.1538911} is given by
\begin{align}
    A_{kij} = \int_{\mathcal{M}}p(t_k,x_k,y)\Psi_{kij}(y)\omega_g(y), \label{Akij}
\end{align}
where $p(t,x,y)$ is the fundamental solution of the heat equation on $(\mathcal{M},g)$. As noted in \citep{steinerberger2017spectral}, for $(t_k,x_k) \in \mathbb{R}_{\geq 0} \times \mathcal{M}$ fixed, we have
\begin{align}
    p(t_k,x_k,y) \sim \left\{\begin{matrix*}[l] t_k^{-d/2} & d_g(x_k,y) \leq t_k^{-1/2}\\0 & \text{otherwise}\end{matrix*}\right. \qquad \text{and} \qquad \int_{M}p(t_k,x_k,y)\omega_g(y) = 1.
\end{align}
Therefore, $p(t_k,x_k,\cdot)$ acts as a local probability measure centered at $x_k$ with scale $t_k^{-1/2}$ (see Eq.~(\ref{ptxy}) in Appendix~\ref{proof:proof21} for a precise form of $p$). We define the scaled local correlation to be the ratio of the local correlation $A_{kij}$ and a factor of $2t_k$.
\begin{thm}
\label{thm:thm2}
Denote the limiting value of the scaled local correlation by $\widetilde{A}_{kij}$,
\begin{align}
    \widetilde{A}_{kij} = \lim_{t_k\rightarrow 0}\frac{A_{kij}}{2t_k} \label{AtildeLimit}
\end{align}
Then $\widetilde{A}_{kij}$ equals the inner product of the gradients of the eigenfunctions $\phi_i$ and $\phi_j$ at $x_k$, that is,
\begin{align}
    \widetilde{A}_{kij}=\nabla\phi_i(x_k)^T\nabla\phi_j(x_k). \label{Atilde}
\end{align}
\end{thm}

\edit{Two proofs are provided in Appendix~\ref{proof:proof21} and~\ref{proof:proof22}. A brief summary is provided below.}


\textit{Proof 1}. In the first proof we choose a sufficiently small $\epsilon_k$ and show that
\begin{align}
    \lim_{t_k\rightarrow0}A_{kij} &= \lim_{t_k\rightarrow 0}\int_{B_{\epsilon_k}(x_k)}G(t_k,x_k,y)\Psi_{kij}(y)\omega_g(y) \label{Akijlim}
\end{align}
where $B_{\epsilon}(x)$ is defined in Eq.~(\ref{Bepsx}) and
\begin{align}
    G(t,x,y) &= \frac{e^{-d_g(x,y)^2/4t}}{(4\pi t)^{d/2}}. \label{Gtxy1}
\end{align}
Then, by using the properties of the exponential map at $x_k$ and applying basic techniques in calculus, we show that $\lim_{t_k\rightarrow 0}A_{kij}/2t_k$ evaluates to $\nabla\phi_i(x_k)^T\nabla\phi_j(x_k)$.


\edit{
\textit{Proof 2}. In the second proof, as in \citep{doi:10.1080/03605302.2014.942739, steinerberger2017spectral},  we used the Feynman-Kac formula,
\begin{align}
    A_{kij} &= [e^{-t_{k}\Delta_g}((\phi_i - \phi_i(x_k))(\phi_j - \phi_j(x_k))](x_k)
\end{align}
and note that
\begin{align}
    \lim_{t_k\rightarrow 0}\frac{A_{kij}}{2t_k} = \left.\frac{1}{2}\frac{\partial A_{kij}}{\partial t_k} \right\vert_{t_k=0} = \frac{-1}{2}\left\{\Delta_g[(\phi_i - \phi_i(x_k)) (\phi_j-\phi_j(x_k))](x_k)\right\}.\label{Atilde_method2}
\end{align}
Then, by applying the formula of the Laplacian of the product of two functions, we show that the above equation equals $\nabla\phi_i(x_k)^T\nabla\phi_j(x_k)$.
}

\subsection{Estimate of $\widetilde{A}_{kij}$ in the Discrete Setting}
\label{subsec:Atilde_esimate}

\edit{To apply Theorems~\ref{thm:thm1} and~\ref{thm:thm2} in practice on data, we need an estimate of $\widetilde{A}_{kij}$ in the discrete setting. There are several ways to obtain this estimate. A generic way is by using the algorithms \citep{doi:10.1080/01621459.2013.827984,aswani2011regression} based on Local Linear Regression (LLR) to estimate the gradient vector $\nabla\phi_i(x_k)$ itself from the values of $\phi_i$ in a neighbor of $x_k$. An alternative approach is to use a finite sum approximation of Eq.~(\ref{Akijlim}) combined with Eq.~(\ref{AtildeLimit}). A third approach is based on the Feynman-Kac formula where we make use of Eq.~(\ref{Atilde_method2}) in the discrete setting. In the following we explain the latter two approaches.}

\subsubsection{Finite sum approximation}
\label{subsubsec:finite}
Let $(x_k)_{k=1}^{n}$ be uniformly distributed points on $(\mathcal{M},g)$. Let $d_e(x_k,x_k')$ be the distance between $x_k$ and $x_{k'}$. The accuracy with which $\widetilde{A}_{kij}$ can be estimated mainly depends on the accuracy of $d_e(\cdot\ ,\ \cdot)$ to the local geodesic distances. For simplicity, we use $d_e(x_k,x_k')$ to be the Euclidean distance $\left\|x_k-x_{k'}\right\|_2$. A more accurate estimate of the local geodesic distances can be computed using the method described in \citep{li2019geodesic}. 


We construct a sparse unnormalized graph Laplacian $L$ using Algo.~\ref{algo:gl}, where the weight matrix $K$ of the graph edges is defined using the Gaussian kernel. The bandwidth of the Gaussian kernel is set using the local scale of the neighborhoods around each point as in self-tuning spectral clustering~\citep{zelnik2005self}. Let $\bm{\phi}_i$ be the $i$th non-trivial eigenvector of $L$ and denote $\phi_i(x_j)$ by $\bm{\phi}_{ij}$.


\begin{algorithm}[H]
\SetAlgoLined
\KwIn{$d_e(x_k,x_{k'})_{k,k'=1}^{n},k_{\textrm{nn}},k_{\textrm{tune}} \text{ where } k_{\textrm{tune}}\leq k_{\textrm{nn}}$}
\KwOut{$L$}
{
    $\mathcal{N}_k \gets $ set of indices of $k_{\textrm{nn}}$ nearest neighbours of $x_k$ based on $d_e(x_k,\cdot)$\;
    $\sigma_{k} \gets d_e(x_k,x_{k^*})$ where $x_{k^*}$ is the $k_{\textrm{tune}}$th nearest neighbor of $x_k$\;
    $K_{kk} \gets 0, K_{kk'} \gets e^{-d_e(x_k,x_{k'})^2/\sigma_{k}\sigma_{k'}}, k' \in \mathcal{N}_k$\;
    $D_{kk} \gets \sum_{k'}K_{kk'}, \ D_{kk'} \gets 0, k \neq k'$\;
    $L \gets D-K$\;
}
\caption{Sparse Unnormalized Graph Laplacian based on \citep{zelnik2005self}}
\label{algo:gl}
\end{algorithm}


We estimate $\widetilde{A}_{kij}$ by evaluating the scaled local correlation $A_{kij}/2t_k$ at a small value of $t_k$. The limiting value of $A_{kij}$ is estimated by substituting a small $t_k$ in the finite sum approximation of the integral in Eq.~(\ref{Akijlim}). The sum is taken on a discrete ball of a small radius $\epsilon_k$ around $x_k$ and is divided by $2t_k$ to obtain an estimate of $\widetilde{A}_{kij}$.


We start by choosing $\epsilon_k$ to be the distance of $k_{\text{lv}}$th nearest neighbor of $x_k$ where $k_{\text{lv}}$ is a hyperparameter with a small integral value (subscript $\text{lv}$ stands for local view).
Thus,
\begin{align}
    \epsilon_k = \text{distance to the } k_{\text{lv}}\text{th nearest neighbor of } x_k. \label{epsk}
\end{align}
Then the limiting value of $t_k$ is given by
\begin{align}
    \sqrt{\text{chi2inv}(p,d)}\sqrt{2 t_k} = \epsilon_k \implies t_k = \frac{1}{2}\frac{\epsilon_k^2}{\text{chi2inv}(p,d)}, \label{tk}
\end{align}
where $\text{chi2inv}$ is the inverse cdf of the chi-squared distribution with $d$ degrees of freedom evaluated at $p$. We take $p$ to be $0.99$ in our experiments. The rationale behind the above choice of $t_k$ is described in Appendix~\ref{rationaletk}.


Now define the discrete ball around $x_k$ as
\begin{align}
    U_k &= \{x_{k'} \ | \ d_e(x_k,x_{k'}) \leq \epsilon_k\}. \label{Uk}
\end{align}

Let $U_k$ denote the $k$th local view of the data in the high dimensional ambient space. For convenience, denote the estimate of $G(t_k,x_k,x_{k'})$ by $G_{kk'}$ where $G$ is as in Eq.~(\ref{Gtxy1}). Then
\begin{align}
    G_{kk'} &= \left\{\begin{matrix*}[l] \frac{\exp(-d_e(x_k,x_{k'})^2/4t_k)}{\sum_{x\in U_k}\exp(-d_e(x_k,x)^2/4t_k)}&,\  x_{k'} \in U_k - \{x_k\}\\
    0&,\ \text{otherwise}. \end{matrix*}\right. \label{Gks}
\end{align}
\edit{Finally, the estimate of $\widetilde{A}_{kij}$ is given by
\begin{align}
    \widetilde{A}_{kij} &= \frac{1}{2t_k} G_k^T((\bm{\phi_i}-\bm{\phi_{ik}})\odot (\bm{\phi_j}-\bm{\phi_{jk}})) \label{Atildek}
\end{align}
where $G_k$ is a column vector containing the $k$th row of the matrix $G$ and $\odot$ represents the Hadamard product.}

\subsubsection{\edit{Estimation based on Feynman-Kac formula}}
This approach to estimate $\widetilde{A}_{kij}$ is simply the discrete analog of Eq.~(\ref{Atilde_method2}),
\begin{align}
    \widetilde{A}_{kij} = \frac{-1}{2}L_k^T((\bm{\phi_i}-\bm{\phi_{ik}})\odot (\bm{\phi_j}-\bm{\phi_{jk}}))\label{Atildek_method2_discrete}
\end{align}
where $L_k$ is a column vector containing the $k$th row of $L$. A variant of this approach which results in better estimates in the noisy case uses a low rank approximation of $L$ using its first few eigenvectors (see Appendix~\ref{sec:supp}).
\begin{rmk}
It is not a coincidence that Eq.~(\ref{Atildek}) and Eq.~(\ref{Atildek_method2_discrete}) look quite similar. In fact, if we take $T$ to be a diagonal matrix with $(t_k)_{k=1}^{n}$ as the diagonal, then the matrix $T^{-1}(I-G)$ approximates $\Delta_g$ in the limit of $(t_k)_{k=1}^{n}$ tending to zero. Replacing $L$ with $T^{-1}(I-G)$ and therefore $L_k$ with $(e_k - G_k)/t_k$ reduces Eq.~(\ref{Atildek_method2_discrete}) to Eq.~(\ref{Atildek}). Here $e_k$ is a column vector with $k$th entry as $1$ and rest zeros. Therefore the two approaches are the same in the limit.
\end{rmk}
\begin{rmk}
The above two approaches can also be generalized to compute the $\nabla f_i(x_k)^T\nabla f_j(x_k)$ for arbitrary $\mathcal{C}^2$ mappings $f_i$ and $f_j$ from $\mathcal{M}$ to $\mathbb{R}$ ( $\nabla f_i(x_k) = \nabla(f_i \ \circ\ \text{exp}_{x_k})(0)$ as per our convention). To achieve this, simply replace $\phi_i$ and $\phi_j$ with $f_i$ and $f_j$ in Eq.~(\ref{Atildek}) and Eq.~(\ref{Atildek_method2_discrete}).
\end{rmk}


\begin{figure}[ht]
    \centering
    \begin{tabular}{cc}
        \begin{tabular}{c} \includegraphics[height=0.18\textwidth,keepaspectratio]{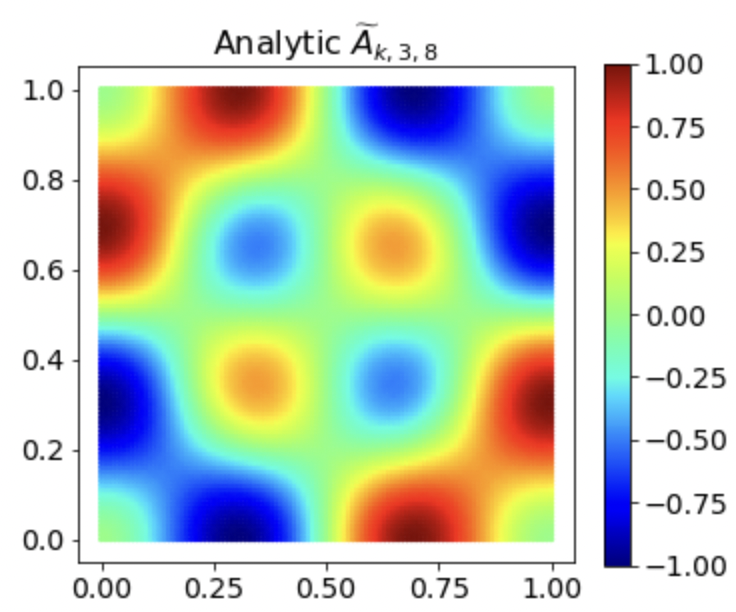}\end{tabular} & \begin{tabular}{c} \includegraphics[height=0.18\textwidth,keepaspectratio]{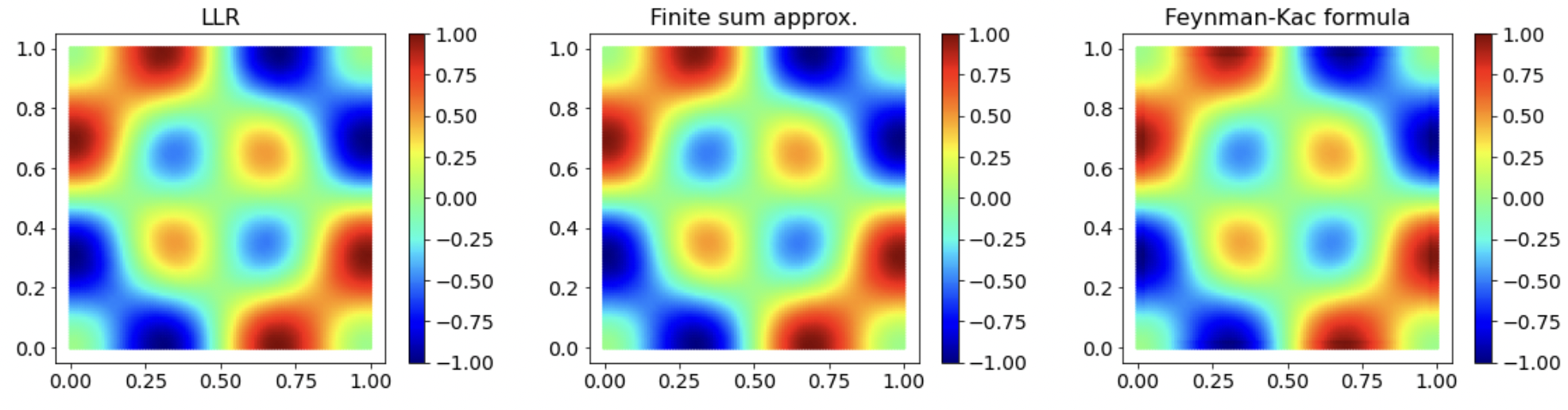}\\
        \includegraphics[height=0.185\textwidth,keepaspectratio]{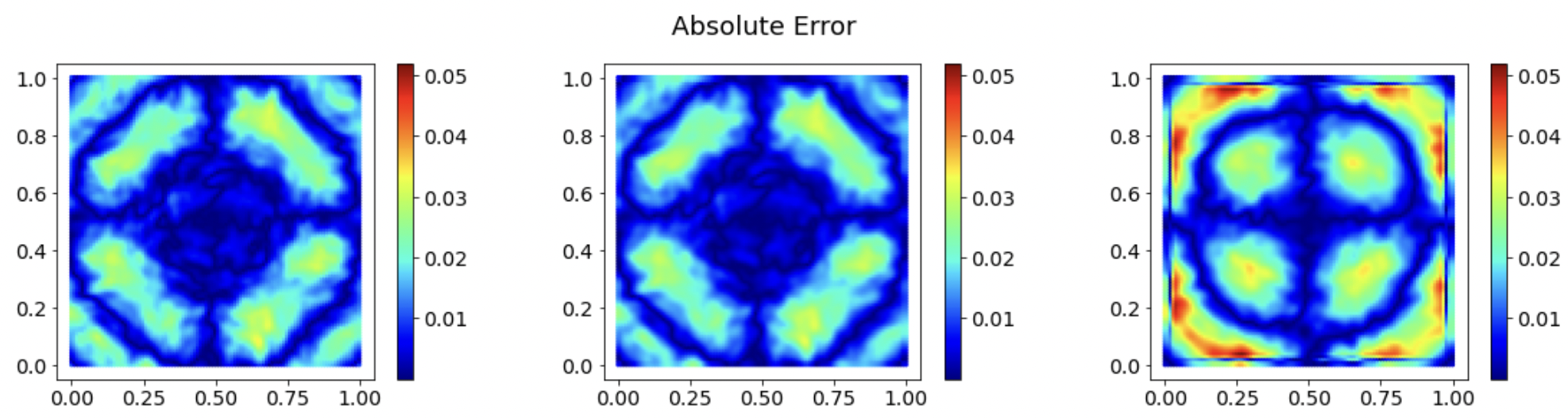}
        \end{tabular}
    \end{tabular}
    \caption{\editt{Comparison of different approaches to estimate $\widetilde{A}_{kij}$ in the discrete setting.}}
    \label{fig:fig21}
\end{figure}

\label{sec:Atilde_compare}
\paragraph{Example.} This example will follow us throughout the paper. Consider a square grid $[0,1]\times [0,1]$ with a spacing of $0.01$ in both $x$ and $y$ direction. With $k_{\textrm{nn}} = 49$, $k_{\textrm{tune}} = 7$ and $d_e(x_k,x_{k'})=\left\|x_k-x_{k'}\right\|_2$ as input to the Algo.~\ref{algo:gl}, we construct the graph Laplacian $L$. Using $k_{\text{lv}} = 25$, $d=2$ and $p=0.99$, we obtain the discrete balls $U_k$ and $t_k$. The $3$rd and $8$th eigenvectors of $L$ 
and the corresponding analytical eigenfunctions are then obtained. The analytical value of $\widetilde{A}_{k38}$ is displayed in Figure~\ref{fig:fig21}, followed by its estimate using LLR \citep{doi:10.1080/01621459.2013.827984}, finite sum approximation and Feynman-Kac formula based approaches. The analytical and the estimated values are normalized by $\max_{k}\widetilde{A}_{kij}$ to bring them to the same scale. The absolute error due to these approaches are shown below the estimates.


Even though, in this example, the Feynman-Kac formulation seem to have a larger error, in our experiments, no single approach seem to be a clear winner across all the examples. This becomes clear in Appendix~\ref{sec:supp} where we provided a comparison of these approaches on a noiseless and a noisy Swiss Roll. The results shown in this paper are based on finite sum approximation to estimate $\widetilde{A}_{kij}$. 

\subsection{\edit{Low Distortion Local Parameterization from Laplacian Eigenvectors}}
\label{subsec:biliplocparam}
We use $\nabla\phi_i \equiv \nabla\phi_i(x_k)$ for brevity. Using the estimates of $\widetilde{A}_{kij}$, we now present an algorithmic construction of low distortion local parameterization $\Phi_{k}$ which maps $U_k$ into $\mathbb{R}^d$. The pseudocode is provided below followed by a full explanation of the steps and a note on the hyperparameters. Before moving forward, it would be helpful for the reader to review the construction procedure in the continuous setting in Section~\ref{subsubsec:jones_method}.

\begin{algorithm}[ht]
\SetAlgoLined
\KwIn{$L, N, k_{\text{lv}}, d, p, (\tau_s,\delta_s)_{s=1}^{d}$}
\KwOut{$(\Phi_k, U_k, \zeta_{kk})_{k=1}^{n}$}
{
    Compute $(\bm{\phi}_i)_{i=1}^{N},\lambda_1\leq\ldots\leq\lambda_N$ by eigendecomposition of $L$\;
    \For{$k\gets1$ \KwTo $n$}{
        Compute $U_k,(\widetilde{A}_{kij})_{i,j=1}^{N}$ (Eq.~(\ref{Uk},~\ref{Atildek}))\;
        Compute $(\gamma_{ki})_{i=1}^{N}$ (Eq.~(\ref{gammaki_discrete}))\;
        $\theta_1 \gets \tau_1$-percentile of $(\widetilde{A}_{kii})_{i=1}^{N}$\;
        Compute $S_k$ (Eq.~(\ref{Sktilde}))\;
        Compute $i_1$ (Eq.~(\ref{i1}))\;
        \For{$s\gets2$ \KwTo $d$}{
            Compute $H^{s}_{kij}$ (Eq.~(\ref{Hdef}))\;
            $\theta_s \gets \tau_s$-percentile of $(H^s_{kii})_{i \in S_k}$\;
            Compute $i_s$ (Eq.~(\ref{is}))\;
        }
        $\Phi_k \gets (\gamma_{ki_{1}}\bm{\phi}_{i_{1}},\ldots,\gamma_{ki_{d}}\bm{\phi}_{i_{d}})$ (Eq.~(\ref{Psik}))\;
        Compute $\zeta_{kk}$ (Eq.~(\ref{zetakj}))\;
    }
}
\caption{BiLipschitz-Local-Parameterization}
\label{algo:blp}
\end{algorithm}

An estimate of $\gamma_{ki}$ is obtained by the discrete analog of Eq.~(\ref{gammai}) and is given by
\begin{align}
    \gamma_{ki} = \text{Root-Mean-Square}(\{\bm{\phi}_{ij}\ |\ x_j\in U_k\})^{-1}. \label{gammaki_discrete}
\end{align}
\paragraph{Step 1. Compute a set $S_k$ of candidate eigenvectors for $\Phi_k$.}
Based on the construction procedure following Theorem~\ref{thm:thm1}, we start by computing a set $S_k$ of candidate eigenvectors to construct $\Phi_k$ of $U_k$. There is no easy way to retrieve the set $S_k$ in the discrete setting as in the procedure. Therefore, we make the natural choice of using the first $N$ nontrivial eigenvectors $(\bm{\phi}_i)_{i=1}^{N}$ of $L$ corresponding to the $N$ smallest eigenvalues $(\lambda_i)_{i=1}^N$, with sufficiently large gradient at $x_k$, as the set $S_k$. The large gradient constraint is required  for the numerical stability of our algorithm. Therefore, we set $S_k$ to be,
\begin{align}
    S_k &= \{i \in \{1,\ldots,N\} \ | \ \left\|\nabla\phi_{i}\right\|^2\geq\theta_1 \} = \{i \in \{1,\ldots,N\}| \ \widetilde{A}_{kii} \geq\theta_1 \}, \label{Sktilde}
\end{align}
where $\theta_1$ is $\tau_1$-percentile of the set $(\widetilde{A}_{kii})_{i=1}^{N}$ and the second equality follows from Eq.~(\ref{Atilde}). Here $N$ and $\tau_1 \in (0,100)$ are hyperparameters.


\paragraph{Step 2. Choose a direction $p_1 \in T_{x_k}\mathcal{M}$.}
The unit norm constraint on $p_1$ is relaxed. This will neither affect the math nor the output of our algorithm. Since $p_1$ can be arbitrary we take $p_1$ to be the gradient of an eigenvector $r_1$, that is $\nabla\phi_{r_1}$. The choice of $r_1$ will determine $\bm{\phi}_{i_1}$. To obtain a low frequency eigenvector, $r_1$ is chosen so that the eigenvalue $\lambda_{r_1}$ is minimal, therefore
\begin{align}
    r_1 &= \argmin_{j\in S_k}\lambda_j. \label{r1}
\end{align}

\paragraph{Step 3. Find $i_1 \in S_k$ such that $\gamma_{ki_1} |\nabla\phi_{i_1}^Tp_1|$ is sufficiently large.} Since $p_1 = \nabla\phi_{r_1}$, using Eq.~(\ref{Atilde}), the formula for $\nabla\phi_i^Tp_1$ becomes
\begin{align}
    \nabla\phi_i^Tp_1 &= \nabla\phi_i^T\nabla\phi_{r_1} = \widetilde{A}_{kir_1}. \label{i_1_formula}
\end{align}
Then we obtain the eigenvector $\bm{\phi}_{i_1}$ so that $\gamma_{ki_1} |\nabla\phi_{i_1}^Tp_1|$ is larger than a certain threshold. We do not know what the value of this threshold would be in the discrete setting. Therefore, we first define the maximum possible value of $\gamma_{ki_1} |\nabla\phi_{i}^Tp_1|$ using Eq.~(\ref{i_1_formula}) as
\begin{align}
    \alpha_1 =\underset{i\in S_k}{\max}\ \gamma_{ki} |\nabla\phi_{i}^Tp_1| = \underset{i\in S_k}{\max}\ \gamma_{ki} |\widetilde{A}_{kir_1}|.
    \label{alpha1}
\end{align}
Then we take the threshold to be $\delta_1\alpha_1$ where $\delta_1 \in (0,1]$ is a hyperparameter. Finally, to obtain a low frequency eigenvector $\bm{\phi}_{i_1}$, we choose $i_1$ 
such that 
\begin{align}
    i_1 &= \argmin_{i \in S_k}\{\lambda_i : \gamma_{ki}|\nabla\phi_{i}^Tp_1| \geq \delta_1\alpha_1\} = \argmin_{i \in S_k}\{\lambda_i : \gamma_{ki} |\widetilde{A}_{kir_1}| \geq \delta_1\alpha_1\}.\label{i1}
\end{align}


After obtaining $\bm{\phi}_{i_1}$, we use a recursive procedure to obtain the $s$-th eigenvector $\bm{\phi}_{i_s}$ where $s \in \{2, \ldots, d\}$ in order.


\paragraph{Step 4. Choose a direction $p_{s} \in T_{x_k}\mathcal{M}$ orthogonal to $\nabla\phi_{i_1},\ldots,\nabla\phi_{i_{s}}$.} Again the unit norm constraint will be relaxed with no change in the output. We are going to take $p_s$ to be the component of $\nabla\phi_{r_s}$ orthogonal to $\nabla\phi_{i_1},\ldots,\nabla\phi_{i_{s}}$ for a carefully chosen $r_s$. For convenience, denote by $V_{s}$ the matrix with $\nabla\phi_{i_1},\ldots,\nabla\phi_{i_{s-1}}$ as columns and let $\mathcal{R}(V_{s})$ be the range of $V_{s}$. Let $\phi_{r_{s}}$ be an eigenvector such that $\nabla\phi_{r_{s}} \not\in \mathcal{R}(V_{s})$. To find such an $r_{s}$, we define
\begin{align}
    H^s_{kij} &= \nabla\phi_i^T (I-V_{s}(V_{s}^TV_{s})^{-1}V_{s}^T) \nabla\phi_{j}\label{Hdef0}\\ 
             &= \widetilde{A}_{kij} - \begin{bmatrix}\widetilde{A}_{kii_1} \ldots \widetilde{A}_{kii_{s-1}}\end{bmatrix}\begin{bmatrix}\widetilde{A}_{ki_1i_1}&\widetilde{A}_{ki_1i_2} &\ldots & \widetilde{A}_{ki_1i_{s-1}}\\\widetilde{A}_{ki_2i_1} & \widetilde{A}_{ki_2i_2} & \ldots & \widetilde{A}_{ki_2i_{s-1}}\\ \vdots & \vdots & \ddots & \vdots\\\widetilde{A}_{ki_{s-1}i_1} & \widetilde{A}_{ki_{s-1}i_2} & \ldots & \widetilde{A}_{ki_{s-1}i_{s-1}}\end{bmatrix}^{-1} \begin{bmatrix}\widetilde{A}_{ki_1j}\\\widetilde{A}_{ki_2j}\\\vdots\\\widetilde{A}_{ki_{s-1}j}\end{bmatrix} \label{Hdef}
\end{align}
Note that $H^s_{kii}$ is the squared norm of the projection of $\nabla\phi_i$ onto the vector space orthogonal to $\mathcal{R}(V_{s})$. Clearly $\nabla\phi_i \not\in \mathcal{R}(V_{s})$ if and only if $H^s_{kii} > 0$. To obtain a low frequency eigenvector $\bm{\phi}_{r_s}$ such that $H^s_{kr_sr_s} > 0$ we choose
\begin{align}
    r_s = \argmin_{i\in S_k}\{\lambda_i: H^s_{kii} \geq \theta_s\} \label{rs}
\end{align}
where $\theta_s$ is the $\tau_s$-percentile of the set $\{H^s_{kii}: i \in S_k\}$ and $\tau_s\in (0,100)$ is a hyperparameter.
Then we take $p_{s}$ to be the component of $\nabla\phi_{r_{s}}$ which is orthogonal to $\mathcal{R}(V_{s})$,
\begin{align}
    p_{s}=(I-V_{s}(V_{s}^TV_{s})^{-1}V_{s}^T)\nabla\phi_{r_{s}}. \label{psp1}
\end{align}
\paragraph{Step 5. Find $i_s \in S_k$ such that $\gamma_{ki_{s}}|\nabla\phi_{i_{s}}^Tp_{s}|$ is sufficiently large.}
Using Eq.~(\ref{Hdef0},~\ref{psp1}), we note that
\begin{align}
    \nabla\phi_{i}^Tp_{s} = H^s_{kir_s}. \label{npips}
\end{align}
To obtain $\bm{\phi}_{i_{s}}$ such that $\gamma_{ki_{s}}|\nabla\phi_{i_{s}}^Tp_{s}|$ is greater than a certain threshold, as in step $3$, we first define the maximum possible value of $\gamma_{ki_{s}}|\nabla\phi_{i}^Tp_{s}|$ using Eq.~(\ref{npips}) as,
\begin{align}
    \alpha_s = \max_{i \in S_k} \gamma_{ki} |\nabla\phi_{i}^Tp_{s}| = \max_{i \in S_k} \gamma_{ki} |H^{s}_{kir_{s}}|. \label{alphas}
\end{align}
Then we take the threshold to be $\delta_s\alpha_s$ where $\delta_{s} \in [0,1]$ is a hyperparameter. Finally, to obtain a low frequency eigenvector $\bm{\phi}_{i_{s}}$ we choose $i_{s}$ such that 
\begin{align}
    i_{s} &= \argmin_{i \in S_k}\{\lambda_i : \gamma_{ki}|\nabla\phi_{i}^Tp_s| \geq \delta_s\alpha_s\} = \argmin_{i \in S_k}\{\lambda_i : \gamma_{ki} |H^{s}_{kir_{s}}| \geq \delta_s\alpha_{s}\}.\label{is}
\end{align}



In the end we obtain a $d$-dimensional parameterization $\Phi_k$ of $U_k$ given by
\begin{align}
    \Phi_k &\equiv (\gamma_{ki_1}\bm{\phi}_{i_1},\ldots,\gamma_{ki_d}\bm{\phi}_{i_d}) \ \text{where}\nonumber\\
    \Phi_{k}(x_{k'})  &= (\gamma_{k{i_1}}\bm{\phi}_{i_{1}k'},\ldots,\gamma_{k{i_d}}\bm{\phi}_{i_{d}k'}) \ \text{and}\label{Psik}\\
    \Phi_k(U_k) &= (\Phi_k(x_{k'}))_{x_{k'} \in U_k}. \nonumber
\end{align}
We call $\Phi_k(U_k)$ the $k$th local view of the data in the $d$-dimensonal embedding space. It is a matrix with $|U_k|$ rows and $d$ columns. Denote the distortion of $\Phi_{k'}$ on $U_k$ by $\zeta_{kk'}$. Using Eq.~(\ref{DistortionPhikUk}) we obtain
\begin{align}
    \zeta_{kk'} &= \text{Distortion}(\Phi_{k'},U_k)\\
    &= \sup_{\substack{x_l, x_{l'} \in U_k\\x_{l}\neq x_{l'}}}\frac{\left\|\Phi_{k'}(x_l)-\Phi_{k'}(x_{l'})\right\|}{d_e(x_l,x_{l'})}\sup_{\substack{x_l, x_{l'} \in U_k\\x_l\neq x_{l'}}}\frac{d_e(x_l,x_{l'})}{\left\|\Phi_{k'}(x_l)-\Phi_{k'}(x_{l'})\right\|}. \label{zetakj}
\end{align}


\paragraph{Postprocessing.} The obtained local parameterizations are  post-processed so as to remove the anomalous parameterizations having unusually high distortion. We replace the local parameterization $\Phi_k$ of $U_k$ by that of a neighbor, $\Phi_{k'}$ where $x_{k'} \in U_k$, if the distortion $\zeta_{kk'}$ produced by $\Phi_{k'}$ on $U_k$ is smaller than the distortion $\zeta_{kk}$ produced by $\Phi_k$ on $U_k$. If $\zeta_{kk'}<\zeta_{kk}$ for multiple $k'$ then we choose the parameterization which produces the least distortion on $U_k$. This procedure is repeated until no replacement is possible. The pseudocode is provided below.

\begin{algorithm}[ht]
\SetAlgoLined
\KwIn{$d_e(x_k,x_{k'})_{k,k'=1}^{n},(I_{k},\Phi_k,\zeta_{kk})_{k=1}^{n}$}
\KwOut{$(\Phi_k,\zeta_{kk})_{k=1}^{n}$}
{
    $N_{\text{replaced}} \gets 1$\;
    \While{$N_{\text{replaced}} > 0$}{
        $N_{\text{replaced}} \gets 0$\;
        $\Phi^{\text{old}}_k \gets \Phi_k$ for all $k\in \{1,\ldots,n\}$\;
        \For{$k\gets1$ \KwTo $n$}{
            Compute $(\zeta_{kk'})_{x_{k'}\in U_k}$ (Eq.~(\ref{zetakj}))\;
            $k^* \gets \argmin_{x_{k'} \in U_k}\zeta_{kk'}$\;   
            \If{$k^* \neq k$}{
                $\ \Phi_{k} \gets \Phi_{k^*}^{\text{old}};\ \ \zeta_{kk} \gets \zeta_{kk^*};\ \ N_{\text{replaced}} \gets N_{\text{replaced}}+1$\;
            }
        }
    }
}
\caption{Postprocess-Local-Parameterization}
\label{algo:post}
\end{algorithm}

\editt{\textbf{A note on hyperparameters $N, (\tau_s, \delta_s)_{s=1}^{d}$}. Generally, $N$ should be small so that the low frequency eigenvectors form the set of candidate eigenvectors. In almost all of our experiments we take $N$ to be $100$. The set of $(\tau_s, \delta_s)_{s=1}^{d}$ is reduced to two hyperparameters, one for all $\tau_s$'s and one for all $\delta_s$'s. As explained above, $\tau_s$ enforces certain vectors to be non-zero and $\delta_s$ enforces certain directional derivatives to be large enough. Therefore, a small value of $\tau_s$ in $(0,100)$ and a large value of $\delta_s$ in $(0,1]$ is suitable. In most of our experiments, we used a value of $50$ for all $\tau_s$ and a value of $0.9$ for all $\delta_s$. Our algorithm is not too sensitive to the values of these hyperparameters. Other values of $N$, $\tau_s$ and $\delta_s$ would also result in the embeddings with high visual quality.}

\begin{figure}[h]
    \begin{tabular}{cc}
        \includegraphics[width=0.45\textwidth,height=0.45\textwidth,keepaspectratio]{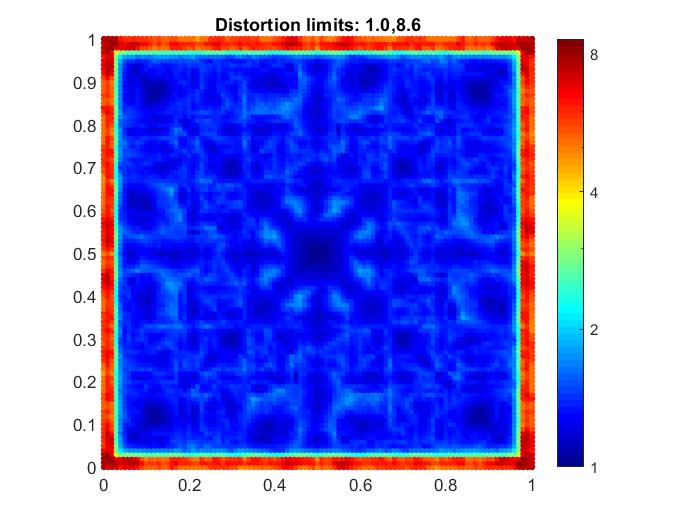}
     & \includegraphics[width=0.45\textwidth,height=0.45\textwidth,keepaspectratio]{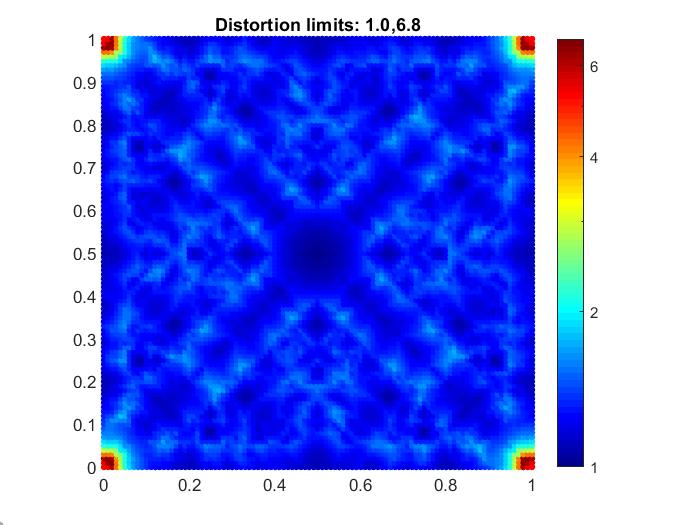}
    \end{tabular}
    \caption{Distortion of the obtained local parameterizations when the points on the boundary are not known (left) versus when they are known apriori (right). Each point $x_k$ is colored by $\zeta_{kk}$ (see Eq. (\ref{zetakj})).}
    \label{fig:fig31}
\end{figure}

\paragraph{Example.} We now build upon the example of the square grid at the end of Section~\ref{subsec:Atilde_esimate}. The values of the additional inputs are $N=100$, $\tau_s = 50$ and $\delta_s = 0.9$ for all $s \in \{1,\ldots,d\}$. Using Algo.~\ref{algo:blp} and~\ref{algo:post} we obtain $10^4$ local views $U_k$ and $\Phi_k(U_k)$ where $|U_k| = 25$ for all $k$. In the left image of Figure~\ref{fig:fig31}, we colored each point $x_k$ with the distortion $\zeta_{kk}$ of the local parameterization $\Phi_k$ on $U_k$. The mapped discrete balls $\Phi_k(U_k)$ for some values of $k$ are also shown in Figure~\ref{fig:local_views} in the Appendix~\ref{sec:supp}.

\begin{rmk}
\label{rmk:highdist}
Note that the parameterizations of the discrete balls close to the boundary have higher distortion. This is because the injectivity radius at the points close to the boundary is low and precisely zero at the points on the boundary. As a result, the size of the balls around these points exceeds the limit beyond which Theorem~\ref{thm:thm1} is applicable.
\end{rmk}

At this point we note the following remark in \citep{jones2007universal}.

\begin{rmk}
\label{rmk:dblmanifold}
As was noted by L. Guibas, when M has a boundary, in the case of Neumann boundary values, one may consider the ``doubled'' manifold, and may apply the result in Theorem~\ref{thm:thm1} for a possibly larger $r_k$.
\end{rmk}

Due to the above remark, assuming that the points on the boundary are known, we computed the distance matrix for the doubled manifold using the method described in \citep{lafondiffusion2004}. Then we recomputed the local parameterizations $\Phi_k$ keeping all other hyperparameters the same as before. In the right image of Figure~\ref{fig:fig31}, we colored each point $x_k$ with the distortion of the updated parameterization $\Phi_k$ on $U_k$. Note the reduction in the distortion of the paramaterizations for the neighborhoods close to the boundary. The distortion is still high near the corners.


\subsection{Time Complexity}
\label{subsec:blp_time}

\edit{The combined worst case time complexity of Algo.~\ref{algo:gl},~\ref{algo:blp} and~\ref{algo:post} is $O(n(N^2(k_{\text{lv}}+d)+k_{\text{lv}}^3 N_{\text{post}}d))$ where $N_{\text{post}}$ is the number of iterations it takes to converge in Algo.~\ref{algo:post} which was observed to be less than $50$ for all the examples in this paper. It took about a minute\footnote{Machine specification: MacOS version 11.4, Apple M1 Chip, $16$GB RAM.} to construct the local views in the above example as well as in all the examples in Section~\ref{sec:compare}.}

\section{Clustering for Intermediate Views}
\label{sec:clustering}
Recall that the discrete balls $U_k$ are the local views of the data in the high dimensional ambient space. In the previous section, we obtained the mappings $\Phi_k$ to construct the local views $\Phi_k(U_k)$ of the data in the $d$-dimensional embedding space. As discussed in Section~\ref{subsec:contrib}, one can use the GPA \citep{fabioprocrustes,gower1975generalized,ten1977orthogonal} to register these local views to recover a global embedding. In practice, too many small local views (high $n$ and small $|U_k|$) result in extremely high computational complexity. Moreover, small overlaps between the local views makes their registration susceptible to errors. Therefore, we perform clustering to obtain $M \ll n$ intermediate views, $\widetilde{U}_m$ and $\widetilde{\Phi}_m(\widetilde{U}_m)$, of the data in the ambient space and the embedding space, respectively. This reduces the time complexity and increases the overlaps between the views, leading to their quick and robust registration.

\subsection{Notation}

Our clustering algorithm is designed so as to ensure low distortion of the parameterizations $\widetilde{\Phi}_{m}$ on $\widetilde{U}_m$. We first describe the notation used and then present the pseudocode followed by a full explanation of the steps. Let $c_k$ be the index of the cluster $x_k$ belongs to. Then the set of points which belong to cluster $m$ is given by
\begin{align}
    \mathcal{C}_m = \{x_k \ | \ c_k = m\}. \label{Cm}
\end{align}
Denote by $c_{U_{k}}$ the set of indices of the neighboring clusters of $x_k$. \edit{The neighboring points of $x_k$ lie in these clusters}, that is,
\begin{align}
    c_{U_k} = \{c_{k'}\ |\ x_{k'} \in U_k\}. \label{cUk}
\end{align}
\edittt{We say that a point $x_k$ lies in the vicinity of a cluster $m$ if $m \in c_{U_k}$. Let $\widetilde{U}_m$ denote the $m$th intermediate view of the data in the ambient space. This constitutes the union of the local views associated with all the points belonging to cluster $m$, that is,
\begin{align}
    \widetilde{U}_m = \bigcup_{k:\ x_k\in \mathcal{C}_m}U_k. \label{Utildek}
\end{align}
Clearly, a larger cluster means a larger intermediate view. In particular, addition of $x_k$ to $\mathcal{C}_m$ grows the intermediate view $\widetilde{U}_m$ to $\widetilde{U}_m \cup U_k$,
\begin{align}
    \mathcal{C}_m \rightarrow \mathcal{C}_m \cup \{x_k\} \implies \widetilde{U}_m \rightarrow \widetilde{U}_m \cup U_k \label{grow}
\end{align}
Let $\widetilde{\Phi}_m$ be the $d$-dimensional parameterization associated with the $m$th cluster. This parameterization maps $\widetilde{U}_m$ to $\widetilde{\Phi}_m(\widetilde{U}_m)$, the $m$th intermediate view of the data in the embedding space. Note that a point $x_k$ generates the local view $U_k$  (see Eq.~(\ref{Uk})) which acts as the domain of the parameterization $\Phi_k$. Similarly, a cluster $\mathcal{C}_m$ obtained through our procedure, generates an intermediate view $\widetilde{U}_m$ (see Eq.~(\ref{Utildek})) which acts as the domain of the parameterization $\widetilde{\Phi}_m$. Overall, our clustering procedure replaces the notion of a local view per an individual point by an intermediate view per a cluster of points.
}

\begin{algorithm}[t]
\SetAlgoLined
\KwIn{$(U_k,\Phi_k)_{k=1}^{n}, \eta_{\text{min}}$}
\KwOut{$(\mathcal{C}_m,\widetilde{U}_m,\widetilde{\Phi}_{m})_{m=1}^{M},(c_k)_{k=1}^{n}$}
{
    Initialize $c_k \gets k$, $\mathcal{C}_m \gets \{x_m\}$, $\widetilde{\Phi}_m \gets \Phi_m$ for all $k,m \in \{1,\ldots,n\}$\;
    \For{$\eta \gets 2$ to $\eta_{\text{min}}$}{
        \edit{Compute $b_{m\leftarrow x_k}$ for all $m,k \in \{1,\ldots,n\}$ (Eq.~(\ref{cUk},~\ref{Utildek},~\ref{costkj}))\;}
    $m,k \gets \argmax_{m',k'} b_{m'\leftarrow x_{k'}};\ \text{bid}^* \gets b_{m\leftarrow x_{k}}$\;
        \While{$\text{bid}^* > 0$}{
            $s \gets c_k;\ \mathcal{C}_s \gets \mathcal{C}_s - x_k;\ c_k \gets m;\ \mathcal{C}_{m} \gets \mathcal{C}_{m} \cup x_k$\;
            \edit{Recompute $b_{m'\leftarrow x_{k'}}$ for all $(m',k') \in \mathcal{S}$ (Eq.~(\ref{SSk}))\;}
            $m,k \gets \argmax_{m',k'} b_{m'\leftarrow x_{k'}};\ \text{bid}^* \gets b_{m\leftarrow x_{k}}$\;
        }
    }
    $M \gets$ the number of non-empty clusters\;
    Remove $\mathcal{C}_m$, $\widetilde{\Phi}_m$ when $|\mathcal{C}_m| = 0$, relabel clusters from $1$ to $M$ and update $c_k$ with new labels\;
    Compute $(\widetilde{U}_m)_{m=1}^{M}$ (Eq.~(\ref{Utildek}))\;
}
\caption{Clustering}
\label{algo:clustering}
\end{algorithm}

\subsection{Low Distortion Clustering}
\label{subsec:low_dist_clustering}
\edittt{Initially, we start with $n$ singleton clusters where the point $x_k$ belongs to the $k$th cluster and the parameterization associated with the $k$th cluster is $\Phi_k$. Thus, $c_k = k$, $\mathcal{C}_m = \{x_m\}$ and $\widetilde{\Phi}_m = \Phi_m$ for all $k,m \in \{1,\ldots,n\}$. This automatically implies that initially $\widetilde{U}_m = U_m$. \textit{The parameterizations associated with the clusters remain the same throughout the procedure}. During the procedure, each cluster $\mathcal{C}_m$ is perceived as an entity which wants to grow the domain $\widetilde{U}_m$ of the associated parameterization $\widetilde{\Phi}_m$ by growing itself (see Eq.~\ref{grow}), while simultaneously keeping the distortion of $\widetilde{\Phi}_m$ on $\widetilde{U}_m$ low (see Eq.~\ref{zetakj}). To achieve that, each cluster $\mathcal{C}_m$ places a \textit{careful} bid $b_{m \leftarrow x_k}$ for each point $x_k$. The global maximum bid is identified and the underlying point $x_k$ is relabelled to the bidding cluster, hence updating $c_k$. With this relabelling, the bidding cluster grows and the source cluster shrinks. This procedure of shrinking and growing clusters is repeated until all \textit{non-empty} clusters are large enough, i.e. have a size at least $\eta_{\text{min}}$, a hyperparameter.
In our experiments, we choose $\eta_{\text{min}}$ from $\{5,10,15,20,25\}$.
We iterate over $\eta$ which varies from $2$ to $\eta_{\text{min}}$.
In the $\eta$-th iteration, we say that the $m$th cluster is \textit{small} if it is non-empty and has a size less than $\eta$, that is, when $|\mathcal{C}_m| \in (0,\eta)$. 
During the iteration, the clusters either shrink or grow  until no small clusters remain.
Therefore, at the end of the $\eta$-th iteration the non-empty clusters are of size at least $\eta$. 
After the last ($\eta_{\text{min}}$th) iteration, each non-empty cluster will have at least $\eta_{min}$ points and the empty clusters are pruned away.}

\begin{figure}[h]
    \centering
    \includegraphics[width=0.75\textwidth,keepaspectratio]{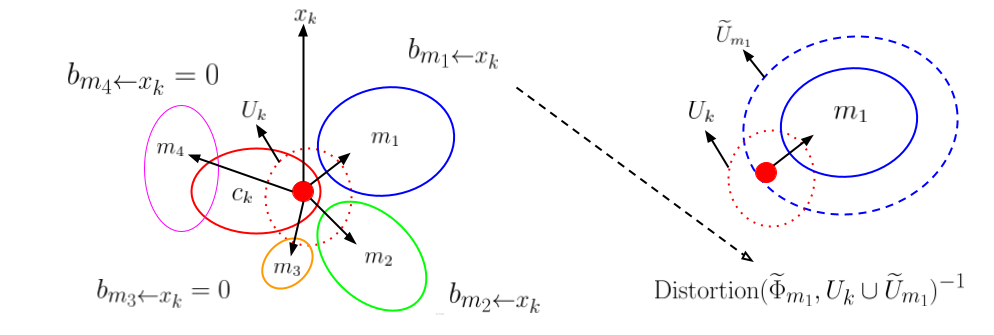}
    \caption{Computation of the bid for a point in a small cluster by the neighboring clusters in $\eta$-th iteration. (left) $x_k$ is a point represented by a small red disc, in a small cluster $c_k$ enclosed by solid red line. The dashed red line enclose $U_k$. \edit{Assume that} the cluster $c_k$ is small so that $|\mathcal{C}_{c_k}| \in (0,\eta)$. Clusters $m_1$, $m_2$, $m_3$ and $m_4$ are enclosed by solid colored lines too. Note that $m_1$, $m_2$ and $m_3$ lie in $c_{U_k}$ (the nonempty overlap between these clusters and $U_k$ indicate that), while $m_4 \not\in c_{U_k}$. Thus, the bid by $m_4$ for $x_k$ is zero. Since the size of cluster $m_3$ is less than the size of cluster $c_k$ i.e. $|\mathcal{C}_{m_3}| < |\mathcal{C}_{c_k}|$, the bid by $m_3$ for $x_k$ is also zero. Since clusters $m_1$ and $m_2$ satisfy all the conditions, the bids by $m_1$ and $m_2$ for $x_k$ are to be computed. (right) The bid $b_{m_1 \leftarrow x_k}$, is given by the inverse of the distortion of $\widetilde{\Phi}_{m_1}$ on $U_k \cup \widetilde{U}_{m_1}$, where the dashed blue line enclose $\widetilde{U}_{m_1}$. If the bid $b_{m_1 \leftarrow x_k}$ is greater (less) than the bid $b_{m_2 \leftarrow x_k}$, then the clustering procedure would favor relabelling of $x_k$ to $m_1$ ($m_2$).}
    \label{fig:clusteringcost}
\end{figure}

\edittt{\paragraph{Bid by cluster $m$ for $x_k$.} In the $\eta$-th iteration, we start by computing the bid $b_{m\leftarrow x_k}$ by each cluster $m$ for each point $x_k$. The bid function is designed so as to satisfy the following conditions. The first two conditions are there to halt the procedure while the last two conditions follow naturally. These conditions are also depicted in Figure~\ref{fig:clusteringcost}.
\begin{enumerate}
    \item No cluster bids for the points in \textit{large} clusters. Since $x_k$ belongs to cluster $c_k$ therefore, if $|\mathcal{C}_{c_k}| > \eta$ then the $b_{m\leftarrow x_k}$ is zero for all $m$.
    \item No cluster bids for a point in another cluster whose size is bigger than its own size. Therefore, if $|\mathcal{C}_m| < |\mathcal{C}_{c_k}|$ then again $b_{m\leftarrow x_k}$ is zero.
    \item A cluster bids for the points in its own vicinity. Therefore, if $m \not\in c_{U_k}$ (see Eq.~\ref{cUk}) then $b_{m\leftarrow x_k}$ is zero.
    \item Recall that a cluster $m$ aims to grow while keeping the distortion of associated parameterization $\widetilde{\Phi}_m$ low on its domain $\widetilde{U}_m$. If the $m$th cluster acquires the point $x_k$, $\widetilde{U}_{m}$ grows due to the addition of $U_k$ to it (see Eq.~(\ref{Utildek})), and so does the distortion of $\widetilde{\Phi}_{m}$ on it. Therefore, to ensure low distortion, the natural bid by $\mathcal{C}_m$ for the point $x_k$, $b_{m\leftarrow x_k}$, is $\text{Distortion}(\widetilde{\Phi}_{m},U_k \cup \widetilde{U}_{m})^{-1}$ (see Eq.~\ref{zetakj}).
\end{enumerate}
Combining the above conditions, we can write the bid by cluster $m$ for the point $x_k$ as,
\begin{align}
    b_{m\leftarrow x_k} = \left\{\begin{matrix*}\text{Distortion}(\widetilde{\Phi}_{m},U_k \cup \widetilde{U}_{m})^{-1} & \text{if }|\mathcal{C}_{c_k}| \in (0,\eta) \land m \in c_{U_k} \land |\mathcal{C}_{m}| \geq |\mathcal{C}_{c_k}| \\0&\text{otherwise}.\end{matrix*}\right. \label{costkj}
\end{align}
In the practical implementation of above equation, $c_{U_k}$ and $\widetilde{U}_{m}$ are computed on the fly using Eq.~(\ref{cUk},~\ref{Utildek}).
}



\edittt{\paragraph{Greedy procedure to grow and shrink clusters.} Given the bids by all the clusters for all the points, we grow and shrink the clusters so that at the end of the current iteration $\eta$, each non-empty cluster has a size at least $\eta$. We start by picking the global maximum bid, say $b_{m \leftarrow x_k}$. Let $x_k$ be in the cluster $s$ (note that $c_k$, the cluster of $x_k$, is $s$ before $x_k$ is relabelled). 
We relabel $c_k$ to $m$, and update the set of points in clusters $s$ and $m$, $\mathcal{C}_s$ and $\mathcal{C}_{m}$, using Eq.~(\ref{Cm}). This implicitly shrinks $\widetilde{U}_s$ and grows $\widetilde{U}_{m}$ (see Eq.~\ref{Utildek}) and affects the bids by clusters $m$ and $s$ or the bids for the points in these clusters. Denote the set of pairs of the indices of all such clusters and the points by
\begin{align}
    \mathcal{S} = \{(m',k') \in \{1,\ldots,n\}^2\ |\ m' \in \{m,s\} \text{ or } x_{k'} \in \mathcal{C}_s \cup \mathcal{C}_m\}. \label{SSk}
\end{align}
Then the bids $b_{m' \leftarrow x_{k'}}$ are recomputed for all $(m',k') \in \mathcal{S}$. It is easy to verify that for all other pairs, neither the conditions nor the distortion in Eq.~(\ref{costkj}) are affected. After this computation, we again pick the global maximum bid and repeat the procedure until the maximum bid becomes zero indicating that no non-empty small cluster remains. This marks the end of the $\eta$-th iteration.}


\edit{\paragraph{Final intermediate views in the ambient and the embedding space.} At the end of the last iteration, all non-empty clusters have at least $\eta_{\text{min}}$ points}. Let $M$ be the number of non-empty clusters. \edit{Using the pigeonhole principle one can show that $M$ would be less than or equal to $n/\eta_{\text{min}}$}. We prune away the empty clusters and relabel the non-empty ones from $1$ to $M$ while updating $c_k$ accordingly. With this, we obtain the clusters $(\mathcal{C}_m)_{m=1}^{M}$ with associated parameterizations $(\widetilde{\Phi}_m)_{m=1}^{M}$. Finally, using Eq.~(\ref{Utildek}), we obtain the $M$ intermediate views $(\widetilde{U}_m)_{m=1}^{M}$ of the data in the ambient space. Then, the intermediate views of the data in the embedding space are given by $(\widetilde{\Phi}_m(\widetilde{U}_m))_{m=1}^{M}$. Note that $\widetilde{\Phi}_m(\widetilde{U}_m)$ is a matrix with $|\widetilde{U}_m|$ rows and $d$ columns (see Eq.~(\ref{Psik})).

\begin{figure}[h]
    \centering
    \begin{tabular}{cc}
        \includegraphics[width=0.45\textwidth,height=0.45\textwidth,keepaspectratio]{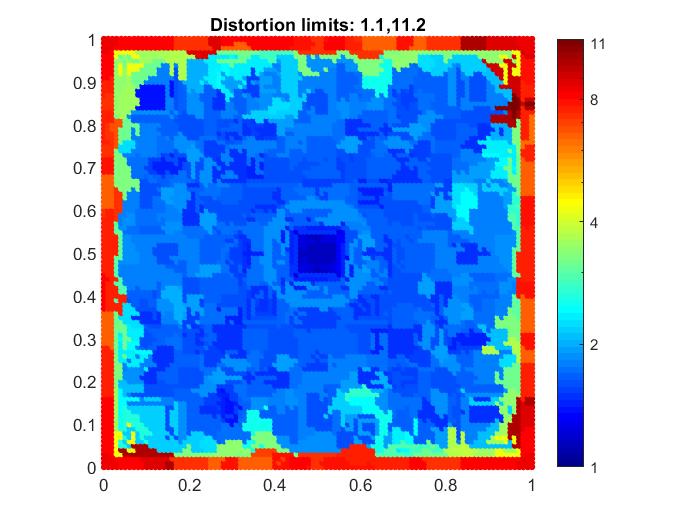} & \includegraphics[width=0.45\textwidth,height=0.45\textwidth,keepaspectratio]{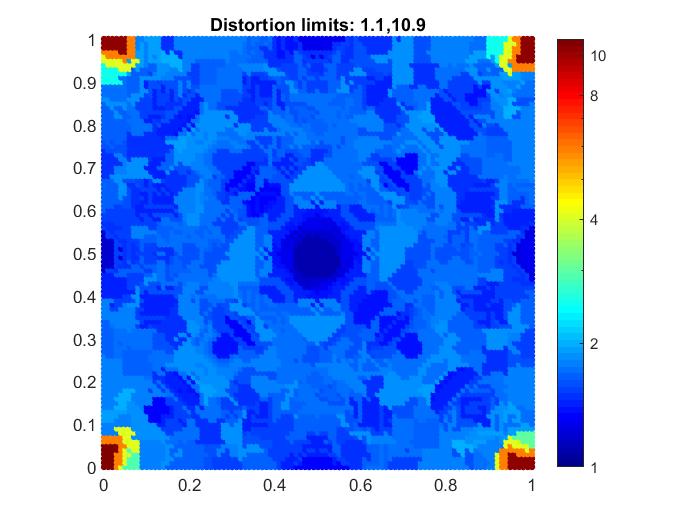}
    \end{tabular}
    \caption{\edit{Each point $x_k$ colored by $\tilde{\zeta}_{c_kc_k}$ when the points on the boundary of the square grid are unknown (left) versus when they are known apriori (right).}}
    \label{fig:fig41}
\end{figure}

\paragraph{Example.} We continue with our example of the square grid which originally contained about $10^4$ points. Therefore, before clustering we had about $10^4$ small local views $U_k$ and $\Phi_k(U_k)$, each containing $25$ points. After clustering with $\eta_{min} = 10$, we obtained $635$ clusters and therefore that many intermediate views $\widetilde{U}_m$ and $\widetilde{\Phi}_m(\widetilde{U}_m)$ with an average size of $79$. When the points on the boundary are known then we obtained $562$ intermediate views with an average size of $90$. Note that there is a trade-off between the size of the intermediate views and the distortion of the parameterizations used to obtain them. For convenience, define $\tilde{\zeta}_{mm}$ to be the distortion of $\widetilde{\Phi}_m$ on $\widetilde{U}_m$ using Eq.~(\ref{zetakj}). Then, as the size of the views are increased (by increasing $\eta_{min}$), the value of $\tilde{\zeta}_{mm}$ would also increase. In Figure~\ref{fig:fig41} we colored the points in cluster $m$, $\mathcal{C}_{m}$, with $\tilde{\zeta}_{mm}$. \edit{In other words, $x_k$ is colored by $\tilde{\zeta}_{c_kc_k}$}. Note the increased distortion in comparison to Figure~\ref{fig:fig31}.

\subsection{Time Complexity}
\label{subsec:cluster_time}
\edit{Our practical implementation of Algo.~\ref{algo:clustering} uses memoization for speed up. It took about a minute to construct intermediate views using in the above example with $n=10^4$, $k_{\text{lv}}=25$, $d=2$ and $\eta_{\text{min}} = 10$, and it took less than $2$ minutes for all the examples in Section~\ref{sec:compare}. It was empirically observed that the time for clustering is linear in $n$, $\eta_{\text{min}}$ and $d$ while it is cubic in $k_{\text{lv}}$.}

\section{Global Embedding using Procrustes Analysis}
\label{sec:ge}
In this section, we present an algorithm based on Procrustes analysis to \edit{align} the intermediate views $\widetilde{\Phi}_m(\widetilde{U}_m)$ and obtain a global embedding. The $M$ views $\widetilde{\Phi}_m(\widetilde{U}_m)$ are transformed by an orthogonal matrix $T_m$ of size $d \times d$, a $d$-dimensional translation vector $v_m$ and a positive scalar $b_m$ as a scaling component. The transformed views are given by $\widetilde{\Phi}^g_m(\widetilde{U}_m)$ such that 
\begin{align}
    \widetilde{\Phi}^g_m(x_k) = b_m\widetilde{\Phi}_m(x_k)T_m + v_m \quad \textrm{for all } x_k \in \widetilde{U}_m.
\end{align}
\edit{First we state a general approach to estimate these parameters, and its limitations in Section~\ref{subsec:general_app}. Then we present an algorithm in Section~\ref{subsec:desc_app} which computes these parameters and a global embedding of the data while addressing the limitations of the general procedure. In Section~\ref{subsec:tear} we describe a simple modification to our algorithm to tear apart closed manifolds. In Appendix~\ref{subsec:ltsa_stitch}, we contrast our global alignment procedure with that of LTSA.}


\begin{figure}[h]
    \centering
    \begin{tabular}{c}
        \includegraphics[height=0.25\textwidth,keepaspectratio]{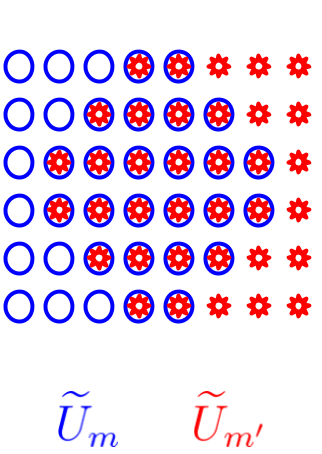} \hspace{30pt} \includegraphics[height=0.25\textwidth,keepaspectratio]{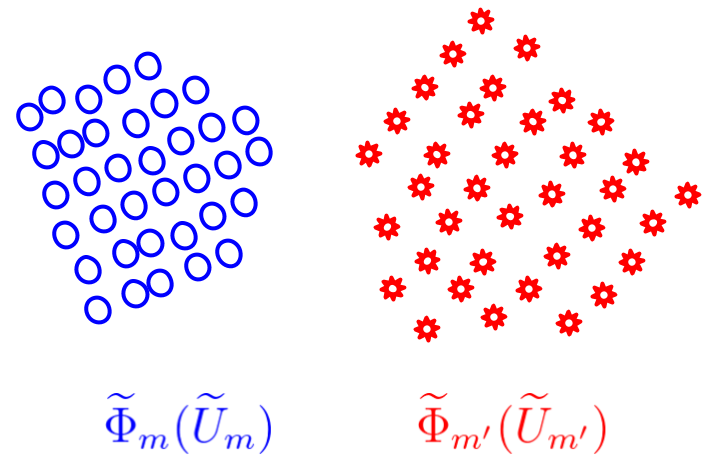} \hspace{30pt} \includegraphics[height=0.25\textwidth,keepaspectratio]{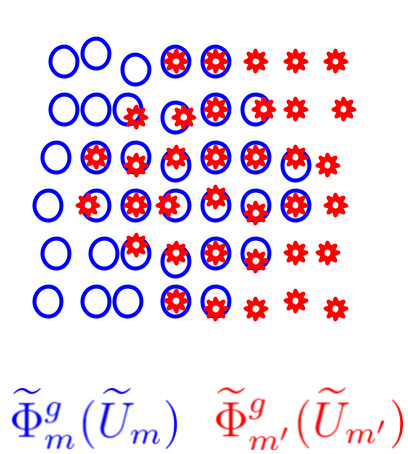}
    \end{tabular}
    \caption{(left) The intermediate views $\widetilde{U}_m$ and $\widetilde{U}_{m'}$ of a $2$d manifold in a possibly high dimensional ambient space. These views trivially align with each other. The red star in blue circles represent their overlap $\widetilde{U}_{mm'}$. (middle) The $m$th and $m'$th intermediate views in the $2$d embedding space. (right) \edit{Transformed} views after aligning $\widetilde{\Phi}_m(\widetilde{U}_{mm'})$ with $\widetilde{\Phi}_{m'}(\widetilde{U}_{mm'})$.}
    \label{fig:fig500}
\end{figure}

\subsection{General Approach for Alignment}
\label{subsec:general_app}

In general, the parameters $(T_m,v_m,b_m)_{m=1}^{M}$ are estimated so that for all $m$ and $m'$, the two transformed views of the overlap between $\widetilde{U}_m$ and $\widetilde{U}_{m'}$, obtained using the parameterizations $\widetilde{\Phi}^g_m$ and $\widetilde{\Phi}^g_{m'}$, align with each other. To be more precise, define the overlap between the $m$th and the $m'$th intermediate views in the ambient space as the set of points which lie in both the views, 
\begin{align}
    \widetilde{U}_{mm'} = \widetilde{U}_m \cap \widetilde{U}_{m'}.\label{Utildekkp}
\end{align}
In the ambient space, the $m$th and the $m'$th views are neighbors if $\widetilde{U}_{mm'}$ is non-empty. As shown in Figure~\ref{fig:fig500} (left), these neighboring views trivially align on the overlap between them. It is natural to ask for a low distortion global embedding of the data. Therefore, we must ensure that the embeddings of $\widetilde{U}_{mm'}$ due to the $m$th and the $m'$th view in the embedding space, also align with each other. Thus, the parameters $(T_m,v_m,b_m)_{m=1}^{M}$ are estimated so that $\widetilde{\Phi}^g_m(\widetilde{U}_{mm'})$ aligns with $\widetilde{\Phi}^g_{m'}(\widetilde{U}_{mm'})$ for all $m$ and $m'$. However, due to the distortion of the parameterizations it is usually not possible to perfectly align the two embeddings (see Figure~\ref{fig:fig500}). We can represent both embeddings of the overlap as matrices with $|\widetilde{U}_{mm'}|$ rows and $d$ columns. \edit{Then we choose the measure of the alignment error to be the squared Frobenius norm of the difference of the two matrices.} The error is trivially zero if $\widetilde{U}_{mm'}$ is empty. \edit{Overall, the parameters are estimated so as to minimize the following alignment error}
\begin{align}
    \mathcal{L}((T_m,v_m,b_m)_{m=1}^{M}) = \frac{1}{2M}\sum_{\substack{m=1\\m'=1}}^{M}\left\|\widetilde{\Phi}^g_m(\widetilde{U}_{mm'})-\widetilde{\Phi}^g_{m'}(\widetilde{U}_{mm'})\right\|^2_F.  \label{alignerr}
\end{align}
In theory, one can start with  \edit{a trivial initialization} of $T_m$, $v_m$ and $b_m$ as $I_d$, $\mathbf{0}$ and $1$, and directly use GPA \citep{fabioprocrustes,gower1975generalized,ten1977orthogonal} to obtain a local minimum of the above alignment error. This approach has two issues.
\begin{enumerate}
    \item Like most optimization algorithms, the rate of convergence to a local minimum and the quality of it depends on the initialization of the parameters. \edit{We empirically observed that with a trivial initialization of the parameters, GPA may take a great amount of time to converge and may also converge to an inferior local minimum.}
    \item Using GPA to align a view with all of its adjacent views would prevent us from tearing apart closed manifolds; as an example see Figure \edit{\ref{fig:gpa1}}.
\end{enumerate} 

\edit{These issues are addressed in subsequent Sections~\ref{subsec:desc_app} and~\ref{subsec:tear}, respectively.}

\subsection{GPA Adaptation for Global Alignment}
\label{subsec:desc_app}

\begin{algorithm}[ht]
\SetAlgoLined
\KwIn{$(x_k,c_k,w)_{k=1}^{n},(\mathcal{C}_m,\widetilde{\Phi}_m,\widetilde{U}_m)_{m=1}^{M}$, $\text{to\_tear}$, $\nu$, $N_r$}
\KwOut{$(T_m,b_m,v_m)_{m=1}^{M}$}
{
    \For{$\text{Iter}\gets 1$ to $N_r+1$}{
        \uIf{$\text{Iter} = 1$}{
            Initialize $T_m \gets I, v_m \gets 0$\;
            Compute $b_{m}$ (Eq.~(\ref{bk}))\;
            Compute $(s_m,p_{s_m})_{m=1}^{M}$ (Eq.~(\ref{s1p1},~\ref{pksk}) in Appendix~\ref{sec:smpm})\;
            \edit{$\mathcal{A} \gets \{s_1\}$ \%The set of already transformed views\;}
        }
        \Else{
            $(s_m)_{m=2}^{M} \gets$ random permutation of $(1,\ldots,M)$ excluding $s_1$\;
        }
        \For{$m\gets 2$ to $M$}{
            $s \gets s_m$, $p \gets p_{s_m}$\;
            (\textbf{Step R1}) $T_{s}, v_{s} \gets $ Procrustes ($\widetilde{\Phi}^g_{p}(\widetilde{U}_{sp})$,$\widetilde{\Phi}^g_{s}(\widetilde{U}_{sp})$, No scaling)\; 
            \uIf{$\text{to\_tear}=$ False}{
                (\textbf{Step R2}) Compute $\mathcal{Z}_s$ (Eq.~(\ref{Zs1}))\;
            }
            \Else{
                (\textbf{Step R2}) Compute $\mathcal{Z}_s$ (Eq.~(\ref{Zs2}))\;
            }
            
            (\textbf{Step R3}) $\mu_{s} \gets $ Centroid of $(\widetilde{\Phi}^g_{m'}(\widetilde{U}_{sm'}))_{m'\in \mathcal{Z}_s}$\; 
            (\textbf{Step R4}) $T_{s}, v_{s} \gets $ Procrustes ($\mu_{s},\widetilde{\Phi}^g_{s}(\cup_{m' \in \mathcal{Z}}U_{sm'})$, No scaling)\;
            \edit{(\textbf{Step R5}) $\mathcal{A} \gets \mathcal{A} \cup \{s\}$\;}
        }
    }
    Compute $(y_k)_{k=1}^{n}$ (Eq.~(\ref{yk})).
}
\caption{Calculate-Global-Embedding}
\label{algo:geinit}
\end{algorithm}


\editt{First we look for a better than trivial initialization of the parameters so that the views are approximately aligned. The idea is to build a rooted tree where nodes represent the intermediate views. This tree is then traversed in a breadth first order starting from the root. As we traverse the tree, the intermediate view associated with a node is aligned with the intermediate view associated with its parent node (and with a few more views), thus giving a better initialization of the parameters. Subsequently, we refine these parameters using a similar procedure involving random order traversal over the intermediate views.}


\begin{figure}[ht]
    \centering
    \begin{tabular}{|>{\raggedright\arraybackslash}M{0.475\textwidth}|>{\raggedright\arraybackslash}M{0.475\textwidth}|}
        \hline
        \makecell{\includegraphics[height=0.225\textwidth,keepaspectratio]{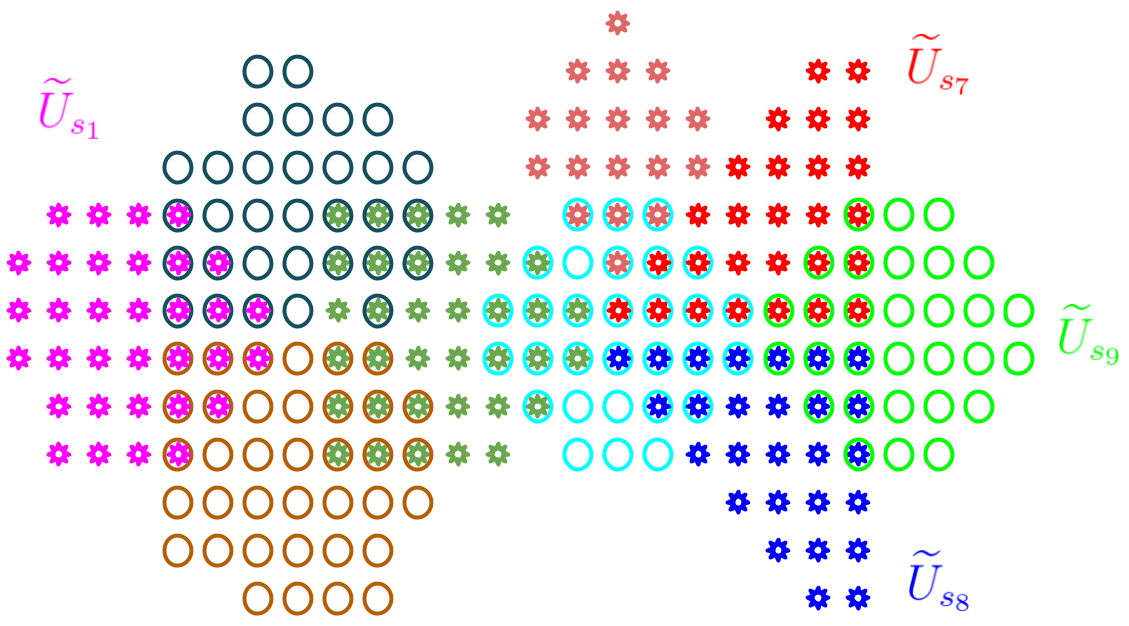}} & \makecell{\includegraphics[height=0.225\textwidth,keepaspectratio]{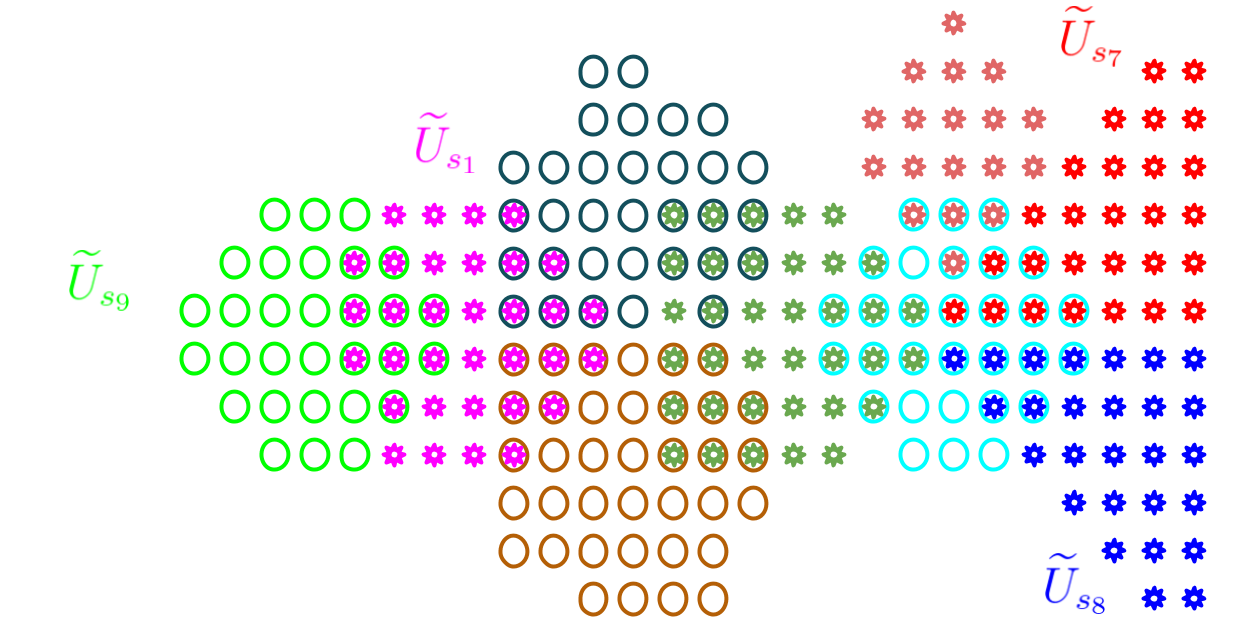}} \\
        \hline
        \footnotesize{(a.1). Nine intermediate views $(\widetilde{U}_{s_m})_{m=1}^{9}$ of a $2d$ manifold with boundary are shown. $\widetilde{U}_{s_9}$ has $\widetilde{U}_{s_7}$ and $\widetilde{U}_{s_8}$ as the neighboring views.} & \footnotesize{(a.2). In combination with (a.1), nine intermediate views $(\widetilde{U}_{s_m})_{m=1}^{9}$ of a closed $2d$ manifold are shown. In addition to $\widetilde{U}_{s_7}$ and $\widetilde{U}_{s_8}$, $\widetilde{U}_{s_9}$ also has $\widetilde{U}_{s_1}$ as the neighboring view.}\\
        \hline
        \makecell{\includegraphics[height=0.225\textwidth,keepaspectratio]{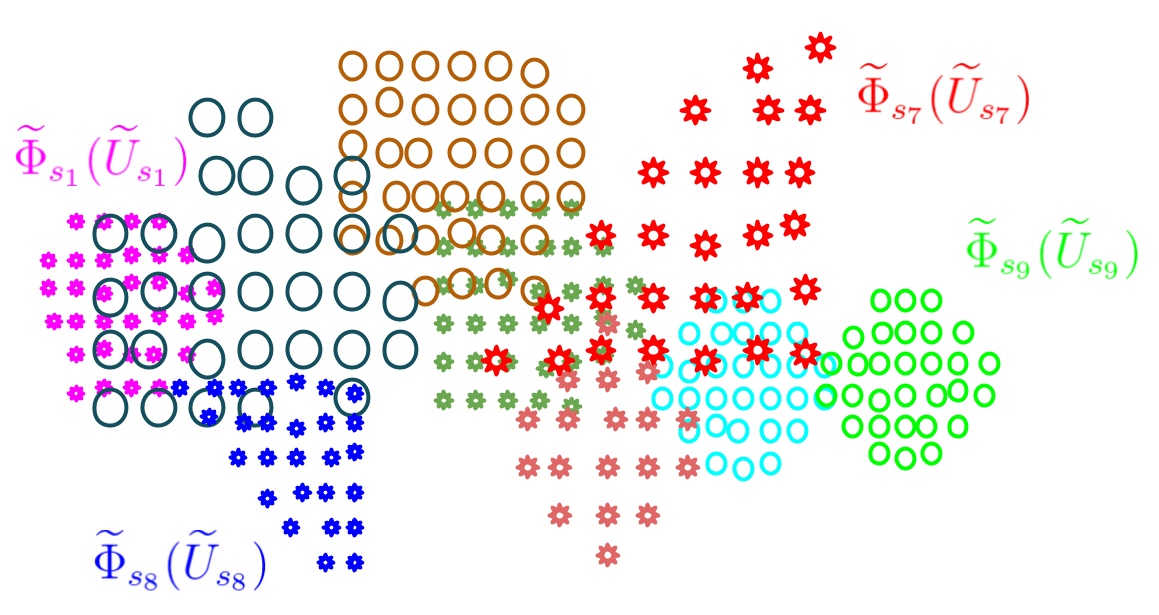}} & \makecell{\includegraphics[height=0.225\textwidth,keepaspectratio]{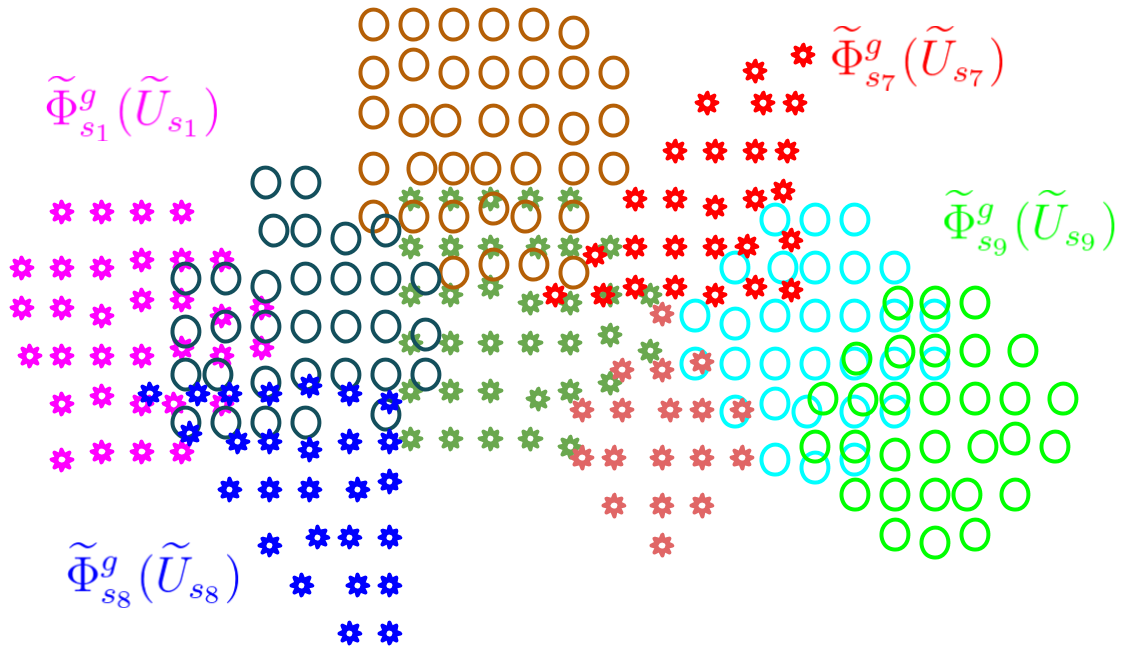}}\\
        \hline
        \footnotesize{(b) The intermediate views $(\widetilde{\Phi}_{s_m}(\widetilde{U}_{s_m}))_{m=1}^{9}$ in the $2$d embedding space, as they were passed as input to Algo.~\ref{algo:geinit}. These views are scrambled in the embedding space and Algo.~\ref{algo:geinit} will move them to the right location.} & \footnotesize{(c) The transformed views after scaling them using $b_m$ as in Eq.~(\ref{bk}).}\\
        \hline
            \end{tabular}
    \caption{An illustration of the intermediate views in the ambient and the embedding space as they are passed as input to Algo.~\ref{algo:geinit} and are scaled using Eq.~(\ref{bk}).}
    \label{fig:ge_algo_new1}
\end{figure}
\paragraph{Initialization ($\text{Iter}=1$, $\text{to\_tear}=$ False).}
In the first outer loop of Algo.~\ref{algo:geinit}, we start with $T_m=I_d$, $v_m$ as the zero vector and compute $b_m$ so as to bring the intermediate views $\widetilde{\Phi}_{m}(\widetilde{U}_m)$ to the same scale as their counterpart $\widetilde{U}_m$ in the ambient space. In turn this brings all the views to similar scale (see Figure~\ref{fig:ge_algo_new1} (c)). We compute the scaling component $b_m$ to be the ratio of the median distance between unique points in $\widetilde{U}_m$ and in $\widetilde{\Phi}_k(\widetilde{U}_m)$, that is,
\begin{align}
    b_m = \frac{\text{median}\left\{d_e(x_k,x_{k'})\ |\ x_k,x_{k'}\in \widetilde{U}_m, x_k \neq x_{k'}\right\}}{\text{median}\left\{\left\|\widetilde{\Phi}_{m}(x_k) - \widetilde{\Phi}_{m}(x_{k'})\right\|_2\ |\ x_{k},x_{k'}\in \widetilde{U}_m, x_k \neq x_{k'}\right\}}. \label{bk}
\end{align}

\begin{figure}[!h]
    \centering
    \begin{tabular}{|>{\raggedright\arraybackslash}M{0.475\textwidth}|>{\raggedright\arraybackslash}M{0.475\textwidth}|}
        \hline  
        \multicolumn{2}{|c|}{\footnotesize{Step R1 ($m=9$)}}\\
        \hline
        \makecell{\includegraphics[height=0.225\textwidth,keepaspectratio]{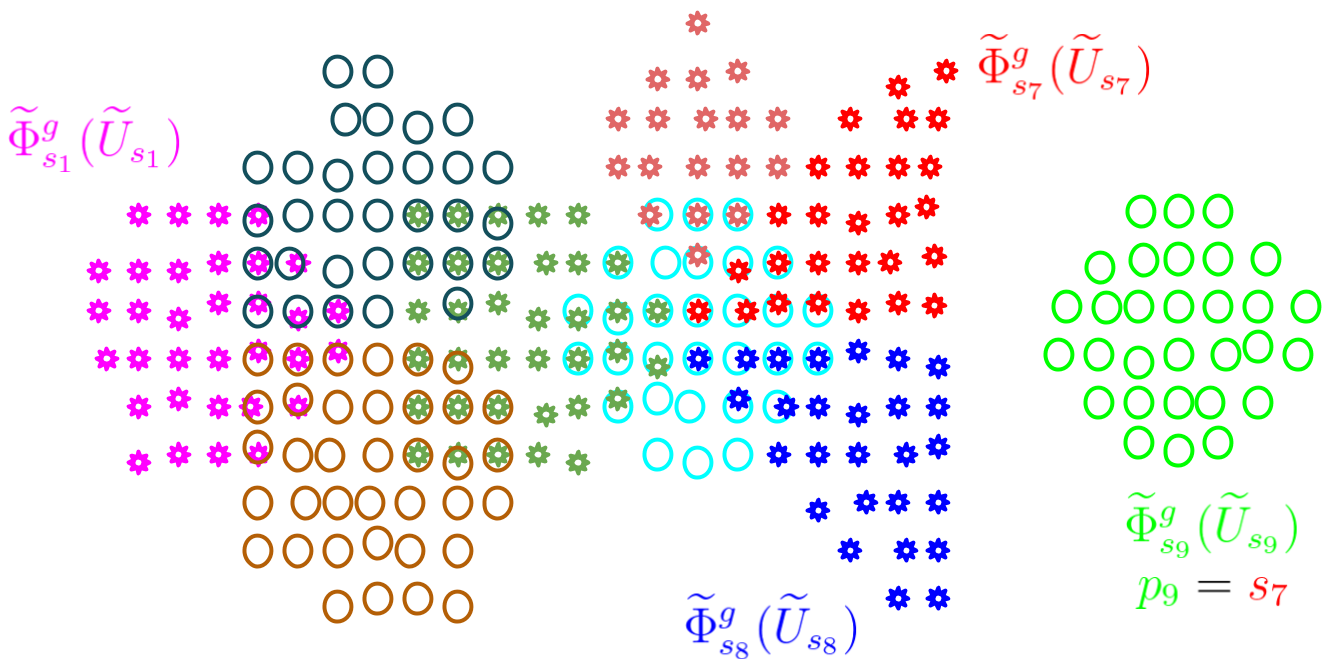}} & \makecell{\includegraphics[height=0.225\textwidth,keepaspectratio]{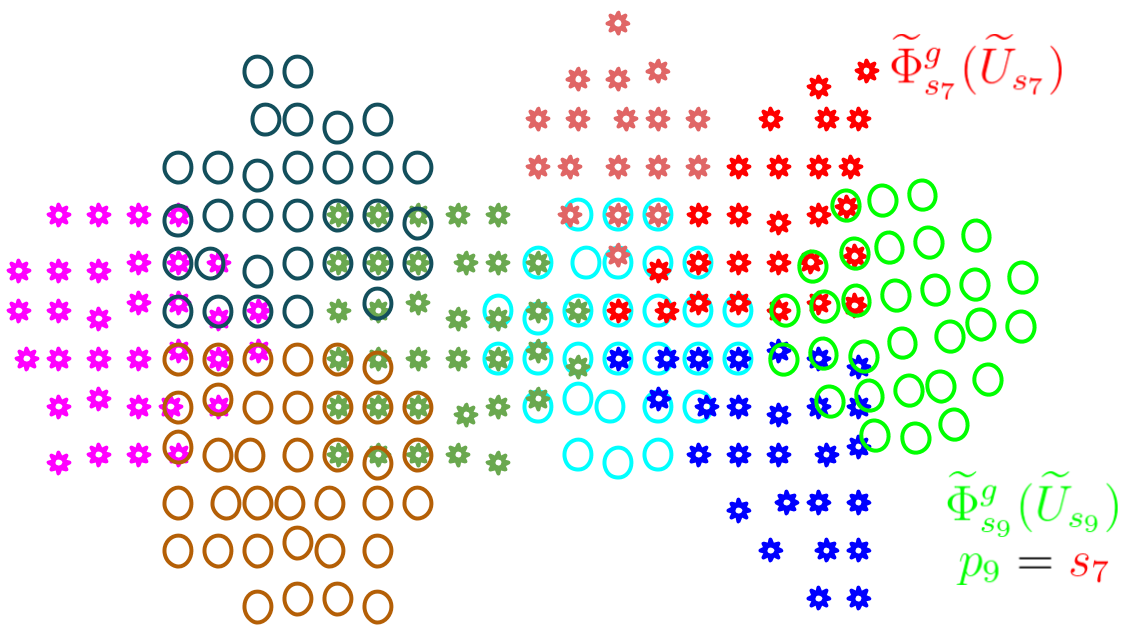}}\\
        \hline
        \footnotesize{(d) The transformed intermediate views $(\widetilde{\Phi}^g_{s_m}(\widetilde{U}_{s_m}))_{m=1}^{9}$ before the start of the iteration $m=9$. The first eight views are approximately aligned and the ninth view is to be aligned. 
        Inaccuracies occur due to distortion.} & \footnotesize{(e) Assuming $p_9 = s_7$, step R1 computed $T_{s_9}$ and $v_{s_9}$ so that $\widetilde{\Phi}^g_{s_9}(\widetilde{U}_{s_9s_7})$ aligns with $\widetilde{\Phi}^g_{s_7}(\widetilde{U}_{s_9s_7})$. The transformed view $\widetilde{\Phi}^g_{s_9}(\widetilde{U}_{s_9})$ is shown. Note that step R1 results in the same output for both cases in Fig.~\ref{fig:ge_algo_new1} (a)}\\
        \hline
        \multicolumn{2}{|c|}{\footnotesize{Step R2 and R3 ($m=9$)}}\\
        \hline
        \makecell{\includegraphics[height=0.225\textwidth,keepaspectratio]{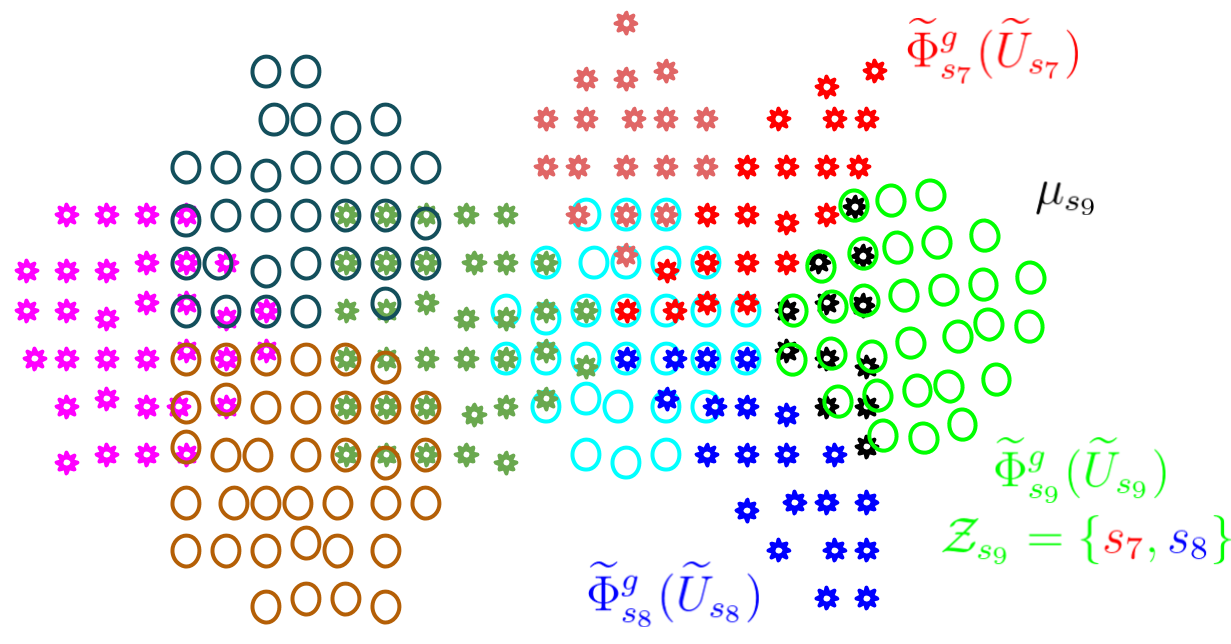}} & \makecell{\includegraphics[height=0.2025\textwidth,keepaspectratio]{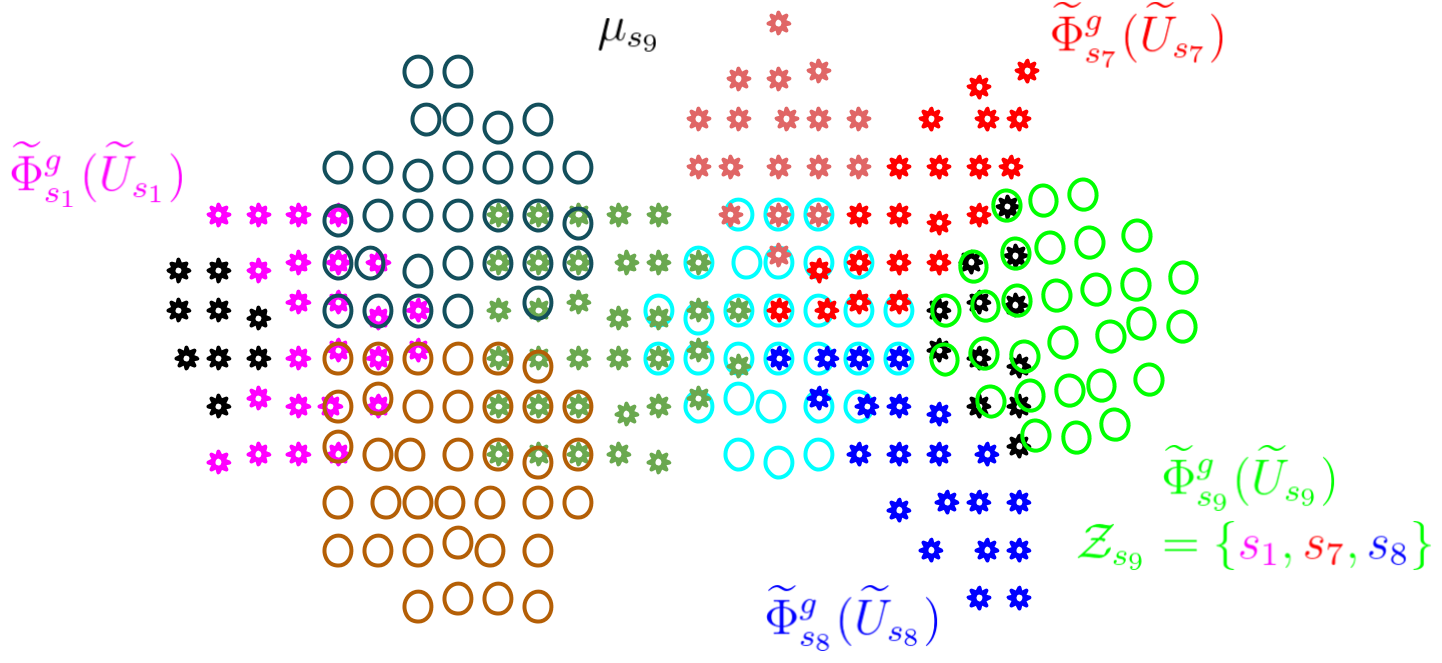}}\\
        \hline
        \footnotesize{(f.1) For a manifold with boundary, $\widetilde{U}_{s_9}$ has non-empty overlaps with $\widetilde{U}_{s_7}$ and $\widetilde{U}_{s_8}$ only. Therefore, step R2 computed $\mathcal{Z}_{s_9} = \{s_7,s_8\}$. The obtained $\mu_{s_9}$ in step R3 is also shown in black.} & \footnotesize{(f.2) For a closed manifold, $\widetilde{U}_{s_9}$ has non-empty overlaps with $\widetilde{U}_{s_1}$, $\widetilde{U}_{s_7}$ and $\widetilde{U}_{s_8}$. Therefore, step R2 computed $\mathcal{Z}_{s_9} = \{s_1,s_7,s_8\}$. The obtained $\mu_{s_9}$ in step R3 is also shown in black.}\\
        \hline
        \multicolumn{2}{|c|}{\footnotesize{Step R4 ($m=9$)}}\\
        \hline
        \makecell{\includegraphics[height=0.225\textwidth,keepaspectratio]{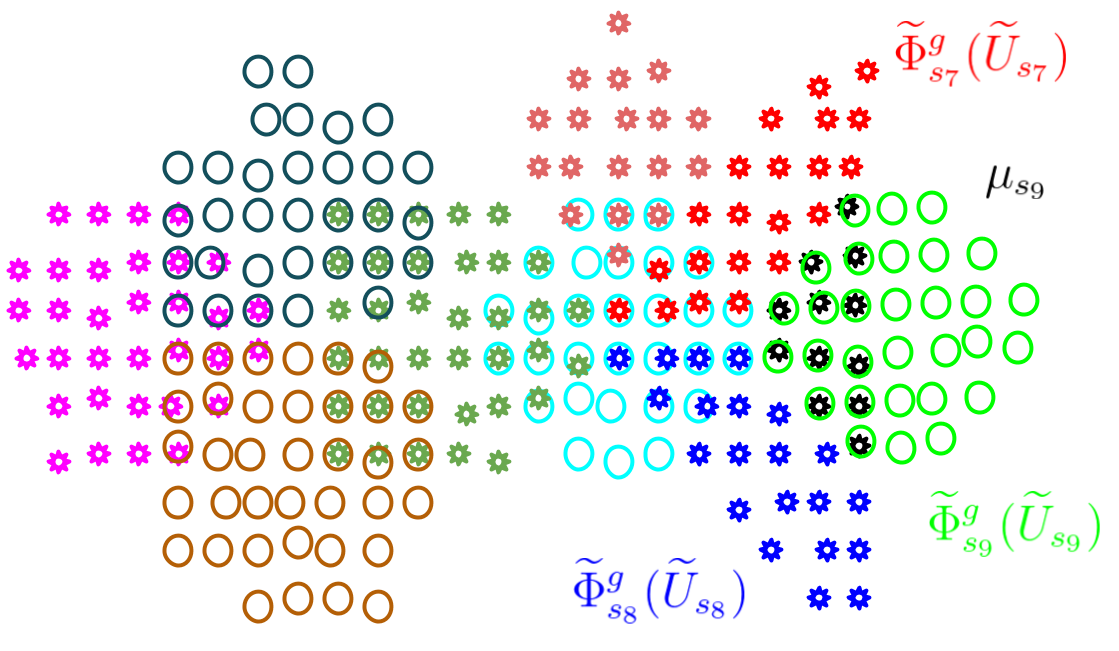}} & \makecell{\includegraphics[height=0.225\textwidth,keepaspectratio]{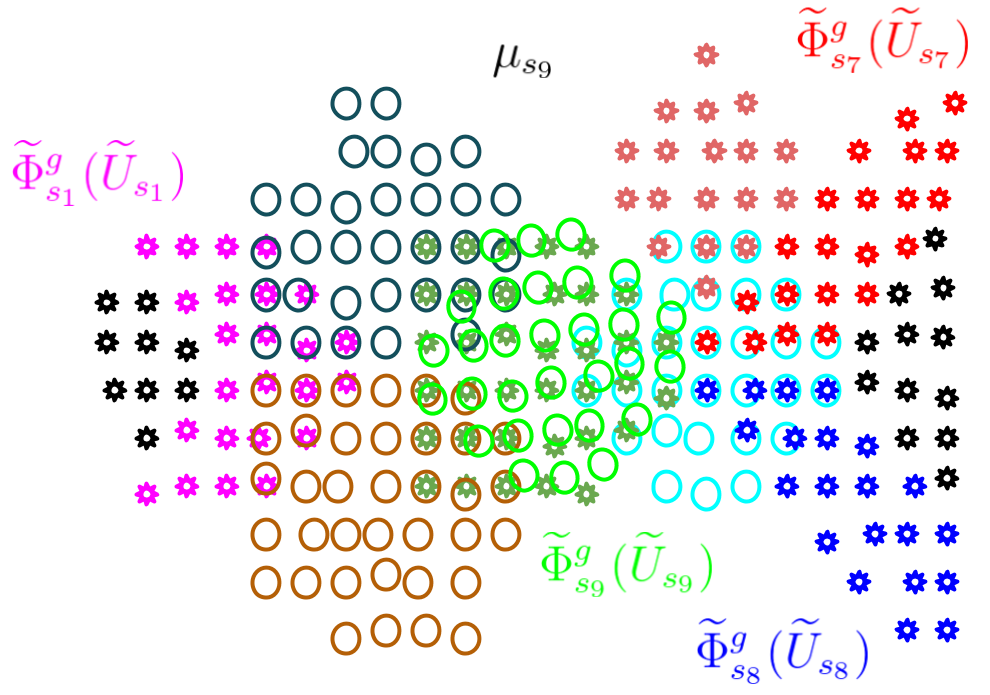}}\\
        \hline
        \footnotesize{(g.1) For a manifold with boundary, step R4 updated $T_{s_9}$ and $v_{s_9}$ so that the view $\widetilde{\Phi}_{s_9}^g(\widetilde{U}_{s_9s_7}\ \cup\ \widetilde{U}_{s_9s_8})$ aligns with $\mu_{s_9}$ in (f.1). The resulting view $\widetilde{\Phi}_{s_9}^g(\widetilde{U}_{s_9})$ is shown.} & \footnotesize{(g.2) For a closed manifold step R4 updated $T_{s_9}$ and $v_{s_9}$ so that view $\widetilde{\Phi}_{s_9}^g(\widetilde{U}_{s_9s_1} \cup \widetilde{U}_{s_9s_7} \cup \widetilde{U}_{s_9s_8})$ aligns with $\mu_{s_9}$ in (f.2). The resulting view $\widetilde{\Phi}_{s_9}^g(\widetilde{U}_{s_9})$ is shown. This is not a desired output as it distorts the global embedding. We resolve this issue in Section~\ref{subsec:tear}.}\\
        \hline
    \end{tabular}
    \caption{An illustration of steps R1 to R4 in Algo.~\ref{algo:geinit}, in continuation of Figure~\ref{fig:ge_algo_new1}.}
    \label{fig:ge_algo_new2}
\end{figure}

\editt{Then we transform the the views in a sequence $(s_m)_{m=1}^{M}$. This sequence corresponds to the breadth first ordering of a tree starting from its root node (which represents $s_1$th view). Let the $p_{s_m}$th view be the parent of the $s_m$th view.  Here $p_{s_m}$ lies in $\{s_1,\ldots,s_{m-1}\}$ and it is a neighboring view of the $s_m$th view in the ambient space, i.e. $\widetilde{U}_{s_mp_{s_m}}$ is non-empty. Details about the computation of these sequences is provided in Appendix~\ref{sec:smpm}. Note that $p_{s_1}$ is not defined and consequently, the first view in the sequence ($s_1$th view) is not transformed, therefore $T_{s_1}$ and $v_{s_1}$ are not updated.  We also define $\mathcal{A}$, initialized with $s_1$, to keep track of visited nodes which also represent the already transformed views. Then we iterate over $m$ which varies from $2$ to $M$. For convenience, denote the current ($m$th) node $s_m$ by $s$ and its parent $p_{s_m}$ by $p$. The following procedure updates $T_{s}$ and $v_{s}$ (refer to Figure~\ref{fig:ge_algo_new1} and \ref{fig:ge_algo_new2} for an illustration of this procedure).}

\paragraph{Step R1.} We compute a temporary value of $T_s$ and $v_s$ by aligning the views $\widetilde{\Phi}^g_{s}(\widetilde{U}_{sp})$ and $\widetilde{\Phi}^g_{p}(\widetilde{U}_{sp})$ of the overlap $\widetilde{U}_{sp}$, using Procrustes analysis \citep{gower2004procrustes} without modifying $b_s$.


\paragraph{Step R2.} Then we identify more views to align the $s$th view with. We compute a subset $\mathcal{Z}_s$ of the set of already visited nodes $\mathcal{A}$ such that $m' \in \mathcal{Z}_s$ if the $s$th view and the $m'$th view are neighbors in the ambient space. \editt{Note that, at this stage, $\mathcal{A}$ is the same as the set $\{s_1,\ldots,s_{m-1}\}$, the indices of the first $m-1$ views.} Therefore,
\begin{align}
    \mathcal{Z}_s = \{m' |\ \widetilde{U}_{sm'} \neq \emptyset\} \cap \mathcal{A}. \label{Zs1}
\end{align}


\paragraph{Step R3.} We then compute the centroid $\mu_{s}$ of the views $(\widetilde{\Phi}^g_{m'}(\widetilde{U}_{sm'}))_{m'\in \mathcal{Z}_s}$. Here $\mu_s$ is a matrix with $d$ columns and the number of rows given by the size of the set $\cup_{m'\in\mathcal{Z}_s}\widetilde{U}_{sm'}$. A point in this set can have multiple embeddings due to multiple parameterizations $(\widetilde{\Phi}^g_{m'})_{m'\in \mathcal{Z}_s}$ depending on the overlaps $(\widetilde{U}_{sm'})_{m'\in \mathcal{Z}_s}$ it lies in. The mean of these embeddings forms a row in $\mu_s$. 


\paragraph{Step R4.} Finally, we update $T_{s}$ and $v_{s}$ by aligning the view $\widetilde{\Phi}^g_{s}(\widetilde{U}_{sm'})$ with $\widetilde{\Phi}^g_{m'}(\widetilde{U}_{sm'})$ for all $m' \in \mathcal{Z}_s$. This alignment is based on the approach in \citep{fabioprocrustes,gower1975generalized} where, using the Procrustes analysis \citep{gower2004procrustes,matlabprocrustes}, the view $\widetilde{\Phi}^g_{s}(\cup_{m' \in \mathcal{Z}_s}\widetilde{U}_{sm'})$ is aligned with the centroid $\mu_{s}$, without modifying $b_s$.

\paragraph{Step R5.} After the $s$th view is transformed, we add it to the set of transformed views $\mathcal{A}$.


\editt{\paragraph{Parameter Refinement ($\text{Iter}\geq 2$, $\text{to\_tear}=$ False).} At the end of the first iteration of the outer loop in Algo.~\ref{algo:geinit}, we have an initialization of $(T_m,b_m,v_m)_{m=1}^{M}$ such that transformed intermediate views are approximately aligned. To further refine these parameters, we iterate over $(s_m)_{m=2}^{M}$ in random order and perform the same five step procedure as above, $N_r$ times. Besides the random-order traversal, the other difference in a refinement iteration is that the set of already visited nodes $\mathcal{A}$, contains all the nodes instead of just the first $m-1$ nodes. This affects the computation of $\mathcal{Z}_s$ (see Eq.~(\ref{Zs1})) in step R2 so that the $s$th intermediate view is now aligned with all those views which are its neighbors in the ambient space. Note that the step R5 is redundant during refinement.} 

\begin{figure}[h!]
    \centering
    \begin{tabular}{M{0.005\textwidth}|M{0.15\textwidth}|M{0.15\textwidth}|M{0.15\textwidth}|M{0.15\textwidth}|M{0.15\textwidth}|}
        & \tiny{Input} & \multicolumn{3}{c}{\tiny{First iteration of the outer loop and stages within inner loop}} & \tiny{End of outer loop}\\
        \hline
        & & \tiny{Before} & \tiny{Half-way} & \tiny{End} & \\
        \hline
        \rotatebox{90}{\tiny{Square}} & \includegraphics[width=0.15\textwidth,height=0.175\textwidth,keepaspectratio]{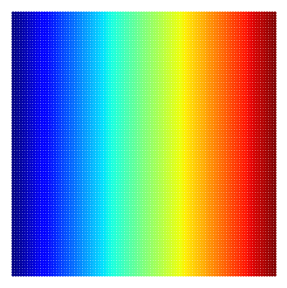} & \includegraphics[width=0.15\textwidth,height=0.175\textwidth,keepaspectratio]{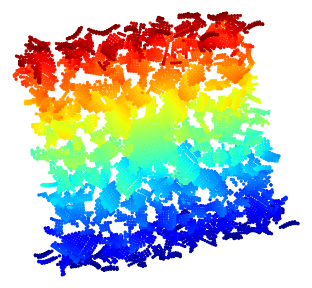} & \includegraphics[width=0.15\textwidth,height=0.175\textwidth,keepaspectratio]{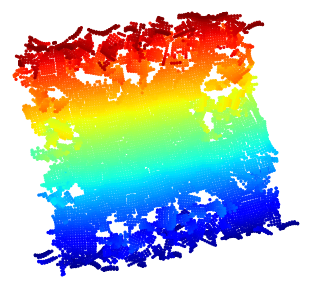} & \includegraphics[width=0.15\textwidth,height=0.175\textwidth,keepaspectratio]{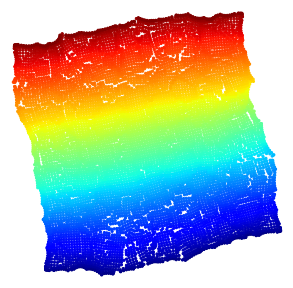} & \includegraphics[width=0.15\textwidth,height=0.175\textwidth,keepaspectratio]{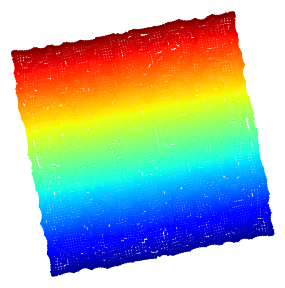} \\
        \hline
        \rotatebox{90}{\tiny{Sphere}} & \includegraphics[width=0.15\textwidth,height=0.175\textwidth,keepaspectratio]{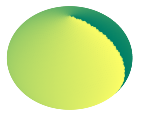} & \includegraphics[width=0.15\textwidth,height=0.175\textwidth,keepaspectratio]{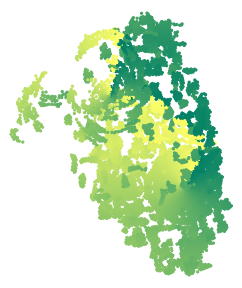} & \includegraphics[width=0.15\textwidth,height=0.175\textwidth,keepaspectratio]{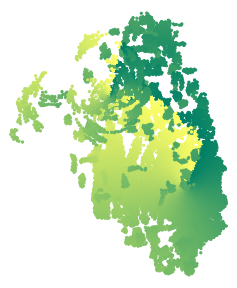} & \includegraphics[width=0.15\textwidth,height=0.175\textwidth,keepaspectratio]{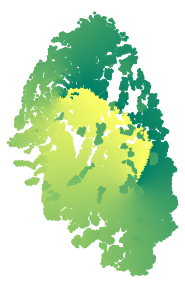} & \includegraphics[width=0.15\textwidth,height=0.175\textwidth,keepaspectratio]{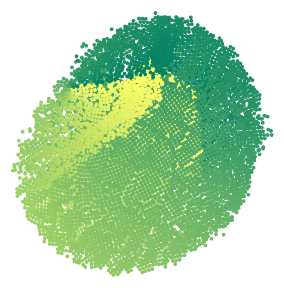} \\
        \hline
    \end{tabular}
    \caption{\editt{$2$d embeddings of a square and a sphere at different stages of Algo.~\ref{algo:geinit}. For illustration purpose, in the plots in the $2$nd and $3$rd columns the translation parameter $v_m$ was manually set for those views which do not lie in the set $\mathcal{A}$.
    Note that the embedding of the sphere is fallacious. The reason and the resolution is provided in Section~\ref{subsec:tear}.}}
    \label{fig:gpa1}
\end{figure}

In the end, we compute the global embedding $y_k$ of $x_k$ by mapping $x_k$ using the transformed parameterization associated with the cluster $c_k$ it belongs to,
\begin{align}
    y_k = \widetilde{\Phi}^g_{c_k}(x_k).\label{yk}
\end{align}
\editt{An illustration of the global embedding at various stages of Algo.~\ref{algo:geinit} is provided in Figure~\ref{fig:gpa1}.}

\subsection{Tearing Closed Manifolds}
\label{subsec:tear}

When the manifold has no boundary, then the step R2 in above section may result in a set $\mathcal{Z}_s$ containing the indices of the views which are neighbors of the $s$th view in the ambient space but are far apart from the transformed $s$th view in the embedding space, obtained right after step R1. For example, as shown in Figure~\ref{fig:ge_algo_new2} (f.2), $s_1 \in \mathcal{Z}_{s_9}$ because the $s_9$th view and the $s_1$th view are neighbors in the ambient space (see Figure~\ref{fig:ge_algo_new1} (a.1, a.2)) but in the embedding space, they are far apart. Due to such indices in $\mathcal{Z}_{s_9}$, the step R3 results in a centroid, which when used in step R4, results in a fallacious estimation of the parameters $T_{s}$ and $v_s$, giving rise to a high distortion embedding. By trying to align with all its neighbors in the ambient space, the $s_9$th view is misaligned with respect to all of them (see Figure~\ref{fig:ge_algo_new2} (g.2)).


\paragraph{Resolution ($\text{to\_tear}=$ True).} We modify the step R2 so as to \edit{introduce a discontinuity by including} the indices of only those views in the set $\mathcal{Z}_s$ which are neighbors of the $s$th view in both the ambient space as well as in the embedding space. We denote the overlap between the $m$th and $m'$th view in the embedding space by $\widetilde{U}^{g}_{mm'}$. There may be multiple heuristics for computing $\widetilde{U}^{g}_{mm'}$ which could work. In the Appendix~\ref{sec:Nmg}, we describe a simple approach based on the already developed machinery in this paper, which uses the hyperparameter $\nu$ provided as input to Algo.~\ref{algo:geinit}. Having obtained $\widetilde{U}^{g}_{mm'}$, we say that the $m$th and the $m'$th intermediate views in the embedding space are neighbors if $\widetilde{U}^g_{mm'}$ is non-empty.


\paragraph{Step R2.} Finally, we compute $\mathcal{Z}_s$ as,
\begin{align}
    \mathcal{Z}_s = \{m'\ |\ \widetilde{U}_{sm'}\neq \emptyset , \; \widetilde{U}^{g}_{sm'}\neq \emptyset \} \cap \mathcal{A}. \label{Zs2}
\end{align}
Note that if it is known apriori that the manifold can be embedded in lower dimension without tearing it apart then we do not require the above modification. In all of our experiments except the one in Section~\ref{subsec:h_d_d}, we do not assume that this information is available.


With this modification, the set $\mathcal{Z}_{s_9}$ in Figure~\ref{fig:ge_algo_new2} (f.2) will not include $s_1$ and therefore the resulting centroid in the step R3 would be the same as the one in Figure~\ref{fig:ge_algo_new2} (f.1). Subsequently, the transformed $s_9$th view would be the one in Figure~\ref{fig:ge_algo_new2} (g.1) rather than Figure~\ref{fig:ge_algo_new2} (g.2).


\paragraph{Gluing instruction for the boundary of the embedding.} Having knowingly torn the manifold apart, we provide at the output, information on the points belonging to the tear and their neighboring points in the ambient space. To encode the ``gluing" instructions along the tear in the form of colors at the output of our algorithm, we recompute $\widetilde{U}^g_{mm'}$. If $\widetilde{U}_{mm'}$ is non-empty but $\widetilde{U}^g_{mm'}$ is empty, then this means that the $m$th and $m'$th views are neighbors in the ambient space but are torn apart in the embedding space. Therefore, we color the global embedding of the points on the overlap \edit{$\widetilde{U}_{mm'}$ which belong to clusters $\mathcal{C}_m$ and $\mathcal{C}_{m'}$} with the same color to indicate that although these points are separated in the embedding space, they are adjacent in the ambient space (see Figures~\ref{fig:fig63},~\ref{fig:fig64} and~\ref{fig:decipher}).

\begin{figure}[h!]
    \centering
    \begin{tabular}{M{0.005\textwidth}|M{0.15\textwidth}|M{0.15\textwidth}|M{0.15\textwidth}|M{0.15\textwidth}|M{0.15\textwidth}|}
        & \tiny{Input} & \multicolumn{3}{c}{\tiny{First iteration of the outer loop and stages within inner loop}}  & \tiny{End of outer loop}\\
        \hline
        & & \tiny{Before} & \tiny{Half-way} & \tiny{End} & \\
        \hline
        \rotatebox{90}{\tiny{Sphere}} & \includegraphics[width=0.15\textwidth,height=0.175\textwidth,keepaspectratio]{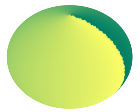} & \includegraphics[width=0.15\textwidth,height=0.175\textwidth,keepaspectratio]{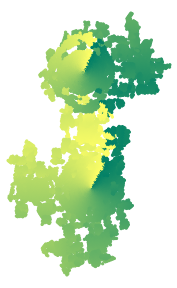} & \includegraphics[width=0.15\textwidth,height=0.175\textwidth,keepaspectratio]{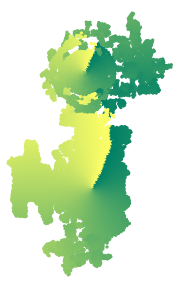} & \includegraphics[width=0.15\textwidth,height=0.175\textwidth,keepaspectratio]{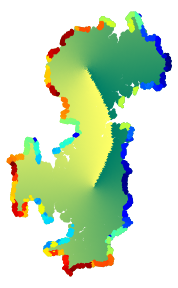} & \includegraphics[width=0.15\textwidth,height=0.175\textwidth,keepaspectratio]{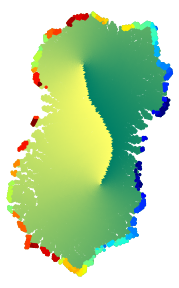} \\
        \hline
    \end{tabular}
    \caption{\editt{$2$d embedding of a sphere at different stages of Algo.~\ref{algo:geinit}. For illustration purpose, in the plots in the $2$nd and $3$rd columns the translation parameter $v_m$ was manually set for those views which do not lie in the set $\mathcal{A}$.}}
    \label{fig:gpa2}
\end{figure}

\editt{An illustration of the global embedding at various stages of Algo.~\ref{algo:geinit} with modified step R2, is provided in Figure~\ref{fig:gpa2}.}

\begin{figure}[h]
    \centering
    \begin{tabular}{cc}
        \includegraphics[width=0.25\textwidth,height=0.25\textwidth,keepaspectratio]{common/Square/final_jet.png} & \includegraphics[width=0.25\textwidth,height=0.25\textwidth,keepaspectratio]{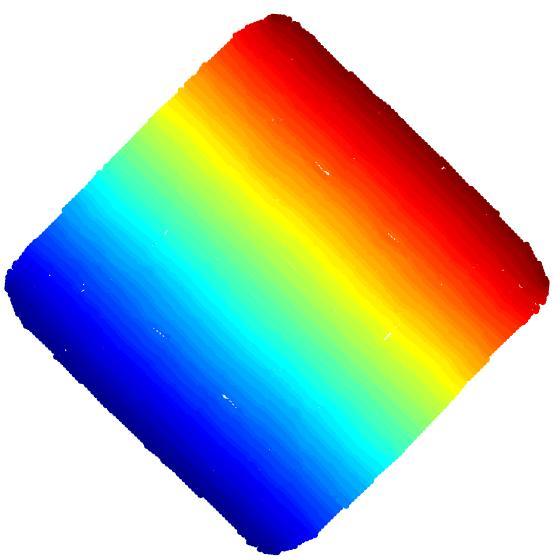}
    \end{tabular}
    \caption{Global embedding of the square grid when the points on the boundary are unknown (left) versus when they are known apriori (right).}
    \label{fig:fig51}
\end{figure}

\paragraph{Example.} The obtained global embeddings of our square grid with $\text{to\_tear}=\text{True}$ and $\nu = 3$, are shown in Figure~\ref{fig:fig51}. Note that the boundary of the obtained embedding is more distorted when the points on the boundary are unknown than when they are known apriori. This is because the intermediate views near the boundary have higher distortion in the former case than in the latter case (see Figure~\ref{fig:fig41}).

\subsection{Time Complexity}
\label{subsec:gpa_time}
\edit{The worst case time complexity of Algo.~\ref{algo:geinit} is $O(N_r nk_{\text{lv}}^2d^2/\eta_{\text{min}})$ when $\text{to\_tear}$ is false. It costs an additional time of $O(N_rn^2 \text{max}(d,k_{\text{lv}}\log n, n/\eta_{\text{min}}^2)))$ when $\text{to\_tear}$ is true. In practice, one refinement step took about $15$ seconds in the above example and between $15$-$20$ seconds for all the examples in Section~\ref{sec:compare}.}

\section{Experimental Results}
\label{sec:compare}
We present experiments to compare LDLE\footnote{The python code is available at \url{https://github.com/chiggum/pyLDLE}} with LTSA \citep{zhang2003nonlinear}, UMAP \citep{mcinnes2018umap}, t-SNE \citep{maaten2008visualizing} and Laplacian eigenmaps \citep{belkin2003laplacian} on several datasets. \edit{First, we compare the embeddings of discretized $2$d manifolds embedded in $\mathbb{R}^2$, $\mathbb{R}^3$ or $\mathbb{R}^4$, containing about $10^4$ points. These manifolds are grouped based on the presence of the boundary and their orientability as in Sections~\ref{subsec:m_w_b},~\ref{subsec:m_wo_b} and~\ref{subsec:n_o_m}. The inputs are shown in the figures themselves except for the flat torus and the Klein bottle, as their $4$D parameterizations cannot be plotted. Therefore, we describe their construction below. A quantitative comparison of the algorithms is provided in Section~\ref{subsubsec:quant_comp}. In Section~\ref{subsubsec:robust} we assess the robustness of these algorithms to the noise in the data. In Section~\ref{subsubsec:sparse} we assess the performance of these algorithms on sparse data. Finally, in Section~\ref{subsec:h_d_d} we compare the embeddings of some high dimensional datasets.}


\textbf{Flat Torus}.
A flat torus is a parallelogram whose opposite sides are identified. In our case, we construct a discrete flat torus using a rectangle with sides $2$ and $0.5$ and embed it in four dimensions as follows,
\begin{align}
    X(\theta_i,\phi_j) &= \frac{1}{4\pi}(4cos(\theta_i),4\sin(\theta_i),\cos(\phi_j),\sin(\phi_j)) \label{ftinput}
\end{align}
where $\theta_i = 0.01i\pi$, $\phi_j = 0.04j\pi$, $i \in \{0,\ldots,199\}$ and $j \in \{0,\ldots,49\}$.


\textbf{Klein bottle}.
A Klein bottle is a non-orientable two dimensional manifold without boundary. We construct a discrete Klein bottle using its $4$D M\"{o}bius tube representation as follows,
\begin{align}
    X(\theta_i,\phi_j) &= (R(\phi_j)\cos\theta_i, R(\phi_j)\sin\theta_i, r\sin\phi_j\cos\frac{\theta_i}{2}, r\sin\phi_j\sin\frac{\theta_i}{2})\\ \label{kbinput}
    R(\phi_j)&=R+r\cos\phi_j
\end{align}
where $\theta_i = i\pi/100$, $\phi_j = j\pi/25$, $i \in \{0,\ldots,199\}$ and $j \in \{0,\ldots,49\}$.

\subsection{Hyperparameters}
\label{subsec:hyperparam}
To embed using LDLE, we use the Euclidean metric and the default values of the hyperparameters and their description are provided in Table~\ref{tab:default}. \edit{Only the value of $\eta_{\text{min}}$ is tuned across all the examples in Sections~\ref{subsec:m_w_b},~\ref{subsec:m_wo_b} and~\ref{subsec:n_o_m} (except for Section~\ref{subsubsec:sparse}), and is provided in Appendix~\ref{sec:hyperparam}. For high dimensional datasets in Section~\ref{subsec:h_d_d}, values of the hyperaparameters which differ from the default values are again provided in Appendix~\ref{sec:hyperparam}.}

\begin{table}[]
    \small
    \centering
    \begin{tabular}{|c|p{0.725\linewidth}|c|}
    \hline
    \makecell{Hyper-\\parameter} & \makecell{Description}  & \makecell{Default\\value}\\
    \hline
    $k_{\text{nn}}$ & No. of nearest neighbors used to construct the graph Laplacian & $49$\\
    \hline
    $k_\text{tune}$ & The nearest neighbor, distance to which is used as a local scaling factor in the construction of graph Laplacian & 7\\
    \hline
    $N$ & No. of nontrivial low frequency Laplacian eigenvectors to consider for the construction of local views in the embedding space & 100\\
    \hline
    $d$ & Intrinsic dimension of the underlying manifold & 2\\
    \hline
    $p$ & Probability mass for computing the bandwidth $t_k$ of the heat kernel & 0.99\\
    \hline
    $k_{\text{lv}}$ & The nearest neighbor, distance to which is used to construct local views in the ambient space & 25\\
    \hline
    $(\tau_s)_{s=1}^{d}$ & Percentiles used to restrict the choice of candidate eigenfunctions & 50\\
    \hline
    $(\delta_s)_{s=1}^{d}$ & Fractions used to restrict the choice of candidate eigenfunctions & 0.9\\
    \hline
    $\eta_{\text{min}}$ & Desired minimum number of points in a cluster & 5\\
    \hline
    $\text{to\_tear}$ & A boolean for whether to tear the manifold or not & True\\
    \hline
    $\nu$ & A relaxation factor to compute the neighborhood graph of the intermediate views in the embedding space & 3\\
    \hline
    $N_r$ & No. of iterations to refine the global embedding & 100\\
    \hline
    \end{tabular}
    \caption{Default values of LDLE hyperparameters.}
    \label{tab:default}
\end{table}


For UMAP, LTSA, t-SNE and Laplacian eigenmaps, we use the Euclidean metric and select the hyperparameters by grid search, choosing the values which result in best visualization quality. \edit{For LTSA, we search for optimal n\_neighbors in $\{5,10,25,50,75,100\}$}. For UMAP, we use $500$ epochs and search for optimal n\_neighbors in $\{25,50,100,200\}$ and min\_dist in $\{0.01,0.1,0.25,0.5\}$. For t-SNE, we use $1000$ iterations and search for optimal perplexity in $\{30,40,50,60\}$ and early exaggeration in $\{2,4,6\}$. For Laplacian eigenmaps, we search for $k_\text{nn}$ in $\{16,25,36,49\}$ and $k_{\text{tune}}$ in $\{3,7,11\}$. The chosen values of the hyperparameters are provided in Appendix~\ref{sec:hyperparam}. We note that the Laplacian eigenmaps fails to correctly embed most of the examples regardless of the choice of the hyperparameters.

\subsection{Manifolds with Boundary}
\label{subsec:m_w_b}
\begin{figure}[h!]
    \centering
    \begin{tabular}{M{0.005\textwidth}|M{0.15\textwidth}|M{0.15\textwidth}|M{0.15\textwidth}|M{0.15\textwidth}|M{0.15\textwidth}|}
        & \tiny{Barbell} & \tiny{Square with two holes} & \tiny{Sphere with a hole} & \tiny{Swiss Roll with a hole} & \tiny{Noisy Swiss Roll}\\
        \hline
        \rotatebox{90}{\tiny{Input}} & \includegraphics[width=0.15\textwidth,height=0.175\textwidth,keepaspectratio]{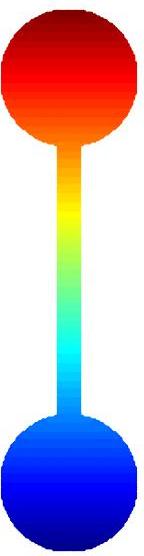} & \includegraphics[width=0.15\textwidth,height=0.175\textwidth,keepaspectratio]{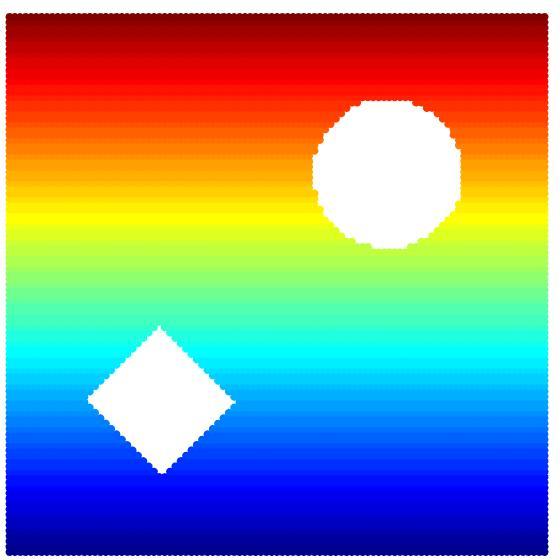} & \includegraphics[width=0.15\textwidth,height=0.175\textwidth,keepaspectratio]{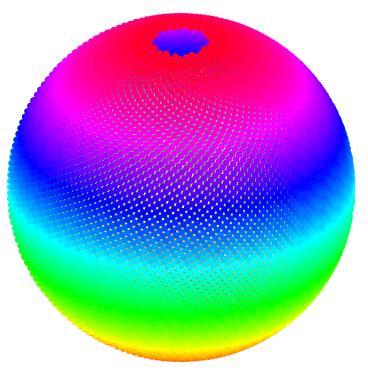} & \includegraphics[width=0.15\textwidth,height=0.175\textwidth,keepaspectratio]{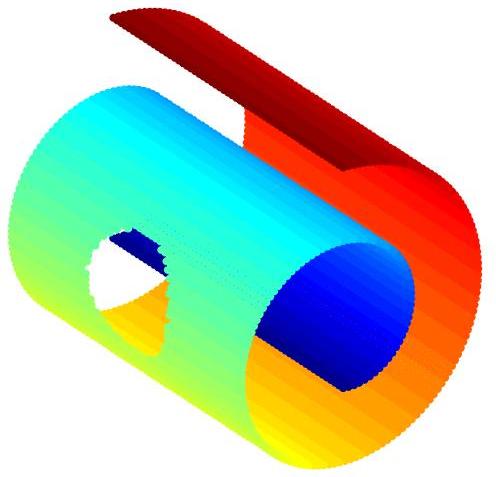} & \includegraphics[width=0.15\textwidth,height=0.175\textwidth,keepaspectratio]{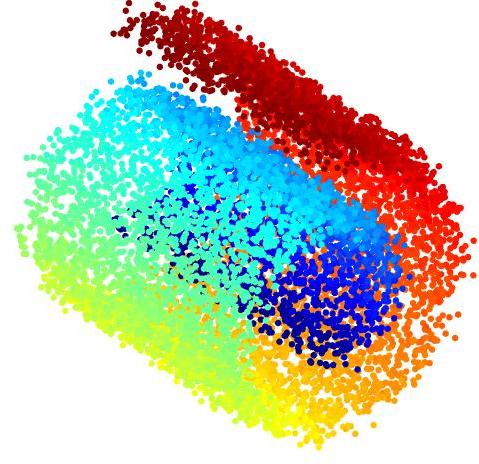} \\
        \hline
        \rotatebox{90}{\tiny{LDLE}} & \rotatebox{90}{\includegraphics[width=0.175\textwidth,height=0.15\textwidth,keepaspectratio]{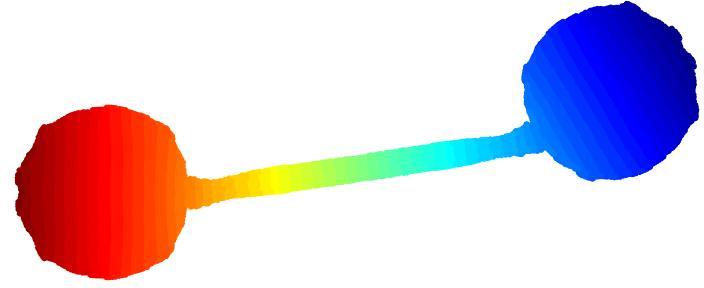}} & \includegraphics[width=0.15\textwidth,height=0.175\textwidth,keepaspectratio]{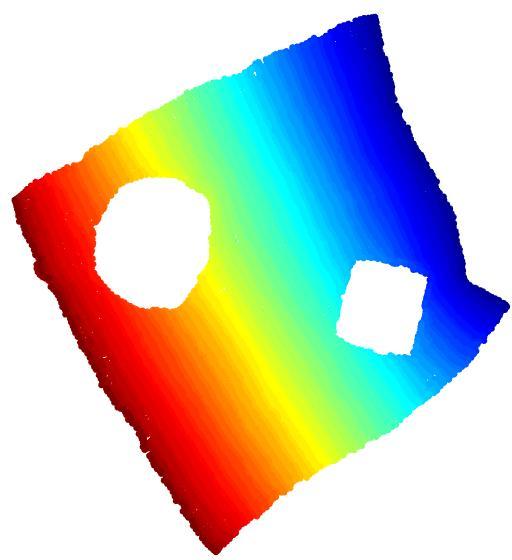} & \includegraphics[width=0.15\textwidth,height=0.175\textwidth,keepaspectratio]{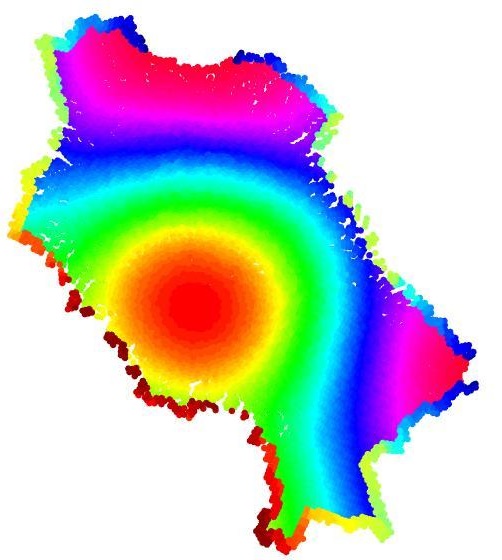} & \includegraphics[width=0.15\textwidth,height=0.175\textwidth,keepaspectratio]{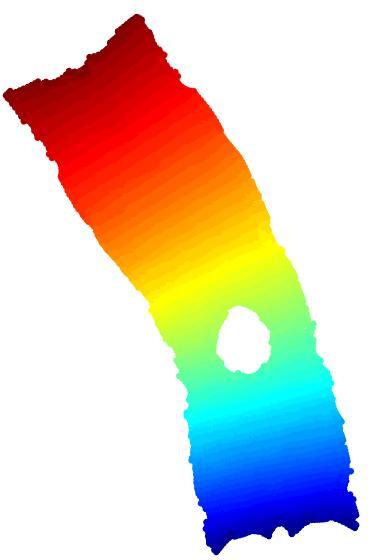} & \includegraphics[width=0.15\textwidth,height=0.175\textwidth,keepaspectratio]{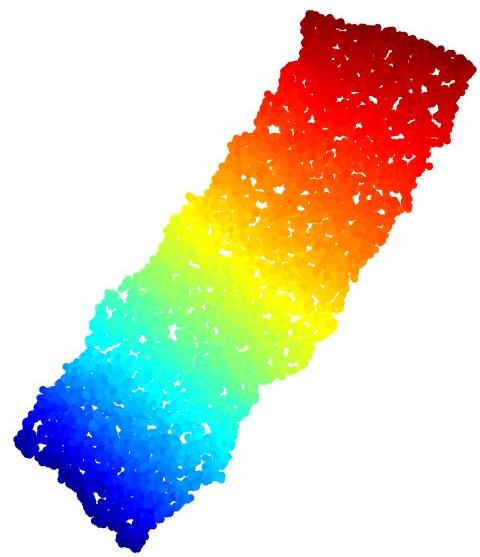} \\
        \hline
        \rotatebox{90}{\tiny{LTSA}} & \rotatebox{90}{\includegraphics[width=0.175\textwidth,height=0.15\textwidth,keepaspectratio]{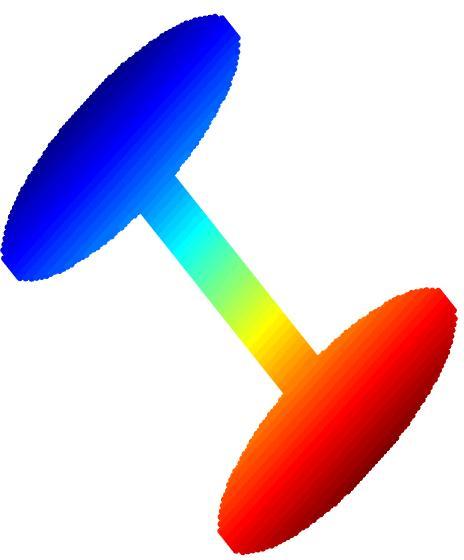}} & \includegraphics[width=0.15\textwidth,height=0.175\textwidth,keepaspectratio]{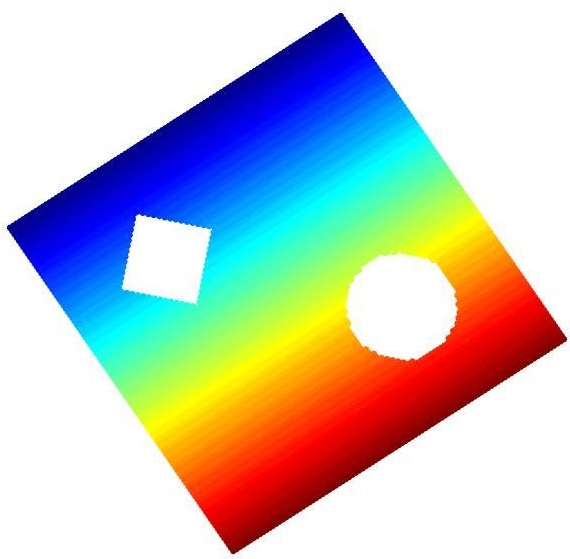} & \includegraphics[width=0.15\textwidth,height=0.175\textwidth,keepaspectratio]{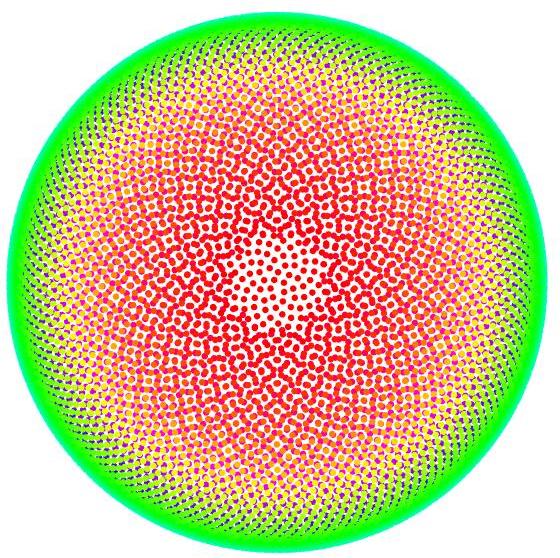} & \includegraphics[width=0.15\textwidth,height=0.175\textwidth,keepaspectratio]{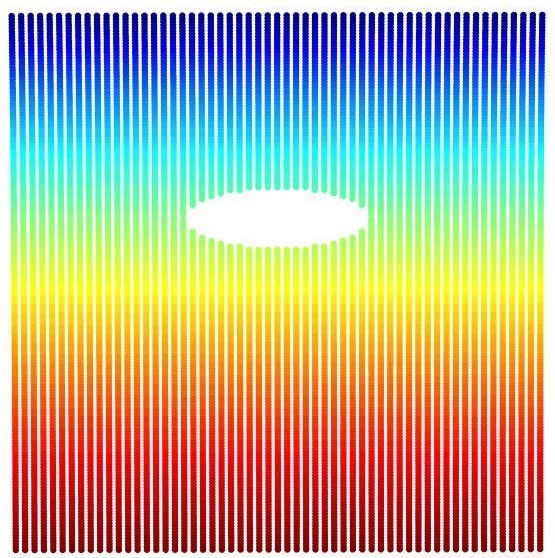} & \includegraphics[width=0.15\textwidth,height=0.175\textwidth,keepaspectratio]{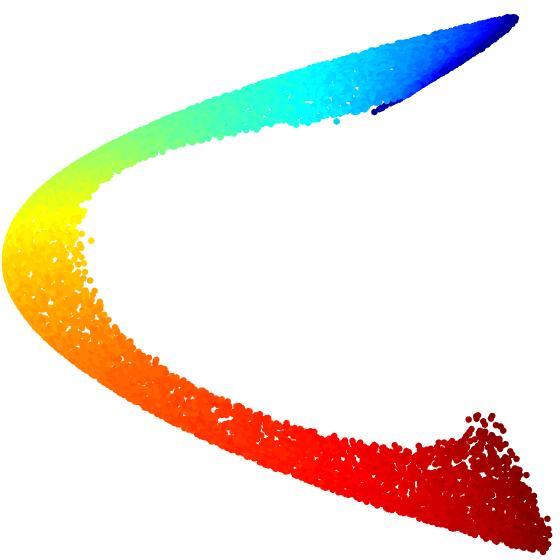} \\
        \hline
        \rotatebox{90}{\tiny{UMAP}} & \includegraphics[width=0.15\textwidth,height=0.175\textwidth,keepaspectratio]{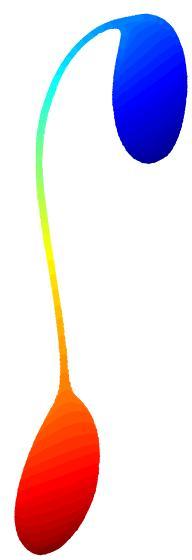} & \includegraphics[width=0.15\textwidth,height=0.175\textwidth,keepaspectratio]{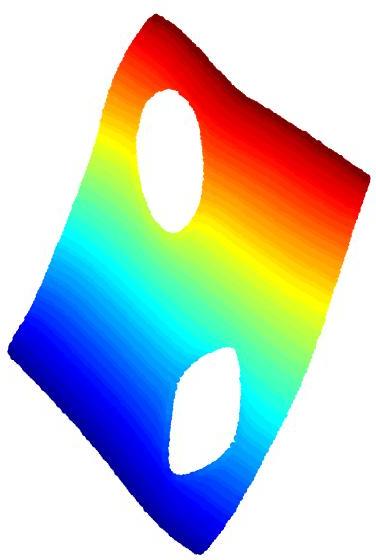} & \includegraphics[width=0.15\textwidth,height=0.175\textwidth,keepaspectratio]{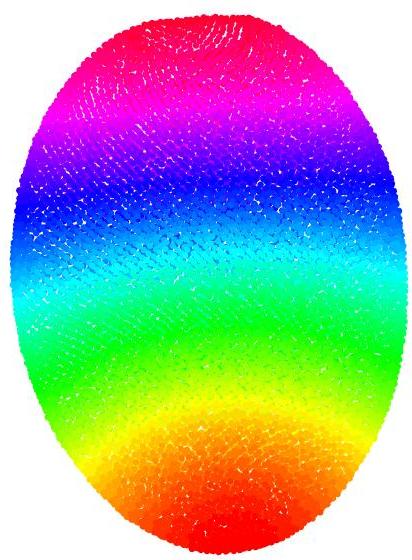} & \includegraphics[width=0.15\textwidth,height=0.175\textwidth,keepaspectratio]{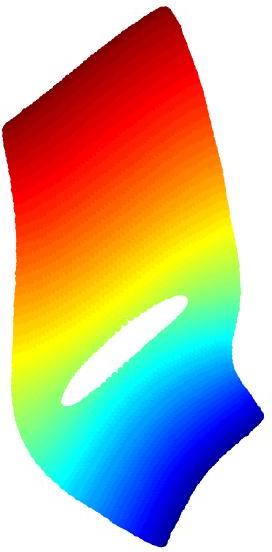} & \includegraphics[width=0.15\textwidth,height=0.175\textwidth,keepaspectratio]{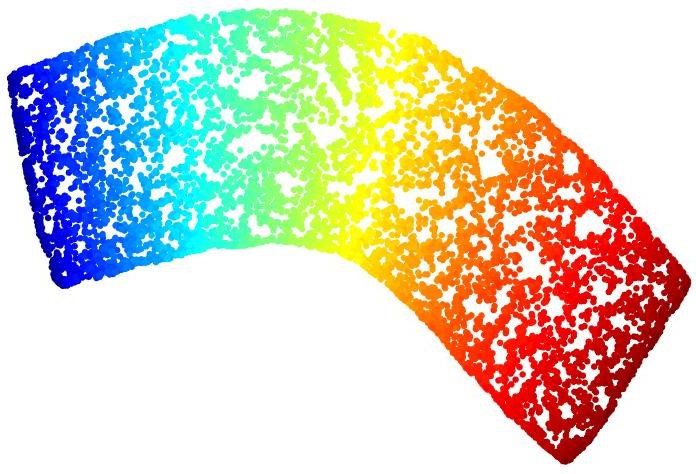}  \\
        \hline
        \rotatebox{90}{\tiny{t-SNE}} & \includegraphics[width=0.15\textwidth,height=0.175\textwidth,keepaspectratio]{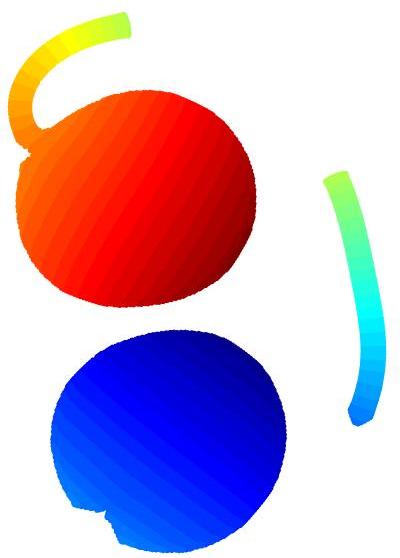} & \includegraphics[width=0.15\textwidth,height=0.175\textwidth,keepaspectratio]{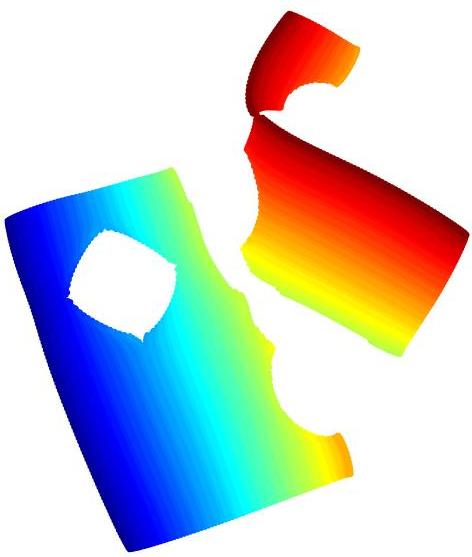} & \rotatebox{90}{\includegraphics[width=0.175\textwidth,height=0.15\textwidth,keepaspectratio]{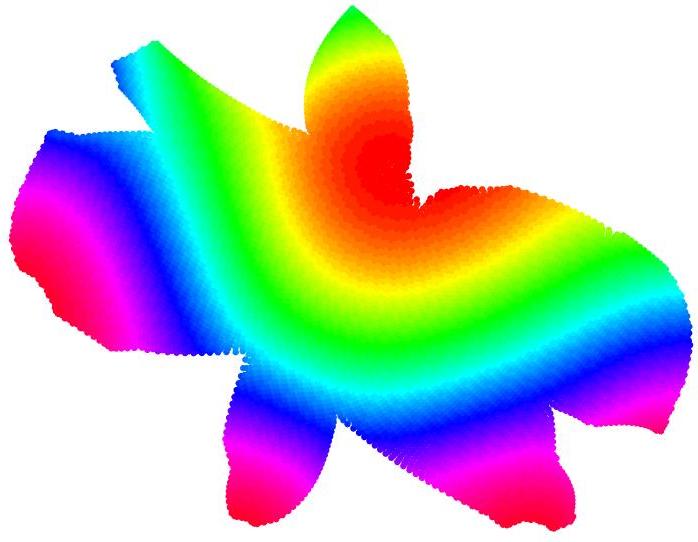}} & \includegraphics[width=0.15\textwidth,height=0.175\textwidth,keepaspectratio]{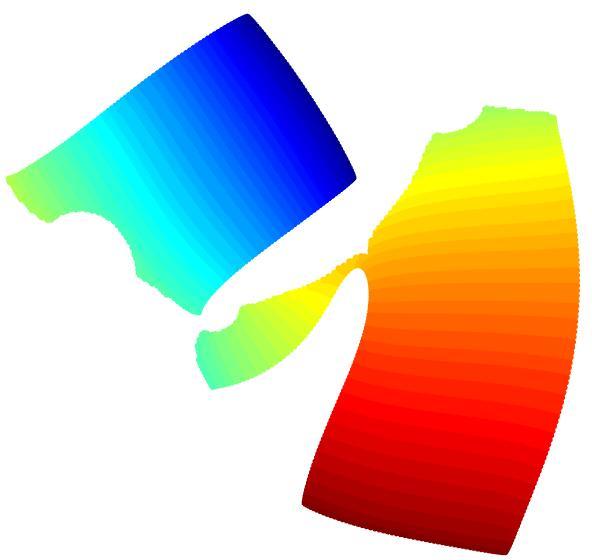} & \includegraphics[width=0.15\textwidth,height=0.175\textwidth,keepaspectratio]{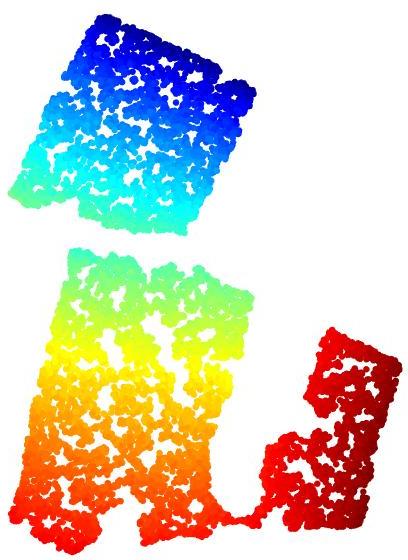}  \\
        \hline
        \rotatebox{90}{\tiny{Laplacian Eigenmaps}} & \includegraphics[width=0.15\textwidth,height=0.175\textwidth,keepaspectratio]{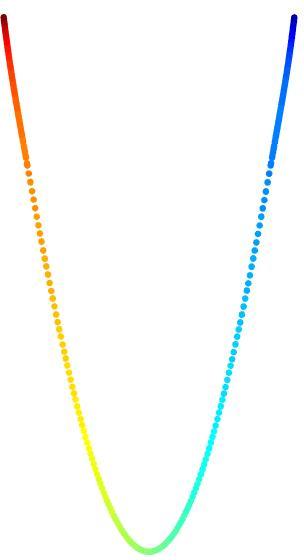} & \includegraphics[width=0.15\textwidth,height=0.175\textwidth,keepaspectratio]{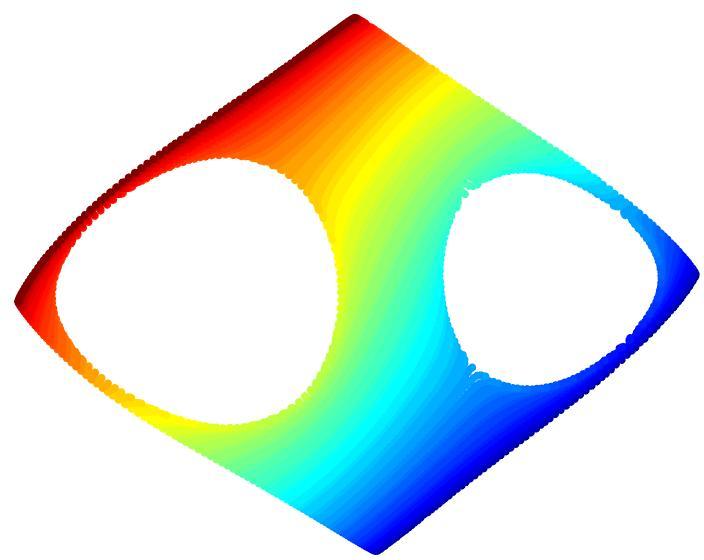} & \includegraphics[width=0.15\textwidth,height=0.175\textwidth,keepaspectratio]{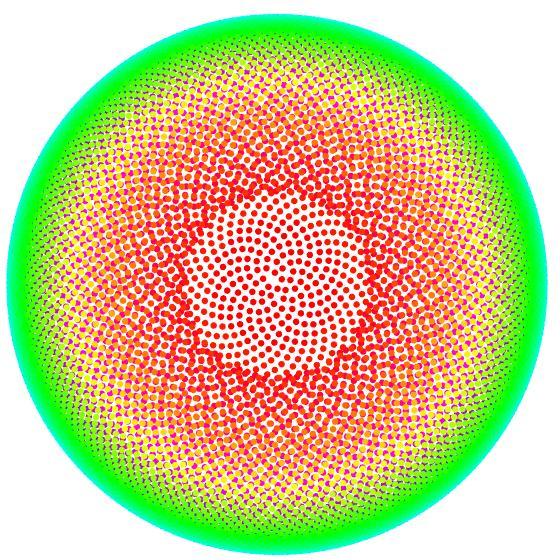} & \includegraphics[width=0.15\textwidth,height=0.175\textwidth,keepaspectratio]{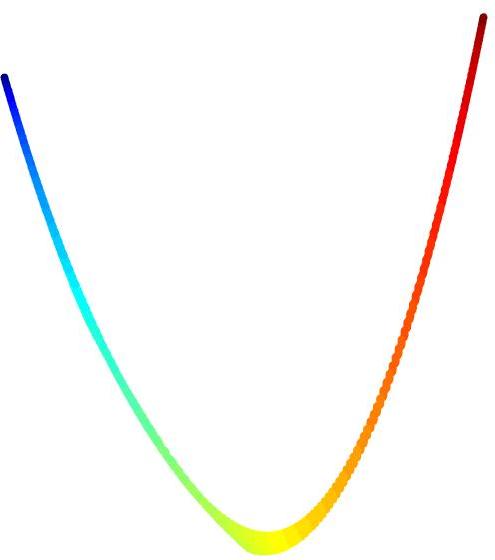} & \includegraphics[width=0.15\textwidth,height=0.175\textwidth,keepaspectratio]{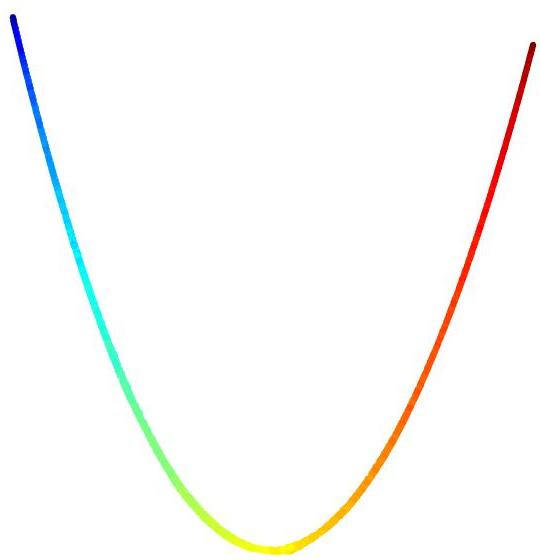} \\
        \hline
    \end{tabular}
    \caption{Embeddings of $2$d manifolds with boundary into $\mathbb{R}^2$. The noisy Swiss Roll is constructed by adding uniform noise in all three dimensions, with support on $[0,0.05]$.}
    \label{fig:fig61}
\end{figure}

\edit{In Figure~\ref{fig:fig61}, we show the $2$d embeddings of $2$d manifolds with boundary, in $\mathbb{R}^2$ or $\mathbb{R}^3$, three of which have holes. To a large extent, LDLE preserved the shape of the holes. 
LTSA perfectly preserved the shape of the holes in the square but deforms it in the Swiss Roll. This is because LTSA embedding does not capture the aspect ratio of the underlying manifold as discussed in Section~\ref{subsec:ltsa_stitch}. UMAP and Laplacian eigenmaps distorted the shape of the holes and the region around them, while t-SNE produced dissected embeddings. For the sphere with a hole which is a curved $2$d manifold with boundary, LTSA, UMAP and Laplacian eigenmaps squeezed it into $\mathbb{R}^2$ while LDLE and t-SNE tore it apart. The correctness of the LDLE embedding is proved in Figure~\ref{fig:decipher}. In the case of noisy swiss roll, LDLE and UMAP produced visually better embeddings in comparison to the other methods.}

\begin{figure}[h!]
    \centering
    \begin{tabular}{M{0.005\textwidth}|M{0.2\textwidth}|M{0.2\textwidth}|M{0.2\textwidth}|}
        & \tiny{Barbell} & \tiny{Square with two holes} & \tiny{Swiss Roll with a hole}\\
        \hline
        \rotatebox{90}{\tiny{LDLE with $\partial\mathcal{M}$ known apriori}} & \includegraphics[width=0.2\textwidth,height=0.25\textwidth,keepaspectratio]{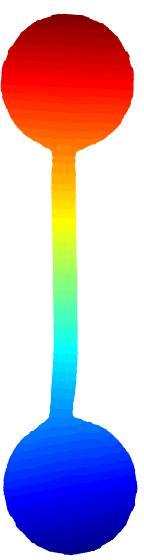} & \includegraphics[width=0.2\textwidth,height=0.25\textwidth,keepaspectratio]{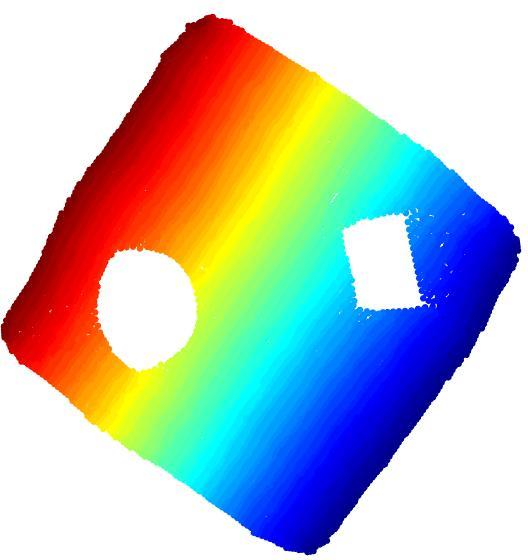} & \rotatebox{90}{\includegraphics[width=0.25\textwidth,height=0.2\textwidth,keepaspectratio]{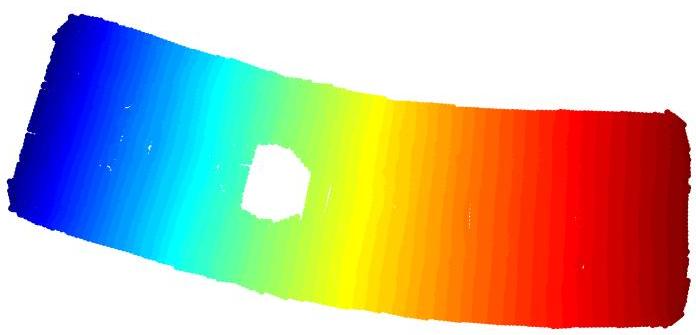}} \\
        \hline
    \end{tabular}
    \caption{LDLE embeddings when the points on the boundary are known apriori.}
    \label{fig:fig62}
\end{figure}

We note that the boundaries of the LDLE embeddings in Figure~\ref{fig:fig61} are usually distorted. The cause of this is explained in Remark~\ref{rmk:highdist}. When the points in the input which lie on the boundary are known apriori then the distortion near the boundary can be reduced using the double manifold as discussed in Remark~\ref{rmk:dblmanifold} and shown in Figure~\ref{fig:fig21}. The obtained LDLE embeddings when the points on the boundary are known, are shown in Figure~\ref{fig:fig62}.

\subsubsection{Quantitative comparison}
\label{subsubsec:quant_comp}
\begin{figure}[h!]  
    \centering
    \footnotesize
    \begin{tabular}{|c|c|c|c|c|}
        \hline
         \makecell{\tiny{Rectangle} \tiny{($4 \times 0.25$)}} & \tiny{barbell} & \makecell{\tiny{Square with two}\\\tiny{holes}} & \makecell{\tiny{Swiss Roll with a}\\\tiny{hole}} &\tiny{Noisy Swiss Roll} \\
         \hline
        \includegraphics[trim=0 0 0 50, clip, width=0.15\textwidth,keepaspectratio]{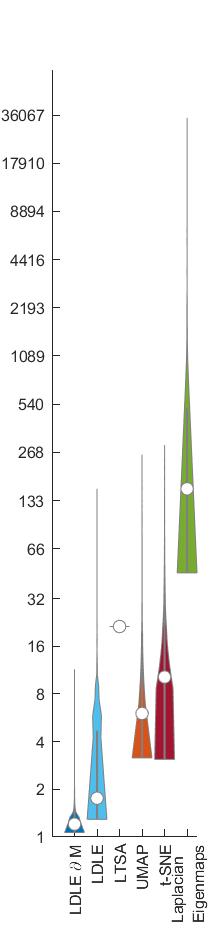}& \includegraphics[trim=0 0 0 50, clip,width=0.15\textwidth,keepaspectratio]{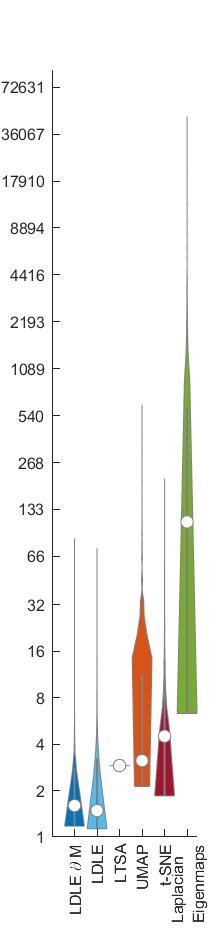}& \includegraphics[trim=0 0 0 50, clip,width=0.15\textwidth,keepaspectratio]{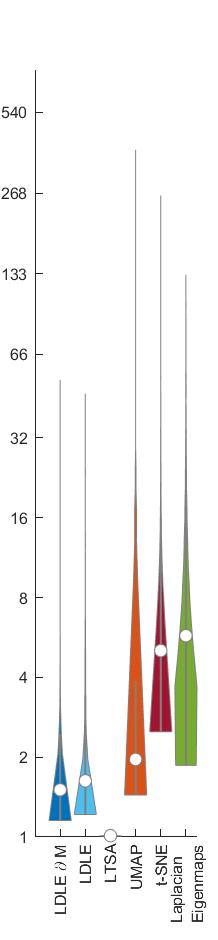}& \includegraphics[trim=0 0 0 50, clip,width=0.15\textwidth,keepaspectratio]{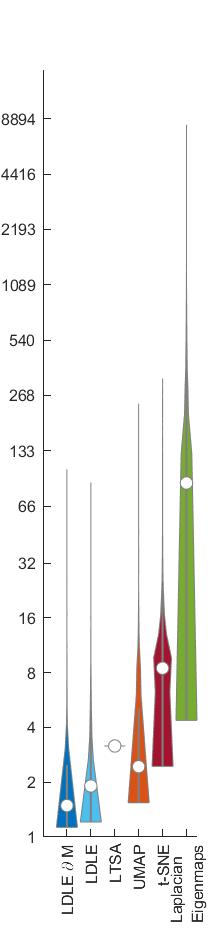}& \includegraphics[trim=0 0 0 50, clip,width=0.1275\textwidth,keepaspectratio]{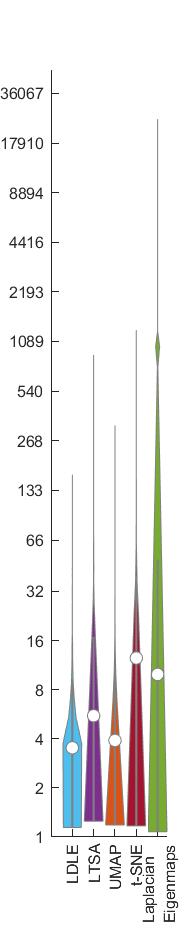}\\
        \hline
    \end{tabular}
    \caption{Violin plots \citep{hintze1998violin, bastian_bechtold_2021_4559847} for the distribution of $\mathcal{D}_k$ (See Eq.~(\ref{Di})). LDLE $\partial M$ means LDLE with boundary known apriori. The white point inside the violin represents the median. The straight line above the end of the violin represents the outliers.}
    \label{fig:fig612}
\end{figure}

\edit{
To compare LDLE with other techniques in a quantitative manner, we compute the distortion $\mathcal{D}_k$ of the embeddings of the geodesics originating from $x_k$ and then plot the distribution of $\mathcal{D}_k$ (see Figure~\ref{fig:fig612}). The procedure to compute $\mathcal{D}_k$ follows. In the discrete setting, we first define the geodesic between two given points as the shortest path between them which in turn is computed by running Dijkstra algorithm on the graph of $5$ nearest neighbors. Here, the distances are measured using the Euclidean metric $d_e$. Denote the number of nodes on the geodesic between $x_{k}$ and $x_{k'}$ by $n_{kk'}$ and the sequence of nodes by $(x_s)_{s=1}^{n_{kk'}}$ where $x_1 = x_k$ and $x_{n_{kk'}} = x_{k'}$. Denote the embedding of $x_k$ by $y_k$. Then the length of the geodesic in the latent space between $x_k$ and $x_{k'}$, and the length of the embedding of the geodesic between $y_k$ and $y_{k'}$ are given by
}
\begin{align}
    L_{kk'} &= \sum_{s=2}^{n_{kk'}}d_e(x_{s},x_{s-1}).
\\
    L^g_{kk'} &= \sum_{s=2}^{n_{kk'}}d_e(y_s, y_{s-1}).
\end{align}
\edit{
Finally, the distortion $\mathcal{D}_k$ of the embeddings of the geodesics originating from $x_k$ is given by the ratio of maximum expansion and minimum contraction, that is,
}
\begin{align}
    \mathcal{D}_{k} &= \sup_{k'}\frac{L^g_{kk'}}{L_{kk'}}/\inf_{k'}\frac{L^g_{kk'}}{L_{kk'}} = \sup_{k'}\frac{L^g_{kk'}}{L_{kk'}}\sup_{k'}\frac{L_{kk'}}{L^g_{kk'}}.\label{Di}
\end{align}
\edit{
A value of $1$ for $\mathcal{D}_k$ means the geodesics originating from $x_k$ have the same length in the input and in the embedding space. If $\mathcal{D}_k = 1$ for all $k$ then the embedding is geometrically, and therefore topologically as well, the same as the input up to scale. Figure~\ref{fig:fig612} shows the distribution of $\mathcal{D}_k$ due to LDLE and other algorithms for various examples. Except for the noisy Swiss Roll, LTSA produced the least maximum distortion. Specifically, for the square with two holes, LTSA produced a distortion of $1$ suggesting its strength on manifolds with unit aspect ratio. In all other examples, LDLE produced the least distortion except for a few outliers. When the boundary is unknown, the points which result in high $\mathcal{D}_k$ are the ones which lie on and near the boundary. When the boundary is known, these are the points which lie on or near the corners (see Figures~\ref{fig:fig21} and~\ref{fig:fig31}). We aim to fix this issue in future work. }

\subsubsection{Robustness to noise}
\label{subsubsec:robust}

\edit{To further analyze the robustness of LDLE under noise we compare the embeddings of the Swiss Roll with Gaussian noise of increasing variance. The resulting embeddings are shown in Figure~\ref{fig:noise}. Note that certain points on LDLE embeddings have a different colormap than the one used for the input. As explained in Section~\ref{subsec:tear}, the points which have the same color under this colormap are adjacent on the manifold but away in the embedding. To be precise, these points lie close to the middle of the gap in the Swiss Roll, creating a bridge between those points which would otherwise be far away on a noiseless Swiss Roll. In a sense, these points cause maximum corruption to the geometry of the underlying noiseless manifold. One can say that these points are have adversarial noise, and LDLE embedding can automatically recognize such points. We will further explore this in future work. LTSA, t-SNE and Laplacian Eigenmaps fail to produce correct embeddings while UMAP embeddings also exhibit high quality.}

\subsubsection{\editt{Sparsity}}
\label{subsubsec:sparse}
A comparison of the embeddings of the Swiss Roll with decreasing resolution and increasing sparsity is provided in Figure~\ref{fig:sparse}. Unlike LTSA and Laplacian Eigenmaps, the embeddings produced by LDLE, UMAP and t-SNE are of high quality. Note that when the resolution is $10$, LDLE embedding of some points have a different colormap. Due to sparsity, certain points on the opposite sides of the gap between the Swiss Roll are neighbors in the ambient space as shown in Figure~\ref{fig:res10issue} in Appendix I. LDLE automatically tore apart these erroneous connections and marked them at the output using a different colormap. A discussion on sample size requirement for LDLE follows.


The distortion of LDLE embeddings directly depend on the distortion of the constructed local parameterizations, which in turn depends on reliable estimates of the graph Laplacian and its eigenvectors. The work in \citep{BELKIN20081289, hein2007graph, trillos2020error, cheng2021eigen} provided conditions on the sample size and the hyperparameters such as the kernel bandwidth, under which the graph Laplacian and its eigenvectors would converge to their continuous counterparts. A similar analysis in the setting of self-tuned kernels used in our approach (see Algo.~\ref{algo:gl}) is also provided in \citep{cheng2020convergence}. These imply that, for a faithful estimation of graph Laplacian and its eigenvectors, the hyperparameter $k_{\text{tune}}$ (see Table~\ref{tab:default}) should be small enough so that the local scaling factors $\sigma_k$ (see Algo.~\ref{algo:gl}) are also small, while the size of the data $n$ should be large enough so that $n \sigma_k^{d+2}/\log(n)$ is sufficiently large for all $k \in \{1,\ldots,n\}$. This suggests that $n$ needs to be exponential in $d$ and inversely related to $\sigma_k$. However, in practice, the data is usually given and therefore $n$ is fixed. So the above mainly states that to obtain accurate estimates, the hyperparameter $k_\text{tune}$ must be decreased. This indeed holds as we had to decrease $k_\text{tune}$ from $7$ to $2$ (see Appendix~\ref{sec:hyperparam}) to produce LDLE embeddings of high quality for increasingly sparse Swiss Roll in Figure~\ref{fig:sparse}.




\subsection{Closed Manifolds}
\label{subsec:m_wo_b}

In Figure~\ref{fig:fig63}, we show the $2$d embeddings of $2$d manifolds without a boundary, a curved torus in $\mathbb{R}^3$ and a flat torus in $\mathbb{R}^4$. LDLE produced similar representation for both the inputs. None of the other methods do that. \editt{The main difference in the LDLE embedding of the two inputs is based on the boundary of the embedding. It is composed of many small line segments for the flat torus, and many small curved segments for the curved torus.} This is clearly because of the difference in the curvature of the two inputs, zero everywhere for the flat torus and non-zero almost everywhere on the curved torus. The mathematical correctness of the LDLE embeddings using the cut and paste argument is shown in Figure~\ref{fig:decipher}. LTSA, UMAP and Laplacian eignemaps squeezed both the manifolds into $\mathbb{R}^2$ while the t-SNE embedding is non-interpretable.

\subsection{Non-Orientable Manifolds}
\label{subsec:n_o_m}

In Figure~\ref{fig:fig64}, we show the $2$d embeddings of non-orientatble $2$d manifolds, a M\"obius strip in $\mathbb{R}^3$ and a Klein bottle in $\mathbb{R}^4$. Laplacian eigenmaps produced incorrect embeddings, t-SNE produced dissected and non-interpretable embeddings and LTSA and UMAP squeezed the inputs into $\mathbb{R}^2$. LDLE produced mathematically correct embeddings by tearing apart both inputs to embed them into $\mathbb{R}^2$ (see Figure~\ref{fig:decipher}).

\subsection{High Dimensional Data}
\label{subsec:h_d_d}

\subsubsection{Synthetic sensor data}
\label{subsubsec:ssd}
In Figure~\ref{fig:sensor}, motivated from \citep{peterfreund2020loca}, we embed a $42$ dimensional synthetic data set representing the signal strength of $42$ transmitters at about $n = 6000$ receiving locations on a toy floor plan. The transmitters and the receivers are distributed uniformly across the floor. Let $(t_{r_k})_{k=1}^{42}$ be the transmitter locations and $r_i$ be the $i$th receiver location. Then the $i$th data point $x_i$ is given by $(e^{-\left\|r_i-t_{r_k}\right\|_2^2})_{k=1}^{42}$. The resulting data set is embedded using and other algorithms into $\mathbb{R}^2$. The hyperparameters resulting in the most visually appealing embeddings were identified for each algorithm and are provided in Table~\ref{tab:hypparam}. The obtained embeddings are shown in Figure~\ref{fig:sensor}. \edit{The shapes of the holes are best preserved by LTSA, then LDLE followed by the other algorithms. The corners of the LDLE embedding are more distorted. The reason for distorted corners is given in Remark~\ref{rmk:highdist}.}

\subsubsection{Face image data}

\edit{In Figure~\ref{fig:face_data}, we show the embedding obtained by applying LDLE on the face image data \citep{tenenbaum2000global} which consists of a sequence of $698$ $64$-by-$64$ pixel images of a face rendered under various pose and lighting conditions. These images are converted to $4096$ dimensional vectors, then projected to $100$ dimensions through PCA  
while retaining about $98\%$ of the variance. These are then embedded using LDLE and other algorithms into $\mathbb{R}^2$. The hyperparameters resulting in the most visually appealing embeddings were identified for each algorithm and are provided in Table~\ref{tab:hypparam_face_yoda}. The resulting embeddings are shown in Figure~\ref{fig:s1_puppets} colored by the pose and lighting of the face. Note that values of the pose and lighting variables for all the images are provided in the dataset itself. We have displayed face images corresponding to few points of the LDLE embeddings as well. Embeddings due to all the techniques except LTSA reasonably capture both the pose and lighting conditions. 
}

\subsubsection{Rotating Yoda-Bulldog dataset}

\edit{
In Figure~\ref{fig:s1_puppets}, we show the $2$d embeddings of the rotating figures dataset presented in \citep{lederman2018learning}. It consists of $8100$ snapshots taken by a camera of a platform with two objects, Yoda and a bull dog, rotating at different frequencies. Therefore, the underlying $2$d parameterization of the data should render a torus. The original images have a dimension of $320 \times 240 \times 3$. In our experiment, we first resize the images to half the original size and then project them to $100$ dimensions through PCA \citep{jolliffe2016principal} while retaining about $98\%$ variance. These are then embedded using LDLE and other algorithms into $\mathbb{R}^2$. The hyperparameters resulting in the most visually appealing embeddings were identified for each algorithm and are provided in Table~\ref{tab:hypparam_face_yoda}. The resulting embeddings are shown in Figure~\ref{fig:s1_puppets} colored by the first dimension of the embedding itself. LTSA and UMAP resulted in a squeezed torus. LDLE tore apart the underlying torus and automatically colored the boundary of the embedding to suggest the gluing instructions. By tracing the color on the boundary we have manually drawn the arrows. Putting these arrows on a piece of paper and using cut and past argument one can establish that the embedding represents a torus (see Figure~\ref{fig:decipher}). The images corresponding to a few points on the boundary are shown. Pairs of images with the same labels represent the two sides of the curve along which LDLE tore apart the torus, and as is evident these pairs are similar.}

\section{Conclusion and Future Work}
\label{sec:conclusion}
We have presented a new bottom-up approach (LDLE) for manifold learning which constructs low-dimensional low distortion local views of the data using the low frequency global eigenvectors of the graph Laplacian, and registers them to obtain a global embedding. Through various examples we demonstrated  that LDLE competes with the other methods in terms of visualization quality. In particular, the embeddings produced by LDLE preserved distances upto a constant scale better than those produced by UMAP, t-SNE, Laplacian Eigenmaps \edit{and for the most part LTSA too}. We also demonstrated that LDLE is robust to the noise in the data and produces fine embeddings even when the data is sparse. We also showed that LDLE can embed closed as well as non-orientable manifolds into their intrinsic dimension, a feature that is missing from the existing techniques. Some of the future directions of our work are as follows.

\begin{figure}[h!]
    \centering
    \begin{tabular}{M{0.005\textwidth}|M{0.225\textwidth}|M{0.225\textwidth}|M{0.225\textwidth}|}
        & $\sigma=0.01$ & $\sigma=0.015$ & $\sigma=0.02$\\
        \hline
        \rotatebox{90}{\tiny{Side view of Swiss Roll}} & \includegraphics[trim=400 160 400 160, clip,width=0.225\textwidth,height=0.225\textwidth,keepaspectratio]{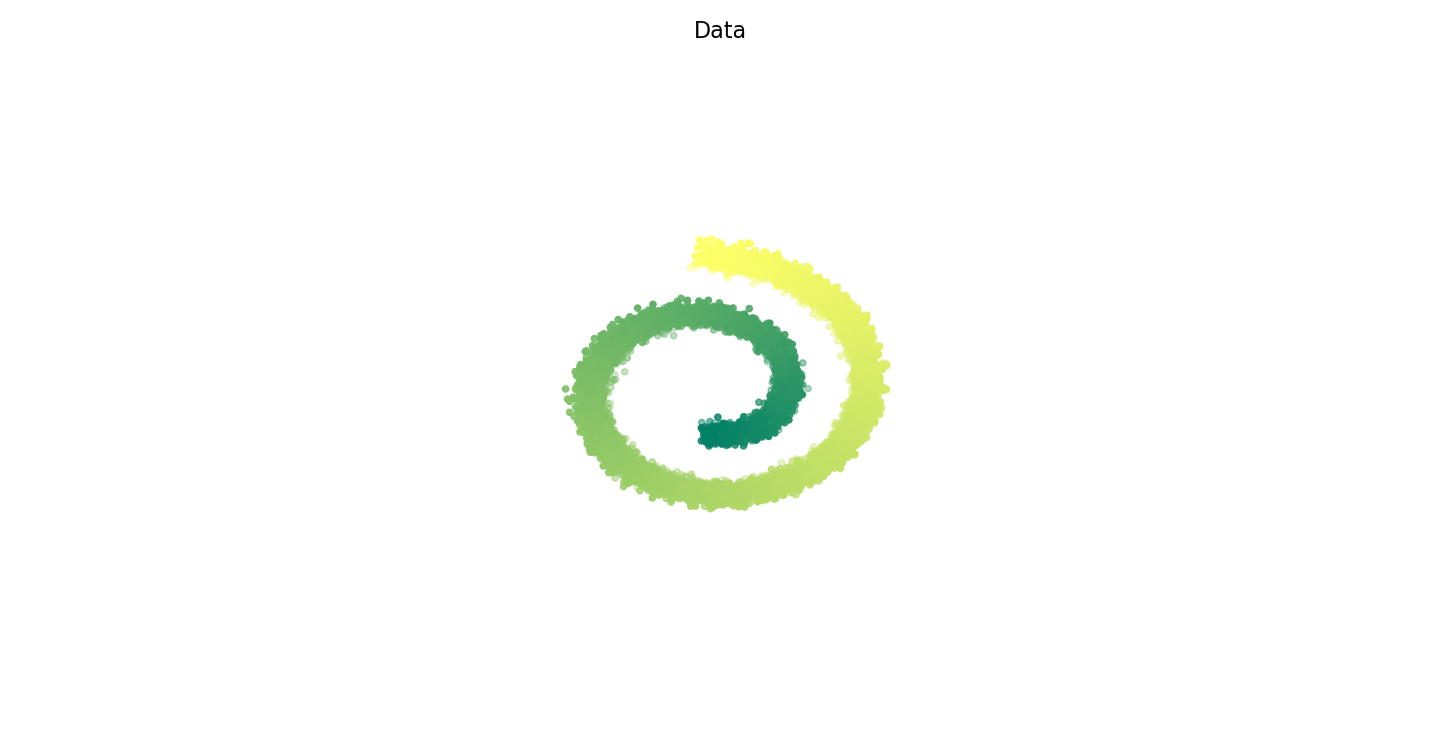} & \includegraphics[trim=400 160 400 160, clip,width=0.225\textwidth,height=0.225\textwidth,keepaspectratio]{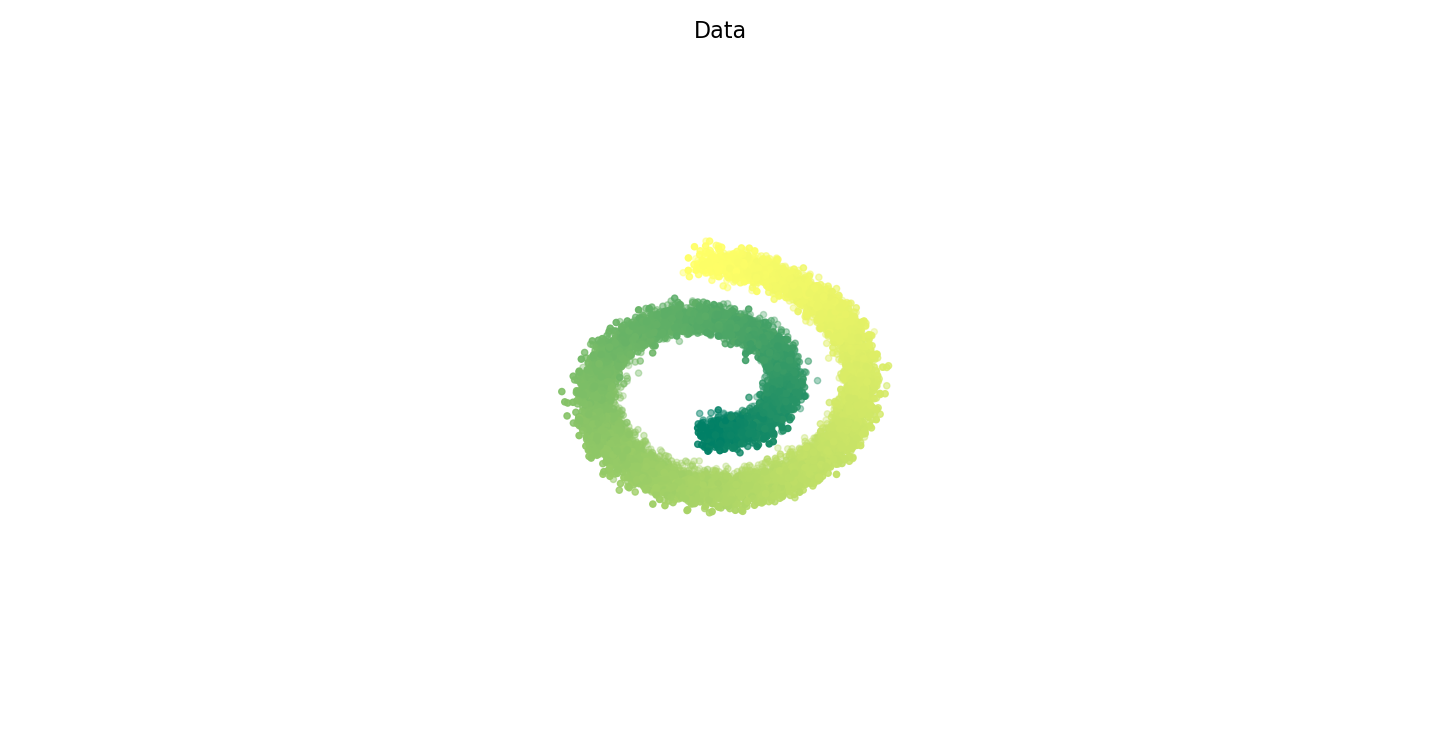} & \includegraphics[trim=400 160 400 160, clip,width=0.225\textwidth,height=0.225\textwidth,keepaspectratio]{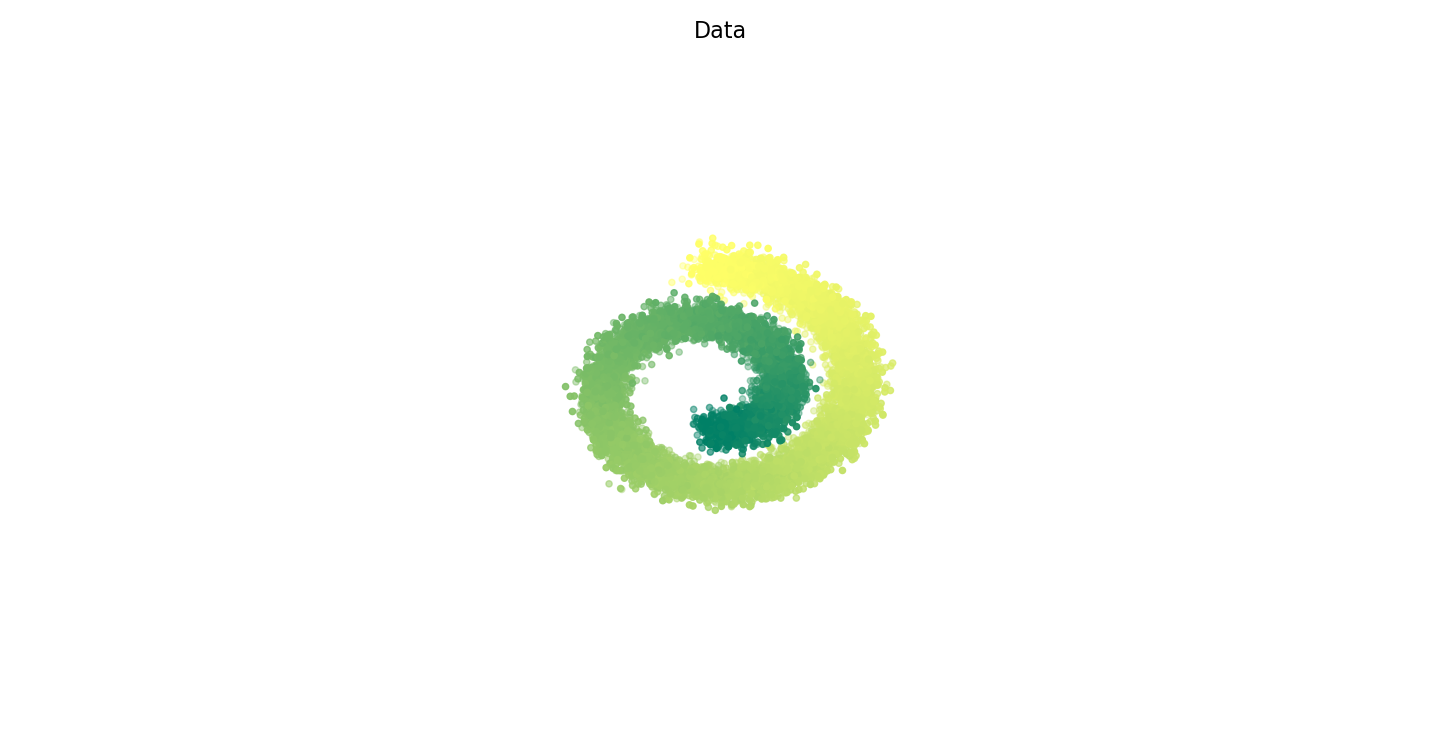}\\
        \hline
        \rotatebox{90}{\tiny{LDLE}} & \includegraphics[width=0.225\textwidth,height=0.225\textwidth,keepaspectratio]{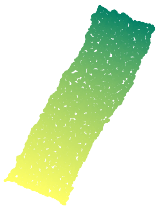} & \includegraphics[width=0.225\textwidth,height=0.225\textwidth,keepaspectratio]{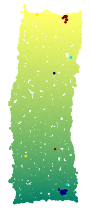} & \includegraphics[width=0.225\textwidth,height=0.225\textwidth,keepaspectratio]{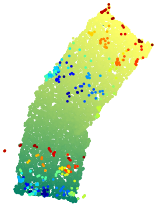}\\
        \hline
        \rotatebox{90}{\tiny{LTSA}} & \includegraphics[trim=70 35 60 40, clip,width=0.225\textwidth,height=0.225\textwidth,keepaspectratio]{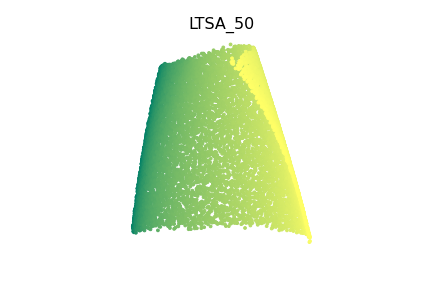} & \includegraphics[trim=70 35 60 40, clip,width=0.225\textwidth,height=0.225\textwidth,keepaspectratio]{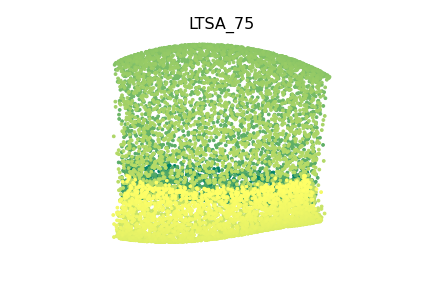} & \includegraphics[trim=70 35 60 40, clip,width=0.225\textwidth,height=0.225\textwidth,keepaspectratio]{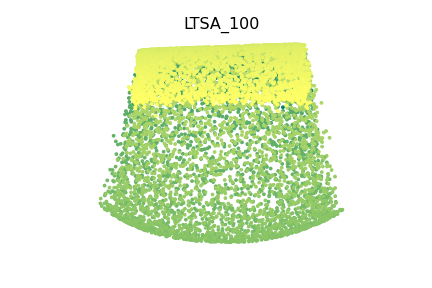}\\
        \hline
        \rotatebox{90}{\tiny{UMAP}} & \includegraphics[trim=100 40 100 50, clip,width=0.225\textwidth,height=0.225\textwidth,keepaspectratio]{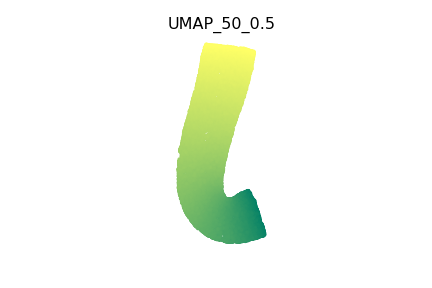} & \includegraphics[trim=100 40 100 50, clip,width=0.225\textwidth,height=0.225\textwidth,keepaspectratio]{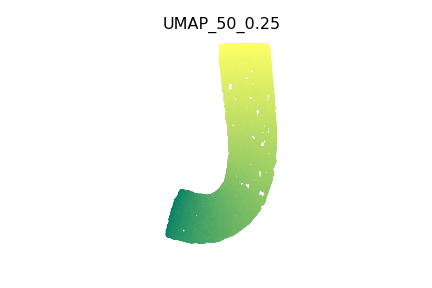} & \includegraphics[trim=100 40 100 50, clip,width=0.225\textwidth,height=0.225\textwidth,keepaspectratio]{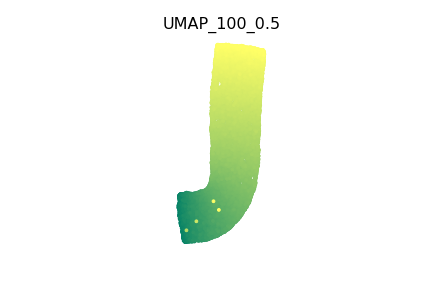}\\
        \hline
        \rotatebox{90}{\tiny{t-SNE}} & \includegraphics[trim=100 35 100 35, clip,width=0.225\textwidth,height=0.225\textwidth,keepaspectratio]{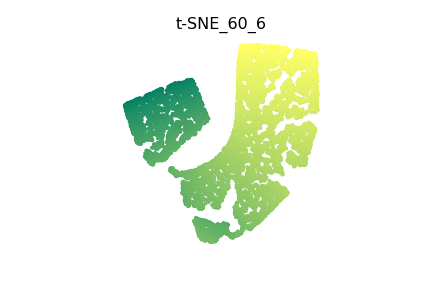} & \includegraphics[trim=100 35 100 35, clip,width=0.225\textwidth,height=0.225\textwidth,keepaspectratio]{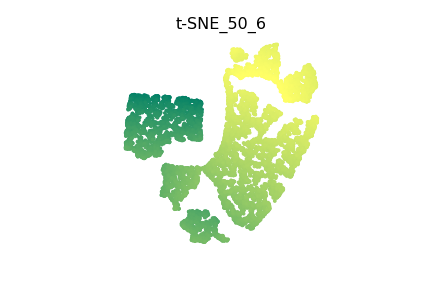} & \includegraphics[trim=100 35 80 35, clip,width=0.225\textwidth,height=0.225\textwidth,keepaspectratio]{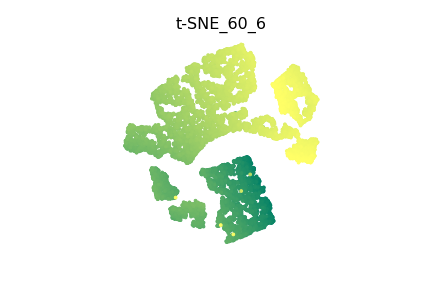}\\
        \hline
        \rotatebox{90}{\tiny{Laplacian Eigenmaps}} & \includegraphics[width=0.225\textwidth,height=0.225\textwidth,keepaspectratio]{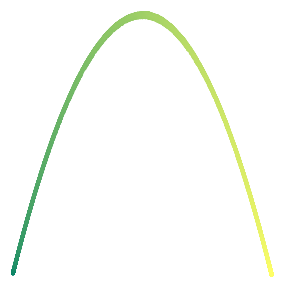} & \includegraphics[width=0.225\textwidth,height=0.225\textwidth,keepaspectratio]{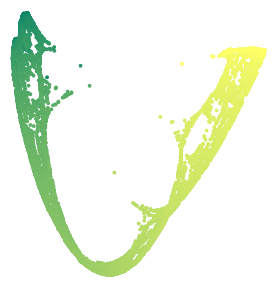} & \includegraphics[width=0.225\textwidth,height=0.225\textwidth,keepaspectratio]{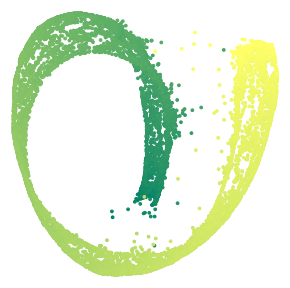}\\
        \hline
    \end{tabular}
    \caption{\editt{Embeddings of the Swiss Roll with additive noise sampled from the Gaussian distribution of zero mean and a variance of $\sigma^2$ (see Section~\ref{subsubsec:robust} for details).}}
    \label{fig:noise}
\end{figure}

\begin{figure}[h!]
    \centering
    \begin{tabular}{M{0.005\textwidth}|M{0.195\textwidth}|M{0.195\textwidth}|M{0.195\textwidth}|M{0.195\textwidth}|}
        & \tiny{RES$=30$ ($n=990$)} & \tiny{RES$=15$ ($n=280$)} & \tiny{RES$=12$ ($n=184$)} & \tiny{RES$=10$ ($n=133$)}\\
        \hline
        \rotatebox{90}{\tiny{Input}} & \includegraphics[width=0.195\textwidth,height=0.195\textwidth,keepaspectratio]{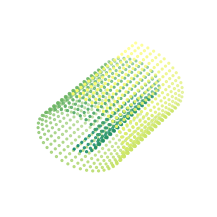} & \includegraphics[width=0.195\textwidth,height=0.195\textwidth,keepaspectratio]{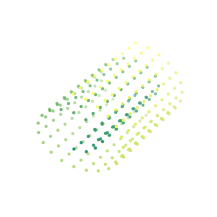} & \includegraphics[width=0.195\textwidth,height=0.195\textwidth,keepaspectratio]{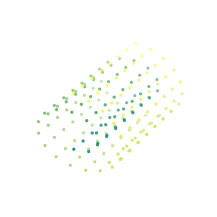} & \includegraphics[width=0.195\textwidth,height=0.195\textwidth,keepaspectratio]{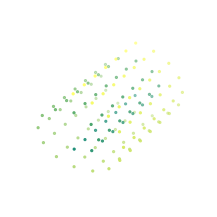}\\
        \hline
        \rotatebox{90}{\tiny{LDLE}} & \includegraphics[width=0.195\textwidth,height=0.195\textwidth,keepaspectratio]{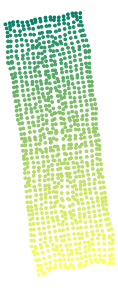} & \includegraphics[width=0.195\textwidth,height=0.195\textwidth,keepaspectratio]{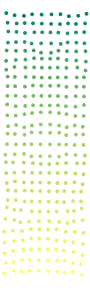} & \includegraphics[width=0.195\textwidth,height=0.195\textwidth,keepaspectratio]{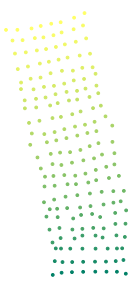} & \includegraphics[width=0.195\textwidth,height=0.195\textwidth,keepaspectratio]{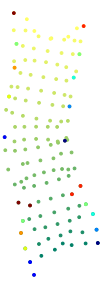}\\
        \hline
        \rotatebox{90}{\tiny{LTSA}} & \includegraphics[width=0.195\textwidth,height=0.195\textwidth,keepaspectratio]{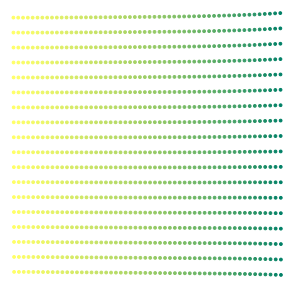} & \includegraphics[width=0.195\textwidth,height=0.195\textwidth,keepaspectratio]{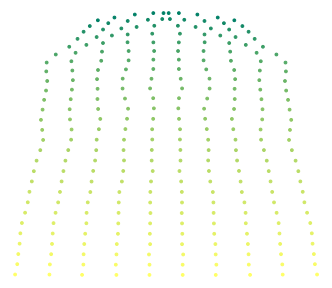} & \includegraphics[width=0.195\textwidth,height=0.195\textwidth,keepaspectratio]{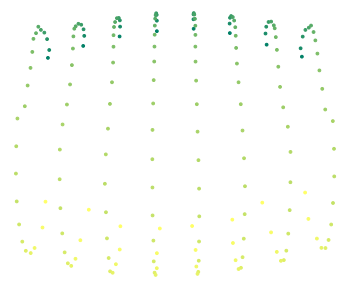} & \includegraphics[width=0.195\textwidth,height=0.195\textwidth,keepaspectratio]{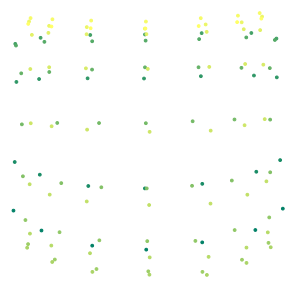}\\
        \hline
        \rotatebox{90}{\tiny{UMAP}} & \includegraphics[width=0.195\textwidth,height=0.195\textwidth,keepaspectratio]{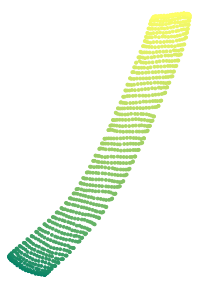} & \includegraphics[width=0.195\textwidth,height=0.195\textwidth,keepaspectratio]{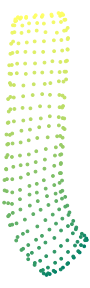} & \includegraphics[width=0.195\textwidth,height=0.195\textwidth,keepaspectratio]{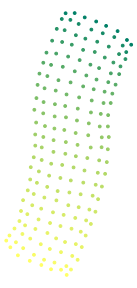} & \includegraphics[width=0.195\textwidth,height=0.195\textwidth,keepaspectratio]{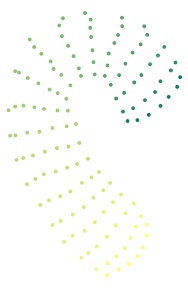}\\
        \hline
        \rotatebox{90}{\tiny{t-SNE}} & \includegraphics[,width=0.195\textwidth,height=0.195\textwidth,keepaspectratio]{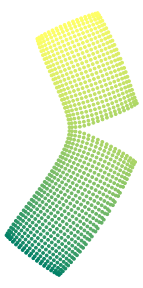} & \includegraphics[width=0.195\textwidth,height=0.195\textwidth,keepaspectratio]{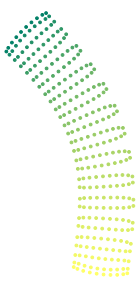} & \includegraphics[width=0.195\textwidth,height=0.195\textwidth,keepaspectratio]{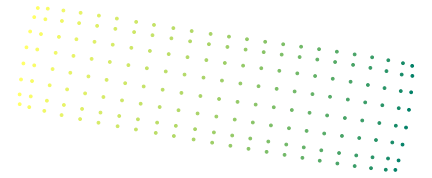} & \includegraphics[width=0.195\textwidth,height=0.195\textwidth,keepaspectratio]{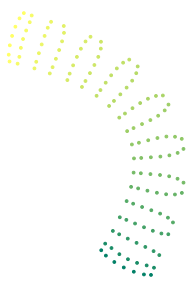}\\
        \hline
        \rotatebox{90}{\tiny{Laplacian Eigenmaps}} & \includegraphics[,width=0.195\textwidth,height=0.195\textwidth,keepaspectratio]{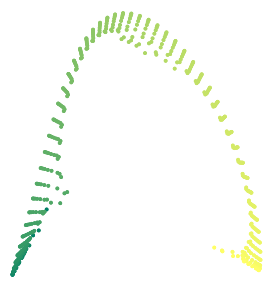} & \includegraphics[width=0.195\textwidth,height=0.195\textwidth,keepaspectratio]{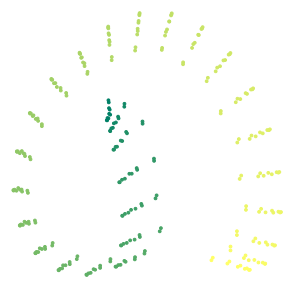} & \includegraphics[width=0.195\textwidth,height=0.195\textwidth,keepaspectratio]{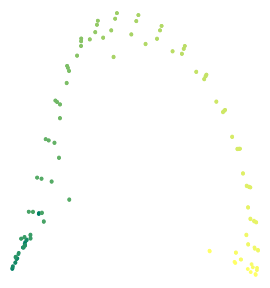} & \includegraphics[width=0.195\textwidth,height=0.195\textwidth,keepaspectratio]{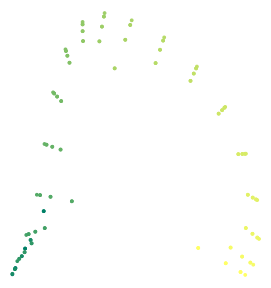}\\
        \hline
    \end{tabular}
    \caption{\editt{Embeddings of the Swiss Roll with decreasing resolution and increasing sparsity (see Section~\ref{subsubsec:sparse} for details). Note that when RES$=7$ ($n=70$) none of the above method produced a correct embedding.}}
    \label{fig:sparse}
\end{figure}

\begin{figure}[h!]
    \centering
    \begin{tabular}{M{0.005\textwidth}|M{0.2\textwidth}|M{0.2\textwidth}|M{0.2\textwidth}|M{0.2\textwidth}|}
        & \multicolumn{2}{c|}{\tiny{Curved torus}} & \multicolumn{2}{c|}{\tiny{Flat torus}}\\
        \hline 
        \rotatebox{90}{\tiny{Input}} & \includegraphics[width=0.2\textwidth,height=0.25\textwidth,keepaspectratio]{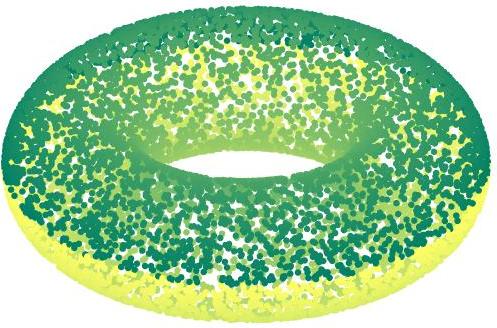} & \includegraphics[width=0.2\textwidth,height=0.25\textwidth,keepaspectratio]{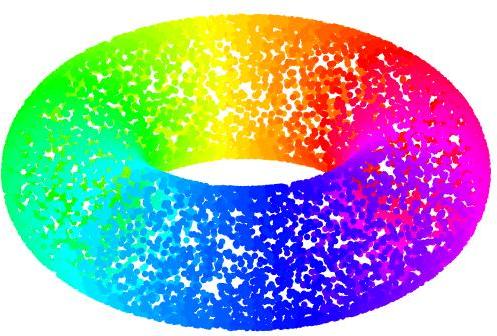}  & \multicolumn{2}{c|}{
        \footnotesize{See Eq.~(\ref{ftinput})}
        }\\
        \hline
        \rotatebox{90}{\tiny{LDLE}} & \includegraphics[width=0.2\textwidth,height=0.25\textwidth,keepaspectratio]{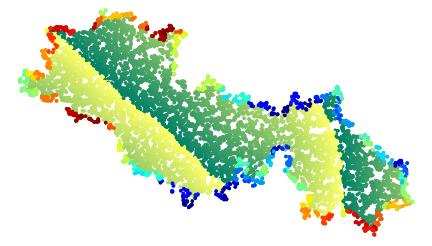} & \includegraphics[width=0.2\textwidth,height=0.25\textwidth,keepaspectratio]{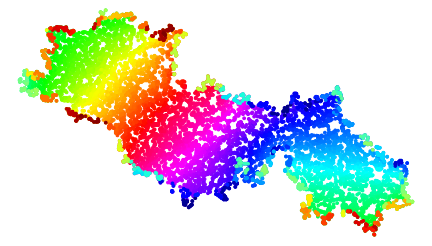} & \includegraphics[width=0.2\textwidth,height=0.25\textwidth,keepaspectratio]{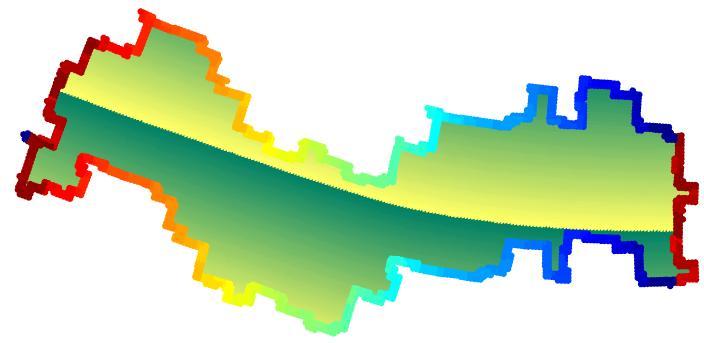} & \includegraphics[width=0.2\textwidth,height=0.25\textwidth,keepaspectratio]{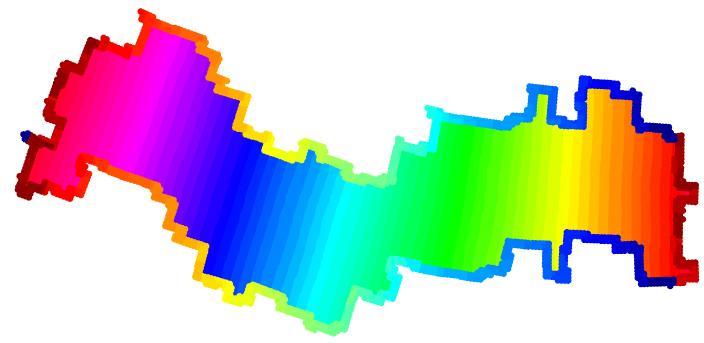} \\
        \hline
        \rotatebox{90}{\tiny{LTSA}} & \includegraphics[width=0.2\textwidth,height=0.25\textwidth,keepaspectratio]{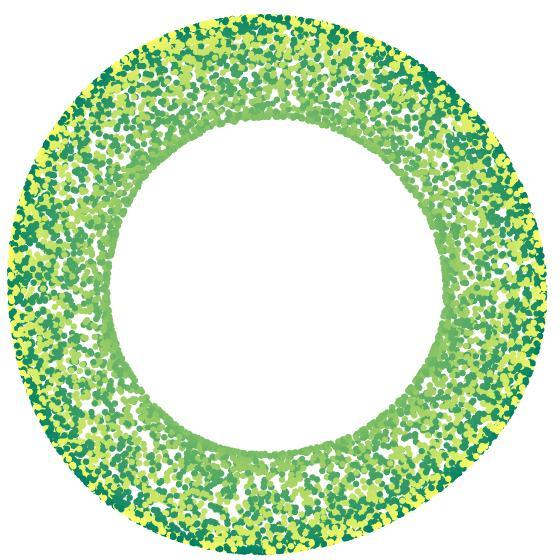} & \includegraphics[width=0.2\textwidth,height=0.25\textwidth,keepaspectratio]{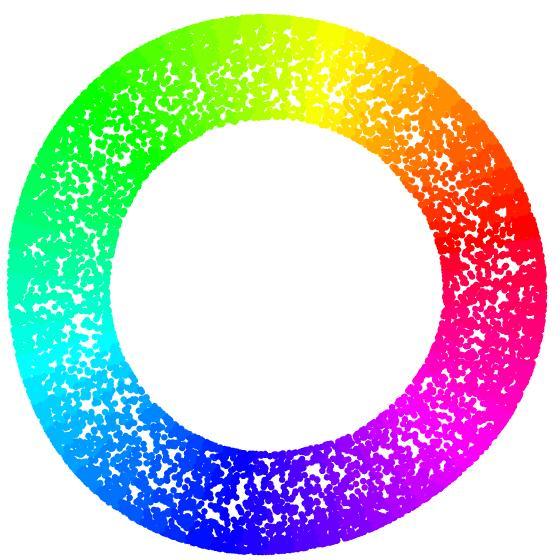} & \includegraphics[width=0.2\textwidth,height=0.25\textwidth,keepaspectratio]{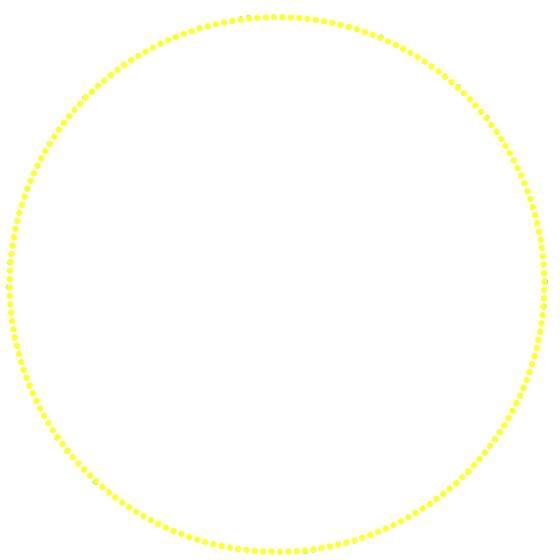} & \includegraphics[width=0.2\textwidth,height=0.25\textwidth,keepaspectratio]{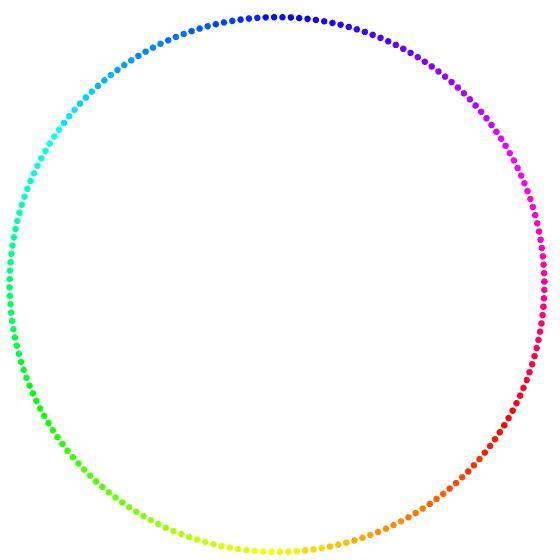} \\
        \hline
        \rotatebox{90}{\tiny{UMAP}} & \includegraphics[width=0.2\textwidth,height=0.25\textwidth,keepaspectratio]{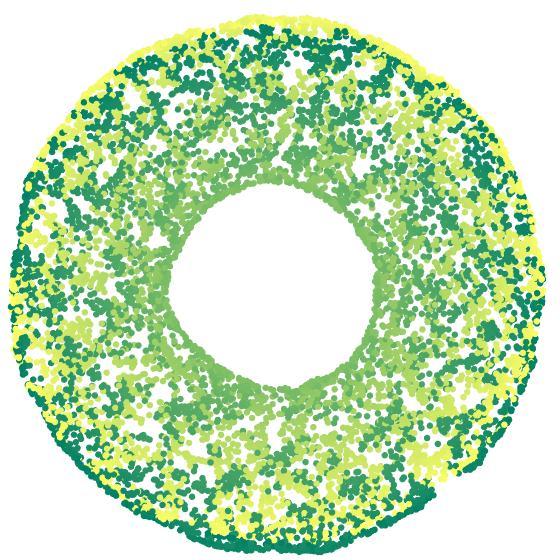} & \includegraphics[width=0.2\textwidth,height=0.25\textwidth,keepaspectratio]{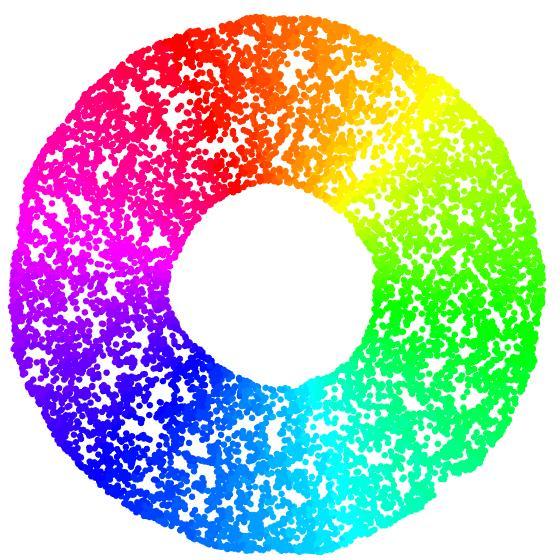} & \includegraphics[width=0.2\textwidth,height=0.25\textwidth,keepaspectratio]{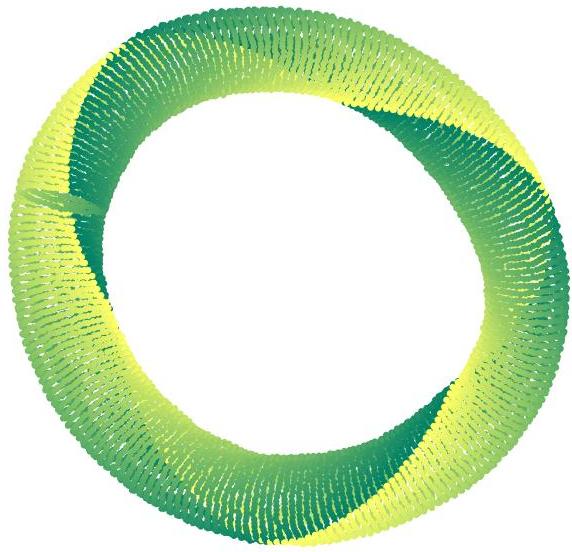} & \includegraphics[width=0.2\textwidth,height=0.25\textwidth,keepaspectratio]{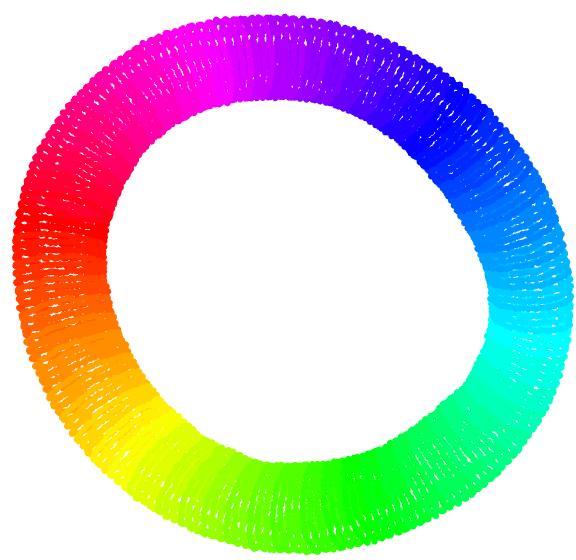} \\
        \hline
        \rotatebox{90}{\tiny{t-SNE}} & \includegraphics[width=0.2\textwidth,height=0.25\textwidth,keepaspectratio]{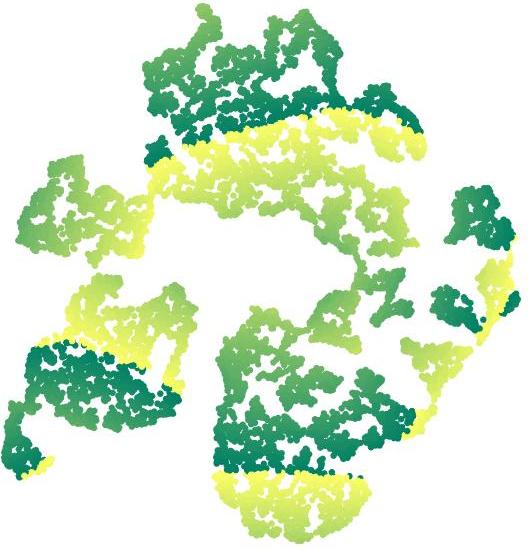} & \includegraphics[width=0.2\textwidth,height=0.25\textwidth,keepaspectratio]{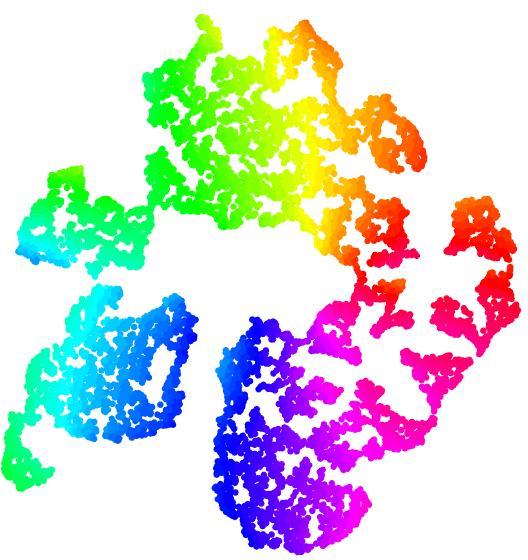} & \includegraphics[width=0.2\textwidth,height=0.25\textwidth,keepaspectratio]{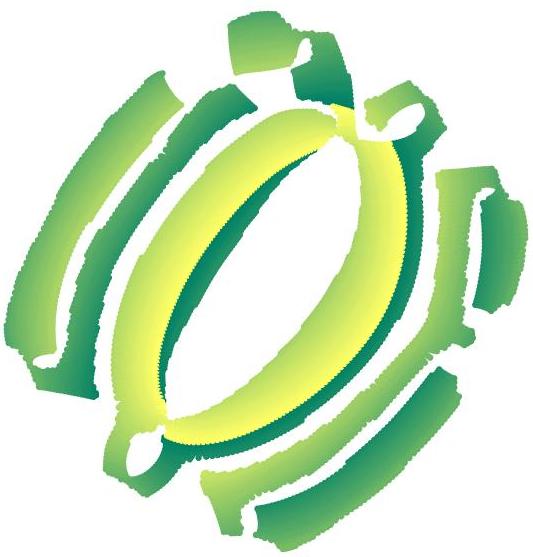} & \includegraphics[width=0.2\textwidth,height=0.25\textwidth,keepaspectratio]{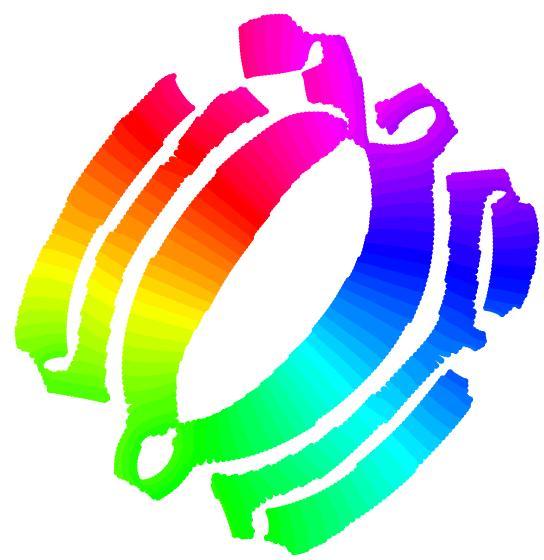} \\
        \hline
        \rotatebox{90}{\tiny{Laplacian Eigenmaps}} & \includegraphics[width=0.2\textwidth,height=0.25\textwidth,keepaspectratio]{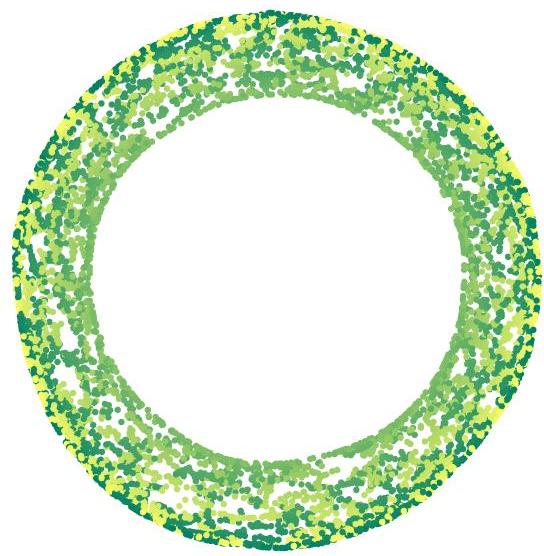} & \includegraphics[width=0.2\textwidth,height=0.25\textwidth,keepaspectratio]{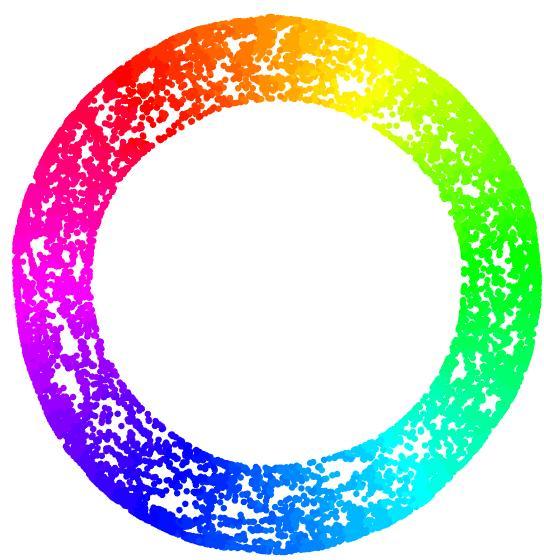} & \includegraphics[width=0.2\textwidth,height=0.25\textwidth,keepaspectratio]{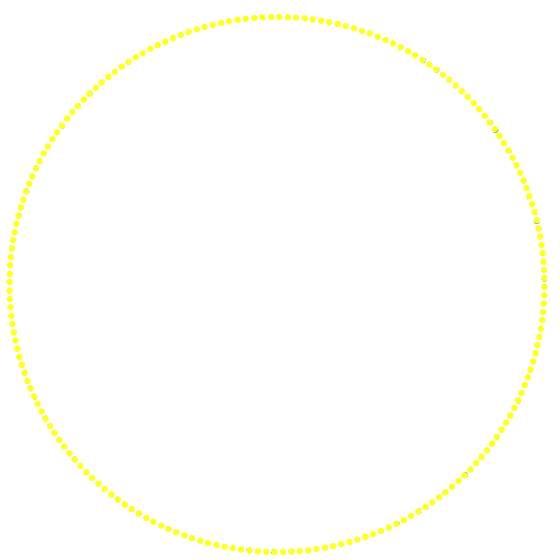} & \includegraphics[width=0.2\textwidth,height=0.25\textwidth,keepaspectratio]{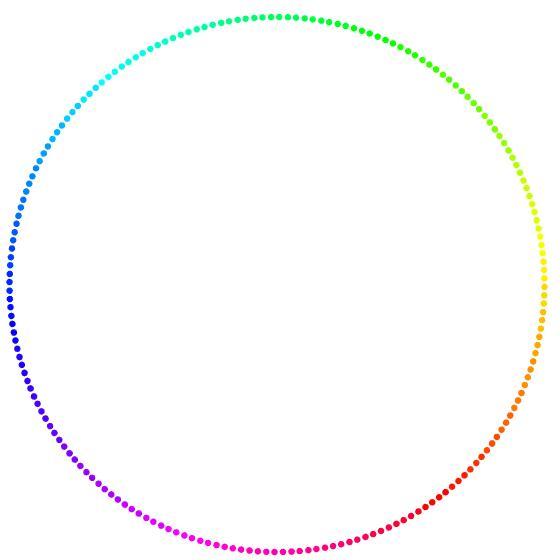} \\
        \hline
    \end{tabular}
    \caption{Embeddings of $2$d manifolds without boundary into $\mathbb{R}^2$. For each manifold, the left and right columns contain the same plots colored by the two parameters of the manifold. A proof of the mathematical correctness of the LDLE embeddings is provided in Figure~\ref{fig:decipher}.}
    \label{fig:fig63}
\end{figure}

\begin{figure}[h!]
    \centering
    \begin{tabular}{M{0.005\textwidth}|M{0.2\textwidth}|M{0.2\textwidth}|M{0.2\textwidth}|M{0.2\textwidth}|}
        & \multicolumn{2}{c|}{\tiny{M\"obius strip}} & \multicolumn{2}{c|}{\tiny{Klein bottle}}\\
        \hline
        \rotatebox{90}{\tiny{Input}} &  \includegraphics[width=0.2\textwidth,height=0.25\textwidth,keepaspectratio]{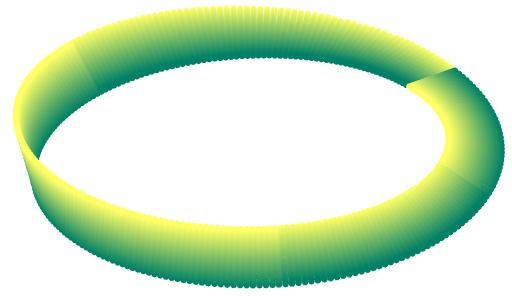} & \includegraphics[width=0.2\textwidth,height=0.25\textwidth,keepaspectratio]{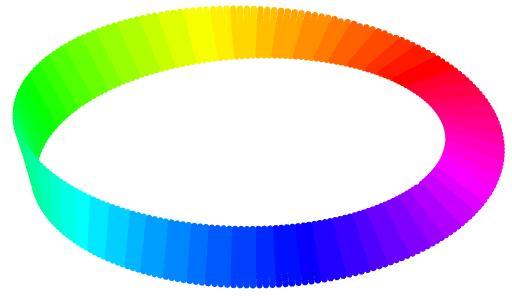} & \multicolumn{2}{c|}{\footnotesize{See Eq.~(\ref{kbinput})}}\\
        \hline
        \rotatebox{90}{\tiny{LDLE}} & \includegraphics[width=0.2\textwidth,height=0.25\textwidth,keepaspectratio]{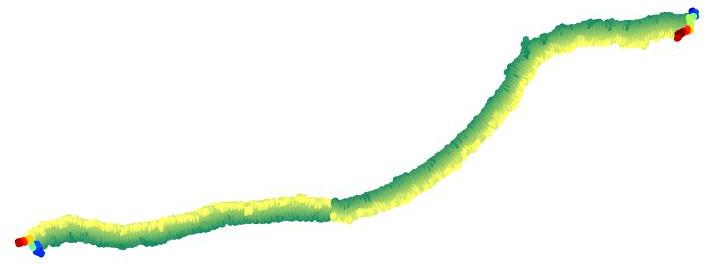} & \includegraphics[width=0.2\textwidth,height=0.25\textwidth,keepaspectratio]{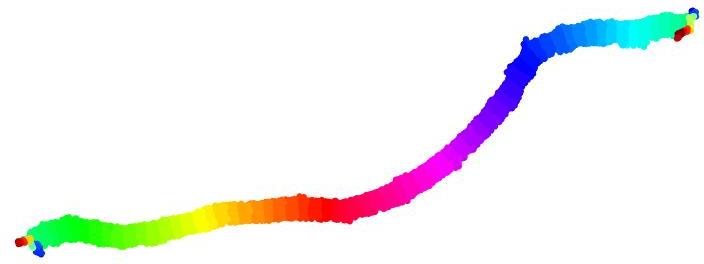} & \includegraphics[width=0.2\textwidth,height=0.25\textwidth,keepaspectratio]{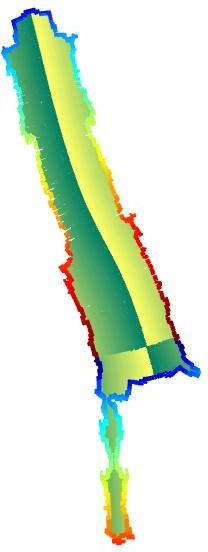} & \includegraphics[width=0.2\textwidth,height=0.25\textwidth,keepaspectratio]{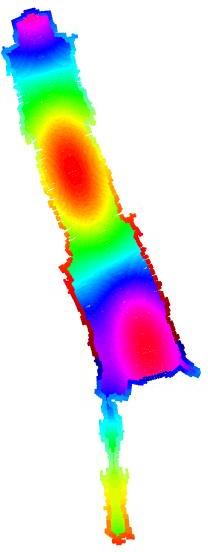} \\
        \hline
        \rotatebox{90}{\tiny{LTSA}} & \includegraphics[width=0.2\textwidth,height=0.25\textwidth,keepaspectratio]{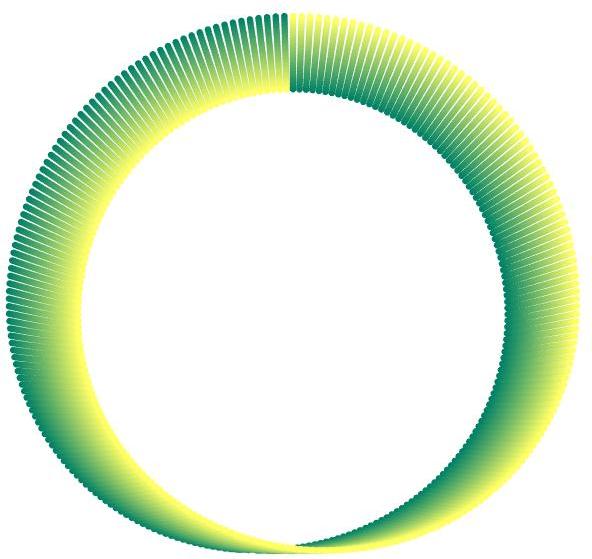} & \includegraphics[width=0.2\textwidth,height=0.25\textwidth,keepaspectratio]{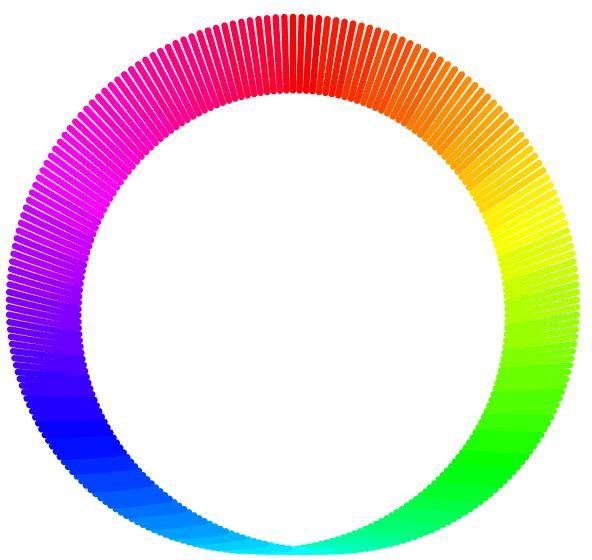} & \includegraphics[width=0.2\textwidth,height=0.25\textwidth,keepaspectratio]{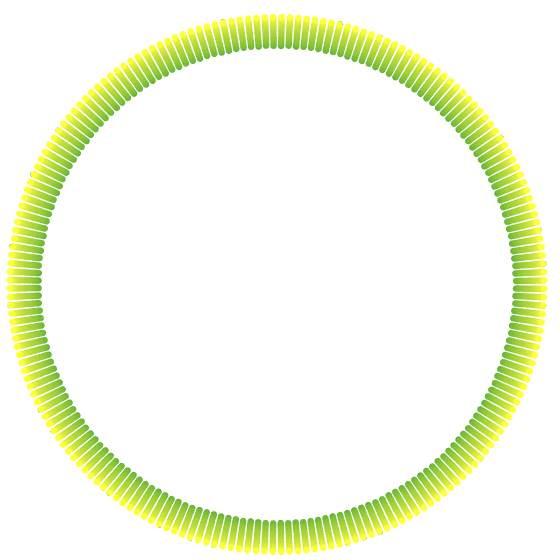} & \includegraphics[width=0.2\textwidth,height=0.25\textwidth,keepaspectratio]{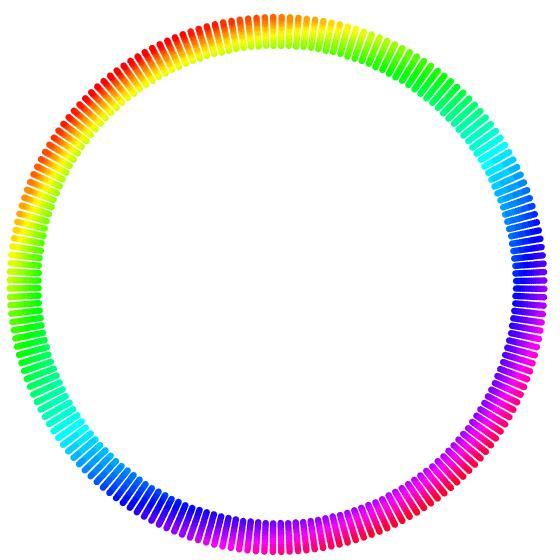} \\
        \hline
        \rotatebox{90}{\tiny{UMAP}} & \includegraphics[width=0.2\textwidth,height=0.25\textwidth,keepaspectratio]{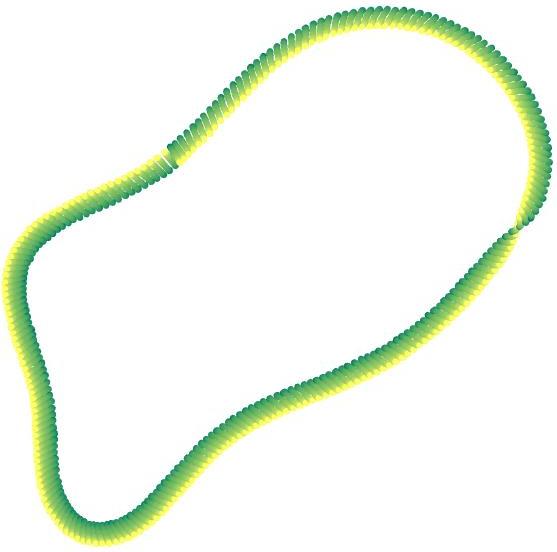} & \includegraphics[width=0.2\textwidth,height=0.25\textwidth,keepaspectratio]{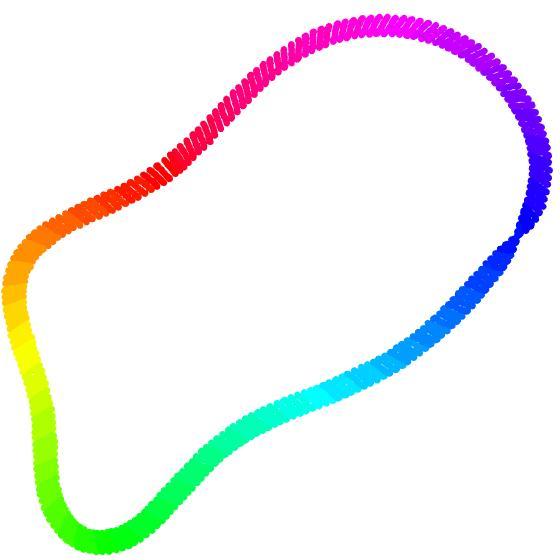} & \includegraphics[width=0.2\textwidth,height=0.25\textwidth,keepaspectratio]{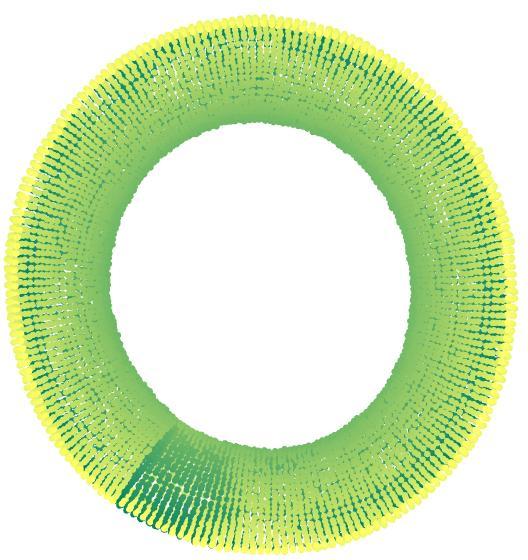} & \includegraphics[width=0.2\textwidth,height=0.25\textwidth,keepaspectratio]{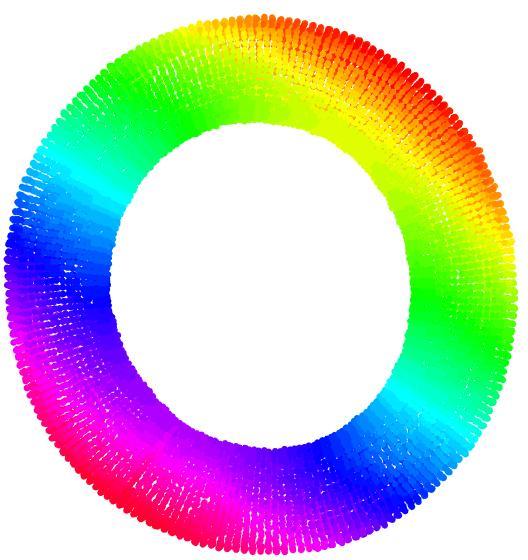} \\
        \hline
        \rotatebox{90}{\tiny{t-SNE}} & \includegraphics[width=0.2\textwidth,height=0.25\textwidth,keepaspectratio]{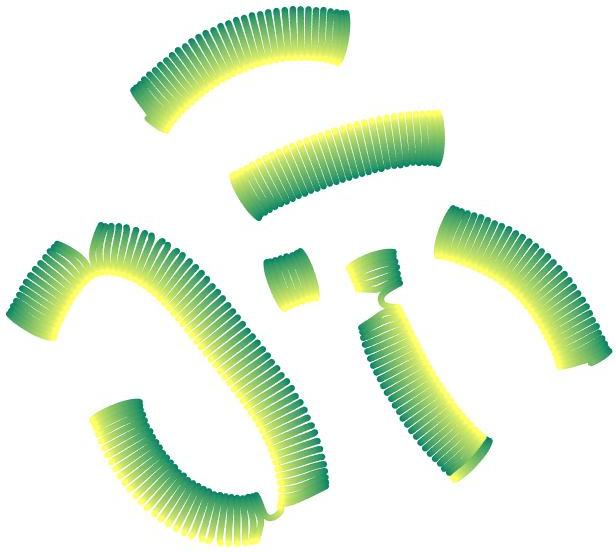} & \includegraphics[width=0.2\textwidth,height=0.25\textwidth,keepaspectratio]{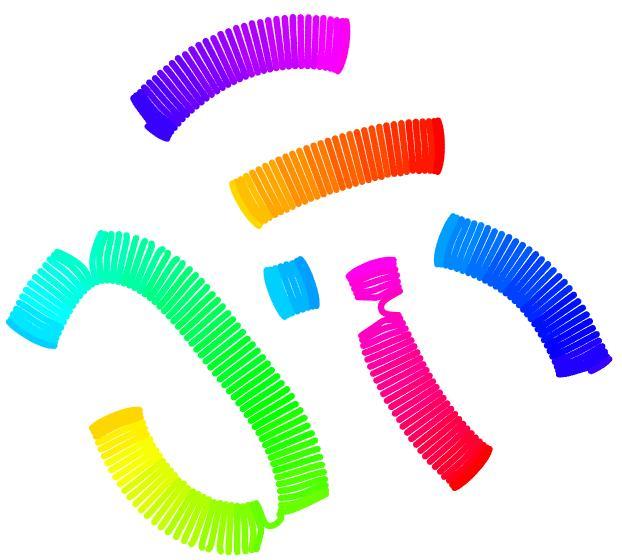} & \includegraphics[width=0.2\textwidth,height=0.25\textwidth,keepaspectratio]{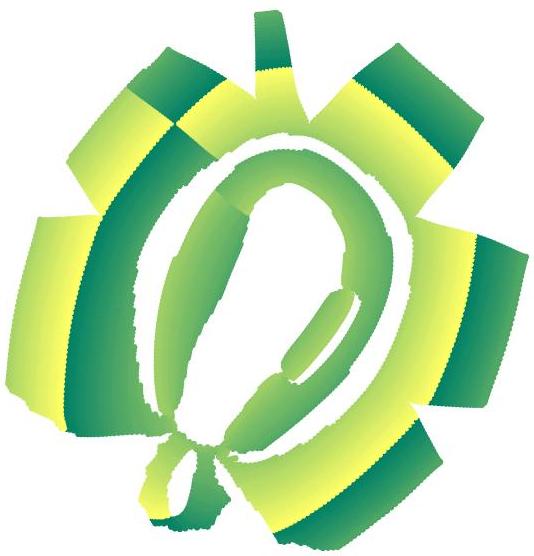} & \includegraphics[width=0.2\textwidth,height=0.25\textwidth,keepaspectratio]{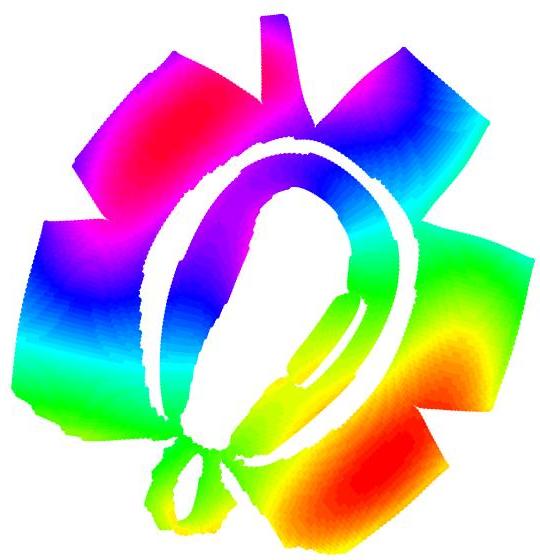} \\
        \hline
        \rotatebox{90}{\tiny{Laplacian Eigenmaps}} & \includegraphics[width=0.2\textwidth,height=0.25\textwidth,keepaspectratio]{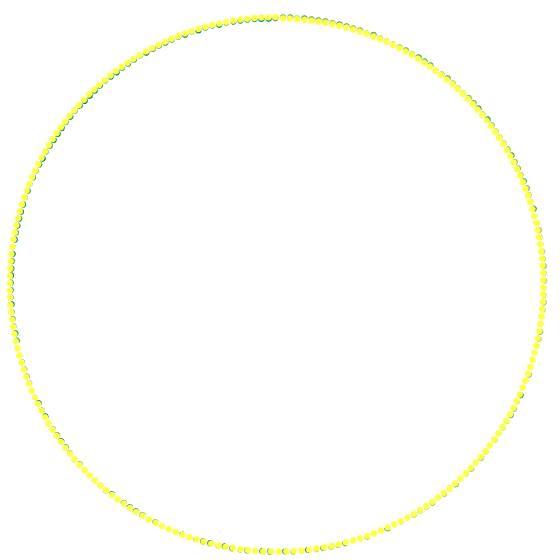} & \includegraphics[width=0.2\textwidth,height=0.25\textwidth,keepaspectratio]{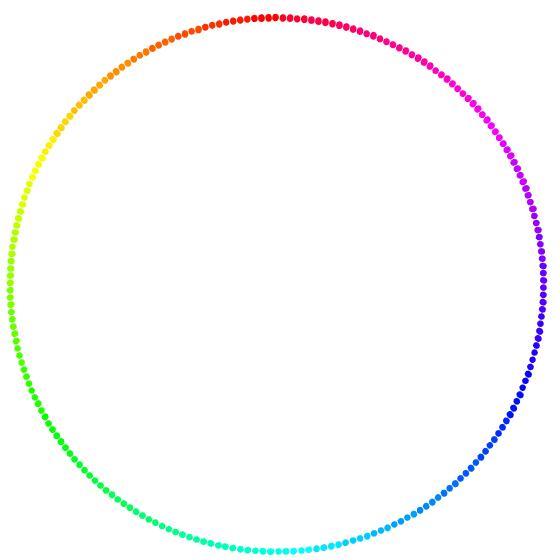} & \includegraphics[width=0.2\textwidth,height=0.25\textwidth,keepaspectratio]{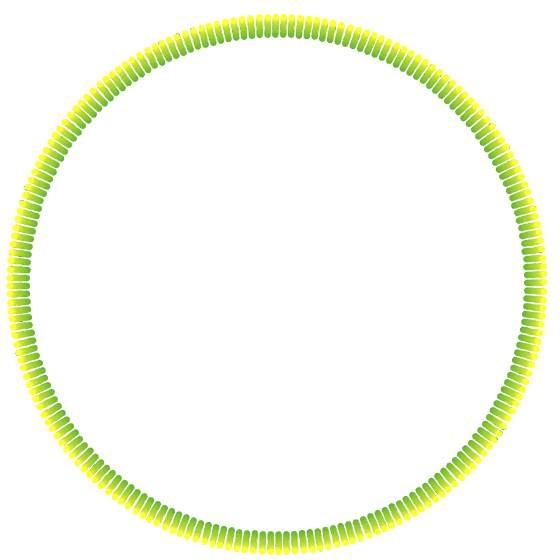} & \includegraphics[width=0.2\textwidth,height=0.25\textwidth,keepaspectratio]{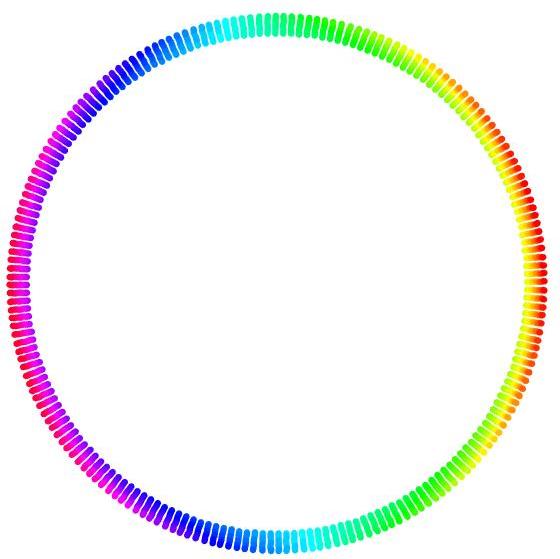} \\
        \hline
    \end{tabular}
    \caption{Embeddings of $2$d non-orientable manifolds into $\mathbb{R}^2$. For each manifold, the left and right columns contain the same plots colored by the two parameters of the manifold. A proof of the mathematical correctness of the LDLE embeddings is provided in Figure~\ref{fig:decipher}.}
    \label{fig:fig64}
\end{figure}

\begin{figure}[h!]
    \centering
    \begin{tabular}{|M{0.125\textwidth}|M{0.125\textwidth}|M{0.125\textwidth}|M{0.125\textwidth}|M{0.125\textwidth}|M{0.125\textwidth}|}
        \tiny{True floor plan} & \tiny{LDLE} & \tiny{LTSA} & \tiny{UMAP} & \tiny{t-SNE} & \tiny{\makecell{Laplacian\\eigenmaps}}\\
        \hline
        \includegraphics[width=0.125\textwidth,keepaspectratio]{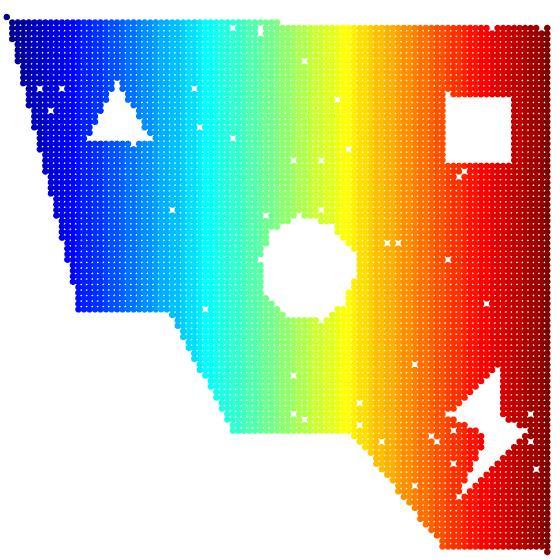} & \includegraphics[width=0.125\textwidth,keepaspectratio]{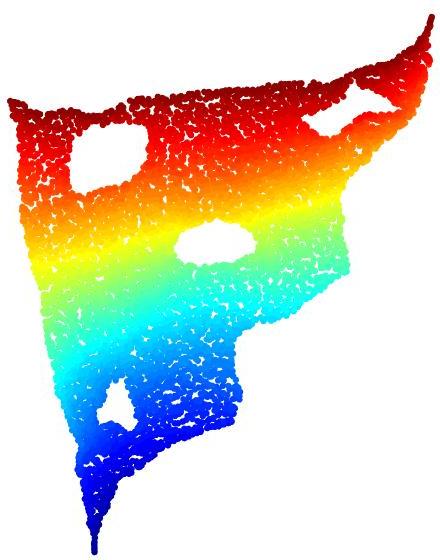} & \includegraphics[width=0.125\textwidth,keepaspectratio]{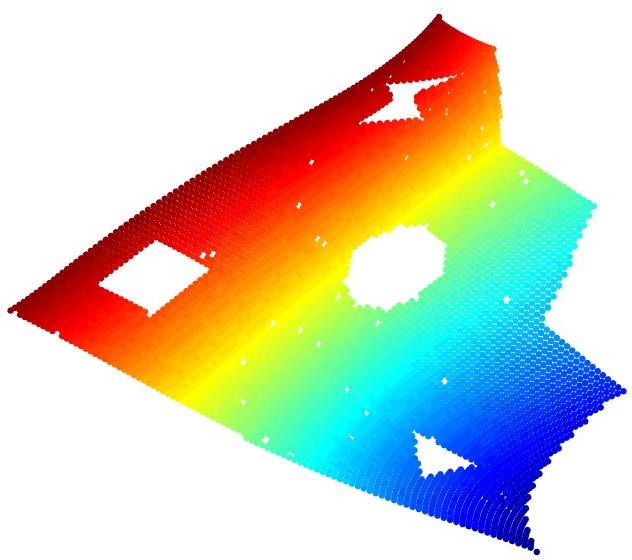} & \includegraphics[width=0.125\textwidth,keepaspectratio]{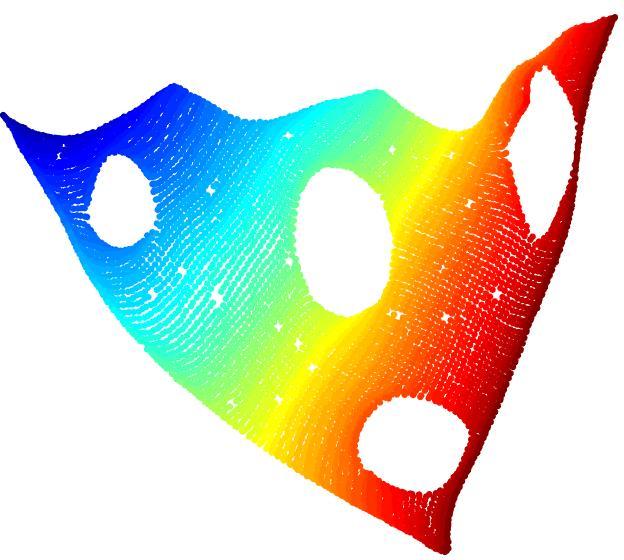} & \includegraphics[width=0.125\textwidth,keepaspectratio]{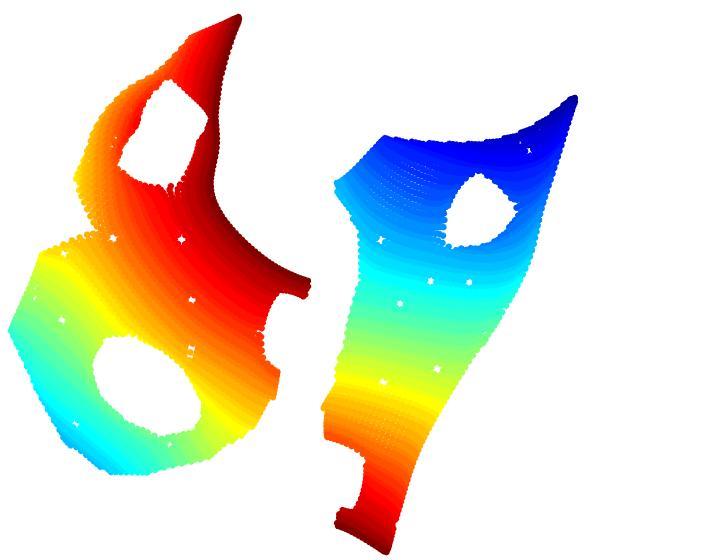} & \includegraphics[width=0.125\textwidth,keepaspectratio]{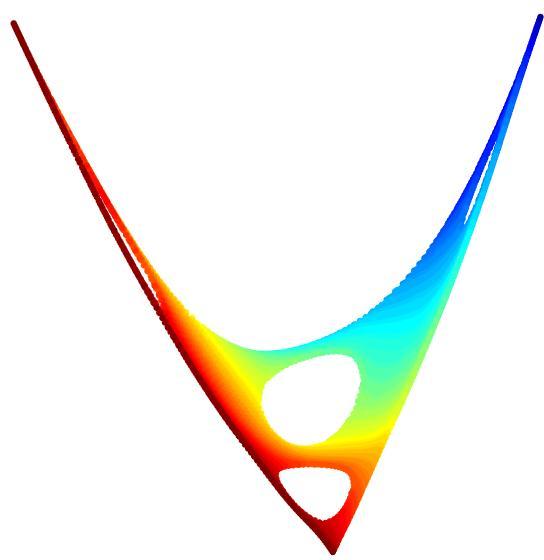} \\
        \hline
    \end{tabular}
    \caption{Embedding of the synthetic sensor data into $\mathbb{R}^2$ (see Section~\ref{subsec:h_d_d} for details).}
    \label{fig:sensor}
\end{figure}

\begin{figure}[h!]
    \centering
    \begin{tabular}{M{0.005\textwidth}|M{0.18\textwidth}|M{0.18\textwidth}|M{0.18\textwidth}|M{0.18\textwidth}|}
        \hline
        & \multicolumn{4}{c|}{\tiny{LDLE}}\\
        \hline
         & \multicolumn{4}{c|}{\includegraphics[width=0.6\textwidth,keepaspectratio]{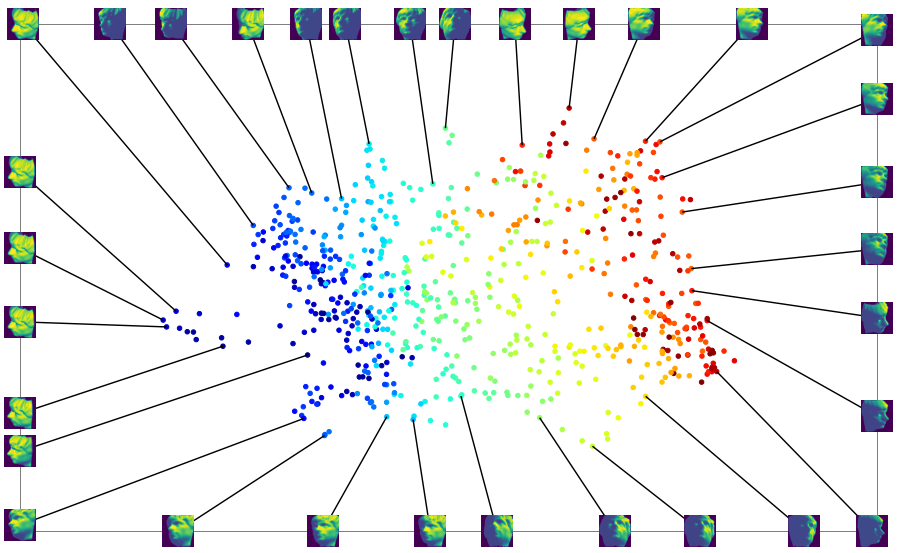}}\\
        \hline
        & \tiny{LDLE} & \tiny{LTSA} & \tiny{UMAP} & \tiny{t-SNE}\\
        \hline
        \rotatebox{90}{\tiny{Pose}} & \includegraphics[width=0.18\textwidth,keepaspectratio]{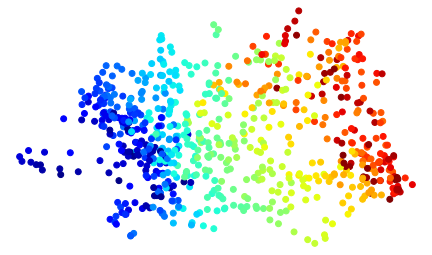} & \rotatebox{90}{\includegraphics[height=0.18\textwidth,keepaspectratio]{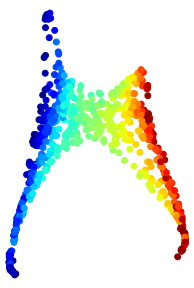}} & \includegraphics[width=0.18\textwidth,keepaspectratio]{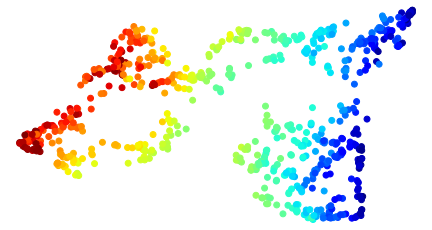} &  \includegraphics[width=0.18\textwidth,keepaspectratio]{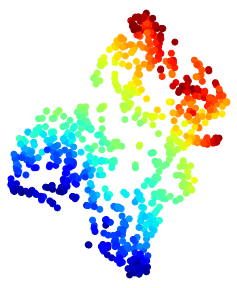}\\
        \hline
        \rotatebox{90}{\tiny{Lighting}} &\includegraphics[width=0.18\textwidth,keepaspectratio]{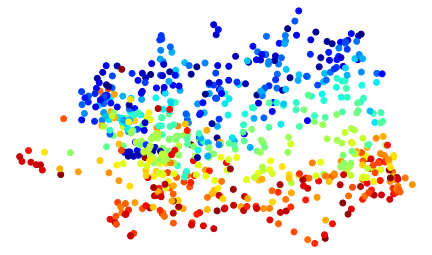} & \rotatebox{90}{\includegraphics[height=0.18\textwidth,keepaspectratio]{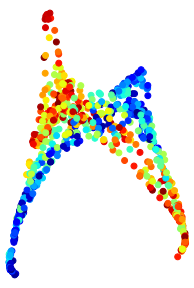}} & \includegraphics[width=0.18\textwidth,keepaspectratio]{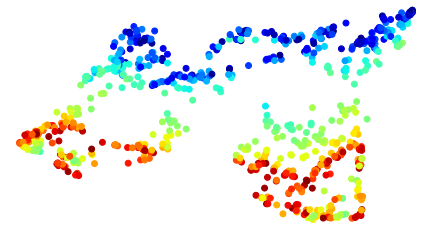} & \includegraphics[width=0.18\textwidth,keepaspectratio]{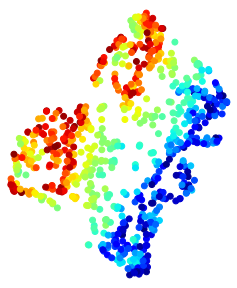}\\
        \hline
    \end{tabular}
    \caption{\editt{Embedding of the face image data set \citep{tenenbaum2000global} into $\mathbb{R}^2$ colored by the pose and lighting conditions (see Section~\ref{subsec:h_d_d} for details).}}
    \label{fig:face_data}
\end{figure}

\begin{figure}[h!]
    \centering
    \begin{tabular}{|M{0.25\textwidth}|M{0.25\textwidth}|M{0.25\textwidth}|}
        \hline
        \multicolumn{3}{|c|}{\tiny{LDLE}}\\
        \hline
        \multicolumn{3}{|c|}{\includegraphics[width=0.65\textwidth,keepaspectratio]{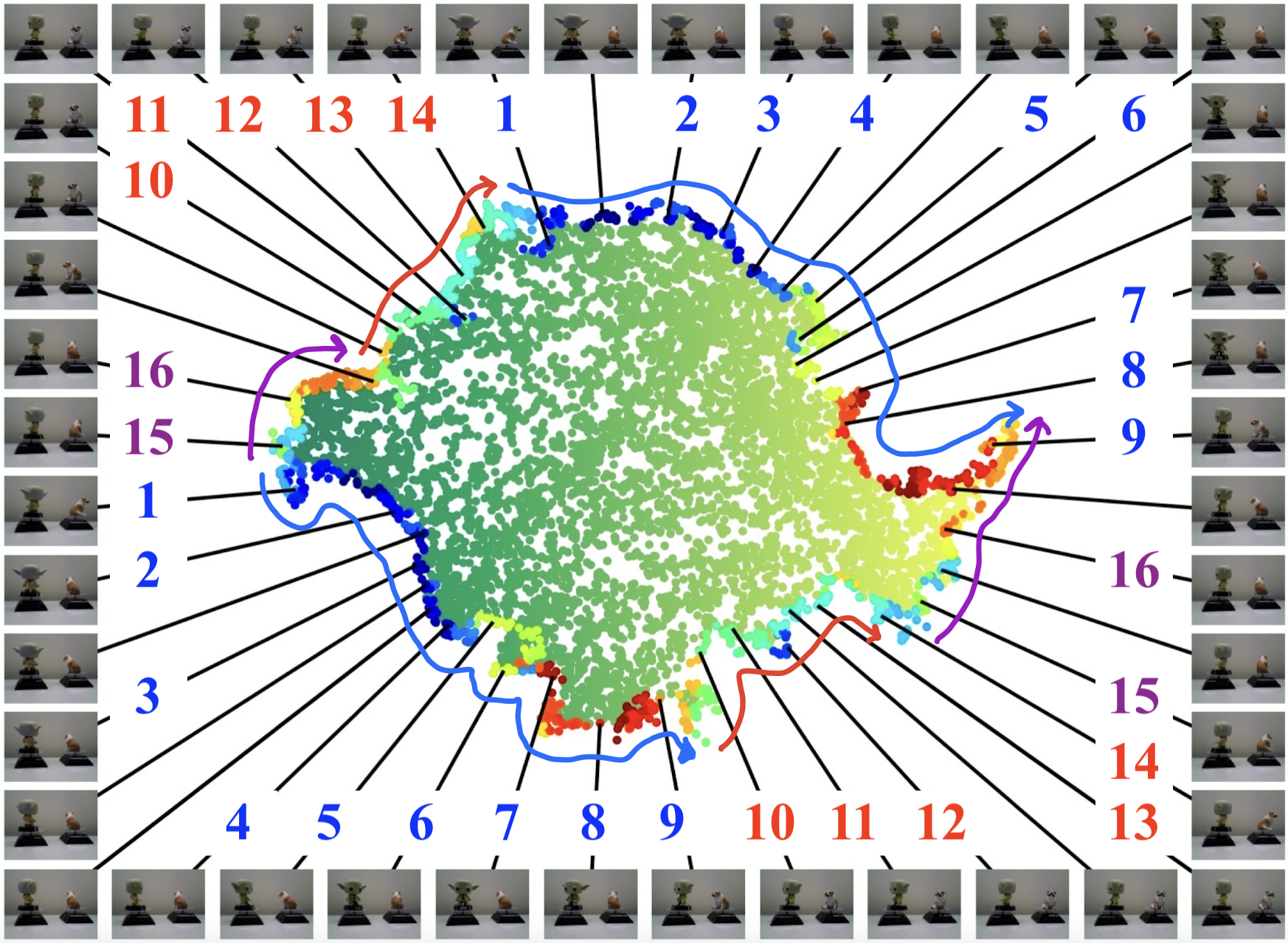}}\\
        \hline
        \tiny{LTSA} & \tiny{UMAP} & \tiny{t-SNE}\\
        \hline
        \includegraphics[width=0.175\textwidth,keepaspectratio]{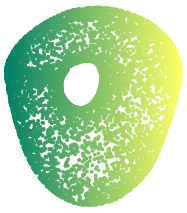} & \includegraphics[width=0.175\textwidth,keepaspectratio]{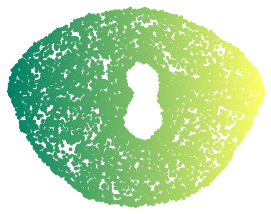} & \includegraphics[width=0.175\textwidth,keepaspectratio]{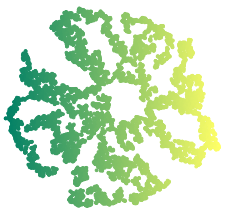}\\
        \hline
    \end{tabular}
    \caption{\editt{Embeddings of snapshots of a platform with two objects, Yoda and a bull dog, each rotating at a different frequency, such that the underlying topology is a torus (see Section~\ref{subsec:h_d_d} for details).}}
    \label{fig:s1_puppets}
\end{figure}

\begin{figure}[h!]
    \centering
    \begin{tabular}{M{0.3\textwidth}M{0.3\textwidth}M{0.3\textwidth}}
        \includegraphics[trim=50 50 50 50, clip, height=0.215\textwidth,keepaspectratio]{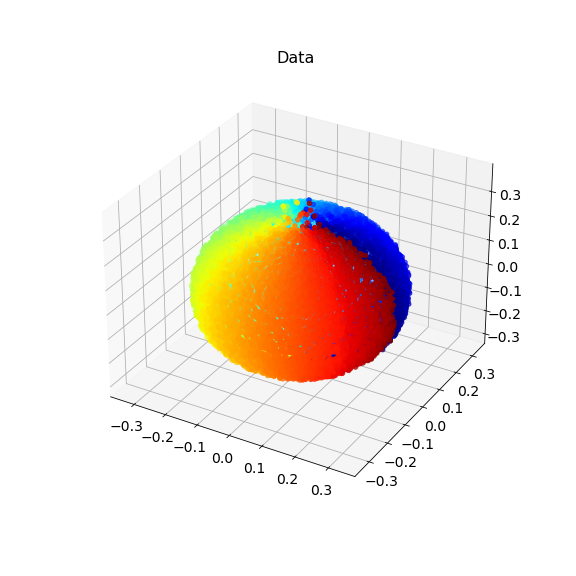} & \includegraphics[trim=60 50 50 50, clip,height=0.25\textwidth,keepaspectratio]{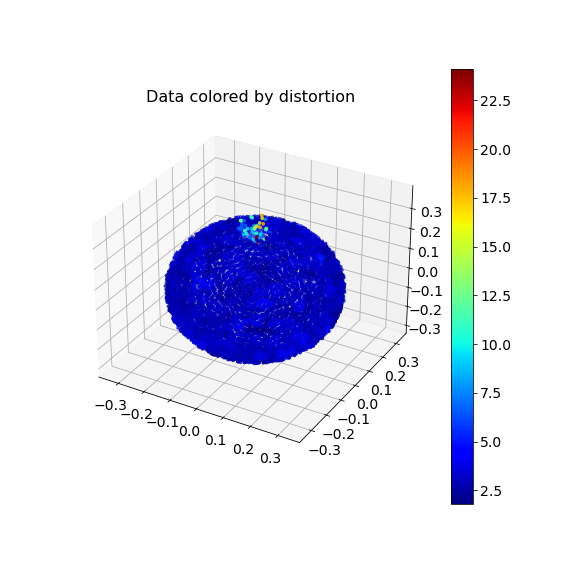} & \includegraphics[trim=35 65 50 50, clip,height=0.25\textwidth,keepaspectratio]{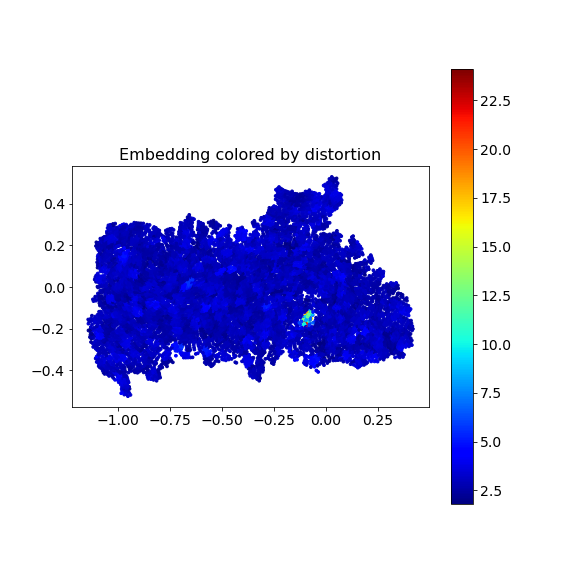}
    \end{tabular}
    \caption{\editt{Local views containing outliers exhibit high distortion. (left) Input data $(x_k)_{k=1}^{n}$. (middle) $x_k$ colored by the distortion $\zeta_{kk}$ of $\Phi_k$ on $U_k$. (right) $y_k$ colored by $\zeta_{kk}$.}}
    \label{fig:anomaly}
\end{figure}

\edit{
\begin{itemize}
    \item It is only natural to expect real world datasets to have boundary and to have many corners. As observed in the experimental results, when the boundary of the manifold is unknown, then the LDLE embedding tends to have distorted boundary. Even when the boundary is known, the embedding has distorted corners. This is caused by high distortion views near the boundary (see Figures~\ref{fig:fig21} and~\ref{fig:fig31}). We aim to fix this issue in our future work. One possible resolution could be based on \citep{berry2017density} which presented a method to approximately calculate the distance of the points from the boundary.
    \item When the data represents a mixture of manifolds, for example, a pair of possibly intersecting spheres or even manifolds of different intrinsic dimensions, it is also natural to expect a manifold learning technique to recover a separate parameterization for each manifold and provide gluing instructions at the output. One way is to perform manifold factorization \citep{pmlr-v130-zhang21j} or multi-manifold clustering \citep{trillos2021large} on the data to recover sets of points representing individual manifolds and then use manifold learning on these separately. We aim to adapt LDLE to achieve this.
    \item The spectrum of the Laplacian has been used in prior work for anomaly detection \citep{7148963,6377228,cheng2018geometry,cheng2020spectral, Mishne2019DN}. Similar to our approach of using a subset of Laplacian eigenvectors to construct low distortion local views in lower dimension, in \citep{Mishne2018biorxiv,cheng2020spectral}, subsets of Laplacian eigenvectors were identified so as to separate small clusters from a large background component. As shown in Figures~\ref{fig:fig21} and~\ref{fig:fig31}, LDLE produced high distortion local views near the boundary and the corners, though these are not outliers. However, if we consider a sphere with outliers (imagine a sphere with noise only at the north pole as in Figure~\ref{fig:anomaly}), then the distortion of the local views containing the outliers is higher than the rest of the views. Therefore, the distortion of the local views can help find anomalies in the data. We aim to further investigate this direction to develop an anomaly detection technique.
    \item Similar to the approach of denoising  a signal by retaining low frequency components, our approach uses low frequency Laplacian eigenvectors to estimate local views. These eigenvectors implicitly capture the global structure of the manifold. Therefore, to construct local views, unlike LTSA which directly relies on the local configuration of data which may be noisy, LDLE relies on the local elements of low frequency global eigenvectors of the Laplacian which are supposed to be robust to the noise. Practical implication of this is shown in Figure~\ref{fig:noise} to some extent while we aim to further investigate the theoretical implications.
\end{itemize}
}

\clearpage

\appendix
\section{First Proof of Theorem~\ref{thm:thm2}}
\label{proof:proof21}
Choose $\epsilon>0$ so that the exponential map $\exp_x:T_x\mathcal{M} \rightarrow \mathcal{M}$ is a well defined diffeomorphism on $\mathcal{B}_{2\epsilon} \subset T_x\mathcal{M}$ where $T_x\mathcal{M}$ is the tangent space to $\mathcal{M}$ at $x$, $\exp_x(0) = x$ and
\begin{align}
    \mathcal{B}_{\epsilon}=\{v\in T_{x}\mathcal{M}\ |\ \left\|v\right\|_2<\epsilon\}. \label{epsball1}
\end{align}
Then using \cite[lem. 48, prop. 50, th. 51]{canzani2013analysis}, for all $y \in B_{\epsilon}(x)$ such that
\begin{align}
    B_{\epsilon}(x)=\{y\in \mathcal{M}\ |\ d_g(x,y)<\epsilon\} \label{epsball}
\end{align}
we have,
\begin{align}
    p(t,x,y) &= G(t,x,y)(u_0(x,y)+tu_1(x,y)+O(t^2)), \label{ptxy}
\end{align}
where
\begin{align}
    G(t,x,y) &= \frac{e^{-d_g(x,y)^2/4t}}{(4\pi t)^{d/2}},\label{Gtxy}\\
    u_0(x,y)&=1+O(\left\|v\right\|^2),\  y=\exp_x(v), v \in T_x\mathcal{M}, \label{u_0}
\end{align}
and for $f \in C(\mathcal{M})$, the following hold
\begin{align}
    f(x)&=\lim_{t\rightarrow 0}\int_{M}p(t,x,y)f(y)\omega_g(y)\\
    &= \lim_{t\rightarrow 0}\int_{B_{\epsilon}(x)}p(t,x,y)f(y)\omega_g(y),\label{MtoepsB}\\
    f(x)&=\lim_{t\rightarrow 0}\int_{B_{\epsilon}(x)}G(t,x,y)f(y)\omega_g(y), \label{Geq}\\
    u_1(x,x)f(x) &= \lim_{t\rightarrow 0}\int_{B_{\epsilon}(x)}G(t,x,y)u_1(x,y)f(y)\omega_g(y). \label{u1eq}
\end{align}
 Using the above equations and the definition of $\Psi_{kij}(y)$ in Eq.~(\ref{Phikij}) and $A_{kij}$ in Eq.~(\ref{Akij}) we compute the limiting value of the scaled local correlation (see Eq.~(\ref{Atilde})),
\begin{align}
    \widetilde{A}_{kij} &=\lim_{t\rightarrow 0}\frac{A_{kij}}{2t} \label{AkijL}\\
    &= \lim_{t\rightarrow 0}\frac{1}{2t}\int_{M}p(t,x_k,y)\Psi_{kij}(y)\omega_g(y).
\end{align}
which will turn out to be the inner product between the gradients of the eigenfunctions $\bm{\phi}_i$ and $\bm{\phi}_j$ at $x_k$.
We start by choosing an $\epsilon_k>0$ so that $\exp_{x_k}$ is a well defined diffeomorphism on $\mathcal{B}_{2\epsilon_k} \subset T_{x_k}\mathcal{M}$. Using Eq.~(\ref{MtoepsB}) we change the region of integration from $\mathcal{M}$ to $B_{\epsilon_k}(x_k)$,
\begin{align}
    \widetilde{A}_{kij} &= \lim_{t_k\rightarrow 0}\frac{1}{2t_k}\int_{B_{\epsilon_k}(x_k)}p(t_k,x_k,y)\Psi_{kij}(y)\omega_g(y). \label{L0}
\end{align}
Substitute $p(t_k,x_k,y)$ from Eq.~(\ref{ptxy}) and simplify using Eq.~(\ref{Geq},~\ref{u1eq}) and the fact that $\Psi_{kij}(x_k)=0$ to get
{
\begin{align}
    \widetilde{A}_{kij} &= \lim_{t_k\rightarrow 0}\frac{1}{2t_k}\int_{B_{\epsilon_k}(x_k)}G(t_k,x_k,y)(u_0(x_k,y)+t_ku_1(x_k,y)+O(t_k^2))\Psi_{kij}(y)\omega_g(y).\nonumber\\
    &= \lim_{t_k\rightarrow 0}\left(\frac{1}{2t_k}\int_{B_{\epsilon_k}(x_k)}G(t_k,x_k,y)u_0(x_k,y)\Psi_{kij}(y)\omega_g(y) +\right.\nonumber\\
    &\qquad\qquad\qquad\qquad\qquad\qquad\qquad\left. \frac{t_ku_1(x_k,x_k)\Psi_{kij}(x_k) + O(t_k^2)\Psi_{kij}(x_k)}{2t_k}\right)\nonumber\\
    &= \lim_{t_k\rightarrow 0}\frac{1}{2t_k}\int_{B_{\epsilon_k}(x_k)}G(t_k,x_k,y)u_0(x_k,y)\Psi_{kij}(y)\omega_g(y) \label{Gu0}.
\end{align}
}%
Replace $y \in B_{\epsilon_k}(x_k)$ by $\exp_{x_{k}}(v)$ where $v \in \mathcal{B}_{\epsilon_k} \subset T_{x_k}\mathcal{M}$ and $\left\|v\right\| = d_g(x_k,y)$. Denote the Jacobian for the change of variable by $J(v)$ i.e. $J(v)=\frac{d}{dv}\exp_{x_k}(v)$. Note that $\exp_{x_k}(0)=x_k$ and $J(0)=I$. Using the Taylor expansion of $\bm{\phi}_{i}$ and $\bm{\phi}_j$ about $0$ we obtain
\begin{align}
    \phi_s(y) = \phi_s(\exp_{x_k}(v)) &= \phi_s(\exp_{x_k}(0)) + \nabla \phi_s(\exp_{x_k}(0))^TJ(0)v + O(\left\|v\right\|^2)\nonumber\\
    &=\phi_s(x_k) + \nabla \phi_s(x_k)^Tv + O(\left\|v\right\|^2),\ s=i,j.
\end{align}
Substituting the above equation in the definition of $\Psi_{kij}(y)$ (see Eq.~(\ref{Phikij})) we get
\begin{align}
    \Psi_{kij}(y)&=\Psi_{kij}(\exp_{x_k}(v))\nonumber\\
    &=v^T\nabla\phi_i\nabla\phi_{j}^Tv + (\nabla\phi_{i}^Tv + \nabla\phi_{j}^Tv)O(\left\|v\right\|^2) + O(\left\|v\right\|^4), \label{phi_taylor}
\end{align}
where $\nabla\phi_s \equiv \nabla\phi_s(x_k), s=i,j$. Now we substitute Eq.~(\ref{phi_taylor},~\ref{Gtxy},~\ref{u_0}) in Eq.~(\ref{Gu0}) while replacing variable $y$ with $\exp_{x_k}(v)$ where $J(v)$ is the Jacobian for the change of variable as before, to get
\begin{align}
    \widetilde{A}_{kij} &=\lim_{t_k\rightarrow 0}\frac{1}{2t_k}\int_{\mathcal{B}_{\epsilon_k}}\frac{e^{-\left\|v\right\|^2/4t_k}}{(4\pi t_k)^{d/2}}(1+O(\left\|v\right\|^2))\Psi_{kij}(\exp_{x_k}(v)) J(v) dv\nonumber\\
    &=L_1+L_2,
\end{align}
where $L_1$ and $L_2$ are the terms obtained by expanding $1+O(\left\|v\right\|^2)$ in the integrand. We will show that $L_2=0$ and $\widetilde{A}_{kij}=L_1 = \nabla\phi_i^T\nabla\phi_j$.
\begin{align}
 L_2 &=\lim_{t_k\rightarrow 0}\frac{1}{2t_k}\int_{\mathcal{B}_{\epsilon_k}}\frac{e^{-\left\|v\right\|^2/4t_k}}{(4\pi t_k)^{d/2}}O(\left\|v\right\|^2)(\operatorname{tr}(\nabla\phi_i\nabla\phi_j^Tvv^T)+ \nonumber\\
 &\qquad\qquad\qquad\qquad\qquad (\nabla\phi_{i}^Tv + \nabla\phi_{j}^Tv)O(\left\|v\right\|^2) + O(\left\|v\right\|^4))J(v) dv\nonumber\\
    &= \lim_{t_k\rightarrow\nonumber 0}\frac{1}{2t_k}(O(t_k^2)+0+0+O(t_k^4))\nonumber\\
    &= 0.
\end{align}
Therefore,
\begin{align}
    \widetilde{A}_{kij} &= L_1 \nonumber\\
    &=\lim_{t_k\rightarrow 0}\frac{1}{2t_k}\int_{\mathcal{B}_{\epsilon_k}}\frac{e^{-\left\|v\right\|^2/4t_k}}{(4\pi t_k)^{d/2}}\Psi_{kij}(\exp_{x_k}(v)) J(v) dv\label{L1y1}\\
    &= \lim_{t_k\rightarrow 0}\frac{1}{2t_k}\int_{\mathcal{B}_{\epsilon_k}}\frac{e^{-\left\|v\right\|^2/4t_k}}{(4\pi t_k)^{d/2}}(v^T\nabla\phi_i\nabla\phi_{j}^Tv +\nonumber\\
    &\qquad\qquad\qquad\qquad (\nabla\phi_{i}^Tv + \nabla\phi_{j}^Tv)O(\left\|v\right\|^2) + O(\left\|v\right\|^4)) J(v) dv\nonumber\\
    &= \lim_{t_k\rightarrow 0}\frac{1}{2t_k}\int_{\mathcal{B}_{\epsilon_k}}\frac{e^{-\left\|v\right\|^2/4t_k}}{(4\pi t_k)^{d/2}}v^T\nabla\phi_i\nabla\phi_{j}^TvJ(v)dv+\frac{0+0+O(t_k^2)}{2t_k} \nonumber\\
    &=\lim_{t_k\rightarrow 0}\frac{1}{2t_k}\int_{\mathcal{B}_{\epsilon_k}}\frac{e^{-\left\|v\right\|^2/4t_k}}{(4\pi t_k)^{d/2}}v^T\nabla\phi_i\nabla\phi_{j}^TvJ(v)dv.
\end{align}
Substitution of $t_k = 0$ leads to the indeterminate form $\frac{0}{0}$. Therefore, we apply L'Hospital's rule and then Leibniz integral rule to get,
{\small
\begin{align}
    \widetilde{A}_{kij}&=\lim_{t_k\rightarrow 0}\frac{1}{2}\int_{\mathcal{B}_{\epsilon_k}}\left(\frac{\left\|v\right\|^2}{4t_k^2}-\frac{d}{2t_k}\right)\frac{e^{-\left\|v\right\|^2/4t_k}}{(4\pi t_k)^{d/2}}v^T\nabla\phi_i\nabla\phi_j^Tv J(v) dv\nonumber\\
     &= \operatorname{tr}\left(\frac{1}{2}\nabla\phi_i\nabla\phi_j^T\lim_{t_k\rightarrow 0}\int_{\mathcal{B}_{\epsilon_k}}\left(\frac{\left\|v\right\|^2}{4t_k^2}-\frac{d}{2t_k}\right)\frac{e^{-\left\|v\right\|^2/4t_k}}{(4\pi t_k)^{d/2}}vv^T J(v) dv\right)\nonumber\\
     &= \operatorname{tr}\left(\frac{1}{2}\nabla\phi_i\nabla\phi_j^T \left(\lim_{t_k\rightarrow 0}\left(\frac{(12+4(d-1))t_k^2}{4t_k^2} - \frac{2t_kd}{2t_k}\right)I+O(t_k)I \right)\right)\nonumber\\
     &= \nabla\phi_i^T\nabla\phi_j.
\end{align}
}%
Finally, note that the Eq.~(\ref{L1y1}) is same as the following equation with $y$ replaced by $\exp_{x_k}(v)$,
\begin{align}
    \widetilde{A}_{kij} &= \lim_{t_k\rightarrow 0}\frac{1}{2t_k}\int_{B_{\epsilon_k}(x_k)}G(t_k,x_k,y)\Psi_{kij}(y)\omega_g(y). \label{L1y}
\end{align}
We used the above equation to estimate $\widetilde{A}_{kij}$ in Section~\ref{subsec:scaledloccorr}.
\qed

\section{\edit{Second Proof of Theorem~\ref{thm:thm2}}}
\label{proof:proof22}
Yet another proof is based on the Feynman-Kac formula \citep{doi:10.1080/03605302.2014.942739, steinerberger2017spectral},
\begin{align}
    A_{kij} &= [e^{-t_{k}\Delta_g}((\phi_i - \phi_i(x_k))(\phi_j - \phi_j(x_k))](x_k).
\end{align}
where 
\begin{align}
    [e^{-t\Delta_g}f](x) = \sum_{i}e^{-\lambda_i t}\langle\phi_i,f\rangle\phi_i(x)
\end{align}
and therefore,
\begin{align}
    \widetilde{A}_{kij} &= \lim_{t_k\rightarrow 0}\frac{A_{kij}}{2t_k} = \left.\frac{1}{2}\frac{\partial A_{kij}}{\partial t_k} \right\vert_{t_k=0}\\
    &= \frac{-1}{2}\left\{\Delta_g[(\phi_i - \phi_i(x_k)) (\phi_j-\phi_j(x_k))](x_k)\right\}\\
    &= \frac{-1}{2}\left\{0+0-2\nabla\phi_i(x_k)^T\nabla\phi_j(x_k)\right\}\\
    &= \nabla\phi_i(x_k)^T\nabla\phi_j(x_k) \label{eq:Atildekij_id}
\end{align}
where we used the fact $\Delta_g (f_if_j)= f_j\Delta_gf_i + f_i\Delta_gf_j - 2\langle \nabla_g f_i (x), \nabla_g f_j(x) \rangle_g$. Note that as per our convention $\nabla\phi_i(x_k) = \nabla(\phi_i \ \circ\ \text{exp}_{x_k})(0)$ and therefore $\langle \nabla_g \phi_i (x), \nabla_g \phi_j(x) \rangle_g = \nabla\phi_i(x_k)^T\nabla\phi_j(x_k)$.

\section{Rationale Behind the Choice of $t_k$ in Eq.~(\ref{tk})}
\label{rationaletk}
Since $|\mathcal{M}|\leq 1$, we note that
\begin{align}
    \epsilon_k \leq \Gamma(d/2+1)^{1/d}/\sqrt{\pi} \label{epskbound}
\end{align}
where the maximum can be achieved when $\mathcal{M}$ is a $d$-dimensional ball of unit volume. Then we take the limiting value of $t_k$ as in Eq.~(\ref{tk}) where $\text{chi2inv}$ is the inverse cdf of the chi-squared distribution with $d$ degrees of freedom evaluated at $p$. Since the covariance matrix of $G(t_k,x,y)$ is $\sqrt{2t_k}I$ (see Eq.~(\ref{Gtxy1})), the above value of $t_k$ ensures $p$ probability mass to lie in $B_{\epsilon_k}(x_k)$. We take $p$ to be $0.99$ in our experiments. Also, using Eq.~(\ref{epskbound}) and Eq.~(\ref{tk}) we have
\begin{align}
    t_k \leq \frac{1}{2\pi}\frac{\Gamma(d/2+1)^{2/d}}{\text{chi2inv}(p,d)} << 1,\ \text{when }p=0.99.
\end{align}
Using the above inequality with $p=0.99$, for $d=2, 10, 100$ and $1000$, the upper bound on $t_k=0.0172, 0.018, 0.0228$ and $0.0268$ respectively. Thus, $t_k$ is indeed a small value close to $0$.

\section{Computation of $(s_m,p_{s_m})_{m=1}^{M}$ in Algo.~\ref{algo:geinit}}
\label{sec:smpm}

Algo.~\ref{algo:geinit} aligns the intermediate views in a sequence. The computation of the sequences $(s_m,p_{s_m})_{m=1}^{M}$ is motivated by the necessary and sufficient conditions for a unique solution to the standard orthogonal Procrustes problem \citep{schonemann1966generalized}. We start by a brief review of a variant of the orthogonal Procrustes problem and then explain how these sequences are computed.

\subsection{A Variant of Orthogonal Procrustes Problem}
Given two matrices $A$ and $B$ of same size with $d$ columns, one asks for an orthogonal matrix $T$ of size $d \times d$ and a $d$-dimensional columns vector $v$ which most closely aligns $A$ to $B$, that is,
\begin{align}
    T,v = \argmin_{\Omega,\omega}\left\| A\Omega + \mathbf{1}_n\omega^T - B\right\|^2_F \text{ such that } \Omega^T\Omega = I. \label{opvariant}
\end{align}
Here $\mathbf{1}_n$ is the $n$-dimensional column vector containing ones. Equating the derivative of the objective with respect to $\omega$ to zero, we obtain the following condition for $\omega$,
\begin{align}
    \omega = \frac{\mathbf{1}_n}{n}^T(A\Omega-B).
\end{align}
Substituting this back in Eq.~(\ref{opvariant}), we reduce the above problem to the standard orthogonal Procrustes problem,
\begin{align}
    T = \argmin_{\Omega}\left\|\overline{A}\Omega-\overline{B}\right\|_F^2
\end{align}
where 
\begin{align}
    \overline{X} = \left(I - \frac{1}{n}\mathbf{1}_n\mathbf{1}_n^T\right)X
\end{align}
for any matrix $X$. This is equivalent to subtracting the mean of the rows in $X$ from each row of $X$.


As proved in \citep{schonemann1966generalized}, the above problem, and therefore the variant, has a unique solution if and only if the square matrix $\overline{A}^T\overline{B}$ has full rank $d$. Denote by $\sigma_d(X)$ the $d$th smallest singular value of $X$. Then $\overline{A}^T\overline{B}$ has full rank if $\sigma_d(\overline{A}^T\overline{B})$ is non-zero, otherwise there exists multiple $T$ which minimize Eq.~(\ref{opvariant}).

\subsection{Computation of $(s_m,p_{s_m})_{m=1}^{M}$}
Here, $s_m$ corresponds to the $s_m$th intermediate view and $p_{s_m}$ corresponds to its parent view. The first view in the sequence corresponds to the largest cluster and it has no parent, that is,
\begin{align}
    s_1 = \argmax_{m=1}^{M}|\mathcal{C}_m| \text{ and } p_{s_1} = \text{none}. \label{s1p1}
\end{align}
For convenience, denote $s_m$ by $s$, $p_{s_m}$ by $p$ and $V_{mm'}$ by $\widetilde{\Phi}^g_{m}(\widetilde{U}_{mm'})$. We choose $s$ and $p$ so that the view $V_{sp}$ can be aligned with the view $V_{ps}$ without any \textit{ambiguity}. In other words, $s$ and $p$ are chosen so that there is a unique solution to the above variant of orthogonal Procrsutes problem (see Eq.~(\ref{opvariant})) with $A$ and $B$ replaced by $V_{sp}$ and $V_{ps}$, respectively. Therefore, an ambiguity (non-uniqueness) would arise when $\sigma_d(\overline{V}_{sp}^T\overline{V}_{ps})$ is zero. We quantify the ambiguity in aligning arbitrary $m$th and the $m'$th intermediate views on their overlap, that is, $V_{mm'}$ and $V_{m'm}$, by
\begin{align}
    W_{mm'} = \sigma_d(\overline{V}_{mm'}^T\overline{V}_{m'm}). \label{Wkk'}
\end{align}
Note that $W_{mm'} = W_{m'm}$. A value of $W_{mm'}$ close to zero means high ambiguity in the alignment of $m$th and $m'$th views. By default, if there is no overlap between $m$th and $m'$th view then $W_{mm'} = W_{m'm} = 0$.


Finally, we compute the sequences $(s_m,p_{s_m})_{m=2}^{M}$ so that $\sum_{m=2}^{M}W_{s_mp_{s_m}}$ is maximized and therefore the net ambiguity is minimized. This is equivalent to obtaining a maximum spanning tree $T$ rooted at $s_1$, of the graph with $M$ nodes and $W$ as the adjacency matrix.
Then $(s_m)_{m=2}^{M}$ is the sequence in which a breadth first search starting from $s_1$ visits the nodes in $T$. And $p_{s_m}$ is the parent of the $s_m$th node in $T$. Thus,
\begin{align}
    (s_m)_{m=2}^{M} = \text{Breadth-First-Search}(T,s_1) \text{ and } p_{s_m} =  \text{parent of } s_m \text{ in } T.\label{pksk}
\end{align}

\section{Computation of $\widetilde{U}^g_{mm'}$ in Eq.~(\ref{Zs2})}
\label{sec:Nmg}

Recall that $\widetilde{U}^g_{mm'}$ is the overlap between the $m$th and $m'$th intermediate views in the embedding space. The idea behind its computation is as follows. We first compute the discrete balls $U^g_{k}$ around each point $y_k$ in the embedding space. These are the analog of $U_k$ around $x_k$ (see Eq.~\ref{Uk}) but in the embedding space, and are given by
\begin{align}
    U^g_k = \{y_{k'}\ |\ d_e(y_k,y_{k'}) < \epsilon^g_k\}. \label{eq:Ugk}
\end{align}
An important point to note here is that while in the ambient space, we used $\epsilon_k$, the distance to the $k_{\text{lv}}$th nearest neighbor, to define a discrete ball around $x_k$, in the embedding space, we must relax $\epsilon_k$ to account for a possibly increased separation between the embedded points. This increase in separation is caused due to the distorted parameterizations. Therefore, to compute discrete balls in the embedding space, we used $\epsilon^g_k$ in Eq.~(\ref{eq:Ugk}), which is the distance to the $\nu k_{\text{lv}}$th nearest neighbor of $y_k$. In all of our experiments, we take $\nu$ to be $3$.

Recall that $c_k$ is the cluster label for the point $x_k$. Using the same label $c_k$ for the point $y_k$, we construct \textit{secondary} intermediate views $\widetilde{U}^g_{m}$ in the embedding space,
\begin{align}
    \widetilde{U}^g_m = \cup_{c_k = m}U^g_k.
\end{align}
Finally, same as the computation of $\widetilde{U}_{mm'}$ in Eq.~(\ref{Utildekkp}), we compute $\widetilde{U}^g_{mm'}$ as the intersection of $\widetilde{U}^g_{m}$ and $\widetilde{U}^g_{m'}$,
\begin{align}
    \widetilde{U}^g_{mm'} = \widetilde{U}^{g}_{m} \cap \widetilde{U}^{g}_{m'}.
\end{align}

\section{\edit{Comparison with the Alignment Procedure in LTSA}}
\label{subsec:ltsa_stitch}
In the following we use the notation developed in this work. LTSA \citep{zhang2003nonlinear} computes the global embedding $Y_m$ of the $m$th intermediate view $\widetilde{U}_m$ so that it respects the local geometry determined by $\widetilde{\Phi}_{m}(\widetilde{U}_m)$. That is,
\begin{align}
    Y_m = \widetilde{\Phi}_{m}(\widetilde{U}_m)L_m + e_mv_m^T + E_m. \label{Ym_ltsa}
\end{align}
Here, $Y = [y_1, y_2, \ldots, y_n]^T$ where $y_i$ is a column vector of length $d$ representing the global embedding of $x_i$, $Y_m$ is a submatrix of $Y$ of size $|\widetilde{U}_m| \times d$ representing the global embeddings of the points in $\widetilde{U}_m$, and $\widetilde{\Phi}_{m}(\widetilde{U}_m)$ is a matrix of size $|\widetilde{U}_m| \times d$ representing the $m$th intermediate view in the embedding space (or in the notation of LTSA, the local embedding of $\widetilde{U}_m$). $e_m$ is a column vector of length $|\widetilde{U}_m|$ containing $1$s. The intermediate view $\widetilde{\Phi}_{m}(\widetilde{U}_m)$ is transformed into the final embedding $Y_m$ through an affine matrix $L_m$ of size $d \times d$ and a translation vector $v_m$ of length $d$. The reconstruction error is captured in the matrix $E_m$. The total reconstruction error is given by,
\begin{align}
    \mathcal{L}'(Y,(L_m,v_m)_{m=1}^{M})&= \sum_{m=1}^{M}\left\|Y_{m} - (\widetilde{\Phi}_{m}(\widetilde{U}_m)L_m + e_mv_m^T)\right\|^2_{F}. \label{ltsa_alignerr}
\end{align}

LTSA estimates $Y$ and $(L_m,v_m)_{m=1}^{M}$ by minimizing the above objective with the constraint $Y^TY = I$. This constraint is the mathematical realization of their assumption that the points are uniformly distributed in the embedding space. Due to this, the obtained global embedding $Y$ does not capture the aspect ratio of the underlying manifold. Also note that due to the overlapping nature of the views $\widetilde{U}_m$, the terms in the above summation are dependent through $Y_m$'s.


Setting aside our adaptation of GPA to tear closed and non-orientable manifolds, our alignment procedure minimizes the error $\mathcal{L}$ in Eq.~(\ref{alignerr}). By introducing the variables $Y$ and $E_m$ as in Eq.~(\ref{Ym_ltsa}), one can deduce that $\mathcal{L}$ is a lower bound of $\mathcal{L}'$ in Eq.~(\ref{ltsa_alignerr}). The main difference in the two alignment procedures is that, while in LTSA, $Y$ is constrained and the transformations are not, in our approach, we restrict the transformations to be rigid. That is, we constrained $L_m$ to be $b_mT_m$ where $b_m$ is a fixed positive scalar as computed in Eq.~(\ref{bk}) and $T_m$ is restricted to be an orthogonal matrix, while there is no constraint on $Y$.



\begin{figure}[h!]  
    \centering
    \footnotesize
    \begin{tabular}{|M{0.14\textwidth}|M{0.13\textwidth}|M{0.14\textwidth}|M{0.14\textwidth}|M{0.14\textwidth}|M{0.14\textwidth}|}
        \tiny{Rectangle} & \tiny{Barbell} & \tiny{Square with two holes} & \tiny{Sphere with a hole} & \tiny{Swiss Roll with a hole} & \tiny{Noisy Swiss Roll}\\
         \hline
         \includegraphics[trim=295 65 250 65, clip, width=0.15\textwidth,keepaspectratio,align=c]{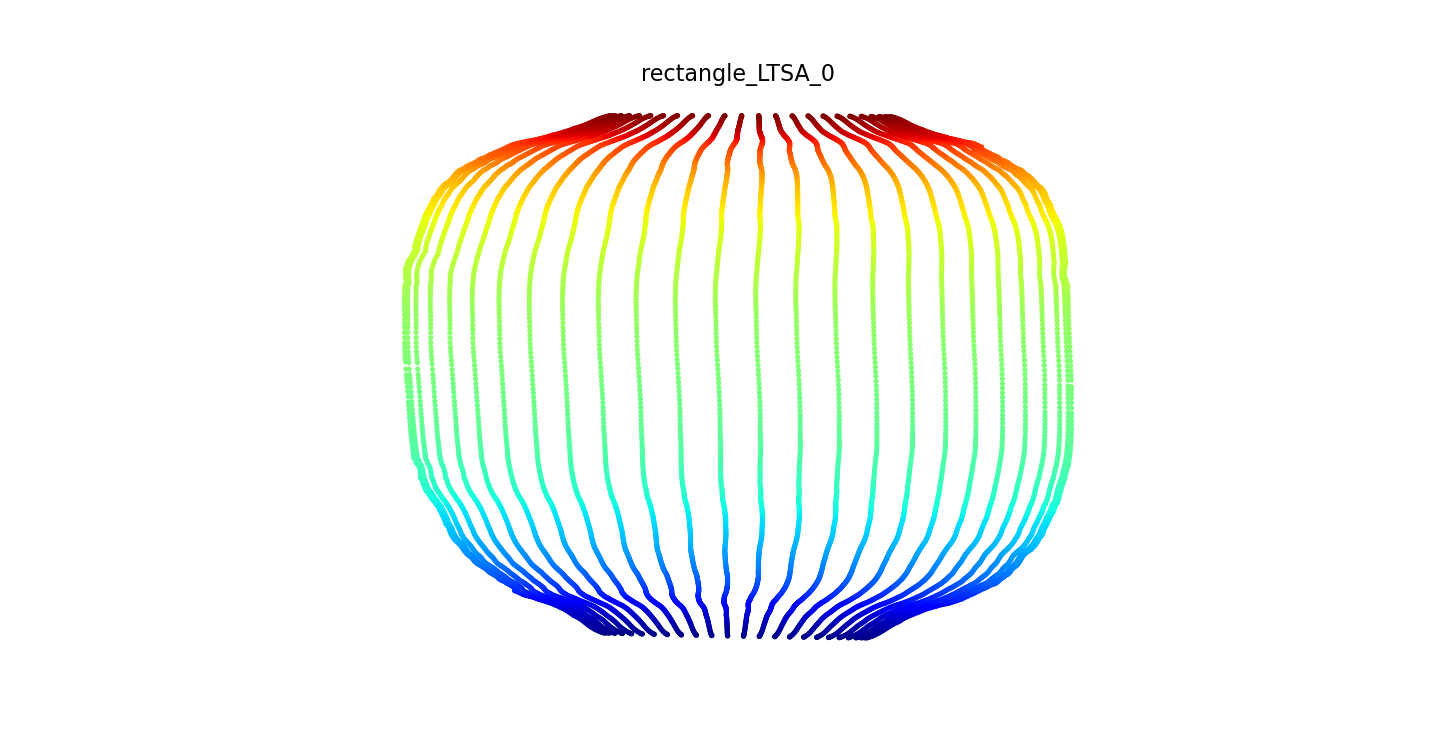} &  \includegraphics[trim=375 65 300 65, clip,width=0.15\textwidth,keepaspectratio,align=c]{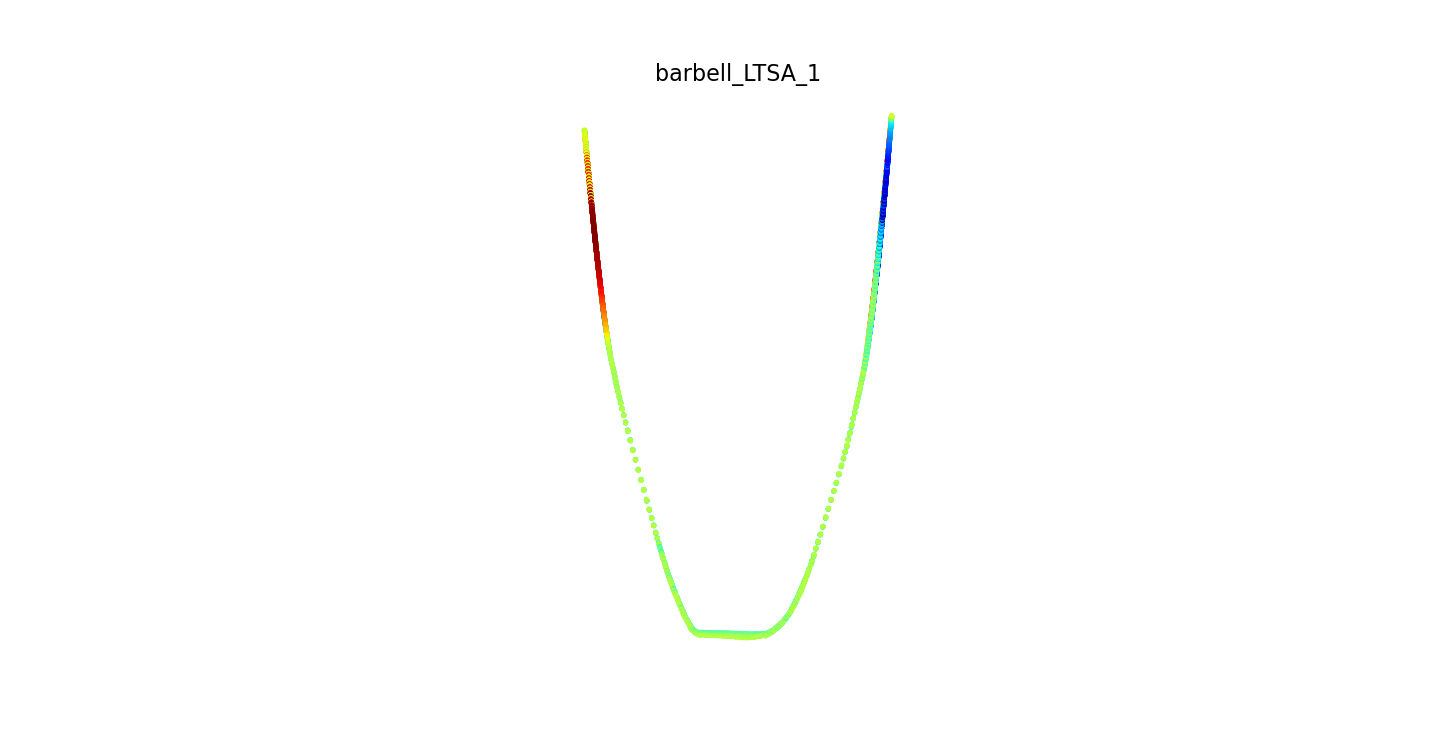}  & \includegraphics[trim=300 65 300 65, clip,width=0.15\textwidth,keepaspectratio,align=c]{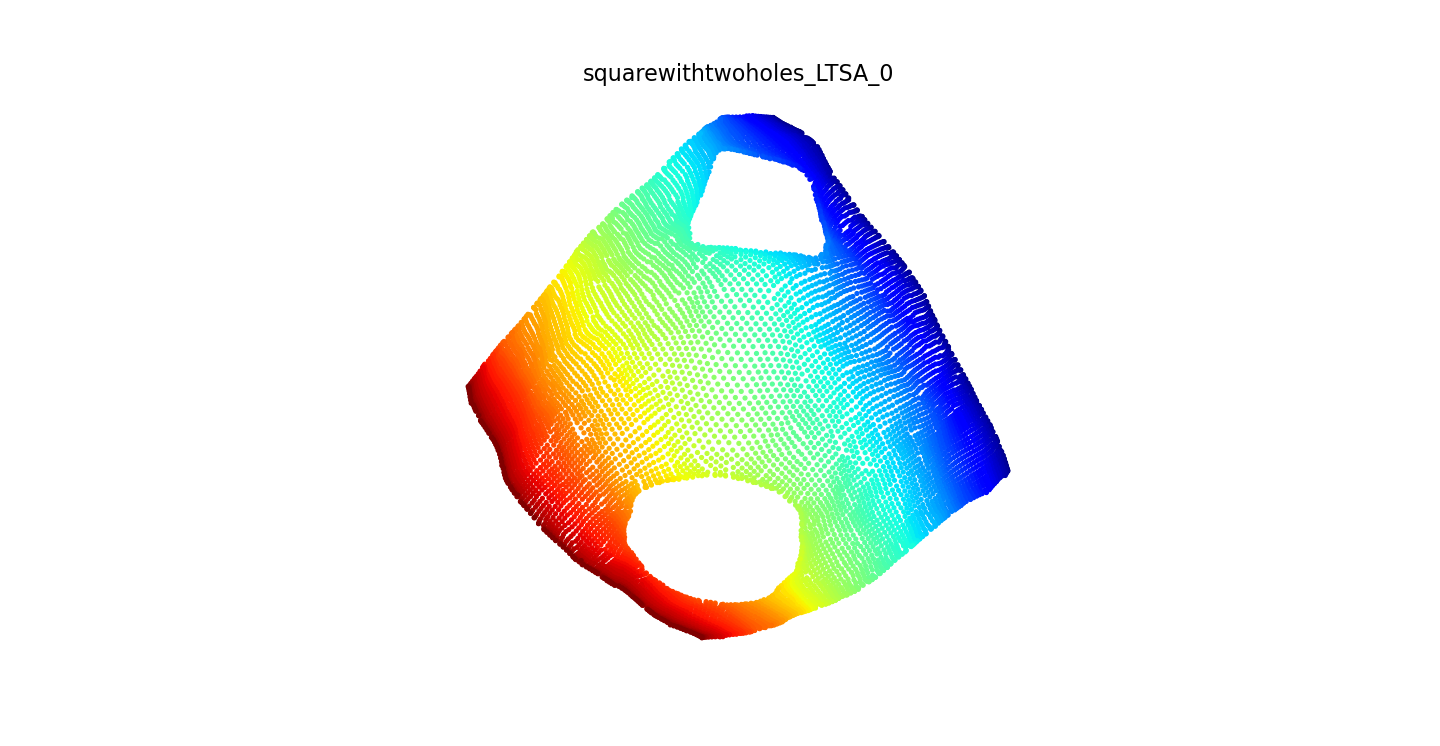}  & \includegraphics[trim=300 65 300 65, clip,width=0.15\textwidth,keepaspectratio,align=c]{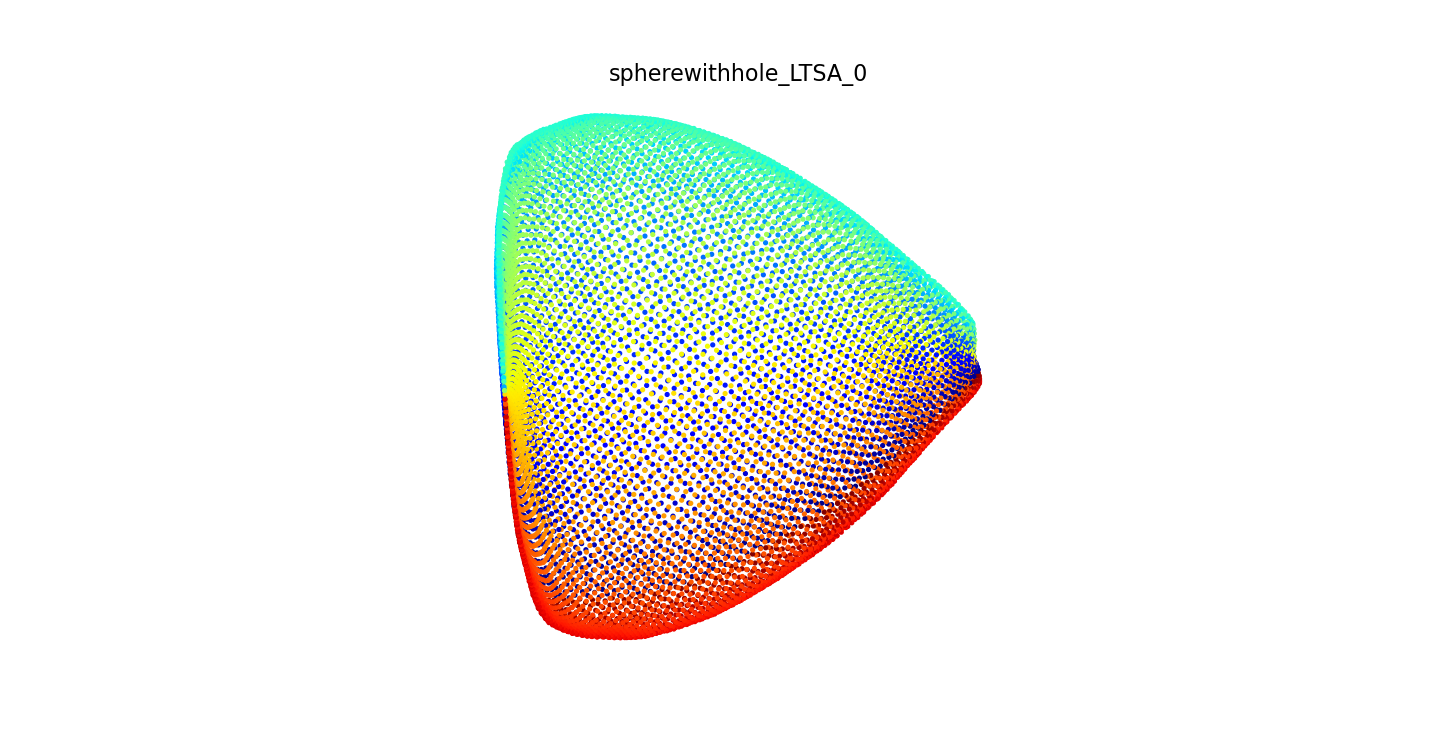} & \includegraphics[trim=355 65 330 65, clip,width=0.15\textwidth,keepaspectratio,align=c]{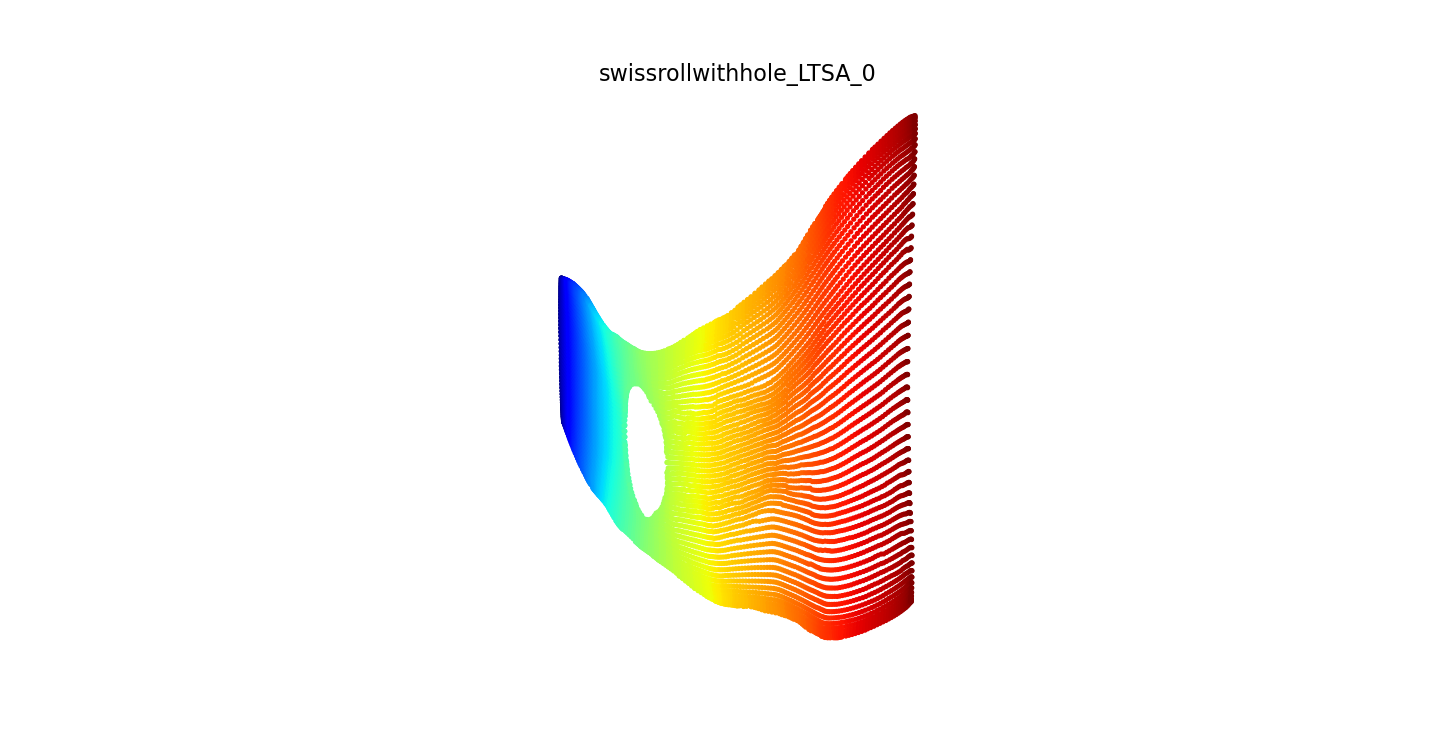}  & \includegraphics[trim=300 65 300 65, clip,width=0.15\textwidth,keepaspectratio,align=c]{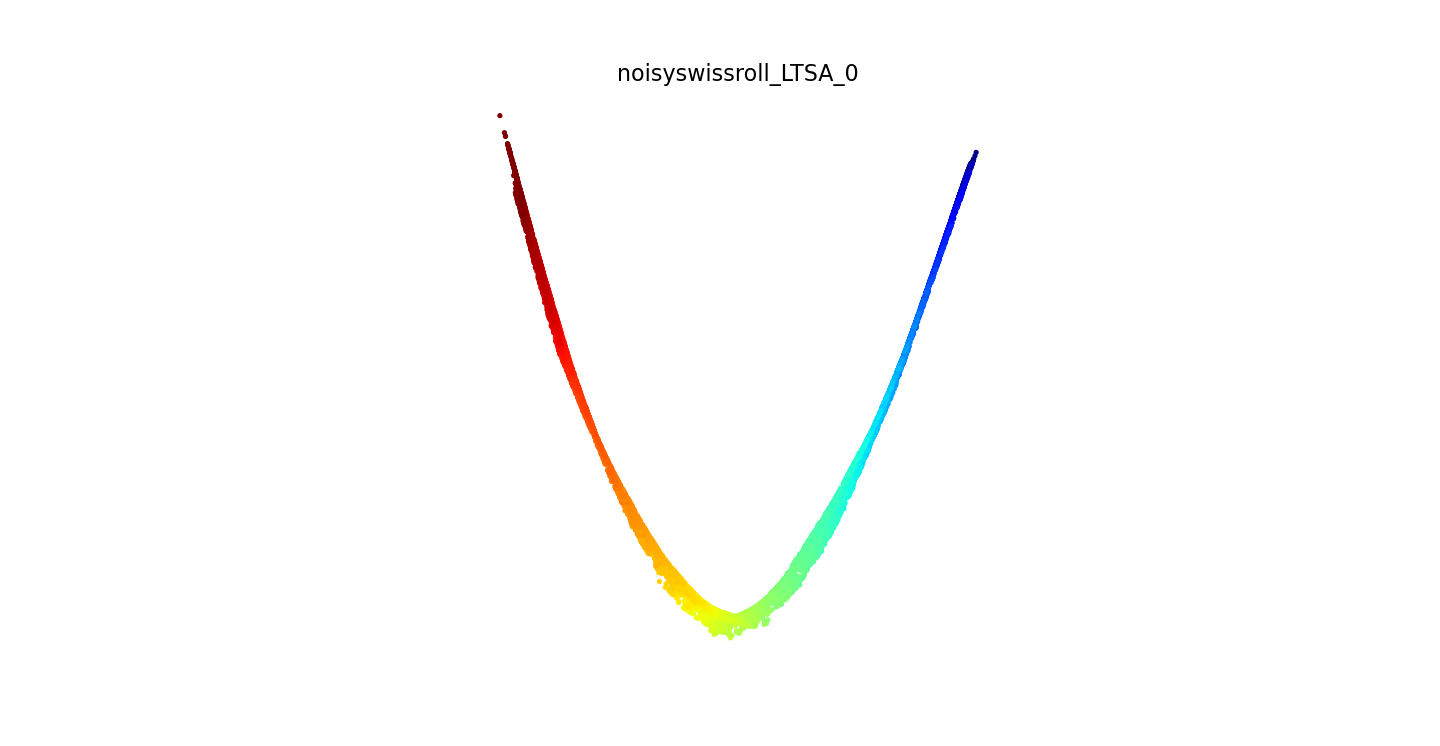}\\
        \hline
    \end{tabular}
    \caption{\editt{Embeddings obtained by using the global alignment procedure in LTSA to align the intermediate views in the embedding space. These views are the result of the clustering step in our algorithm.}}
    \label{fig:stitch_w_ltsa}
\end{figure}

From a practical standpoint, when the tearing of manifolds is not needed, one can use either procedure to align the intermediate views and obtain a global embedding. However, as shown in the Figure~\ref{fig:stitch_w_ltsa}, the embeddings produced by aligning our intermediate views using the alignment procedure in LTSA, are visually incorrect. The high distortion views near the boundary must be at cause here (see Figure~\ref{fig:fig41}). Since our alignment procedure works well on the same views as shown in Section~\ref{subsec:m_w_b}, this suggests that, compared to LTSA, our alignment procedure is more robust to the high distortion views. For similar reasons, one would expect LTSA to be less robust to the noisy data. This is indeed true as depicted in Figure~\ref{fig:noise}.

One advantage of using LTSA is the efficiency. LTSA reduces the optimal $Y$ to be the eigenvectors of a certain matrix leading to a fast algorithm. Our constraint does not allow such simplification and therefore we developed an iterative procedure by adapting GPA \citep{fabioprocrustes,gower1975generalized,ten1977orthogonal}. This procedure is slower than that in LTSA. We aim to improve the run-time in the subsequent versions of our code.

\section{Hyperparameters}
\label{sec:hyperparam}

\begin{table}[H]
    \renewcommand{\arraystretch}{2}
    \tiny
    \centering
    \begin{tabular}{|c|c|c|c|c|c|c|c|c|c|c|c|c|c|}
        \hline
        \makecell{\rotatebox[origin=c]{90}{Input}\\\\Algorithm} & Hyperparameters &
        \rotatebox[origin=c]{90}{Rectangle} &
        \rotatebox[origin=c]{90}{Barbell} &
        \rotatebox[origin=c]{90}{\makecell{Square with\\two holes}} &
        \rotatebox[origin=c]{90}{\makecell{Sphere with\\a hole}} &
        \rotatebox[origin=c]{90}{\makecell{Swissroll\\with a hole}} &
        \rotatebox[origin=c]{90}{\makecell{Noisy\\swissroll}} &
        \rotatebox[origin=c]{90}{Sphere} &
        \rotatebox[origin=c]{90}{\makecell{Curved\\torus}} &
        \rotatebox[origin=c]{90}{Flat torus} &
        \rotatebox[origin=c]{90}{\makecell{M\"obius strip}} &
        \rotatebox[origin=c]{90}{Klein Bottle}&
        \rotatebox[origin=c]{90}{\makecell{42-dim signal\\ strength data}}\\
        \hline
        LDLE & $\eta_{\text{min}}$ & 5 & 5 & 10 & 5 & 20 & 15 & 5 & 18 & 10 & 10 & 5 & 5\\
        \hline
        LTSA & n\_neighbors & 75 & 25 & 10 & 5 & 5 & 50 & 5 & 25 & 25 & 75 & 25 & 50\\
        \hline
        \multirow{2}{*}{UMAP} & n\_neighbors & 200 & 200 & 200 & 200 & 200 & 200 & 200 & 200 & 200 & 200 & 200 & 50 \\\cline{2-14}
        & min\_dist & 0.1 & 0.05 & 0.5 & 0.5 & 0.25 & 0.05 & 0.5 & 0.25 & 0.5 & 0.05 & 0.5 & 0.25\\\cline{2-14}
        \hline
        \multirow{2}{*}{t-SNE} & perplexity & 50 & 40 & 50 & 50 & 50 & 60 & 60 & 60 &60  & 60 & 50 & 60\\\cline{2-14}
        & exaggeration & 4 & 6 & 6 & 4 & 4 & 4 & 4 & 4 & 6 & 4 & 6 & 4\\\cline{2-14}
        \hline
        \multirow{2}{*}{\makecell{Laplacian\\Eigenmaps}} & $k_{\text{nn}}$ & - & - & 16 & - & - & - & - & - & - & - & - & 16\\\cline{2-14}
        & $k_{\text{tune}}$ & - & - & 7 & - & - & - & - & - & - & - & - & 7\\\cline{2-14}
        \hline
    \end{tabular}
    \caption{Hyperparameters used in the algorithms for the examples in Sections~\ref{subsec:m_w_b},~\ref{subsec:m_wo_b},~\ref{subsec:n_o_m} and~\ref{subsubsec:ssd}. For Laplacian eigenmaps, in all the examples except for square with two holes, all the searched values of the hyperparameters result in similar plots.}
    \label{tab:hypparam}
\end{table}

\begin{table}[H]
    \renewcommand{\arraystretch}{2}
    \tiny
    \centering
    \begin{tabular}{|c|c|c|c|c|}
        \hline
        \makecell{\rotatebox[origin=c]{90}{Noise}\\\\Algorithm} & Hyperparameters &
        $\sigma=0.01$&
        $\sigma=0.015$ &
        $\sigma=0.02$\\
        \hline
        LDLE & $\eta_{\text{min}}$ & 5 & 15 & 10\\
        \hline
        LTSA & n\_neighbors & 50 & 75 & 100\\
        \hline
        \multirow{2}{*}{UMAP} & n\_neighbors & 50 & 50 & 100 \\\cline{2-5}
        & min\_dist & 0.5 & 0.25 & 0.5\\\cline{2-5}
        \hline
        \multirow{2}{*}{t-SNE} & perplexity & 60 & 50 & 60 \\\cline{2-5}
        & exaggeration & 6 & 6 & 6\\\cline{2-5}
        \hline
    \end{tabular}
    \caption{Hyperparameters used in the algorithms for the Swiss Roll with increasing Gaussian noise (see Figure~\ref{fig:noise})}
    \label{tab:noise_hyperparam}
\end{table}

\begin{table}[H]
    \renewcommand{\arraystretch}{2}
    \tiny
    \centering
    \begin{tabular}{|c|c|c|c|c|c|}
        \hline
        \makecell{\rotatebox[origin=c]{90}{Resolution}\\\\Algorithm} & Hyperparameters &
        RES $=30$&RES $=15$&RES $=12$&RES $=10$\\
        \hline
        \multirow{4}{*}{LDLE} & $\eta_{\text{min}}$ & 3 & 3 & 3&3\\\cline{2-6}
        & $k_\text{tune}$ & 7 & 2 & 2 & 2\\\cline{2-6}
        & $N$ & 100 & 25 & 25 & 25\\\cline{2-6}
        & $k_{\text{lv}}$ & 7 & 4 & 4& 4\\\cline{2-6}
        \hline
        LTSA & n\_neighbors & 5 & 4 & 5 & 10\\
        \hline
        \multirow{2}{*}{UMAP} & n\_neighbors & 25 & 25 & 10 & 5 \\\cline{2-6}
        & min\_dist & 0.01 & 0.01 & 0.5 & 0.5\\\cline{2-5}
        \hline
        \multirow{2}{*}{t-SNE} & perplexity & 10 & 5 & 5 & 5 \\\cline{2-6}
        & exaggeration & 4 & 2 & 4 & 2\\\cline{2-6}
        \hline
    \end{tabular}
    \caption{Hyperparameters used in the algorithms for the Swiss Roll with increasing sparsity (see Figure~\ref{fig:sparse})}
    \label{tab:sparse_hyperparam}
\end{table}

\begin{table}[H]
    \renewcommand{\arraystretch}{2}
    \tiny
    \centering
    \begin{tabular}{|c|c|c|}
    \hline
    Method & \multicolumn{2}{c|}{Hyperparameters}\\
    \hline
    & face image data & Yoda-bulldog data\\
    \hline
    LDLE & \makecell{$N=25$, $k_{\text{lv}} = 12$, $\tau_s=5$, $\delta_s=0.25$ for all $s \in \{1,2\}$,\\ $\eta_{\text{min}}=4$, $\text{to\_tear}=$ False} & $N=25$, $\tau_s=10$, $\delta_s=0.5$ for all $s \in \{1,2\}$, $\eta_{\text{min}}=10$\\
    \hline
    LTSA & $\text{n\_neighbors}=10$ &$\text{n\_neighbors}=10$\\
    \hline
    UMAP & $\text{n\_neighbors}=50$, $\text{min\_dist}=0.01$ & $\text{n\_neighbors}=50$, $\text{min\_dist}=0.01$\\
    \hline
    t-SNE & $\text{perplexity}=60$, $\text{early\_exaggeration}=2$ & $\text{perplexity}=60$, $\text{early\_exaggeration}=2$\\
    \hline
    \end{tabular}
    \caption{Hyperparameters used in the algorithms for the face image data \citep{tenenbaum2000global} (see Figure~\ref{fig:face_data}) and the Yoda-bulldog dataset \citep{lederman2018learning} (see Figure~\ref{fig:s1_puppets}).}
    \label{tab:hypparam_face_yoda}
\end{table}



\section{Supplementary Figures}
\label{sec:supp}

\begin{figure}[h!]
    \centering
    \begin{tabular}{c}
        \includegraphics[ height=0.25\textwidth,keepaspectratio]{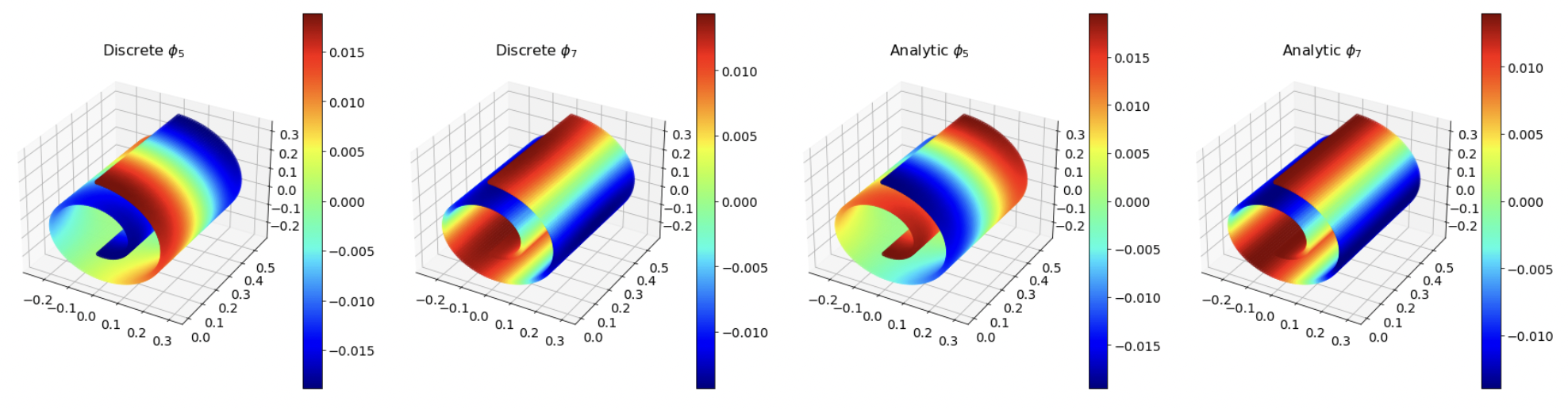} \\
        \includegraphics[height=0.25\textwidth,keepaspectratio]{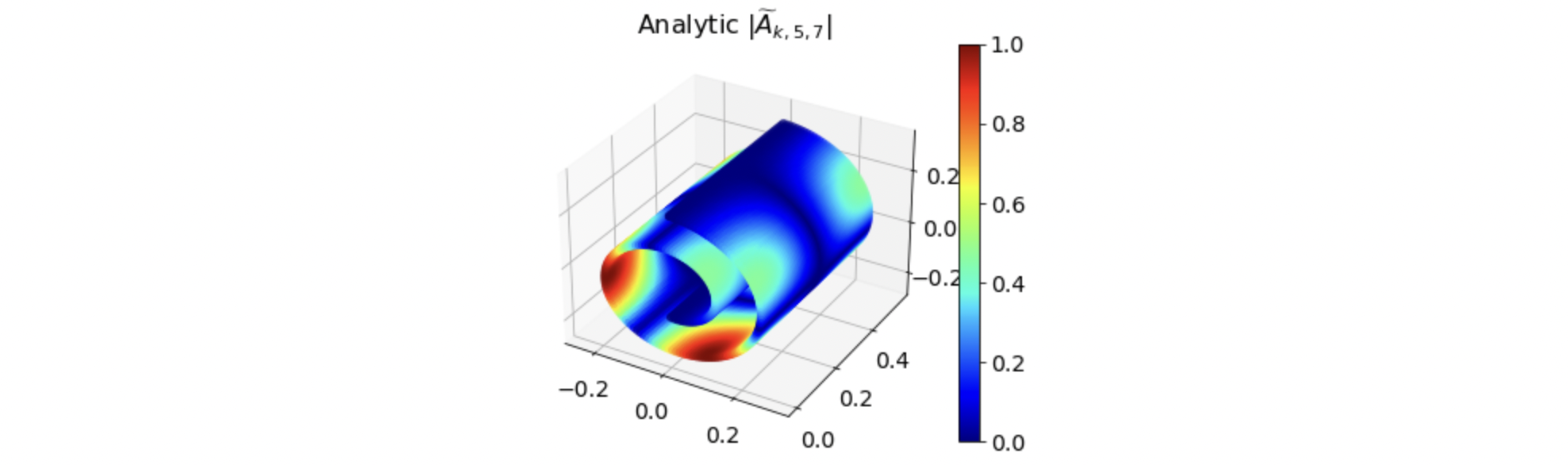} \\
        \includegraphics[height=0.5\textwidth,keepaspectratio]{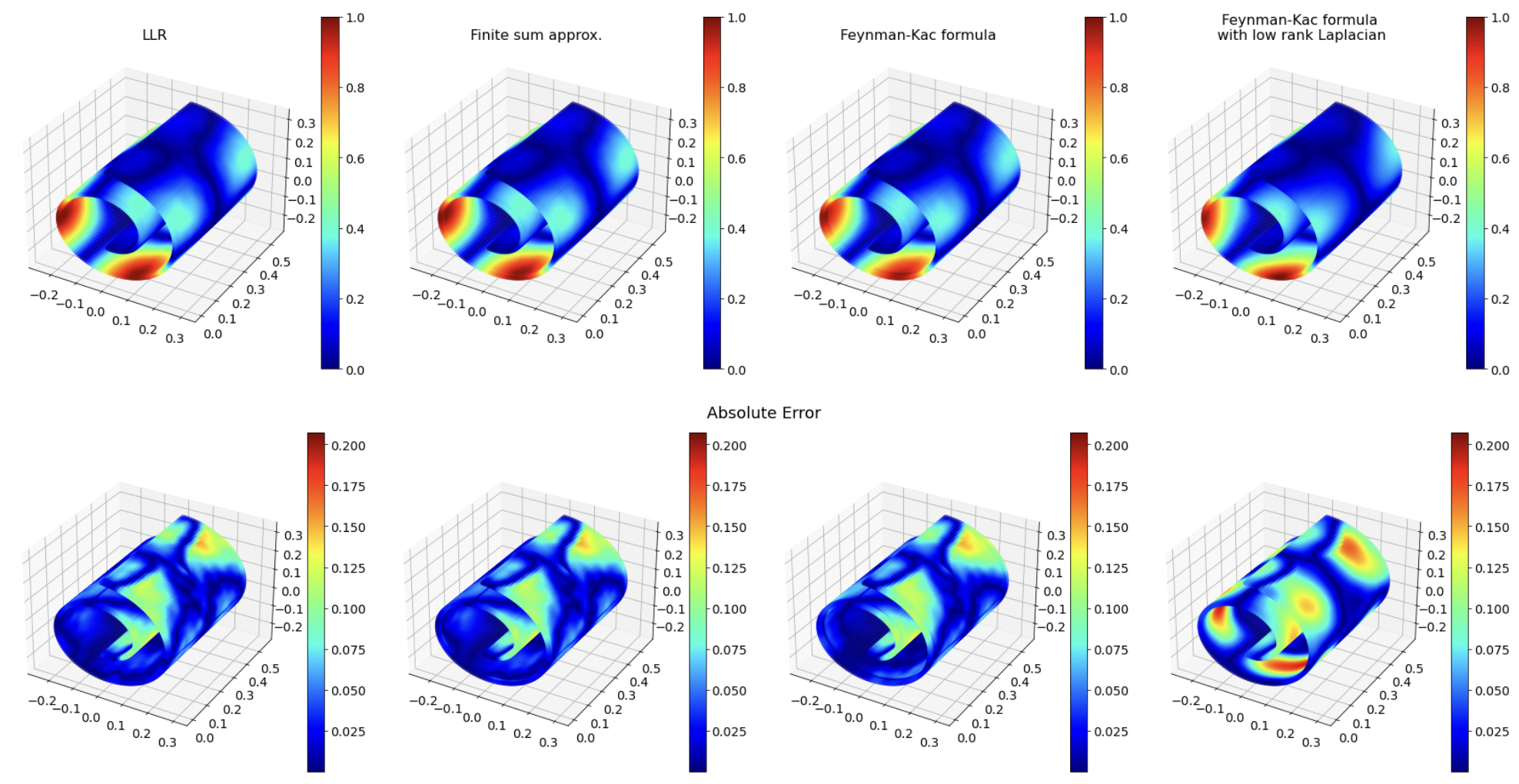}
    \end{tabular}
    \caption{\editt{Comparison of different techniques to estimate $\widetilde{A}_{kij}$ on a Swiss Roll with no noise, where $i=5$ and $j=7$. (first row) Analytical eigenfunctions and the obtained discrete eigenvectors are shown. (second row) Analytical value of $|\widetilde{A}_{kij}|$ is shown. Note that LDLE depends on the absolute values of $\widetilde{A}_{kij}$. (third row) Estimation of $|\widetilde{A}_{kij}|$ are shown due to Local Linear Regression based approach \citep{doi:10.1080/01621459.2013.827984}, finite sum approximation and Feynman-Kac formula based approaches as described in Section~\ref{subsec:Atilde_esimate} and a variant of the latter which uses low rank (of $100$) approximation of the graph Laplacian in Eq.~(\ref{Atildek_method2_discrete}). (fourth row) Absolute difference between the estimates and the analytical value. LLR, finite sum approx. and Feynman-Kac formula based approaches seem to perform slightly better.}}
    \label{fig:Atilde_compare_1}
\end{figure}

\begin{figure}[h!]
    \centering
    \begin{tabular}{c}
        \includegraphics[ height=0.25\textwidth,keepaspectratio]{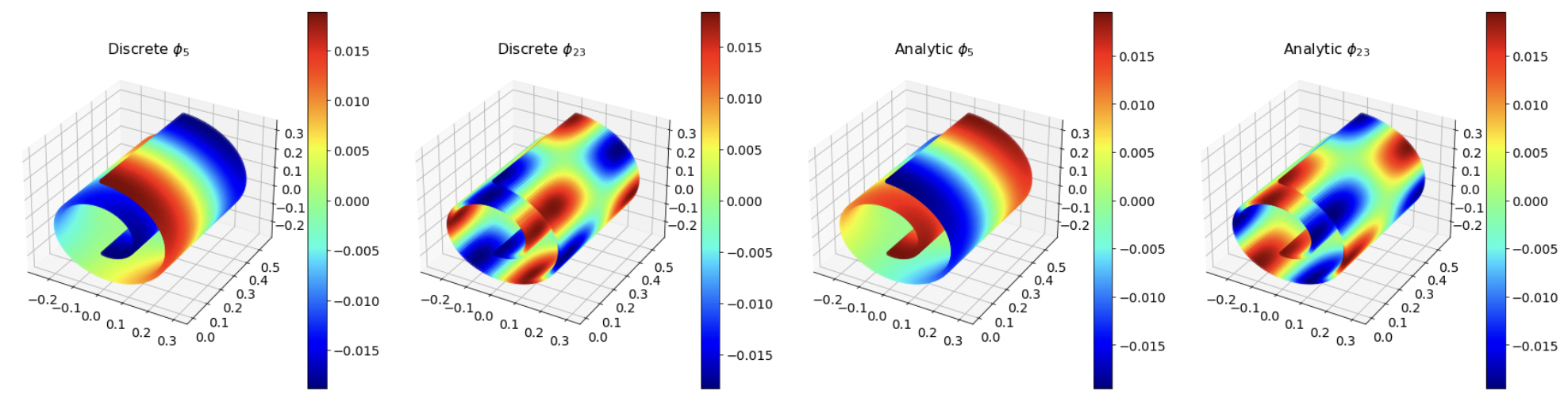} \\
        \includegraphics[height=0.25\textwidth,keepaspectratio]{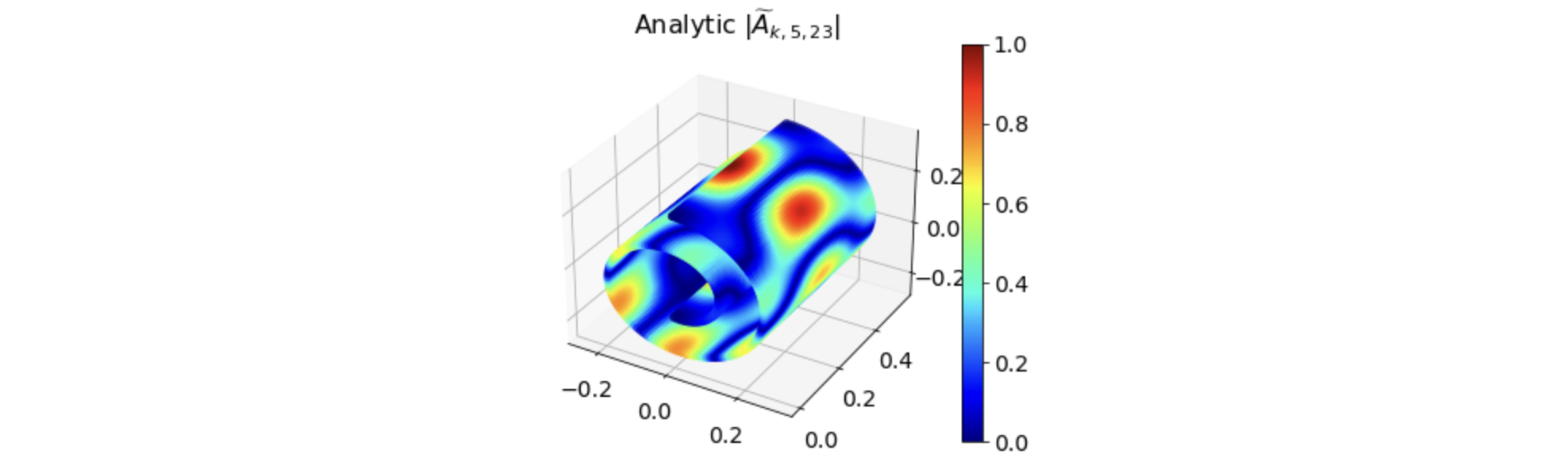} \\
        \includegraphics[height=0.5\textwidth,keepaspectratio]{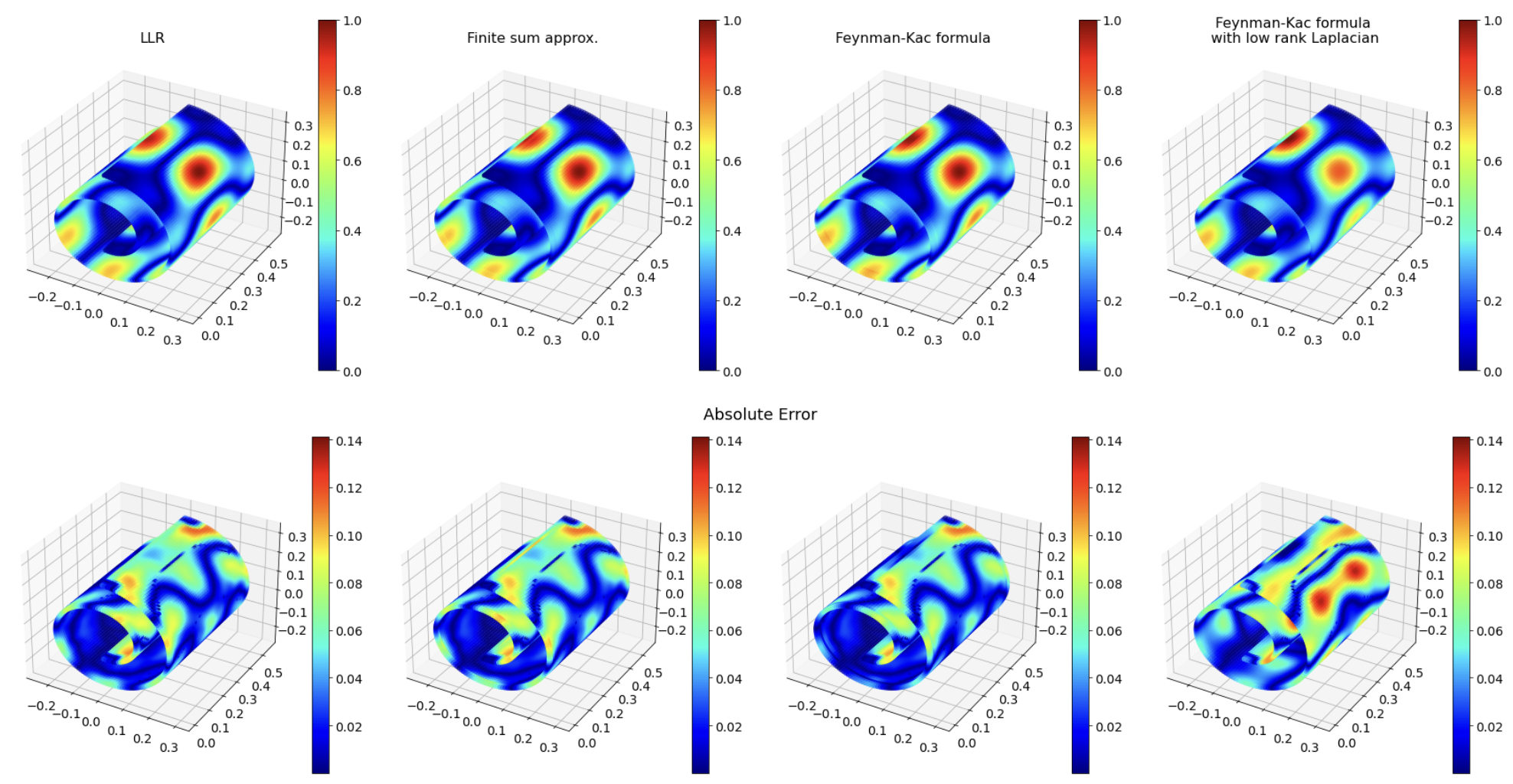}
    \end{tabular}
    \caption{\editt{Comparison of different techniques to estimate $\widetilde{A}_{kij}$ on a Swiss Roll with no noise, where $i=5$ and $j=23$. (first row) Analytical eigenfunctions and the obtained discrete eigenvectors are shown. (second row) Analytical value of $|\widetilde{A}_{kij}|$ is shown. Note that LDLE depends on the absolute values of $\widetilde{A}_{kij}$. (third row) Estimation of $|\widetilde{A}_{kij}|$ are shown due to Local Linear Regression based approach \citep{doi:10.1080/01621459.2013.827984}, finite sum approximation and Feynman-Kac formula based approaches as described in Section~\ref{subsec:Atilde_esimate} and a variant of the latter which uses low rank (of $100$) approximation of the graph Laplacian in Eq.~(\ref{Atildek_method2_discrete}). (fourth row) Absolute difference between the estimates and the analytical value. LLR, finite sum approx. and Feynman-Kac formula based approaches seem to perform slightly better.}}
    \label{fig:Atilde_compare_2}
\end{figure}

\begin{figure}[h!]
    \centering
    \begin{tabular}{c}
        \includegraphics[ height=0.25\textwidth,keepaspectratio]{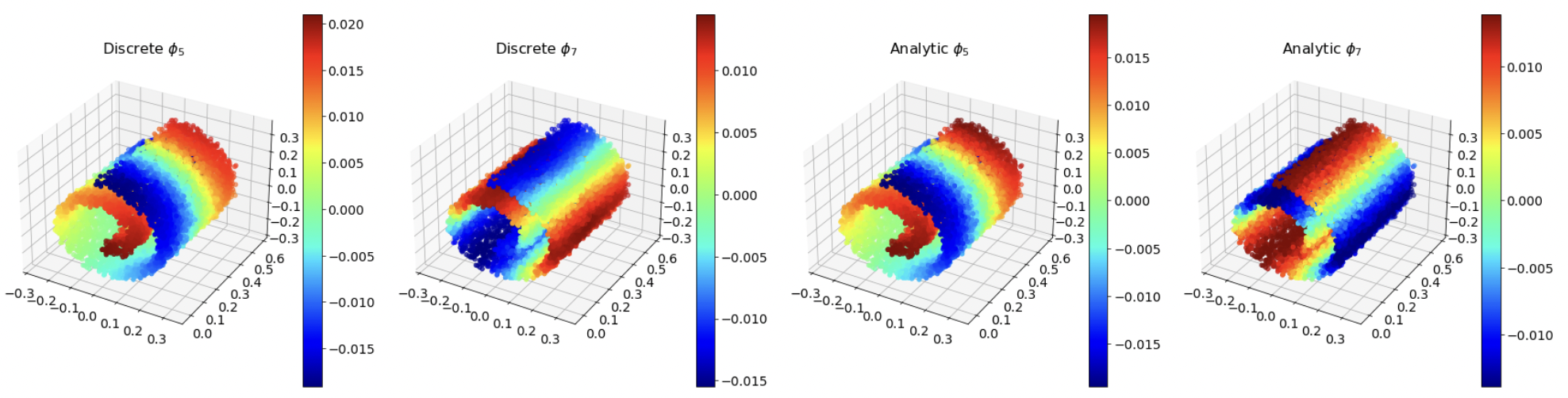} \\
        \includegraphics[height=0.25\textwidth,keepaspectratio]{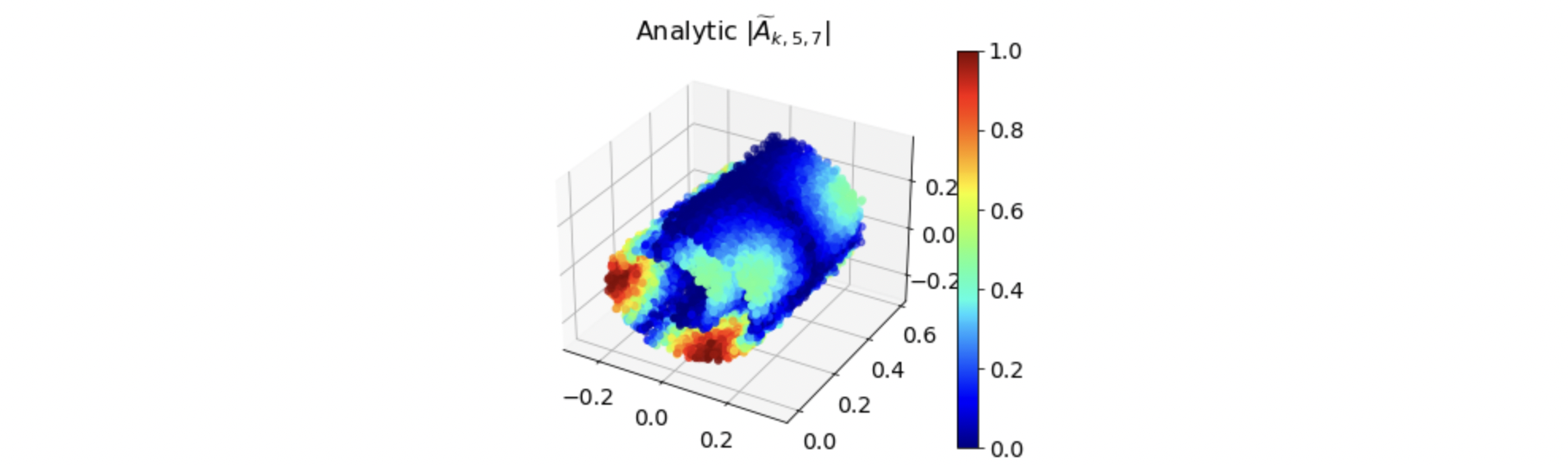} \\
        \includegraphics[height=0.5\textwidth,keepaspectratio]{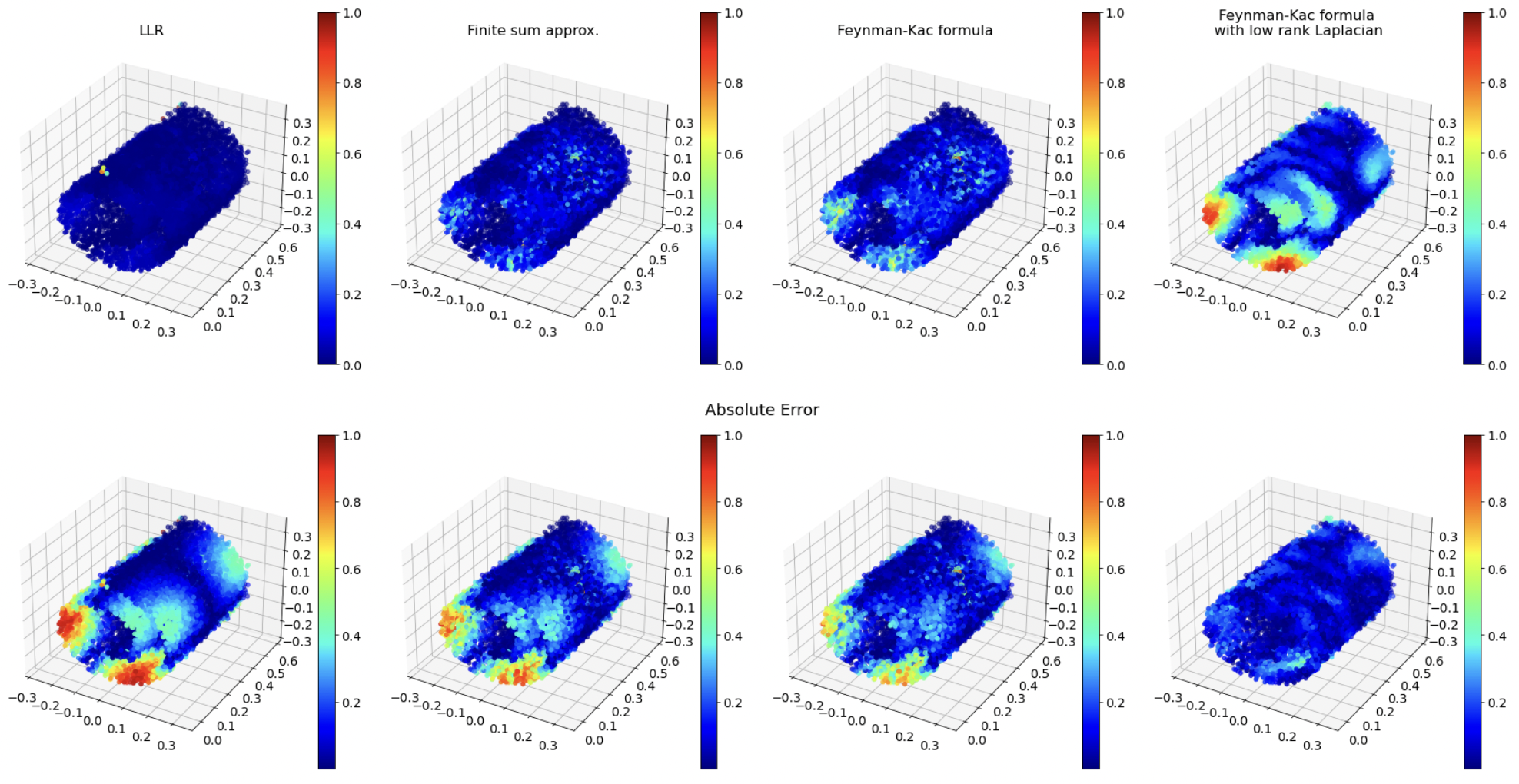}
    \end{tabular}
    \caption{\editt{Comparison of different techniques to estimate $\widetilde{A}_{kij}$ on a Swiss Roll with Gaussian noise of variance $10^{-4}$, where $i=5$ and $j=7$. (first row) Analytical eigenfunctions obtained for the noiseless version of the Swiss Roll, and the obtained discrete eigenvectors are shown. (second row) Analytical value of $|\widetilde{A}_{kij}|$ is shown. Note that LDLE depends on the absolute values of $\widetilde{A}_{kij}$. (third row) Estimation of $|\widetilde{A}_{kij}|$ are shown due to Local Linear Regression based approach \citep{doi:10.1080/01621459.2013.827984}, finite sum approximation and Feynman-Kac formula based approaches as described in Section~\ref{subsec:Atilde_esimate} and a variant of the latter which uses low rank (of $100$) approximation of the graph Laplacian in Eq.~(\ref{Atildek_method2_discrete}). (fourth row) Absolute difference between the estimates and the analytical value.  The Feynman-Kac formula based approach which uses low rank approximation of $L$ seem to perform the best while the LLR based approach produced high error.}}
    \label{fig:Atilde_compare_3}
\end{figure}

\begin{figure}[h!]
    \centering
    \begin{tabular}{c}
        \includegraphics[ height=0.25\textwidth,keepaspectratio]{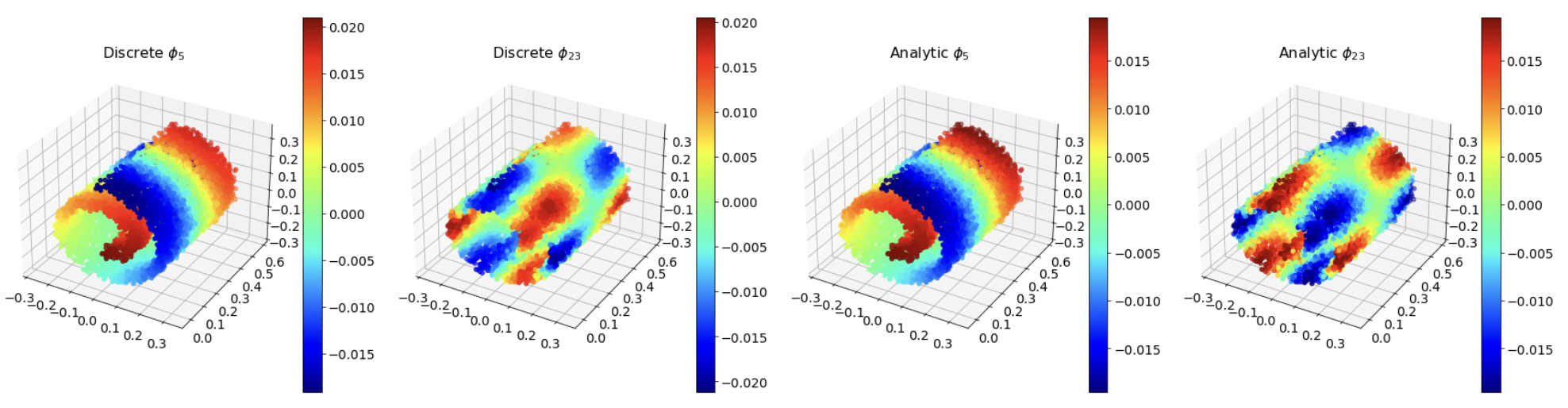} \\
        \includegraphics[height=0.25\textwidth,keepaspectratio]{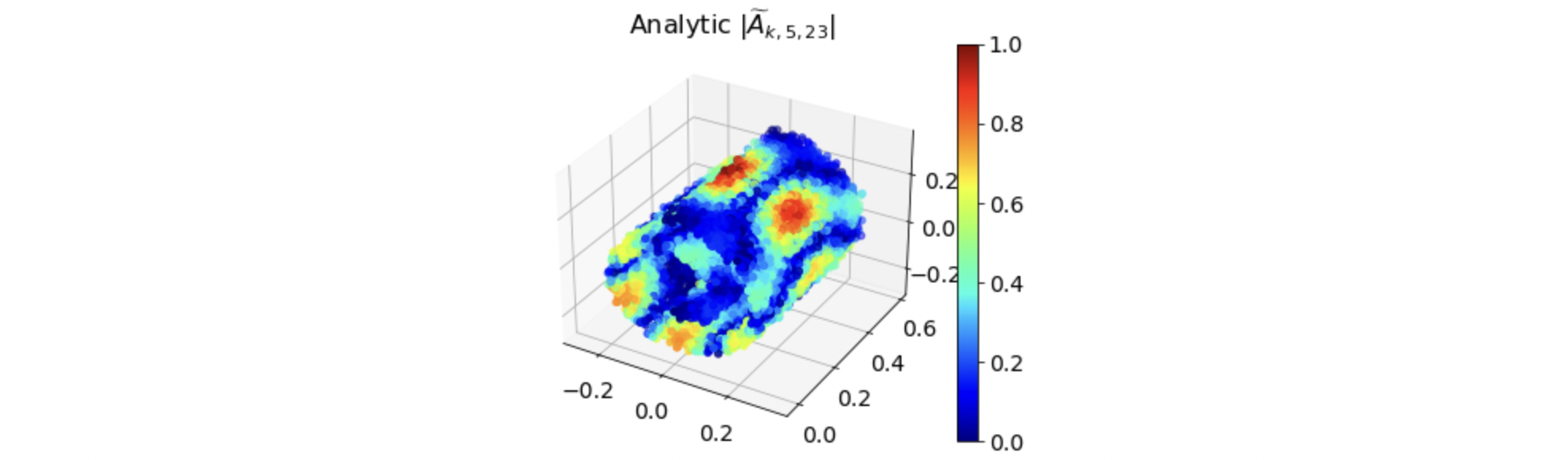} \\
        \includegraphics[height=0.5\textwidth,keepaspectratio]{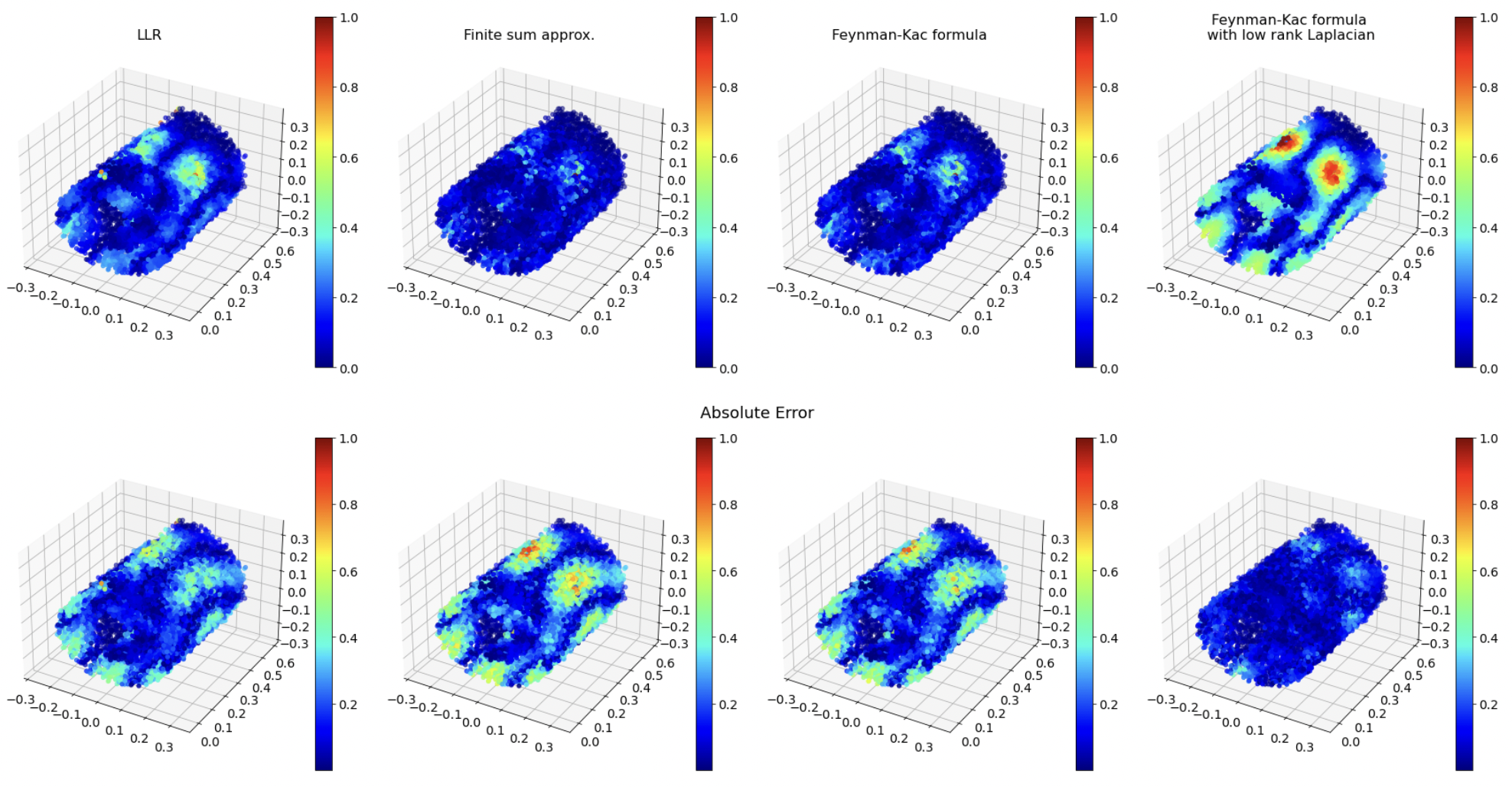}
    \end{tabular}
    \caption{\editt{Comparison of different techniques to estimate $\widetilde{A}_{kij}$ on a Swiss Roll with Gaussian noise of variance $10^{-4}$, where $i=5$ and $j=23$. (first row) Analytical eigenfunctions obtained for the noiseless version of the Swiss Roll, and the obtained discrete eigenvectors are shown. (second row) Analytical value of $|\widetilde{A}_{kij}|$ is shown. Note that LDLE depends on the absolute values of $\widetilde{A}_{kij}$. (third row) Estimation of $|\widetilde{A}_{kij}|$ are shown due to Local Linear Regression based approach \citep{doi:10.1080/01621459.2013.827984}, finite sum approximation and Feynman-Kac formula based approaches as described in Section~\ref{subsec:Atilde_esimate} and a variant of the latter which uses low rank (of $100$) approximation of the graph Laplacian in Eq.~(\ref{Atildek_method2_discrete}). (fourth row) Absolute difference between the estimates and the analytical value. The Feynman-Kac formula based approach which uses low rank approximation of $L$ seem to perform the best while the errors due to other three approaches are somewhat similar.}}
    \label{fig:Atilde_compare_4}
\end{figure}

\begin{figure}[!h]
    \begin{tabular}{c}
        \includegraphics[trim=0 300 0 300, clip,width=1\textwidth,height=0.5\textwidth,keepaspectratio]{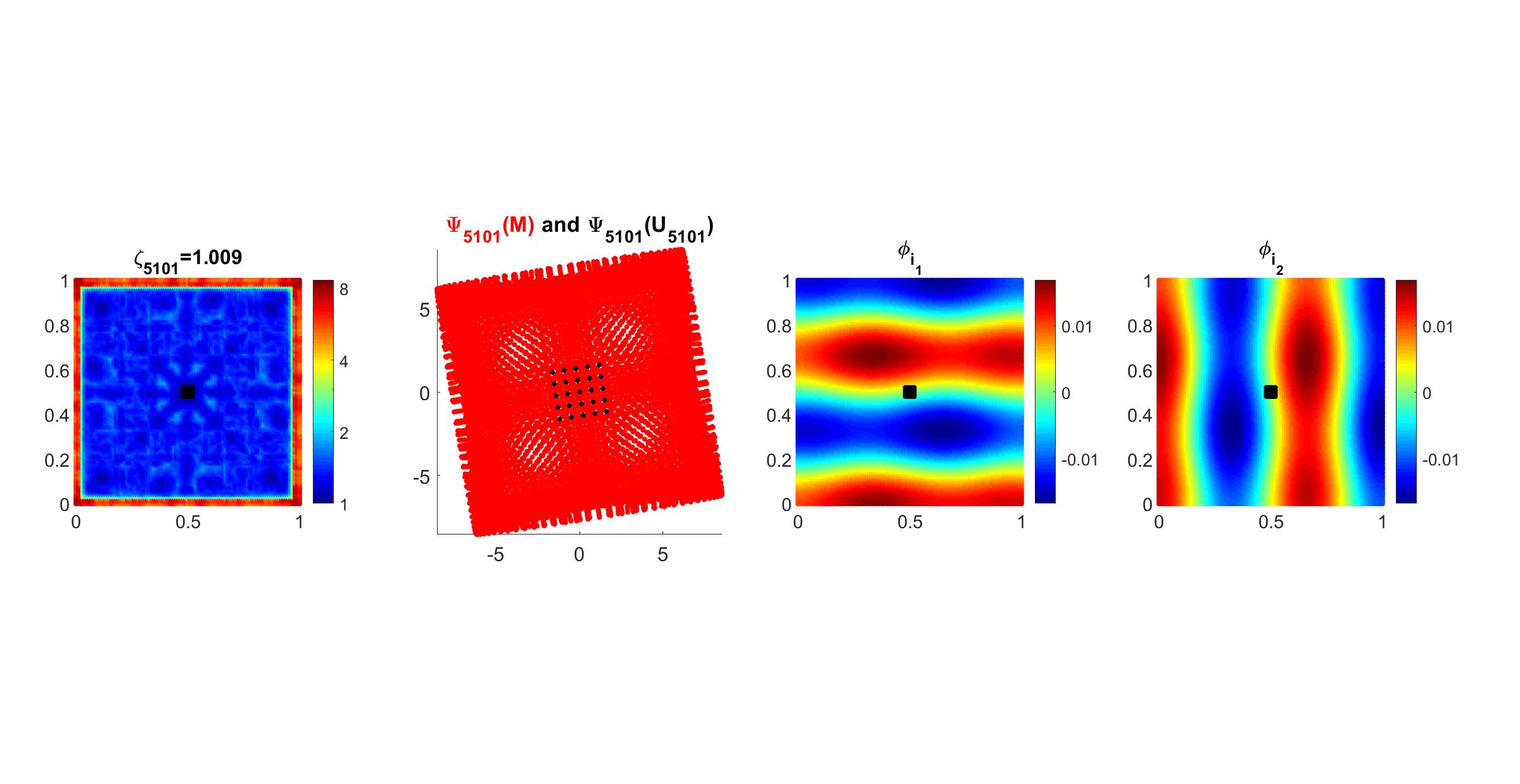}\\
        \includegraphics[trim=0 300 0 300, clip,width=1\textwidth,height=0.5\textwidth,keepaspectratio]{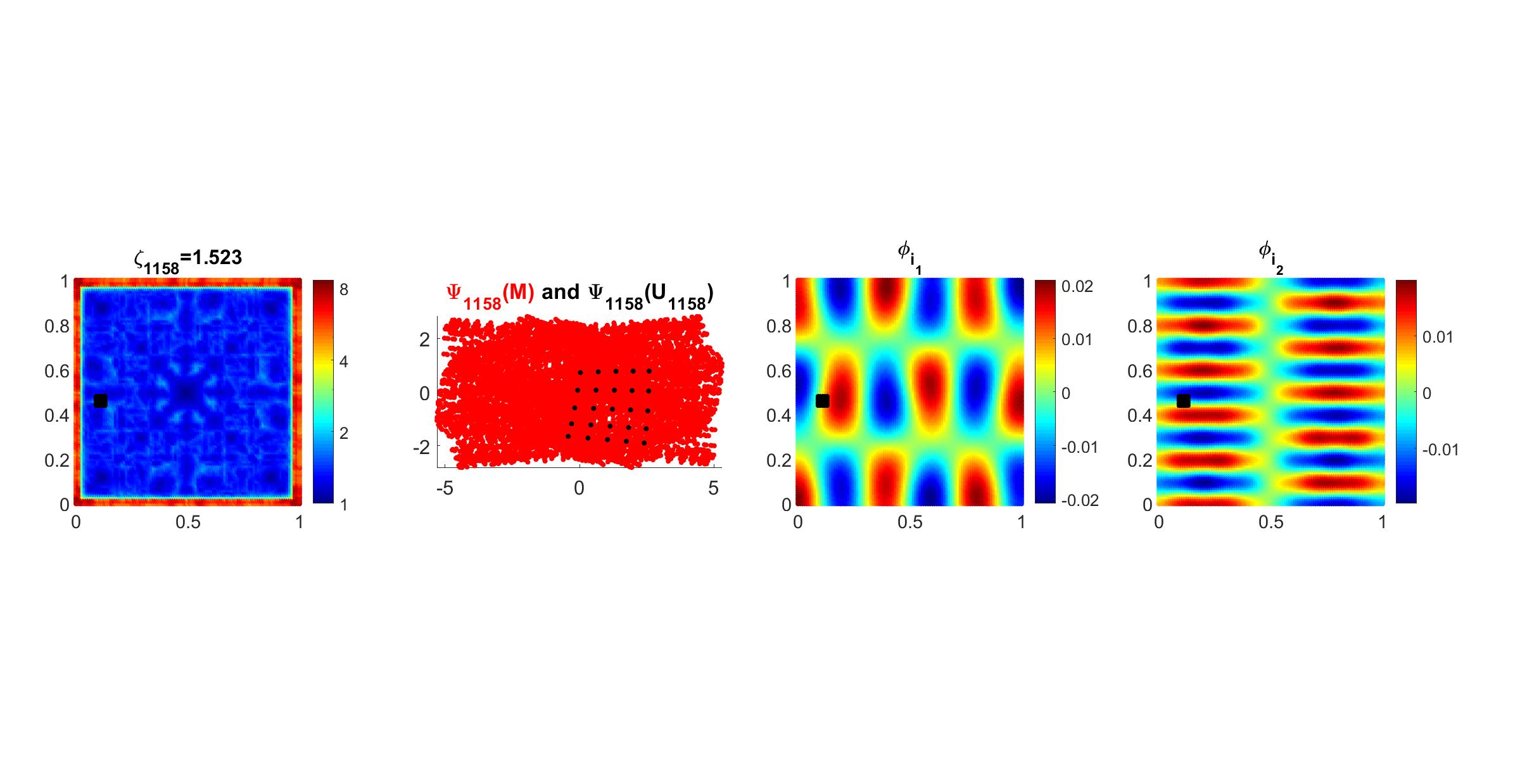}\\
        \includegraphics[trim=0 300 0 300, clip,width=1\textwidth,height=0.5\textwidth,keepaspectratio]{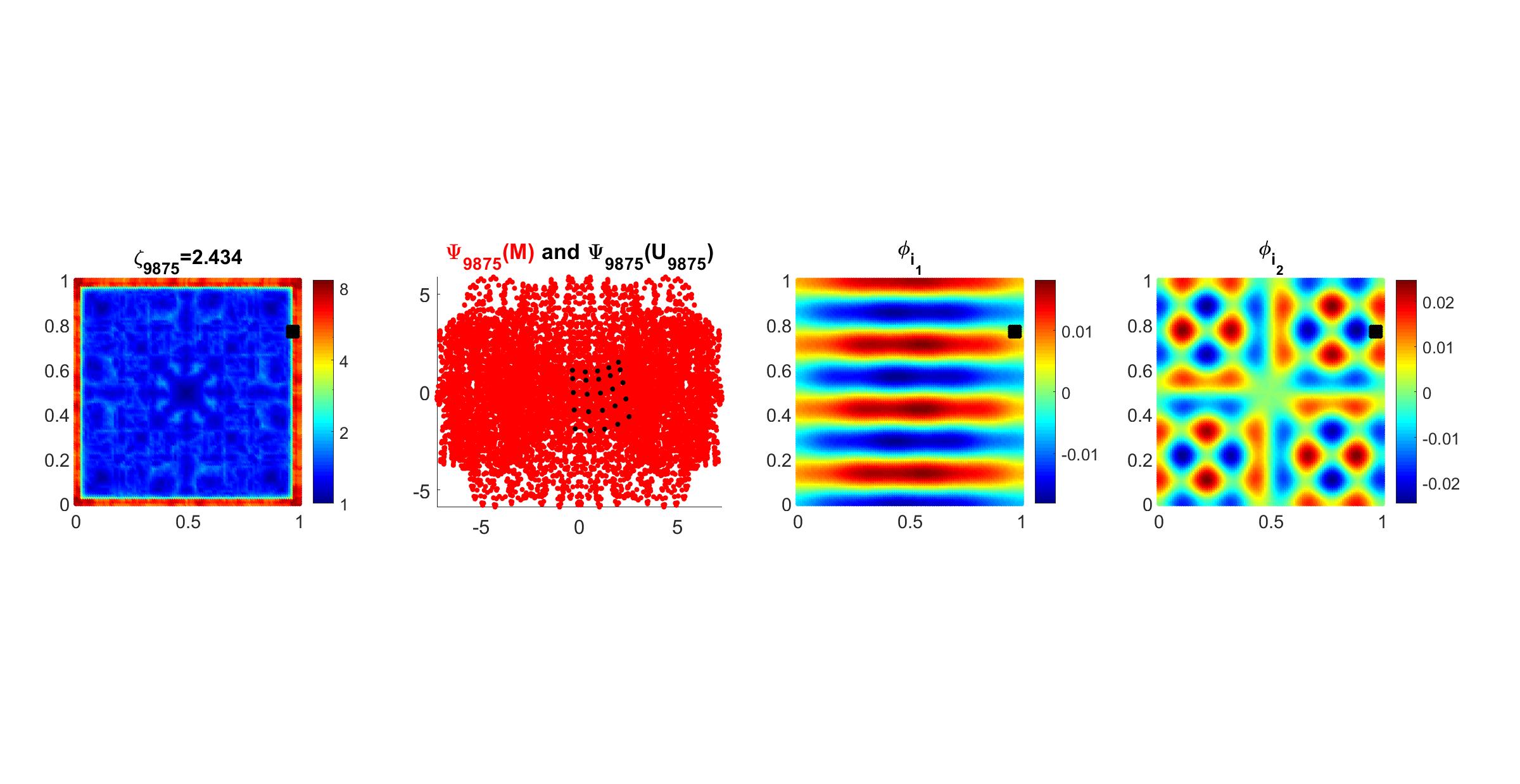}\\
        \includegraphics[trim=0 300 0 300, clip,width=1\textwidth,height=0.5\textwidth,keepaspectratio]{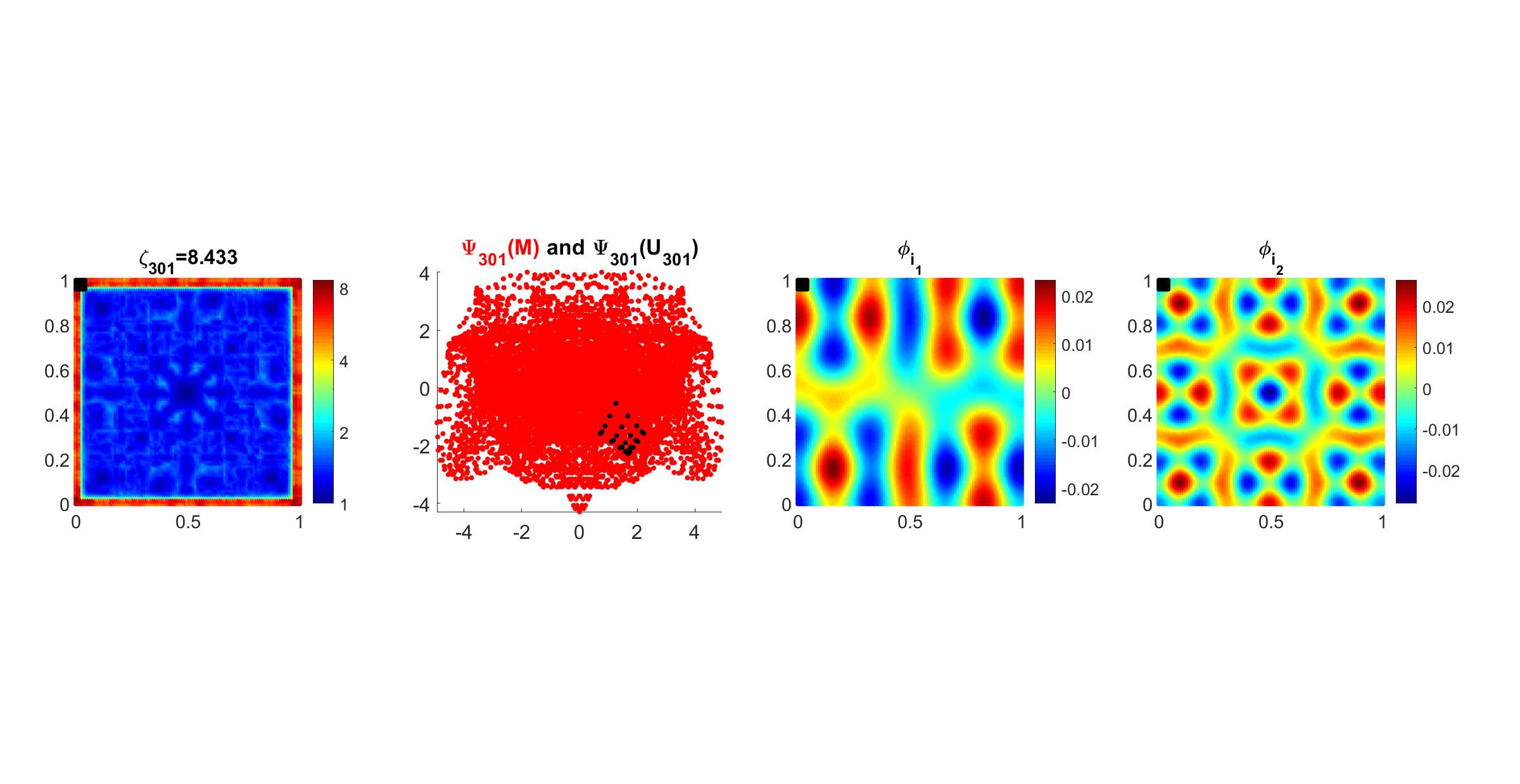}\\
    \end{tabular}
    \caption{(first column) Input square grid is shown. The points $x_k$ are colored by the distortion $\zeta_{kk}$ of the obtained local parameterizations $\Phi_k$ on the neighborhood $U_k$ surrounding them. A local view $U_{k_0}$ around $x_{k_0}$ for a fixed $k_0$ is also shown in black. (second column) The corresponding local view in the embedding space $\Phi_{k_0}(U_{k_0})$ is shown in black. Although of no significance to our algorithm, for visualization purpose, the embedding of the square due to $\Phi_{k_0}$, $\Phi_{k_0}(M)$, is shown in red. (third and fourth columns) The eigenvectors $\bm{\phi}_{i_1}$ and $\bm{\phi}_{i_2}$ chosen for the construction of $\Phi_{k_0}$ are shown. Points in $U_{k_0}$ are again colored in black. Note that the gradient of these eigenvectors are close to being orthogonal in the vicinity of $U_{k_0}$ and in particular, at $x_{k_0}$.}
    \label{fig:local_views}
\end{figure}

\begin{figure}[!h]
    \centering
    \begin{tabular}{M{0.005\textwidth}|c|c|}
        & \tiny{LDLE with arrows} &\tiny{Derived cut and paste diagrams}\\
        \hline 
        \rotatebox[origin=l]{90}{\tiny{Sphere}} & \includegraphics[width=0.2\textwidth,height=0.2\textwidth,keepaspectratio,align=c]{fig2/sphere/arrow2.png}  &  \includegraphics[width=0.5\textwidth,height=0.2\textwidth,keepaspectratio,align=c]{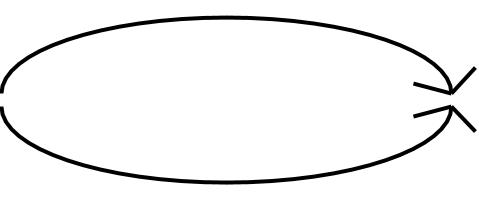}
        \\
        \hline
        \rotatebox[origin=l]{90}{\tiny{Sphere with a hole}} & \includegraphics[width=0.2\textwidth,height=0.2\textwidth,keepaspectratio,align=c]{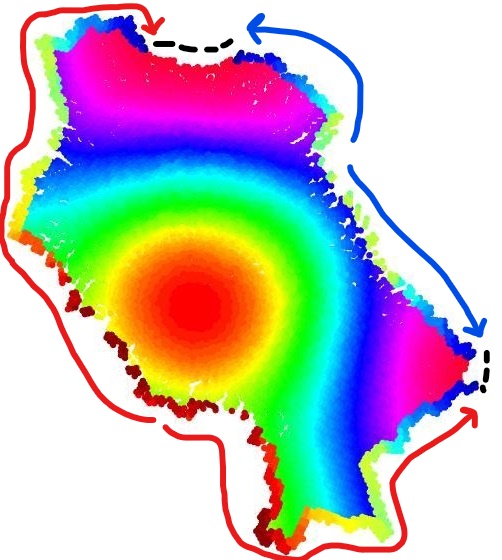}   &  \includegraphics[width=0.5\textwidth,height=0.2\textwidth,keepaspectratio,align=c]{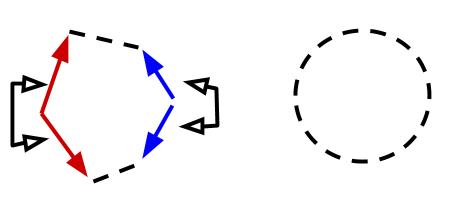}
        \\
        \hline
        \rotatebox[origin=l]{90}{\tiny{Curved torus}} & \includegraphics[width=0.2\textwidth,height=0.2\textwidth,keepaspectratio,align=c]{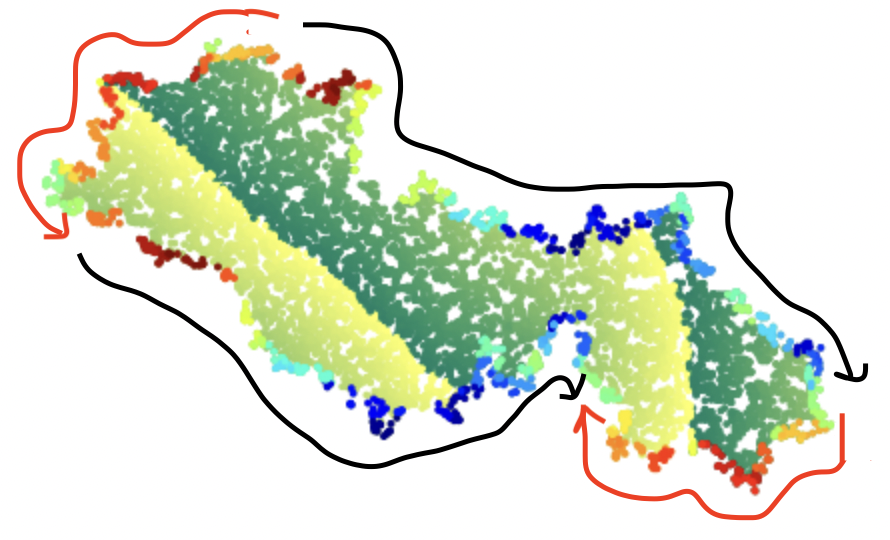} &  \includegraphics[width=0.5\textwidth,height=0.2\textwidth,keepaspectratio,align=c]{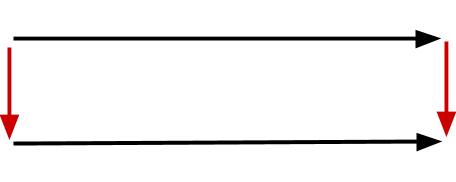}\\
        \hline
        \rotatebox[origin=l]{90}{\tiny{Flat torus}} & \includegraphics[width=0.2\textwidth,height=0.2\textwidth,keepaspectratio,align=c]{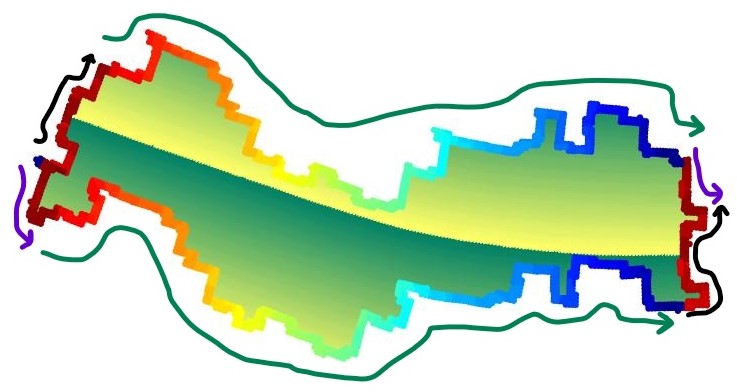} &  \includegraphics[width=0.5\textwidth,height=0.2\textwidth,keepaspectratio,align=c]{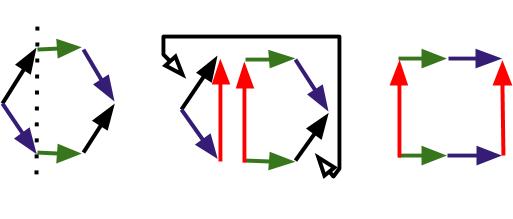}\\
        \hline
        \rotatebox[origin=l]{90}{\tiny{M\"obius strip}} & \includegraphics[width=0.2\textwidth,height=0.2\textwidth,keepaspectratio,align=c]{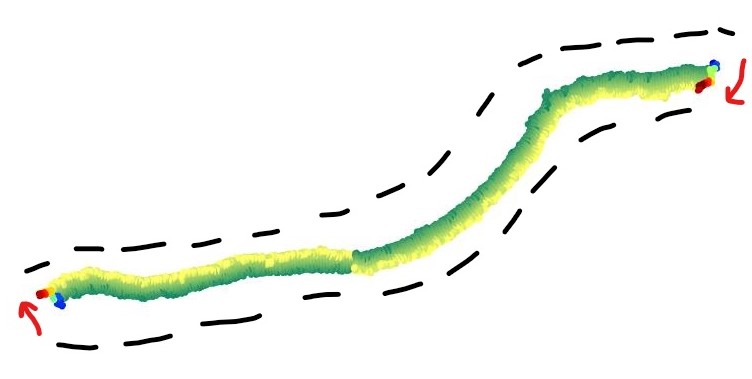} &  \includegraphics[width=0.5\textwidth,height=0.2\textwidth,keepaspectratio,align=c]{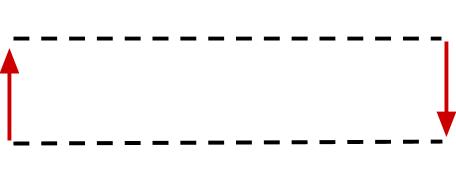}\\
        \hline
        \rotatebox[origin=l]{90}{\tiny{Klein bottle}} & \includegraphics[width=0.2\textwidth,height=0.2\textwidth,keepaspectratio,align=c]{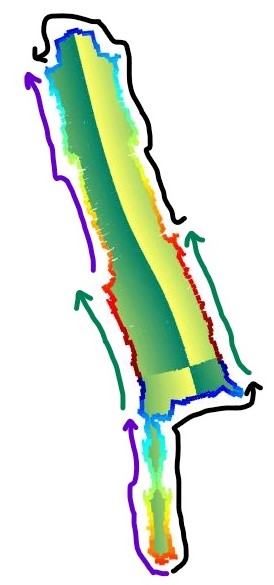} &  \includegraphics[width=0.5\textwidth,height=0.2\textwidth,keepaspectratio,align=c]{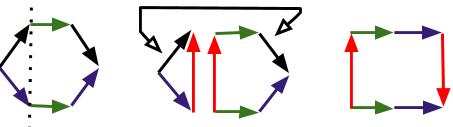}\\
        \hline
    \end{tabular}
    \caption{(Left) LDLE embedding with arrows drawn by tracing the colored boundary. (Right) Derived cut and paste diagrams to prove the correctness of the embedding. Pieces of the boundary represented by filled arrows of the same color are to be stitched together. Pieces of the boundary represented by black dashed lines are not to be stitched. Dotted lines and shallow arrows represent cut and paste instructions, respectively.}
    \label{fig:decipher}
\end{figure}

\begin{figure}[h!]  
    \centering
    \includegraphics[width=0.5\textwidth,keepaspectratio,align=c]{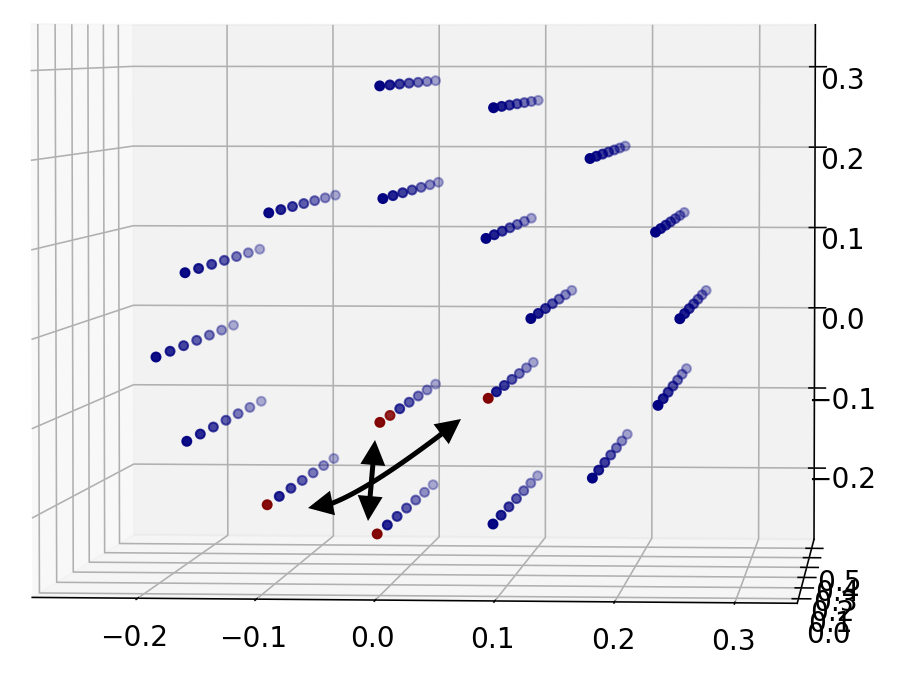}
    \caption{\editt{In Figure~\ref{fig:sparse}, for the case when RES $=10$, certain points on the opposite sides of the gap between the Swiss Roll are neighbors in the ambient space. These points are shown in red.}}
    \label{fig:res10issue}
\end{figure}

\clearpage

\begin{thebibliography}{10}

\bibitem{aswani2011regression}
Anil Aswani, Peter Bickel, Claire Tomlin, et~al.
\newblock Regression on manifolds: Estimation of the exterior derivative.
\newblock {\em The Annals of Statistics}, 39(1):48--81, 2011.

\bibitem{bastian_bechtold_2021_4559847}
Bastian Bechtold, Patrick Fletcher, Seamus Holden, and Srinivas
  Gorur-Shandilya.
\newblock {\em {bastibe Violinplot-Matlab: A Good Starting Point}}.
\newblock github, 2021.

\bibitem{belkin2003laplacian}
Mikhail Belkin and Partha Niyogi.
\newblock {Laplacian eigenmaps for dimensionality reduction and data
  representation}.
\newblock {\em Neural computation}, 15(6):1373--1396, 2003.

\bibitem{BELKIN20081289}
Mikhail Belkin and Partha Niyogi.
\newblock {Towards a theoretical foundation for Laplacian-based manifold
  methods}.
\newblock {\em Journal of Computer and System Sciences}, 74(8):1289--1308,
  2008.
\newblock Learning Theory 2005.

\bibitem{berry2017density}
Tyrus Berry and Timothy Sauer.
\newblock {Density estimation on manifolds with boundary}.
\newblock {\em Computational Statistics \& Data Analysis}, 107:1--17, 2017.

\bibitem{DBLP:journals/corr/BlauM16}
Yochai Blau and Tomer Michaeli.
\newblock {Non-Redundant Spectral Dimensionality Reduction}.
\newblock {\em CoRR}, abs/1612.03412, 2016.

\bibitem{canzani2013analysis}
Yaiza Canzani.
\newblock Analysis on manifolds via the laplacian.
\newblock {\em Lecture Notes available at: http://www. math. harvard.
  edu/canzani/docs/Laplacian. pdf}, 2013.

\bibitem{NEURIPS2019_6a10bbd4}
Yu-Chia Chen and Marina Meila.
\newblock {Selecting the independent coordinates of manifolds with large aspect
  ratios}.
\newblock In {\em Advances in Neural Information Processing Systems},
  volume~32, pages 1088--1097. Curran Associates, Inc., 2019.

\bibitem{doi:10.1080/01621459.2013.827984}
Ming-yen Cheng and Hau-tieng Wu.
\newblock {Local Linear Regression on Manifolds and Its Geometric
  Interpretation}.
\newblock {\em {Journal of the American Statistical Association}},
  108(504):1421--1434, 2013.

\bibitem{cheng2020spectral}
Xiuyuan Cheng and Gal Mishne.
\newblock Spectral embedding norm: Looking deep into the spectrum of the graph
  laplacian.
\newblock {\em SIAM Journal on Imaging Sciences}, 13(2):1015--1048, 2020.

\bibitem{cheng2018geometry}
Xiuyuan Cheng, Gal Mishne, and Stefan Steinerberger.
\newblock The geometry of nodal sets and outlier detection.
\newblock {\em Journal of Number Theory}, 185:48--64, 2018.

\bibitem{cheng2020convergence}
Xiuyuan Cheng and Hau-Tieng Wu.
\newblock {Convergence of graph laplacian with knn self-tuned kernels}.
\newblock {\em arXiv preprint arXiv:2011.01479}, 2020.

\bibitem{cheng2021eigen}
Xiuyuan Cheng and Nan Wu.
\newblock {Eigen-convergence of Gaussian kernelized graph Laplacian by manifold
  heat interpolation}.
\newblock {\em arXiv preprint arXiv:2101.09875}, 2021.

\bibitem{NEURIPS2018_4c5bcfec}
Leena Chennuru~Vankadara and Ulrike von Luxburg.
\newblock {Measures of distortion for machine learning}.
\newblock In {S. Bengio and H. Wallach and H. Larochelle and K. Grauman and N.
  Cesa-Bianchi and R. Garnett}, editor, {\em {Advances in Neural Information
  Processing Systems}}, volume~31. {Curran Associates, Inc.}, 2018.

\bibitem{7148963}
A.~{Cloninger} and W.~{Czaja}.
\newblock {Eigenvector localization on data-dependent graphs}.
\newblock In {\em {2015 International Conference on Sampling Theory and
  Applications (SampTA)}}, pages 608--612, 2015.

\bibitem{doi:10.1080/10586458.2018.1538911}
Alexander Cloninger and Stefan Steinerberger.
\newblock {On the Dual Geometry of Laplacian Eigenfunctions}.
\newblock {\em Experimental Mathematics}, 0(0):1--11, 2018.

\bibitem{coifman2006diffusion}
Ronald~R Coifman and St{\'e}phane Lafon.
\newblock {Diffusion maps}.
\newblock {\em Applied and computational harmonic analysis}, 21(1):5--30, 2006.

\bibitem{fabioprocrustes}
Fabio Crosilla and Alberto Beinat.
\newblock {Use of generalised Procrustes analysis for the photogrammetric block
  adjustment by independent models}.
\newblock {\em ISPRS Journal of Photogrammetry and Remote Sensing},
  56:195--209, 04 2002.

\bibitem{DSILVA2018759}
Carmeline~J. Dsilva, Ronen Talmon, Ronald~R. Coifman, and Ioannis~G.
  Kevrekidis.
\newblock {Parsimonious representation of nonlinear dynamical systems through
  manifold learning: A chemotaxis case study}.
\newblock {\em Applied and Computational Harmonic Analysis}, 44(3):759 -- 773,
  2018.

\bibitem{gower1975generalized}
John~C Gower.
\newblock {Generalized procrustes analysis}.
\newblock {\em Psychometrika}, 40(1):33--51, 1975.

\bibitem{gower2004procrustes}
John~C Gower, Garmt~B Dijksterhuis, et~al.
\newblock {\em {Procrustes problems}}, volume~30.
\newblock Oxford University Press on Demand, 2004.

\bibitem{hein2007graph}
Matthias Hein, Jean-Yves Audibert, and Ulrike~von Luxburg.
\newblock {Graph laplacians and their convergence on random neighborhood
  graphs.}
\newblock {\em Journal of Machine Learning Research}, 8(6), 2007.

\bibitem{hintze1998violin}
Jerry~L Hintze and Ray~D Nelson.
\newblock {Violin plots: a box plot-density trace synergism}.
\newblock {\em The American Statistician}, 52(2):181--184, 1998.

\bibitem{jolliffe2016principal}
Ian~T Jolliffe and Jorge Cadima.
\newblock {Principal component analysis: a review and recent developments}.
\newblock {\em Philosophical Transactions of the Royal Society A: Mathematical,
  Physical and Engineering Sciences}, 374(2065):20150202, 2016.

\bibitem{jones2007universal}
Peter~W Jones, Mauro Maggioni, and Raanan Schul.
\newblock {Universal local parametrizations via heat kernels and eigenfunctions
  of the Laplacian}.
\newblock {\em arXiv preprint arXiv:0709.1975}, 2007.

\bibitem{kobak2021initialization}
Dmitry Kobak and George~C Linderman.
\newblock {Initialization is critical for preserving global data structure in
  both t-SNE and UMAP}.
\newblock {\em {Nature biotechnology}}, 39(2):156--157, 2021.

\bibitem{lafondiffusion2004}
S.S. Lafon.
\newblock {Diffusion Maps and Geometric Harmonics}.
\newblock {\em PhD Thesis}, page~45, 2004.

\bibitem{lederman2018learning}
Roy~R Lederman and Ronen Talmon.
\newblock {Learning the geometry of common latent variables using
  alternating-diffusion}.
\newblock {\em Applied and Computational Harmonic Analysis}, 44(3):509--536,
  2018.

\bibitem{li2019geodesic}
Didong Li and David~B Dunson.
\newblock {Geodesic distance estimation with spherelets}.
\newblock {\em arXiv preprint arXiv:1907.00296}, 2019.

\bibitem{maaten2008visualizing}
Laurens van~der Maaten and Geoffrey Hinton.
\newblock {Visualizing data using t-SNE}.
\newblock {\em Journal of machine learning research}, 9(Nov):2579--2605, 2008.

\bibitem{matlabprocrustes}
MATLAB.
\newblock {\em {Procrustes analysis, Statistics and Machine Learning Toolbox}}.
\newblock The MathWorks, Natick, MA, USA, 2018.

\bibitem{mcinnes2018umap}
Leland McInnes, John Healy, and James Melville.
\newblock {Umap: Uniform manifold approximation and projection for dimension
  reduction}.
\newblock {\em arXiv preprint arXiv:1802.03426}, 2018.

\bibitem{6377228}
G.~{Mishne} and I.~{Cohen}.
\newblock Multiscale anomaly detection using diffusion maps.
\newblock {\em IEEE Journal of Selected Topics in Signal Processing},
  7(1):111--123, 2013.

\bibitem{Mishne2018biorxiv}
Gal Mishne, Ronald~R. Coifman, Maria Lavzin, and Jackie Schiller.
\newblock Automated cellular structure extraction in biological images with
  applications to calcium imaging data.
\newblock {\em bioRxiv}, 2018.

\bibitem{Mishne2019DN}
Gal Mishne, Uri Shaham, Alexander Cloninger, and Israel Cohen.
\newblock Diffusion nets.
\newblock {\em Applied and Computational Harmonic Analysis}, 47(2):259 -- 285,
  2019.

\bibitem{peterfreund2020loca}
Erez Peterfreund, Ofir Lindenbaum, Felix Dietrich, Tom Bertalan, Matan Gavish,
  Ioannis~G Kevrekidis, and Ronald~R Coifman.
\newblock {LOCA: LOcal Conformal Autoencoder for standardized data
  coordinates}.
\newblock {\em arXiv preprint arXiv:2004.07234}, 2020.

\bibitem{roweis2000nonlinear}
Sam~T Roweis and Lawrence~K Saul.
\newblock {Nonlinear dimensionality reduction by locally linear embedding}.
\newblock {\em science}, 290(5500):2323--2326, 2000.

\bibitem{8450808}
N.~{Saito}.
\newblock {How Can We Naturally Order and Organize Graph Laplacian
  Eigenvectors?}
\newblock In {\em 2018 IEEE Statistical Signal Processing Workshop (SSP)},
  pages 483--487, 2018.

\bibitem{schonemann1966generalized}
Peter~H Sch{\"o}nemann.
\newblock A generalized solution of the orthogonal procrustes problem.
\newblock {\em Psychometrika}, 31(1):1--10, 1966.

\bibitem{singer2011orientability}
Amit Singer and Hau-tieng Wu.
\newblock {Orientability and diffusion maps}.
\newblock {\em Applied and computational harmonic analysis}, 31(1):44--58,
  2011.

\bibitem{doi:10.1080/03605302.2014.942739}
Stefan Steinerberger.
\newblock {Lower Bounds on Nodal Sets of Eigenfunctions via the Heat Flow}.
\newblock {\em Communications in Partial Differential Equations},
  39(12):2240--2261, 2014.

\bibitem{steinerberger2017spectral}
Stefan Steinerberger.
\newblock {On the spectral resolution of products of Laplacian eigenfunctions}.
\newblock {\em arXiv preprint arXiv:1711.09826}, 2017.

\bibitem{ten1977orthogonal}
Jos~MF Ten~Berge.
\newblock {Orthogonal Procrustes rotation for two or more matrices}.
\newblock {\em Psychometrika}, 42(2):267--276, 1977.

\bibitem{tenenbaum2000global}
Joshua~B Tenenbaum, Vin De~Silva, and John~C Langford.
\newblock {A global geometric framework for nonlinear dimensionality
  reduction}.
\newblock {\em science}, 290(5500):2319--2323, 2000.

\bibitem{trillos2020error}
Nicol{\'a}s~Garc{\'\i}a Trillos, Moritz Gerlach, Matthias Hein, and Dejan
  Slep{\v{c}}ev.
\newblock {Error estimates for spectral convergence of the graph Laplacian on
  random geometric graphs toward the Laplace--Beltrami operator}.
\newblock {\em Foundations of Computational Mathematics}, 20(4):827--887, 2020.

\bibitem{trillos2021large}
Nicolas~Garcia Trillos, Pengfei He, and Chenghui Li.
\newblock Large sample spectral analysis of graph-based multi-manifold
  clustering, 2021.

\bibitem{zelnik2005self}
Lihi Zelnik-Manor and Pietro Perona.
\newblock {Self-tuning spectral clustering}.
\newblock In {\em Advances in neural information processing systems}, pages
  1601--1608, 2005.

\bibitem{pmlr-v130-zhang21j}
Sharon Zhang, Amit Moscovich, and Amit Singer.
\newblock {Product Manifold Learning }.
\newblock In Arindam Banerjee and Kenji Fukumizu, editors, {\em Proceedings of
  The 24th International Conference on Artificial Intelligence and Statistics},
  volume 130 of {\em Proceedings of Machine Learning Research}, pages
  3241--3249. PMLR, 13--15 Apr 2021.

\bibitem{zhang2003nonlinear}
Zhenyue Zhang and Hongyuan Zha.
\newblock {Nonlinear dimension reduction via local tangent space alignment}.
\newblock In {\em {International Conference on Intelligent Data Engineering and
  Automated Learning}}, pages 477--481. Springer, 2003.

\end{thebibliography}

\end{document}